\documentclass[12pt]{article}
\usepackage{bm,subeqnarray,eqsection,indent,cite,url}
\usepackage{latexsym}
\usepackage{amsfonts}          
\usepackage{amssymb}           

\begin{document}

\date{May 30, 2011 \\[1mm] revised November 27, 2012}

\title{\vspace*{-2.5cm}
    \hspace*{-6mm}\hbox{Algebraic/combinatorial proofs of Cayley-type identities} \\
       for derivatives of determinants and pfaffians}

\author{
  {\small Sergio Caracciolo}                       \\[-1.7mm]
  {\small\it Dipartimento di Fisica and INFN}      \\[-1.7mm]
  {\small\it Universit\`a degli Studi di Milano}   \\[-1.7mm]
  {\small\it via Celoria 16}                       \\[-1.7mm]
  {\small\it I-20133 Milano, ITALY}                \\[-1.7mm]
  {\small\tt Sergio.Caracciolo@mi.infn.it}         \\[-1.7mm]
  {\protect\makebox[5in]{\quad}}
   \\
  {\small Alan D.~Sokal\thanks{Also at Department of Mathematics,
           University College London, London WC1E 6BT, England.}}  \\[-1.7mm]
  {\small\it Department of Physics}       \\[-1.7mm]
  {\small\it New York University}         \\[-1.7mm]
  {\small\it 4 Washington Place}          \\[-1.7mm]
  {\small\it New York, NY 10003 USA}      \\[-1.7mm]
  {\small\tt sokal@nyu.edu}           \\[-1.7mm]
  {\protect\makebox[5in]{\quad}}  
   \\
  {\small Andrea Sportiello\thanks{Address after 1 October 2012:
         LIPN, UMR CNRS 7030, Universit\'e Paris-Nord,
         99 avenue Jean-Baptiste Cl\'ement, 93430 Villetaneuse, FRANCE}}
                                                   \\[-1.7mm]
  {\small\it Dipartimento di Fisica and INFN}      \\[-1.7mm]
  {\small\it Universit\`a degli Studi di Milano}   \\[-1.7mm]
  {\small\it via Celoria 16}                       \\[-1.7mm]
  {\small\it I-20133 Milano, ITALY}                \\[-1.7mm]
  {\small\tt Andrea.Sportiello@mi.infn.it}         \\[-4mm]
  {\protect\makebox[5in]{\quad}}
}

\maketitle
\thispagestyle{empty}   

\vspace{-7mm}

\begin{abstract}
The classic Cayley identity states that
$$\det(\partial) \, (\det X)^s  \;=\; s(s+1) \cdots (s+n-1) \, (\det X)^{s-1}$$
where $X=(x_{ij})$ is an $n \times n$ matrix of indeterminates
and $\partial=(\partial/\partial x_{ij})$ is the corresponding matrix of
partial derivatives.
In this paper we present straightforward algebraic/combinatorial proofs
of a variety of Cayley-type identities, both old and new.
The most powerful of these proofs employ Grassmann algebra
(= exterior algebra) and Grassmann--Berezin integration.
Among the new identities proven here are
a pair of ``diagonal-parametrized'' Cayley identities,
a pair of ``Laplacian-parametrized'' Cayley identities,
and the ``product-parametrized'' and ``border-parametrized''
rectangular Cayley identities.
\end{abstract}


\bigskip
\noindent
{\bf Mathematics Subject Classification (MSC) codes:}
05A19 (Primary);
05E15, 05E99, 11S90, 13A50, 13N10, 14F10,
15A15, 15A23, 15A24, 15A33, 15A72, 15A75, 16S32,
20G05, 20G20, 32C38, 43A85, 81T18, 82B20 (Secondary).

\bigskip
\bigskip

\newtheorem{defin}{Definition}[section]
\newtheorem{definition}[defin]{Definition}
\newtheorem{prop}[defin]{Proposition}
\newtheorem{proposition}[defin]{Proposition}
\newtheorem{lem}[defin]{Lemma}
\newtheorem{lemma}[defin]{Lemma}
\newtheorem{guess}[defin]{Conjecture}
\newtheorem{ques}[defin]{Question}
\newtheorem{question}[defin]{Question}
\newtheorem{prob}[defin]{Problem}
\newtheorem{problem}[defin]{Problem}
\newtheorem{thm}[defin]{Theorem}
\newtheorem{theorem}[defin]{Theorem}
\newtheorem{cor}[defin]{Corollary}
\newtheorem{corollary}[defin]{Corollary}
\newtheorem{conj}[defin]{Conjecture}
\newtheorem{conjecture}[defin]{Conjecture}
\newtheorem{examp}[defin]{Example}
\newtheorem{example}[defin]{Example}
\newtheorem{claim}[defin]{Claim}

\renewcommand{\theenumi}{\alph{enumi}}
\renewcommand{\labelenumi}{(\theenumi)}
\def\prf{\par\noindent{\bf Proof.\enspace}\rm}
\def\rmk{\par\medskip\noindent{\bf Remark.\enspace}\rm}

\newcommand{\be}{\begin{equation}}
\newcommand{\ee}{\end{equation}}
\newcommand{\<}{\langle}
\renewcommand{\>}{\rangle}
\newcommand{\widebar}{\overline}
\def\reff#1{(\protect\ref{#1})}
\def\spose#1{\hbox to 0pt{#1\hss}}
\def\ltapprox{\mathrel{\spose{\lower 3pt\hbox{$\mathchar"218$}}
 \raise 2.0pt\hbox{$\mathchar"13C$}}}
\def\gtapprox{\mathrel{\spose{\lower 3pt\hbox{$\mathchar"218$}}
 \raise 2.0pt\hbox{$\mathchar"13E$}}}
\def\textprime{${}^\prime$}
\def\proof{\par\medskip\noindent{\sc Proof.\ }}
\newcommand{\qed}{\quad $\Box$ \medskip \medskip}
\def\firstproof{\par\medskip\noindent{\sc First proof.\ }}
\def\secondproof{\par\medskip\noindent{\sc Second proof.\ }}
\def\proofof#1{\bigskip\noindent{\sc Proof of #1.\ }}
\def\alternateproofof#1{\bigskip\noindent{\sc Alternate proof of #1.\ }}
\def\half{ {1 \over 2} }
\def\third{ {1 \over 3} }
\def\twothird{ {2 \over 3} }
\def\smfrac#1#2{{\textstyle{#1\over #2}}}
\def\smsmfrac#1#2{{\scriptstyle{#1\over #2}}}
\def\smhalf{ {\smfrac{1}{2}} }
\def\smsmhalf{ {\smsmfrac{1}{2}} }
\newcommand{\real}{\mathop{\rm Re}\nolimits}
\renewcommand{\Re}{\mathop{\rm Re}\nolimits}
\newcommand{\imag}{\mathop{\rm Im}\nolimits}
\renewcommand{\Im}{\mathop{\rm Im}\nolimits}
\newcommand{\sgn}{\mathop{\rm sgn}\nolimits}
\newcommand\supp{\mathop{\rm supp}\nolimits}
\newcommand{\diag}{\mathop{\rm diag}\nolimits}
\newcommand{\pf}{\mathop{\rm pf}\nolimits}
\newcommand{\hf}{\mathop{\rm hf}\nolimits}
\newcommand{\tr}{\mathop{\rm tr}\nolimits}
\newcommand{\per}{\mathop{\rm per}\nolimits}
\newcommand{\adj}{\mathop{\rm adj}\nolimits}
\newcommand{\Res}{\mathop{\rm Res}\nolimits}
\def\hboxscript#1{ {\hbox{\scriptsize\em #1}} }
\def\hboxrm#1{ {\hbox{\scriptsize\rm #1}} }

\newcommand{\restrict}{\upharpoonright}
\renewcommand{\emptyset}{\varnothing}

\def\Z{{\mathbb Z}}
\def\ZZ{{\mathbb Z}}
\def\R{{\mathbb R}}
\def\RR{{\mathbb R}}
\def\C{{\mathbb C}}
\def\CC{{\mathbb C}}
\def\N{{\mathbb N}}
\def\NN{{\mathbb N}}
\def\Q{{\mathbb Q}}
\def\QQ{{\mathbb Q}}

\newcommand{\scra}{{\mathcal{A}}}
\newcommand{\scrb}{{\mathcal{B}}}
\newcommand{\scrc}{{\mathcal{C}}}
\newcommand{\scrd}{{\mathcal{D}}}
\newcommand{\scre}{{\mathcal{E}}}
\newcommand{\scrf}{{\mathcal{F}}}
\newcommand{\scrg}{{\mathcal{G}}}
\newcommand{\scrh}{{\mathcal{H}}}
\newcommand{\scri}{{\mathcal{I}}}
\newcommand{\scrj}{{\mathcal{J}}}
\newcommand{\scrk}{{\mathcal{K}}}
\newcommand{\scrl}{{\mathcal{L}}}
\newcommand{\scrm}{{\mathcal{M}}}
\newcommand{\scrn}{{\mathcal{N}}}
\newcommand{\scro}{{\mathcal{O}}}
\newcommand{\scrp}{{\mathcal{P}}}
\newcommand{\scrq}{{\mathcal{Q}}}
\newcommand{\scrr}{{\mathcal{R}}}
\newcommand{\scrs}{{\mathcal{S}}}
\newcommand{\scrt}{{\mathcal{T}}}
\newcommand{\scru}{{\mathcal{U}}}
\newcommand{\scrv}{{\mathcal{V}}}
\newcommand{\scrw}{{\mathcal{W}}}
\newcommand{\scrx}{{\mathcal{X}}}
\newcommand{\scry}{{\mathcal{Y}}}
\newcommand{\scrz}{{\mathcal{Z}}}

\newcommand{\bgamma}{{\boldsymbol{\gamma}}}
\newcommand{\bsigma}{{\boldsymbol{\sigma}}}
\newcommand{\balpha}{{\boldsymbol{\alpha}}}
\newcommand{\bbeta}{{\boldsymbol{\beta}}}
\renewcommand{\pmod}[1]{\;({\rm mod}\:#1)}
\def\psibar{{\bar{\psi}}}
\def\etabar{{\bar{\eta}}}
\def\zetabar{{\bar{\zeta}}}
\def\chibar{{\bar{\chi}}}
\def\xibar{{\bar{\xi}}}
\def\lambdabar{{\bar{\lambda}}}
\def\mubar{{\bar{\mu}}}
\def\varphibar{{\bar{\varphi}}}
\def\phibar{{\bar{\phi}}}
\def\cz{\overline{z}}
\newcommand{\psietavert}{{\biggl( \!\! \begin{array}{cc} \psi \\ \eta
                                    \end{array} \!\! \biggr)}}
\newcommand{\psibaretabarvert}{{\biggl( \!\! \begin{array}{cc} \psibar \\ \etabar
                                    \end{array} \!\! \biggr)}}
\newcommand{\psietaprimevert}{{\biggl( \!\! \begin{array}{cc} \psi' \\ \eta'
                                    \end{array} \!\! \biggr)}}
\newcommand{\psibarzerovert}{{\biggl( \!\! \begin{array}{cc} \psibar \\ 0
                                        \end{array} \!\! \biggr)}}

\newcommand{\ch}  {\chi}
\newcommand{\al}  {\alpha}
\newcommand{\ta}  {\tau}
\newcommand{\si}  {\sigma}
\newcommand{\ep}  {\epsilon}
\newcommand{\Si}  {\Sigma}
\newcommand{\lp}  {\left(}
\newcommand{\rp}  {\right)}
\newcommand{\Br}  {\overline}
\newcommand{\la}{\langle}
\newcommand{\ra}{\rangle}

\newcommand{\XRL}{ {X^{\hbox{\rm\scriptsize row-Lap}}} }
\newcommand{\partialRL}{ {\partial^{\hbox{\rm\scriptsize row-Lap}}} }
\newcommand{\XCL}{ {X^{\hbox{\rm\scriptsize col-Lap}}} }
\newcommand{\partialCL}{ {\partial^{\hbox{\rm\scriptsize col-Lap}}} }
\newcommand{\XSL}{ {X^{\hbox{\rm\scriptsize sym-Lap}}} }
\newcommand{\partialSL}{ {\partial^{\hbox{\rm\scriptsize sym-Lap}}} }

\newcommand{\de}{\partial}
\newcommand{\ddotsinverse}
  {\textrm{\raisebox{-1.2mm}{$\cdot$}\,$\cdot$\,\raisebox{1.2mm}{$\cdot$}}}

%
%
\newcommand{\ofo}{ {{}_1 \! F_1} }
\newcommand{\tfo}{ {{}_2 \! F_1} }


\newenvironment{sarray}{
	  \textfont0=\scriptfont0
	  \scriptfont0=\scriptscriptfont0
	  \textfont1=\scriptfont1
	  \scriptfont1=\scriptscriptfont1
	  \textfont2=\scriptfont2
	  \scriptfont2=\scriptscriptfont2
	  \textfont3=\scriptfont3
	  \scriptfont3=\scriptscriptfont3
	\renewcommand{\arraystretch}{0.7}
	\begin{array}{l}}{\end{array}}

\newenvironment{scarray}{
	  \textfont0=\scriptfont0
	  \scriptfont0=\scriptscriptfont0
	  \textfont1=\scriptfont1
	  \scriptfont1=\scriptscriptfont1
	  \textfont2=\scriptfont2
	  \scriptfont2=\scriptscriptfont2
	  \textfont3=\scriptfont3
	  \scriptfont3=\scriptscriptfont3
	\renewcommand{\arraystretch}{0.7}
	\begin{array}{c}}{\end{array}}

\newcommand{\bydef}{:=}
\newcommand{\defby}{=:}
\newcommand{\eqdef} {\stackrel{\rm def}{=}}
\newcommand{\implies}{\Longrightarrow}
\newcommand{\binom}[2]{\left(#1\atop#2\right)}
\newcommand{\bigpartial}{\displaystyle\partial}

%
\newcommand{\ef}[1]{\, #1}     

\newcommand{\Reof}[1]{\mathfrak{Re}(#1)}
\newcommand{\Imof}[1]{\mathfrak{Im}(#1)}
\newcommand{\eval}[1]{\left\langle {#1} \right\rangle}
\newcommand{\leval}[1]{\langle {#1} \rangle}
\newcommand{\reval}[1]{\overline{#1}}

\newcommand{\sspan}{\mathrm{span}} 
\newcommand{\kker}{\mathrm{ker}}
\newcommand{\rrank}{\mathrm{rank}}

\newcommand{\bigast}[1]{\underset{#1}{\textrm{{\huge $\ast$}}}}

\newcommand{\dx}[1] {\mathrm{d}{#1}}
\newcommand{\dede}[1]{\frac{\partial}{\partial #1}}
\newcommand{\deenne}[2]{\frac{\partial^#2}{\partial #1 ^#2}}
\newcommand{\vett}[1]{#1}

\newcommand{\tinyfrac}[2] {\genfrac{}{}{}{1}{#1}{#2} }
\newcommand{\Lfrac}[2] {\genfrac{}{}{}{0}{#1}{#2} }

\newtheorem{ansatz}{Ansatz}[section]
\newtheorem{theor}{Theorem}
\newtheorem{coroll}{Corollary}

\tableofcontents
\clearpage

\section{Introduction}

Let $X=(x_{ij})$ be an $n \times n$ matrix of indeterminates, and let
$\partial=(\partial/\partial x_{ij})$ be the corresponding matrix of
partial derivatives.
The following beautiful identity is conventionally\footnote{
   But erroneously:  see Section~\ref{sec.historical} below.
}
attributed to Arthur Cayley (1821--1895):
\be
   \det(\partial) \, (\det X)^s  \;=\;
   s(s+1) \cdots (s+n-1) \, (\det X)^{s-1}
   \;.
 \label{eq.intro.1}
\ee
[When $n=1$ this is of course the elementary formula $(d/dx) x^s = s x^{s-1}$.]
A generalization of \reff{eq.intro.1} to arbitrary minors also holds,
and is sometimes\footnote{
   Also erroneously:  see again Section~\ref{sec.historical}.
}
 attributed to Alfredo Capelli (1855--1910):
if $I,J \subseteq \{1,\ldots,n\}$ with $|I| = |J| = k$, then
\be
   \det(\partial_{IJ}) \, (\det X)^s  \;=\;
    s(s+1) \cdots (s+k-1) \, (\det X)^{s-1} \,
    \epsilon(I,J) \, (\det X_{I^c J^c})
 \label{eq.intro.2}
\ee
where $\epsilon(I,J) = (-1)^{\sum_{i \in I} i + \sum_{j \in J} j}$.
Analogous identities for symmetric and antisymmetric matrices were proved
by G{\aa}rding \cite{Garding_48} in 1948
and Shimura \cite{Shimura_84} in 1984, respectively.

Although these identities are essentially algebraic or combinatorial in nature,
the simplest proofs currently available in the literature are analytic,
exploiting Laplace-type integral representations for $(\det X)^s$
\cite{Shimura_84,Faraut_94}.
Indeed, most of the existing algebraic/combinatorial proofs
\cite{Grace_03,Turnbull_28,Garding_48,Canfield_81,Stein_82,Williamson_87,%
Dolgachev_03}
are somewhat difficult to follow,
partly because of old-fashioned notation.\footnote{
   Among the exceptions are \cite[Theorem~1.3 and Lemma~2.12]{Canfield_81}
   and \cite[Lemma~2.1]{Dolgachev_03}.
}
In this paper we would like to give straightforward algebraic/combinatorial
proofs of a variety of Cayley-type identities, some of which are known
and others of which are new.  The most powerful of these proofs
employ Grassmann algebra (= exterior algebra)
and Grassmann--Berezin integration.\footnote{
   We stress that our proofs are {\em not}\/ ``combinatorial''
   in the narrow sense of exhibiting bijections or exploiting double-counting.
   Rather, our proofs are based on straightforward algebraic manipulation
   combined with some elementary arguments of enumerative combinatorics
   (enumeration of permutations by number of cycles,
    identities involving binomial coefficients, etc.).
   For this reason we have opted to describe our methods
   using the rather awkward adjective ``algebraic/combinatorial''.
}

Nowadays, identities like \reff{eq.intro.1} are best understood as
calculations of Bernstein--Sato type \cite{Bjork_79,Coutinho_95,Krause_00}
for special polynomials.
To see what is at issue,
let $P(x_1,\ldots,x_n) \not\equiv 0$ be a polynomial in $n$ variables with
coefficients in a field $K$ of characteristic 0.
Then Bernstein \cite{Bernstein_72} proved in 1972 that there exist
a polynomial-coefficient partial differential operator
$Q(s,x,\partial/\partial x)$ and a polynomial $b(s) \not\equiv 0$
(both with coefficients in $K$)
satisfying
\be
   Q(s,x,\partial/\partial x) \, P(x)^s
   \;=\;
   b(s) \, P(x)^{s-1}
   \;.
 \label{eq.intro.bernstein}
\ee
We call any pair $(Q,b)$ satisfying \reff{eq.intro.bernstein}
a  {\em Bernstein--Sato pair}\/ for $P$.

The set of all $b$ for which there exists a $Q$
satisfying \reff{eq.intro.bernstein}
is easily seen to be an ideal in the polynomial ring $K[s]$.
By Bernstein's theorem this ideal is nontrivial,
so it is generated by a unique monic polynomial $b(s)$,
called the {\em Bernstein--Sato polynomial}\/ (or {\em b-function}\/)
of $P$.%
\footnote{
    In the literature on Bernstein--Sato equations it is customary
    to shift our $s$ by 1, i.e.~write
    $Q(s,x,\partial/\partial x) P(x)^{s+1} = b(s) P(x)^{s}$,
   so that the usual Bernstein--Sato polynomial is our $b(s+1)$. We
   choose here the slightly unconventional notation 
   (\ref{eq.intro.bernstein}) because it seems better adapted to the
   Cayley identity (\ref{eq.intro.1}).
 \label{footnote.shifted}
}
Cayley-type identities thus provide Bernstein--Sato pairs for certain
polynomials $P$ arising from determinants.

Bernstein--Sato pairs are especially useful in treating
the problem of analytically continuing the distribution $P_\Omega^s$,
which can be posed as follows
\cite{Gelfand_54,Gelfand_64,Bernstein_69,Atiyah_70,Bernstein_72,%
Bjork_79,Coutinho_95,Igusa_00,Krause_00}:
Let $P(x_1,\ldots,x_n) \not\equiv 0$ be a polynomial with real coefficients,
and let $\Omega \subseteq \R^n$ be an open set such that
$P \ge 0$ on $\Omega$ and $P=0$ on $\partial\Omega$.
Then, for any complex number $s$ satisfying $\real s > 0$,
the function $P^s$ is well-defined on $\Omega$ and polynomially bounded,
and thus defines a tempered distribution $P_\Omega^s \in \scrs'(\R^n)$
by the formula
\be
   \< P_\Omega^s, \varphi \>
   \;=\;
   \int\limits_\Omega P(x)^s \, \varphi(x) \, dx
\ee
for any test function $\varphi \in \scrs(\R^n)$.
Furthermore, the function $s \mapsto \< P_\Omega^s, \varphi \>$
is analytic on the half-plane $\real s > 0$,
with complex derivative given by
\be
   {d \over ds} \,
   \< P_\Omega^s, \varphi \>
   \;=\;
   \int\limits_\Omega P(x)^s \, (\log P(x)) \, \varphi(x) \, dx
   \;.
\ee
Thus $P_\Omega^s$ is a distribution-valued analytic
function of $s$ on the right half-plane.
We want to know whether $P_\Omega^s$ can be analytically continued
to the whole complex plane as a meromorphic function of $s$.
This problem was first posed by I.M.~Gel'fand \cite{Gelfand_54}
at the 1954 International Congress of Mathematicians.
It was answered affirmatively in 1969 independently by
Bernstein and S.I.~Gel'fand \cite{Bernstein_69} and Atiyah \cite{Atiyah_70},
using deep results from algebraic geometry
(Hironaka's resolution of singularities \cite{Hironaka_64}).
A few years later, Bernstein \cite{Bernstein_72} produced a
much simpler proof based on using
the Bernstein--Sato equation \reff{eq.intro.bernstein}
to extend $P_\Omega^s$ successively to half-planes
$\real s > -1$, $\real s > -2$, etc.
See e.g.\ \cite[sections~7.1 and~7.3]{Bjork_79}
or \cite[section~5.3]{Igusa_00} for details.

The special case in which $P$ is a determinant
of a symmetric or hermitian matrix
(and $\Omega$ is e.g.\ the cone of positive-definite matrices)
has been studied by several authors
\cite{Rais_72,Blekher_82,Ricci_86,Faraut_94,Blind_97,Muro_99,Muro_01,Muro_03};
it plays a central role in the theory of Riesz distributions
on Euclidean Jordan algebras (or equivalently on symmetric cones)
\cite[Chapter~VII]{Faraut_94}.
This case is also useful in quantum field theory in studying the analytic
continuation of Feynman integrals to ``complex space-time dimension''
\cite{Speer_76,Blekher_82,Etingof_99}.
In an analogous way, the parametrized symmetric Cayley identity
(Theorem~\ref{thm.para.sym.cayley} below)
will play a key role in studying the analytic continuation of integrals
over products of spheres $S^{N-1} \subset \R^N$ to ``complex dimension $N$''
\cite{Caracciolo_et_al_in_prep},
with the aim of giving a rigorous nonperturbative formulation
of the correspondence found in \cite{Caracciolo_04,CSS_hyperforests}
between spanning forests and the $N$-vector model in
statistical mechanics at $N=-1$.
This latter application was, in fact,
our original motivation for studying Cayley-type identities.
The original Cayley identity \reff{eq.intro.1}
was also rediscovered by Creutz \cite{Creutz_78}
and used by him to compute certain invariant integrals over $SU(n)$
that arise in lattice gauge theory.

Many of the polynomials $P$ treated here can also be understood
as relative invariants of prehomogeneous vector spaces
\cite{Igusa_00,Kimura_03}.
When applicable, this connection allows the immediate identification
of a suitable operator $Q(\partial/\partial x)$
--- namely, the dual of $P$ itself ---
and provides a general proof that the corresponding $b(s)$
satisfies $\deg b = \deg P$
and is indeed (up to a constant factor)
the Bernstein--Sato polynomial of $P$.\footnote{
   See \cite[Corollary~6.1.1 and Theorem~6.1.1]{Igusa_00}
   \cite[Proposition~2.22]{Kimura_03}
   for the first two points,
   and \cite[Theorem~6.3.2]{Igusa_00} for the third.
}
Furthermore, this approach sometimes allows the explicit calculation
of $b(s)$ by means of microlocal calculus
\cite{Sato_80,Kimura_82,Ozeki_80,Ozeki_90,Ukai_03,%
Wakatsuki_04,Sugiyama_05,Sato_06,Sugiyama_11}
or other methods \cite{Sugiyama_02}.\footnote{
   We are grateful to Nero Budur for explaining to us
   the connection between our results and
   the theory of prehomogeneous vector spaces.
}

The purpose of the present paper is to give straightforward
(and we hope elegant) algebraic/combinatorial proofs
of a variety of Cayley-type identities, both old and new.
Since our main aim is to illustrate proof techniques that may be useful
in other contexts, we shall give, wherever possible,
several alternate proofs of each result.
One purpose of this paper is, in fact, to make propaganda
among mathematicians for the power of Grassmann--Berezin integration
as a tool for proving algebraic or combinatorial identities.
Among the new results in this paper are
the ``diagonal-parametrized'' Cayley identities
(Theorems~\ref{thm.para.cayley} and \ref{thm.para.sym.cayley}),
the ``Laplacian-parametrized'' Cayley identities
(Theorems~\ref{thm.Lap.cayley} and \ref{thm.Lap.sym.cayley}),
and the ``product-parametrized'' and ``border-parametrized''
rectangular Cayley identities
(Theorems~\ref{thm.productcayley} and \ref{thm.borderedcayley}).
We also give an elementary (though rather intricate)
proof of the multi-matrix rectangular Cayley identity
(Theorem~\ref{thm.multirectcayley})
that was proven recently by Sugiyama \cite{Sugiyama_11}.

The plan of this paper is as follows:
In Section~\ref{sec.statement} we state the identities to be proven
and briefly discuss their interpretation.
In Section~\ref{sec.elementary} we give elementary algebraic/combinatorial proofs
of the three basic Cayley-type identities
(ordinary, symmetric and antisymmetric).
In Section~\ref{sec.grassmann.1} we give very simple proofs
of these same identities, based on representing $(\det X)^s$ as a
fermionic or bosonic Gaussian integral.\footnote{
   The proofs based on representing $(\det X)^s$ as a
   bosonic Gaussian integral are closely related to the
   existing analytic proofs \cite{Shimura_84,Faraut_94}.
   The proofs using {\em fermionic}\/ Gaussian integrals,
   by contrast, are really {\em algebraic/combinatorial}\/ proofs,
   as fermionic ``integration'' is a purely algebraic/combinatorial
   construction (see Appendix~\ref{app.grassmann.4} below).
}
In Section~\ref{sec.grassmann.2} we give alternate
(and arguably even simpler) proofs,
based on representing $\det(\partial)$ as a fermionic Gaussian integral;
this method is very powerful and
allows us to prove also the (considerably more difficult)
``rectangular Cayley identities''.
In Section~\ref{sec.param} we prove the diagonal-parametrized
Cayley identities, in Section~\ref{sec.Lap} we prove the
Laplacian-parametrized Cayley identities,
and in Section~\ref{sec.border} we prove the product-parametrized
and border-parametrized rectangular Cayley identities.
Finally, in Section~\ref{sec.conj.minimality} we formulate some conjectures
concerning the minimality of our Bernstein--Sato pairs.
In Appendix~\ref{app.grassmann} we provide a brief introduction to
Grassmann algebra and Grassmann--Berezin (fermionic) integration;
we hope that this appendix will prove useful to mathematicians
seeking a mathematically rigorous (but conceptually elementary)
presentation of this powerful algebraic/combinatorial tool.\footnote{
     Unfortunately, most of the existing presentations of
     Grassmann--Berezin integration are aimed at physicists
     (see e.g.\ \cite[Chapter~1]{Zinn-Justin} for an excellent treatment)
     and may not meet mathematicians' standards of precision and rigor,
     or else are aimed at applications to differential geometry
     and hence involve heavier theoretical machinery than is
     needed for combinatorial applications.
     One exception is the brief summary given by
     Abdesselam \cite[Section~2]{Abdesselam_03};
     our presentation can be viewed as an enlargement of his.
  }
In Appendix~\ref{app.identities} we collect some formulae
that will be needed in the proofs.

We have tried hard to write this paper in a ``modular'' fashion,
so that the reader can skip around according to his/her interests
without having to read the whole thing.
Indeed, after a brief perusal of Section~\ref{sec.statement},
the reader can proceed directly to Section~\ref{sec.elementary},
\ref{sec.grassmann.1} or \ref{sec.grassmann.2} as desired,
consulting Appendices~\ref{app.grassmann} and \ref{app.identities} as needed.

\section{Statement of main results}  \label{sec.statement}

{\bf Notation:}
We write $[n] = \{1,\ldots,n\}$.
Give a matrix $A$, we denote its transpose by $A^{\rm T}$.
For an invertible square matrix $A$, we use the shorthand $A^{\rm -T}$
for $(A^{-1})^{\rm T}$ (which is also equal to
$(A^{\rm T})^{-1}$).
If $A = (a_{ij})_{i,j=1}^n$ is an $n \times n$ matrix
and $I,J \subseteq [n]$,
we denote by $A_{IJ}$ the submatrix of $A$
corresponding to the rows $I$ and the columns $J$,
all kept in their original order.
We write $I^c$ to denote the complement of $I$ in $[n]$.
We define $\epsilon(I) = (-1)^{|I| (|I|-1)/2} (-1)^{\sum_{i \in I} i}$;
it is the sign of the permutation
that takes the sequence $1 \cdots n$ into $I I^c$
when the sets $I$ and $I^c$ are each written in increasing order.
We also define $\epsilon(I,J) = \epsilon(I) \epsilon(J)$;
it is the sign of the permutation that takes $I I^c$ into $J J^c$.
In particular, if $|I| = |J|$, we have
$\epsilon(I,J) = (-1)^{\sum_{i \in I} i + \sum_{j \in J} j}$.

\subsection{Ordinary, symmetric, antisymmetric and hermitian Cayley identities}

The basic Cayley-type identity is the following:

\begin{theorem}[ordinary Cayley identity]
 \label{thm.cayley}
Let $X=(x_{ij})$ be an $n \times n$ matrix of indeterminates, and let
$\partial=(\partial/\partial x_{ij})$ be the corresponding matrix of
partial derivatives.  Then
\be
   \det(\partial) \, (\det X)^s  \;=\;  s(s+1) \cdots (s+n-1) \, (\det X)^{s-1}
   \;.
 \label{eq.cayley.1}
\ee
More generally, if $I,J \subseteq [n]$ with $|I| = |J| = k$, then
\be
   \det(\partial_{IJ}) \, (\det X)^s  \;=\;
    s(s+1) \cdots (s+k-1) \, (\det X)^{s-1} \,
    \epsilon(I,J) \, (\det X_{I^c J^c})
   \;.
 \label{eq.cayley.2}
\ee
\end{theorem}

\smallskip

{\bf Remark.}  Since $\det(\partial)$ and $\det(\partial_{IJ})$
are constant-coefficient differential operators,
the matrix $X$ can everywhere be replaced by $X+A$
for any fixed matrix $A$,
and the identities \reff{eq.cayley.1}/\reff{eq.cayley.2} remain valid.

\bigskip

We consider next a version of the Cayley identity for {\em symmetric}\/
matrices $X^{\rm sym}=(x_{ij})$.  What this means is that only the
variables $(x_{ij})_{1 \le i \le j \le n}$ are taken as
independent indeterminates; then $x_{ij}$ for $i > j$ is regarded
as a synonym for $x_{ji}$.

\begin{theorem}[symmetric Cayley identity]
 \label{thm.symcayley}
Let $X^{\rm sym}=(x_{ij})$ be an $n \times n$ symmetric matrix of
indeterminates,
and let $\partial^{\rm sym}$ be the matrix whose elements are
\be
   (\partial^{\rm sym})_{ij}
    \;=\;  \smhalf (1 + \delta_{ij})
                   {\bigpartial \over \bigpartial x_{ij}}
    \;=\;  \cases{ \partial / \partial x_{ii}  & if $i=j$  \cr
                   \noalign{\vskip 4pt}
                   \smhalf \partial / \partial x_{ij}   & if $i < j$ \cr
                   \noalign{\vskip 4pt}
                   \smhalf \partial / \partial x_{ji}   & if $i > j$ \cr
                 }
 \label{def.partial.sym}
\ee
Then
\be
   \det(\partial^{\rm sym}) \, (\det X^{\rm sym})^s
   \;=\;  s(s+\smhalf) \cdots \left(s+ {n-1 \over 2} \right) \,
          (\det X^{\rm sym})^{s-1}
   \;.
 \label{eq.symcayley.1}
\ee
More generally, if $I,J \subseteq [n]$ with $|I| = |J| = k$, then
\be
   \det(\partial^{\rm sym}_{IJ}) \, (\det X^{\rm sym})^s  \;=\;
      s(s+\smhalf) \cdots \left(s+ {k-1 \over 2} \right) \,
      (\det X^{\rm sym})^{s-1} \, \epsilon(I,J) \, (\det X^{\rm sym}_{I^c J^c})
   \;.
 \label{eq.symcayley.2}
\ee
\end{theorem}

{\bf Remarks.}  1.  The matrix $X^{\rm sym}$ can everywhere
be replaced by $X^{\rm sym}+A$ for any fixed {\em symmetric}\/ matrix $A$.

2. If we prefer to work over the integers rather than the
rationals, it suffices to multiply $\partial^{\rm sym}$ by 2
and correspondingly multiply the right-hand side of
\reff{eq.symcayley.1} [resp.\ \reff{eq.symcayley.2}]
by $2^n$ (resp.\ $2^k$).

\bigskip

Next let us state a version of the Cayley identity for
{\em antisymmetric}\/ matrices $X^{\rm antisym}=(x_{ij})$.\footnote{
   More precisely, {\em alternating}\/ matrices,
   i.e.\ matrices satisfying $x_{ij} = -x_{ji}$ {\em and}\/ $x_{ii} = 0$.
   The latter identity is a consequence of the former whenever the
   underlying ring of coefficients is an integral domain of
   characteristic $\neq 2$ (so that $2x=0$ implies $x=0$),
   but not in general otherwise.
   See e.g.\ \cite[sections~XIII.6 and XV.9]{Lang_02}.
   In this paper we use the term ``antisymmetric'' to denote
   $x_{ij} = -x_{ji}$ {\em and}\/ $x_{ii} = 0$.
}
Here only the variables $(x_{ij})_{1 \le i < j \le n}$
are taken as independent indeterminates;
then $x_{ij}$ for $i > j$ is regarded as a synonym for $-x_{ji}$,
and $x_{ii}$ is regarded as a synonym for 0.
As befits antisymmetric matrices, the corresponding identity
involves pfaffians in place of determinants:

\begin{theorem}[antisymmetric Cayley identity]
 \label{thm.antisymcayley}
Let $X^{\rm antisym}=(x_{ij})$ be a $2m \times 2m$ antisymmetric matrix of
indeterminates,
and let $\partial^{\rm antisym}$ be the corresponding matrix
of partial derivatives, i.e.
\be
   (\partial^{\rm antisym})_{ij}
    \;=\;  \cases{ 0     & if $i=j$  \cr
                   \noalign{\vskip 4pt}
                   \partial / \partial x_{ij}   & if $i < j$ \cr
                   \noalign{\vskip 4pt}
                   -\, \partial / \partial x_{ji}   & if $i > j$ \cr
                 }
\ee
Then
\be
    \pf(\partial^{\rm antisym}) \, (\pf X^{\rm antisym})^s  \;=\;
        s(s+2) \cdots (s+2m-2) \, (\pf X^{\rm antisym})^{s-1}
   \;.
 \label{eq.antisymcayley.1}
\ee
More generally, if $I \subseteq [2m]$ with $|I| = 2k$, then
\be
    \pf(\partial^{\rm antisym}_{II}) \, (\pf X^{\rm antisym})^s  \;=\;
        s(s+2) \cdots (s+2k-2) \, (\pf X^{\rm antisym})^{s-1}
                               \, \epsilon(I) 
                               \, (\pf X^{\rm antisym}_{I^c I^c})
   \;.
 \label{eq.antisymcayley.2}
\ee
\end{theorem}

As an immediate corollary we get a result for antisymmetric determinants:

\begin{corollary}[antisymmetric Cayley identity for determinants]
 \label{cor.antisymcayley}
Let $X^{\rm antisym}=(x_{ij})$ be a $2m \times 2m$ antisymmetric matrix of
indeterminates,
and let $\partial^{\rm antisym}$ be the corresponding matrix of
partial derivatives.  Then
\be
    \det(\partial^{\rm antisym}) \, (\det X^{\rm antisym})^s  \;=\;
        (2s-1)(2s) \cdots (2s+2m-2) \, (\det X^{\rm antisym})^{s-1}
   \;.
 \label{eq.antisymcayley.3}
\ee
More generally, if $I \subseteq [2m]$ with $|I| = 2k$, then
\be
    \det(\partial^{\rm antisym}_{II}) \, (\det X^{\rm antisym})^s  \;=\;
        (2s-1)(2s) \cdots (2s+2k-2) \, (\det X^{\rm antisym})^{s-1}
                               \, (\det X^{\rm antisym}_{I^c I^c})
   \;.
 \label{eq.antisymcayley.4}
\ee
\end{corollary}

\noindent
Please note that in the antisymmetric case we are able at present
to handle only principal minors,
i.e.\ we have been unable to find a general formula for
$ \det(\partial^{\rm antisym}_{IJ}) \, (\det X^{\rm antisym})^s$
when $I \neq J$.

\bigskip

Next we state a version of the Cayley identity for
{\em ``hermitian''}\/ matrices $Z^{\rm herm}=(z_{ij})$.\footnote{
   We put ``hermitian'' in quotation marks because the variables
   $(x_{ij})$ and $(y_{ij})$ in this identity
   are neither real nor complex numbers,
   but are simply indeterminates.
}
By this we mean the following:
We introduce indeterminates $(x_{ij})_{1 \le i \le j \le n}$
and $(y_{ij})_{1 \le i < j \le n}$,
and define matrices $X^{\rm sym}$ and $Y^{\rm antisym}$
as before;
we then set $Z^{\rm herm}= X^{\rm sym} + i Y^{\rm antisym}$
and $\partial^{\rm herm} =
     \partial_X^{\rm sym} - (i/2) \partial_Y^{\rm antisym}$.
In terms of the usual complex derivatives
\be
   {\partial \over \partial z_{ij}}  \;=\;
   \half \left( {\partial \over \partial x_{ij}}
                \,-\,  i \, {\partial \over \partial y_{ij}} \right)
   ,\qquad
   {\partial \over \partial \bar{z}_{ij}}  \;=\;
   \half \left( {\partial \over \partial x_{ij}}
                \,+\,  i \, {\partial \over \partial y_{ij}} \right)
\ee
for $i < j$,
this can be written as
\be
   (\partial^{\rm herm})_{ij}
    \;=\;  \cases{ \partial / \partial x_{ii}  & if $i=j$  \cr
                   \noalign{\vskip 4pt}
                   \partial / \partial z_{ij}   & if $i < j$ \cr
                   \noalign{\vskip 4pt}
                   \partial / \partial \bar{z}_{ji}   & if $i > j$ \cr
                 }
 \label{def.partial.herm}
\ee
We then have:

\begin{theorem}[``hermitian'' Cayley identity]
 \label{thm.hermcayley}
Let $Z^{\rm herm}$ and $\partial^{\rm herm}$ be defined as above.
Then
\be
   \det(\partial^{\rm herm}) \, (\det Z^{\rm herm})^s  
\;=\;  s(s+1) \cdots (s+n-1) \, (\det Z^{\rm herm})^{s-1}
   \;.
 \label{eq.hermcayley.1}
\ee
More generally, if $I,J \subseteq [n]$ with $|I| = |J| = k$, then
\be
   \det(\partial^{\rm herm}_{IJ}) \, (\det Z^{\rm herm})^s  \;=\;
    s(s+1) \cdots (s+k-1) \, (\det Z^{\rm herm})^{s-1} \,
    \epsilon(I,J) \, (\det Z^{\rm herm}_{I^c J^c})
   \;.
 \label{eq.hermcayley.2}
\ee
\end{theorem}

The resemblance of this theorem to the
{\em ordinary}\/ Cayley identity \reff{eq.cayley.1}/\reff{eq.cayley.2}
is no accident;
indeed, the two identities are immediately interderivable.
To see this, it suffices to notice that the action of derivatives
on indeterminates is identical in the two cases:
\begin{subeqnarray}
\partial_{ij} X_{i'j'} & = & \delta_{i,i'} \delta_{j,j'}   \\[1mm]
\partial^{\rm herm}_{ij} Z^{\rm herm}_{i'j'}
    & = & \delta_{i,i'} \delta_{j,j'}
 \label{eq.derivs.herm}
\end{subeqnarray}
and that this relation completely determines the expressions
\reff{eq.cayley.1}/\reff{eq.cayley.2} and
\reff{eq.hermcayley.1}/\reff{eq.hermcayley.2}.
For this reason, we need not consider further
the ``hermitian'' Cayley identity.\footnote{
   Let us remark that, by contrast, in the analytic proofs
   \cite{Shimura_84,Faraut_94} it is convenient to consider
   the hermitian case instead of the ordinary case,
   as $(\det Z)^s$ for a positive-definite complex hermitian matrix $Z$
   has a simple Laplace-type integral representation.
}

\subsection{Rectangular Cayley identities}  \label{sec.statement.rectangular}


Let us now formulate some Cayley-type identities for
{\em rectangular}\/ matrices of size $m \times n$ with $m \le n$.
These identities are somewhat more complicated than the preceding ones,
because the matrices appearing in the determinants (or pfaffians)
are {\em quadratic}\/ (or of higher order)
rather than linear in the indeterminates.

\begin{theorem}[two-matrix rectangular Cayley identity]
 \label{thm.tmrectcayley}
\ Let $X=(x_{ij})$ and
\linebreak
$Y=(y_{ij})$ be $m \times n$ matrices of indeterminates
with $m \le n$,
and let $\partial_X=(\partial/\partial x_{ij})$
and $\partial_Y=(\partial/\partial y_{ij})$
be the corresponding matrices of partial derivatives.  Then
\be
   \det(\partial_X \partial_Y^{\rm T}) \, \det (X Y^{\rm T})^s
   \;=\;  \left( \prod\limits_{j=0}^{m-1} (s+j)(s+n-m+j) \right)
          \, \det (X Y^{\rm T})^{s-1}
   \;.
 \label{eq.tmrectcayley.1}
\ee
More generally, if $I,J \subseteq [m]$ with $|I| = |J| = k$, then
\begin{eqnarray}
   \det[(\partial_X \partial_Y^{\rm T})_{IJ}] \, \det (X Y^{\rm T})^s
   & = &
   \left( \prod\limits_{j=0}^{k-1} (s+j)(s+n-m+j) \right)
          \, \det (X Y^{\rm T})^{s-1}
      \,\times  \nonumber \\
   & & \qquad
   \epsilon(I,J) \det[(X Y^{\rm T})_{I^c J^c}]
   \;.
 \label{eq.tmrectcayley.2}
\end{eqnarray}
\end{theorem}

\noindent
If $m=1$, Theorem~\ref{thm.tmrectcayley} reduces to the easily-derived formula
$(\nabla_{\bf x} \cdot \nabla_{\bf y}) ({\bf x} \cdot {\bf y})^s
 = s (s+n-1) ({\bf x} \cdot {\bf y})^{s-1}$
for ${\bf x} = (x_1,\ldots,x_n)$ and ${\bf y} = (y_1,\ldots,y_n)$.
If $m=n$, Theorem~\ref{thm.tmrectcayley} can be derived by
separate applications of Theorem~\ref{thm.cayley} to $X$ and $Y$.
In other cases it appears to be new.

\begin{theorem}[one-matrix rectangular symmetric Cayley identity]
 \label{thm.rectcayley}
\quad\break
Let $X=(x_{ij})$ be an $m \times n$ matrix of indeterminates with $m \le n$,
and let $\partial=(\partial/\partial x_{ij})$ be the corresponding matrix of
partial derivatives.  Then
\be
   \det(\partial \partial^{\rm T}) \, \det (X X^{\rm T})^s
   \;=\;  \left( \prod\limits_{j=0}^{m-1} (2s+j)(2s+n-m-1+j) \right)
          \, \det (X X^{\rm T})^{s-1}
   \;.
 \label{eq.rectcayley.1}
\ee
More generally, if $I,J \subseteq [m]$ with $|I| = |J| = k$, then
\begin{eqnarray}
   \det[(\partial \partial^{\rm T})_{IJ}] \, \det (X X^{\rm T})^s
   & = &
   \left( \prod\limits_{j=0}^{k-1} (2s+j)(2s+n-m-1+j) \right)
          \, \det (X X^{\rm T})^{s-1}
      \,\times  \nonumber \\
   & & \qquad
   \epsilon(I,J) \det[(X X^{\rm T})_{I^c J^c}]
   \;.
 \label{eq.rectcayley.2}
\end{eqnarray}
\end{theorem}

\noindent
If $m=1$, Theorem~\ref{thm.rectcayley} reduces to the well-known formula
$\Delta ({\bf x}^2)^s = 2s(2s+n-2) ({\bf x}^2)^{s-1}$
for ${\bf x} = (x_1,\ldots,x_n)$.
If $m=n$, Theorem~\ref{thm.rectcayley} can be derived from
two applications of Theorem~\ref{thm.cayley}.
The general case has been proven recently by several authors,
using analytic methods.\footnote{
   See Faraut--Kor\'anyi \cite[section~XVI.4]{Faraut_94},
   Kh\'ekalo \cite{Khekalo_01,Khekalo_05}
   and Rubin \cite{Rubin_06}.
}
We call Theorem~\ref{thm.rectcayley} a ``symmetric'' identity
because the matrices $X X^{\rm T}$ and $\partial \partial^{\rm T}$
that appear in it are symmetric by construction.


Finally, here is an analogue of Theorem~\ref{thm.rectcayley}
that involves matrices that are {\em antisymmetric}\/
rather than symmetric by construction;
in place of the identity matrix we use
the standard $2n \times 2n$ symplectic form
\be
J \;=\;
\left(
\begin{array}{cc|cc|c}
 0 & 1 & \multicolumn{3}{c}{} \\
-1 & 0 & \multicolumn{3}{c}{} \\
\cline{1-4}
 &  & 0 & 1 & \\
 &  & -1 & 0 & \\
\cline{3-4}
\multicolumn{4}{c}{} & \ddots
\end{array}
\right)  \;.
  \label{def.J.intro}
\ee
Not surprisingly, this formula involves pfaffians rather than determinants:

\begin{theorem}[one-matrix rectangular antisymmetric Cayley identity]
 \label{thm.antisymrectcayley}
\quad\break
Let $X=(x_{ij})$ be a $2m \times 2n$ matrix of indeterminates with $m \le n$,
and let $\partial=(\partial/\partial x_{ij})$ be the corresponding matrix of
partial derivatives.  Then
\be
   \pf(\partial J \partial^{\rm T}) \, \pf (X J X^{\rm T})^s
   \;=\;  \left( \prod\limits_{j=0}^{m-1} (s+2j)(s+2n-2m+1+2j) \right)
          \, \pf (X J X^{\rm T})^{s-1}
   \;.
 \label{eq.antisymrectcayley.1}
\ee
More generally, if $I \subseteq [2m]$ with $|I| = 2k$, then
\begin{eqnarray}
   \pf[(\partial J \partial^{\rm T})_{II}] \, \pf (X J X^{\rm T})^s
   & = &
   \left( \prod\limits_{j=0}^{k-1} (s+2j)(s+2n-2m+1+2j) \right)
      \,\times  \nonumber \\
   & & \quad
          \pf (X J X^{\rm T})^{s-1} \,
   \epsilon(I) \pf[(X J X^{\rm T})_{I^c I^c}]
   \;.  \qquad
 \label{eq.antisymrectcayley.2}
\end{eqnarray}
\end{theorem}

\noindent
If $m=n$, \reff{eq.antisymrectcayley.1}
reduces to the ordinary Cayley identity \reff{eq.cayley.1}
of size $2m$ since $\pf(AJ A^{\rm T}) = \det A$.

These three rectangular Cayley identities are roughly analogous
to the ordinary, symmetric and antisymmetric Cayley identities,
respectively, but their proofs are more intricate.

Finally, here is a generalization of
Theorems~\ref{thm.cayley} and \ref{thm.tmrectcayley}
to an arbitrary number~$\ell$ of rectangular matrices
of arbitrary compatible sizes:

\begin{theorem}[multi-matrix rectangular Cayley identity]
 \label{thm.multirectcayley}
Fix integers $\ell \ge 1$ and $n_1,\ldots,n_{\ell} \ge 0$
and write $n_{\ell+1} = n_1$.
For $1 \le \alpha \le \ell$,
let $X^{(\alpha)}$ be an $n_\alpha \times n_{\alpha+1}$
matrix of indeterminates, and let $\partial^{(\alpha)}$
be the corresponding matrix of partial derivatives.
Then
\be
   \det(\partial^{(1)} \cdots \partial^{(\ell)}) \,
      \det(X^{(1)} \cdots X^{(\ell)})^s
   \;=\;  \left( \prod\limits_{\alpha=1}^\ell
                 \prod\limits_{j=0}^{n_1 -1} (s+n_\alpha-n_1 +j) \right)
          \, \det (X^{(1)} \cdots X^{(\ell)})^{s-1}
   \;.
 \label{eq.multirectcayley.1}
\ee
More generally, if $I,J \subseteq [n_1]$ with $|I| = |J| = k$, then
\begin{eqnarray}
   &  &
   \det[(\partial^{(1)} \cdots \partial^{(\ell)})_{IJ}] \,
      \det(X^{(1)} \cdots X^{(\ell)})^s
   \;=\;
   \left( \prod\limits_{\alpha=1}^\ell
          \prod\limits_{j=0}^{k-1} (s+n_\alpha-n_1 +j) \right)
      \times
      \nonumber \\
   & &
   \qquad\qquad
   \det (X^{(1)} \cdots X^{(\ell)})^{s-1} \,
   \epsilon(I,J) \det[(X^{(1)} \cdots X^{(\ell)})_{I^c J^c}]
   \;.
 \label{eq.multirectcayley.2}
\end{eqnarray}
\end{theorem}


We expect that there will exist ``symmetric'' and ``antisymmetric''
variants of Theorem~\ref{thm.multirectcayley},
for matrix products of the form
$X^{(1)} \cdots X^{(\ell)} X^{(\ell){\rm T}} \cdots X^{(1){\rm T}}$
or
\linebreak
$X^{(1)} \cdots X^{(\ell)} X^{(\ell+1)} X^{(\ell){\rm T}} \cdots X^{(1){\rm T}}$
(with $X^{(\ell+1)}$ square and symmetric)
for the symmetric case, and
$X^{(1)} \cdots X^{(\ell)} J X^{(\ell){\rm T}} \cdots X^{(1){\rm T}}$
or
$X^{(1)} \cdots X^{(\ell)} X^{(\ell+1)} X^{(\ell){\rm T}} \cdots X^{(1){\rm T}}$
(with $X^{(\ell+1)}$ square and antisymmetric)
for the antisymmetric case
(and analogous matrices of differential operators).
But all good things must come to an end,
and for lack of time we have chosen not to pursue this direction.

\subsection{Diagonal-parametrized Cayley identities} \label{sec.statement.para}

We would now like to formulate analogues of the ordinary and symmetric
Cayley identities in which the diagonal elements of the matrix $X$
are treated as parameters (i.e., not differentiated with respect to)
and only the off-diagonal elements are treated as variables.

\begin{theorem}[diagonal-parametrized ordinary Cayley identity]
 \label{thm.para.cayley}
Let $X=(x_{ij})$ be an $n \times n$ matrix of indeterminates,
let $\balpha=(\alpha_1,\ldots,\alpha_n)$ and $\bbeta=(\beta_1,\ldots,\beta_n)$
be arbitrary numbers (or symbols),
and let $D_{\balpha,\bbeta,s}$ be the matrix of differential operators
defined by
\be
   (D_{\balpha,\bbeta,s})_{ij}  \;=\;
   \cases{  s
            - \alpha_i
              \sum\limits_{k \neq i} x_{ik} \,
                     {\bigpartial \over \bigpartial x_{ik}}
            - (1-\alpha_i)
              \sum\limits_{l \neq i} x_{li} \,
                     {\bigpartial \over \bigpartial x_{li}}
                    &  if $i = j$  \cr
            \noalign{\vskip 4pt}
            x_{ii}^{\beta_i} x_{jj}^{1-\beta_j}
                     {\bigpartial \over \bigpartial x_{ij}}
                    &  if $i \neq j$  \cr
         }
 \label{def.para.cayley}
\ee
Then
\be
   \det(D_{\balpha,\bbeta,s}) \, (\det X)^s  \;=\;
      s(s+1) \cdots (s+n-1) \,
      \Biggl( \prod\limits_{i=1}^n x_{ii} \Biggr) \,
      (\det X)^{s-1}
   \;.
 \label{eq.para.cayley1}
\ee
More generally, if $I,J \subseteq [n]$ with $|I| = |J| = k$, we have
\begin{eqnarray}
   & &
   \det((D_{\balpha,\bbeta,s})_{IJ}) \, (\det X)^s  \;=\;
       \nonumber \\
   & & \qquad
      s(s+1) \cdots (s+k-1) \, (\det X)^{s-1} \,
      \Biggl( \prod\limits_{i \in I} x_{ii}^{\beta_i} \Biggr) \,
      \Biggl( \prod\limits_{j \in J} x_{jj}^{1-\beta_j} \Biggr) \,
      \epsilon(I,J) \, (\det X_{I^c J^c})
      \;.
       \nonumber \\
 \label{eq.para.cayley2}
\end{eqnarray}
\end{theorem}

Please note that although the elements of the matrix $D_{\balpha,\bbeta,s}$
belong to a {\em non-commutative}\/ algebra of differential operators,
this matrix has the special property that each of its elements
commutes with all the elements not in its own row or column;
therefore, the determinant is well-defined without any special
ordering prescriptions.

\bigskip

{\bf Remarks.}
1.  Expressions involving $x_{ii}^{\beta_i}$ and $x_{jj}^{1-\beta_j}$
can be understood as follows:  we work in the Weyl algebra
generated by $\{x_{ij}\}_{i \neq j}$ and  $\{\partial_{ij}\}_{i \neq j}$,
augmented by $s$ in the usual way (see Section~\ref{subsec.alg})
as well as by the central elements
$y_i = x_{ii}^{\beta_i}$ and $z_i = x_{ii}^{1-\beta_i}$,
it being understood that $x_{ii}$ is a shorthand for $y_i z_i$.

2.  It is easy to see why \reff{eq.para.cayley1} holds with a
right-hand side that is independent of the choice of $\bbeta$.
Indeed, the operators $D_{\balpha,\bbeta,s}$ and $D_{\balpha,\bbeta',s}$
are related by the similarity transformation
\be
D_{\balpha,\bbeta',s}
  \;=\;
\mathrm{diag}(x_{ii}^{\beta'_i - \beta_i})
\,
D_{\balpha,\bbeta,s}
\,
\mathrm{diag}(x_{ii}^{\beta_i - \beta'_i})
\ee
where the quantities $x_{ii}^{\pm (\beta_i - \beta'_i)}$
commute with all entries in all matrices, so that
$\det(D_{\balpha,\bbeta,s}) = \det(D_{\balpha,\bbeta',s})$.
Similar reasoning explains why the right-hand side of \reff{eq.para.cayley2}
depends on $\bbeta$ in the way it does.

3.  It should be stressed that \reff{eq.para.cayley1} provides
a {\em non-minimal}\/ Bernstein--Sato pair.
In fact, a lower-order
Bernstein--Sato pair can be obtained from \reff{eq.para.cayley2}
by taking $I=J = [n] \setminus \{i_0\}$ for any fixed $i_0 \in [n]$:
\be
   \det((D_{\balpha,\bbeta,s})_{\{i_0\}^c \{i_0\}^c}) \, (\det X)^s  \;=\;
      s(s+1) \cdots (s+n-2) \,
      \Biggl( \prod\limits_{i=1}^n x_{ii} \Biggr) \,
      (\det X)^{s-1}
   \;.
 \label{eq.para.cayley3}
\ee


\begin{theorem}[diagonal-parametrized symmetric Cayley identity]
 \label{thm.para.sym.cayley}
\quad\break
Let $X^{\rm sym}=(x_{ij})$ be an $n \times n$ symmetric matrix of
indeterminates,
let $\bbeta=(\beta_1,\ldots,\beta_n)$ be arbitrary numbers (or symbols),
and let $D^{\rm sym}_{\bbeta,s}$ be the matrix of differential operators
defined by
\be
   (D^{\rm sym}_{\bbeta,s})_{ij}  \;=\;
   \cases{  s
            - \smhalf
              \sum\limits_{k > i} x_{ik} \,
                  {\bigpartial \over \bigpartial x_{ik}}
            - \smhalf
              \sum\limits_{l < i} x_{li} \,
                  {\bigpartial \over \bigpartial x_{li}}
                    &  if $i = j$  \cr
            \noalign{\vskip 4pt}
            \smhalf
            x_{ii}^{\beta_i} x_{jj}^{1-\beta_j} \,
                  {\bigpartial \over \bigpartial x_{ij}}
                    &  if $i < j$  \cr
            \noalign{\vskip 4pt}
            \smhalf
            x_{ii}^{\beta_i} x_{jj}^{1-\beta_j} \,
                  {\bigpartial \over \bigpartial x_{ji}}
                    &  if $i > j$  \cr
         }
 \label{def.para.sym.cayley}
\ee
Then
\be
   \det(D^{\rm sym}_{\bbeta,s}) \, (\det X^{\rm sym})^s
   \;=\;
   s(s+\smhalf) \cdots \left(s+ {n-1 \over 2} \right) \,
      \Biggl( \prod\limits_{i=1}^n x_{ii} \Biggr) \,
      (\det X^{\rm sym})^{s-1}
   \;.
 \label{eq.para.sym.cayley1}
\ee
More generally, if $I,J \subseteq [n]$ with $|I| = |J| = k$, we have
\begin{eqnarray}
   & &
   \det((D^{\rm sym}_{\bbeta,s})_{IJ}) \, (\det X^{\rm sym})^s  \;=\;
       \nonumber \\
   & & \quad
      s(s+\smhalf) \cdots \left(s+ {k-1 \over 2} \right) \,
        (\det X^{\rm sym})^{s-1} \,
      \Biggl( \prod\limits_{i \in I} x_{ii}^{\beta_i} \Biggr) \,
      \Biggl( \prod\limits_{j \in J} x_{jj}^{1-\beta_j} \Biggr) \,
      \epsilon(I,J) \, (\det X^{\rm sym}_{I^c J^c})
      \;.
       \nonumber \\
 \label{eq.para.sym.cayley2}
\end{eqnarray}
\end{theorem}

Please note that in the symmetric case we are forced to take
$\alpha_i = \smhalf$ for all $i$.
Note also that \reff{eq.para.sym.cayley1} provides
a non-minimal Bernstein--Sato pair,
and that a lower-order
pair can be obtained from \reff{eq.para.sym.cayley2}
by taking $I=J = [n] \setminus \{i_0\}$ for any fixed $i_0 \in [n]$:
\be
   \det((D^{\rm sym}_{\bbeta,s})_{\{i_0\}^c \{i_0\}^c}) \, (\det X^{\rm sym})^s
   \;=\;
   s(s+\smhalf) \cdots \left(s+ {n-2 \over 2} \right) \,
      \Biggl( \prod\limits_{i=1}^n x_{ii} \Biggr) \,
      (\det X^{\rm sym})^{s-1}
   \;.
 \label{eq.para.sym.cayley3}
\ee

%

\subsection{Laplacian-parametrized Cayley identities} \label{sec.statement.Lap}

In the preceding subsection we treated the off-diagonal elements
$\{x_{ij}\}_{i \neq j}$ as indeterminates
and the diagonal elements $x_{ii}$ as parameters.
Here we again treat the off-diagonal elements as indeterminates,
but now we use the {\em row sums}\/ $t_i = \sum\limits_{j=1}^n x_{ij}$
as the parameters.
(This way of writing a matrix arises in the matrix-tree theorem
 \cite{Chaiken_82,Moon_94,Abdesselam_03}.)
In other words, we define the row-Laplacian matrix
with off-diagonal elements $\{x_{ij}\}_{i \neq j}$
and row sums 0,
\be
   (\XRL)_{ij}
    \;=\;  \cases{ x_{ij}     & if $i \neq j$  \cr
                   \noalign{\vskip 6pt}
                   - \sum\limits_{k \neq i} x_{ik}  & if $i = j$ \cr
                 }
\ee
and the diagonal matrix $T = \diag(t_i)$,
and we then study $\det(T+\XRL)$ where the $t_i$ are treated as parameters.
We shall need the matrix of differential operators
\be
   (\partialRL)_{ij}
    \;=\;  \cases{ \partial/\partial x_{ij}     & if $i \neq j$  \cr
                   \noalign{\vskip 6pt}
                   0                            & if $i = j$ \cr
                 }
\ee

\begin{theorem}[Laplacian-parametrized ordinary Cayley identity]
 \label{thm.Lap.cayley}
 \label{thm.cayley.rowlapl.pre}
\quad\break
Let $\XRL$, $T$ and $\partialRL$ be $n \times n$ matrices defined as above.
Let $U$ be the $n \times n$ matrix with all entries equal to 1.
Then
\begin{eqnarray}
  & &   [\det(U + \partialRL) - \det (\partialRL)]
             \: (\det (T + \XRL))^s
     \nonumber \\[2mm]
  & &
  \qquad =\;
 \Big(\sum_i t_i \Big)\, s(s+1) \cdots (s+n-2) \, (\det (T + \XRL))^{s-1}
   \;. \qquad
 \label{eq.cayley.1.rowlapl.pre.INTRO}
\end{eqnarray}
\end{theorem}


Let us recall that,
by the matrix-tree theorem \cite{Chaiken_82,Moon_94,Abdesselam_03},
$\det(T + \XRL)$ is the generating polynomial of
rooted directed spanning forests on the vertex set $[n]$,
with weight
$\prod_{i \in R} t_i \prod_{\vec{ij} \in E(F)} (-x_{ij})$ 
for a (rooted directed) forest with roots at the vertices $i \in R$ and
edges $\vec{ij} \in E(F)$, directed towards the roots.
In particular, the term linear in $t_i$
--- whose coefficient is the principal minor of order $n-1$,
$\det (\XRL)_{\{i\}^c \{i\}^c}$ ---
enumerates the directed spanning trees rooted at $i$.
Taking this limit in \reff{eq.cayley.1.rowlapl.pre.INTRO}, we obtain:

\begin{corollary}[Cayley identity for the directed-spanning-tree polynomial]
 \label{cor.cayley.rowlapl}
For each $i \in [n]$, we have
\begin{eqnarray}
 & &
   [ \det(U + \partialRL) - \det (\partialRL) ]
      \: (\det (\XRL)_{\{i\}^c \{i\}^c})^s
   \nonumber \\[2mm]
 & &
 \qquad\qquad =\;
   s(s+1) \cdots (s+n-2) \, (\det (\XRL)_{\{i\}^c \{i\}^c})^{s-1}
   \;.
\end{eqnarray}
\end{corollary}

We also have an analogous identity for symmetric Laplacian-parametrized
matrices.  What this means is that we introduce indeterminates
$\{x_{ij}\}_{1 \le i < j \le n}$ and regard $x_{ij}$ for $i>j$
as a synonym for $x_{ji}$;
we then define the symmetric Laplacian matrix
\be
   (\XSL)_{ij}
    \;=\;  \cases{ x_{ij}     & if $i < j$  \cr
                   \noalign{\vskip 6pt}
                   x_{ji}     & if $i > j$  \cr
                   \noalign{\vskip 6pt}
                   - \sum\limits_{k \neq i} x_{ik}  & if $i = j$ \cr
                 }
\ee
and the corresponding matrix of partial derivatives
\be
   (\partialSL)_{ij}
    \;=\;  \cases{ \partial/\partial x_{ij}     & if $i < j$  \cr
                   \noalign{\vskip 6pt}
                   \partial/\partial x_{ji}     & if $i > j$  \cr
                   \noalign{\vskip 6pt}
                   0                            & if $i = j$ \cr
                 }
\ee
We then have:

\begin{theorem}[Laplacian-parametrized symmetric Cayley identity]
 \label{thm.Lap.sym.cayley}
 \label{thm.cayley.lapl.pre}
\quad\break
Let $\XSL$, $T$ and $\partialSL$ be $n \times n$ matrices defined as above.
Let $U$ be the $n \times n$ matrix with all entries equal to 1.
Then
\begin{eqnarray}
 & &  [\det(U + \partialSL) - \det (\partialSL)]
             \: (\det (T + \XSL))^s
     \nonumber \\[2mm]
 & &
 \qquad =\;
 \Big(\sum_i t_i \Big)\, 2s(2s+1) \cdots (2s+n-2) \, (\det (T + \XSL))^{s-1}
   \;. \qquad
 \label{eq.thm.Lap.sym.cayley}
\end{eqnarray}
\end{theorem}

Similarly, $\det(T + \XSL)$ is the generating polynomial of
rooted (undirected) spanning forests on the vertex set $[n]$,
with weight
$\prod_{i \in R} t_i \prod_{ij \in E(F)} (-x_{ij})$ 
for a (rooted undirected) forest with roots at the vertices $i \in R$ and
edges $ij \in E(F)$.
In particular, the term linear in $t_i$
--- whose coefficient is $\det (\XSL)_{\{i\}^c \{i\}^c}$ ---
is independent of $i$ and enumerates the spanning trees.
We thus obtain:

\begin{corollary}[Cayley identity for the spanning-tree polynomial]
 \label{cor.cayley.symlapl}
For each $i \in [n]$, we have
\begin{eqnarray}
  & &
   [ \det(U + \partialSL) - \det (\partialSL) ]
      \: (\det (\XSL)_{\{i\}^c \{i\}^c})^s
  \nonumber \\[2mm]
  & &
  \qquad\quad =\;
   2s(2s+1) \cdots (2s+n-2) \, (\det (\XSL)_{\{i\}^c \{i\}^c})^{s-1}
   \;.  \qquad
\end{eqnarray}
\end{corollary}

\noindent
See \cite{Sokal_bcc2005} for many interesting additional properties
of the spanning-tree polynomial.

\subsection{Product-parametrized and border-parametrized
   rectangular Cayley identities}

Here we will present two curious Cayley identities for rectangular matrices
that are much simpler than the identities presented in
Section~\ref{sec.statement.rectangular}, because the indeterminates
occur linearly rather than quadratically in the argument of the determinant.

So let $X=(x_{ij})$ be an $m \times n$ matrix of indeterminates with $m \le n$,
and let $\partial=(\partial/\partial x_{ij})$ be the corresponding matrix of
partial derivatives.  One easy way to obtain square ($m \times m$)
matrices from $X$ and $\partial$ is to right-multiply them
by $n \times m$ matrices $A$ and $B$, respectively.
We then have the following Cayley-type identity,
in which $A$ and $B$ occur only as {\em parameters}\/:

\begin{theorem}[product-parametrized rectangular Cayley identity]
  \label{thm.cayley.XA2m}
  \label{thm.productcayley}
\quad\break
Let $X=(x_{ij})$ be an $m \times n$ matrix of indeterminates with 
$m \leq n$, and let
$\partial=(\partial/\partial x_{ij})$ be the corresponding matrix of
partial derivatives.  Let $A=(a_{ij})$ and $B=(b_{ij})$ be
$n \times m$ matrices of constants.
Then
\be
   \det(\partial B) 
\, (\det X A)^s  \;=\;  
\det(A^{\rm T} B) \;
s(s+1) \cdots (s+m-1) \, (\det X A)^{s-1}
   \;.
 \label{eq.cayley.XA2m}
\ee
More generally, if $I,J \subseteq [m]$ with $|I| = |J| = k$, then
\be
   \det((\partial B)_{IJ})
\, (\det X A)^s  \;=\; 
\det(MA) \;
s(s+1) \cdots (s+k-1) \, (\det X A)^{s-1}
   \;.
 \label{eq.cayley.XA2m.2}
\ee
where $M$ is an $m \times n$ matrix defined as follows:
if $I=\{i_1,\ldots,i_k\}$ and $J=\{j_1,\ldots,j_k\}$ in increasing order, then
\be
M_{\alpha\beta}
\;=\;
  \cases{  x_{\alpha \beta}  & if $\alpha \notin I$  \cr
          \noalign{\vskip 2mm}
           b_{\beta j_h}     & if $\alpha = i_h$  \cr
        }
  \label{eq.defMinXa2m}
\ee
\end{theorem}

Note that if $n=m$, then \reff{eq.cayley.XA2m} reduces to
the ordinary Cayley identity \reff{eq.cayley.1} multiplied on both sides
by $(\det B) (\det A)^s$,
while \reff{eq.cayley.XA2m.2} specialized to $A=B=I_m$
gives the all-minors identity \reff{eq.cayley.2}.
In the general case $m \le n$,
\reff{eq.cayley.XA2m.2} reduces to \reff{eq.cayley.XA2m} when $I=J=[m]$
(since we then have $M=B^{\rm T}$),
while \reff{eq.cayley.XA2m.2} reduces to the trivial identity
$(\det X A)^s  = (\det X A)^s$ when $I = J = \emptyset$
(since we then have $M=X$).

\bigskip

A second easy way to complete $X$ and $\partial$ to
square ($n \times n$) matrices is to adjoin $n-m$ rows of constants
at the bottom:
\be
   \widehat{X}  \;=\; \left( \! \begin{array}{c}
                                    X \\
                                    \hline
                                    A
                                \end{array} \! \right)
   ,\qquad
   \widehat{\partial}  \;=\; \left( \! \begin{array}{c}
                                    \partial \\
                                    \hline
                                    B
                                \end{array} \! \right)
 \label{def.Xhat}
\ee
where $A$ and $B$ are $(n-m) \times n$ matrices of constants.
We can think of $\widehat{X}$ and $\widehat{\partial}$
as ``bordered'' matrices obtained by filling out $X$ and $\partial$.
We then have the following Cayley-type identity,
in which $A$ and $B$ again occur only as parameters:

\begin{theorem}[border-parametrized rectangular Cayley identity]
   \label{thm.borderedcayley}
Let $X=(x_{ij})$ be an $m \times n$ matrix of indeterminates with $m \le n$,
let $\partial=(\partial/\partial x_{ij})$ be the corresponding matrix of
partial derivatives,
and let $A$ and $B$ be $(n-m) \times n$ matrices of constants.
Define $\widehat{X}$ and $\widehat{\partial}$ as in \reff{def.Xhat}.
Then
\be
   \det(\widehat{\partial}) \, (\det \widehat{X})^s
   \;=\;
   \det(A B^{\rm T}) \,
     s(s+1) \cdots (s+m-1) \, (\det \widehat{X})^{s-1}
   \;.
 \label{eq.borderedcayley.1}
\ee
\end{theorem}

\noindent
When $m=n$ this formula reduces to the ordinary Cayley identity
\reff{eq.cayley.1}, if we make the convention that the determinant
of an empty matrix is 1.

Please note that, by Laplace expansion, $\det \widehat{X}$ is
a linear combination of $m \times m$ minors of $X$
(and likewise for $\det \widehat{\partial}$).
Indeed, when $n-m=1$ one can obtain {\em all}\/ such linear
combinations in this way, by suitable choice of the row vector $A$;
for $n-m \ge 2$ one obtains in general a subset of such linear combinations.
It is striking that the form of the identity ---
and in particular the polynomial $b(s)$ occurring in it ---
does not depend on the choice of the matrices $A$ and $B$.

In Section~\ref{sec.border} we will prove the
``product-parametrized'' and ``border-parametrized'' identities
and then explain the close relationship between them.

\subsection{Historical remarks}   \label{sec.historical}

As noted in the Introduction, the identity \reff{eq.intro.1}
is conventionally attributed to Arthur Cayley (1821--1895);
the generalization \reff{eq.intro.2} to arbitrary minors
is sometimes attributed to Alfredo Capelli (1855--1910).
The trouble is, neither \reff{eq.intro.1} nor \reff{eq.intro.2}
occurs anywhere --- as far as we can tell ---
in the {\em Collected Papers}\/ of Cayley \cite{Cayley_collected}.
Nor are we able to find these formulae in any of the relevant works of Capelli
\cite{Capelli_1882,Capelli_1887,Capelli_1888,Capelli_1890,Capelli_02}.
The operator $\Omega = \det(\partial)$ was indeed introduced
by Cayley on the second page of his famous 1846 paper on invariants
\cite{Cayley_1846};  it became known as Cayley's $\Omega$-process
and went on to play an important role in classical invariant theory
(see e.g.\ \cite{Weyl_46,Schur_68,Fulton_91,Olver_99,Dolgachev_03}).
But we strongly doubt that Cayley ever knew \reff{eq.intro.1}.

A detailed history of \reff{eq.intro.1} and \reff{eq.intro.2}
will be presented elsewhere \cite{cayley_history}.
Suffice it to say that the special case for $2 \times 2$ matrices
appears already in the 1872 book of Alfred Clebsch (1833--1872)
on the invariant theory of binary forms \cite[p.~20]{Clebsch_1872}.
But even for $n=3$, the first unambiguous statement of which we are aware
appears in an 1890 paper of Giulio Vivanti (1859--1949)
\cite{Vivanti_1890}.\footnote{
   This paper is also cited in Muir's massive annotated bibliography
   of work on the theory of determinants \cite[vol.~4, p.~479]{Muir}.
   Indeed, it was thanks to Muir that we discovered Vivanti's paper.
   Malek Abdesselam has also drawn our attention to the papers
   of Clebsch (1861) \cite[pp.~7--14]{Clebsch_1861}
   and Gordan (1872) \cite[pp.~107--116]{Gordan_1872},
   where formulae closely related to \reff{eq.intro.1} can be found.
%
}
And amazingly, in this very first paper,
Vivanti proves not only the basic ``Cayley'' identity \reff{eq.intro.1}
but also the generalization \reff{eq.intro.2} for minors,
for completely general $n$ and $k$
albeit only in the case $I=J$ (i.e., principal minors).
In fact, his inductive method of proof works {\em only}\/
because he is handling the ``all-principal-minors'' version;
it would not work for the ``simple'' identity \reff{eq.intro.1} alone.

Proofs (by direct computation) of the ``Cayley'' identity \reff{eq.intro.1}
for $n=3$ can be found in the early-twentieth-century books
of Grace and Young (1903) \cite{Grace_03},
Weitzenbock (1923) \cite{Weitzenbock_23}
and Turnbull (1928) \cite{Turnbull_28}.
Weitzenbock also states \reff{eq.intro.1} without proof for general $n$,
saying that it is obtained ``by completely analogous calculation''
\cite[p.~16]{Weitzenbock_23};
similarly, Turnbull states both \reff{eq.intro.1} and \reff{eq.intro.2}
for general $n$ and leaves them as exercises for the reader
\cite[pp.~114--116]{Turnbull_28}.
(Unfortunately, Turnbull's old-fashioned notation
is very difficult to follow.)
We are not convinced that the extension from $n=3$ to general $n$
is quite so trivial as these authors imply.
We will, in any case, provide an elementary algebraic/combinatorial proof
of \reff{eq.intro.1}/\reff{eq.intro.2}
in Section~\ref{sec.elementary.1} below.


The symmetric analogues \reff{eq.symcayley.1}/\reff{eq.symcayley.2}
are due to G{\aa}rding in 1948 \cite{Garding_48};
see also \cite[Proposition~VII.1.4]{Faraut_94}.
The antisymmetric analogue \reff{eq.antisymcayley.1}
is due to Shimura \cite{Shimura_84} in 1984;
see again \cite[Proposition~VII.1.4]{Faraut_94},
where the quaternionic hermitian determinant
(Moore determinant \cite{Aslaksen_96})
is equivalent to a pfaffian,
and see also \cite[Corollary~3.13]{Kinoshita_02}.


Shimura \cite{Shimura_84} also gives generalizations of all these formulae
in which $\det(\partial)$ is replaced by
other homogeneous differential operators.
Similarly,
Rubenthaler and Schiffmann \cite[especially Section~5]{Rubenthaler_87}
and Faraut and Kor\'anyi \cite[Proposition~VII.1.6]{Faraut_94}
give generalizations in which both $\det(\partial)$ and $\det X$
are replaced by suitable products of leading principal minors.
All these proofs are analytic in nature.
Elegant algebraic/combinatorial proofs of identities in which $\det X$
is replaced by a product of minors
have been given by Canfield, Williamson and Evans
\cite[Theorem~1.3 and Lemma~2.12]{Canfield_81}
\cite[Theorem~4.1]{Williamson_87}.\footnote{
   See also Turnbull \cite{Turnbull_49} for a similar result,
   but expressed in difficult-to-follow notation.
}

The one-matrix rectangular symmetric Cayley identity
\reff{eq.rectcayley.1}
has been proven recently by
Faraut and Kor\'anyi \cite[section~XVI.4]{Faraut_94},
Kh\'ekalo \cite{Khekalo_01,Khekalo_05} and Rubin \cite{Rubin_06},
using analytic methods.

Many of the polynomials $P$ treated here can also be understood
as relative invariants of prehomogeneous vector spaces
\cite{Igusa_00,Kimura_03};
and the corresponding $b$-functions have been computed in that context,
mostly by means of microlocal calculus
\cite{Sato_80,Kimura_82,Sugiyama_11}.
See \cite[Appendix]{Kimura_03} for a table of all irreducible reduced
prehomogeneous vector spaces and their relative invariants and $b$-functions.
All but one of the ``generic'' cases in this table
(i.e., those involving matrices with one or more arbitrary dimensions)
correspond to identities treated here:
namely, the ordinary, symmetric and antisymmetric Cayley identities
[\reff{eq.cayley.1}, \reff{eq.symcayley.1} and \reff{eq.antisymcayley.1}]
correspond to cases (1), (2) and (3), respectively,
while the one-matrix rectangular symmetric and antisymmetric Cayley identities
[\reff{eq.rectcayley.1} and \reff{eq.antisymrectcayley.1}]
correspond to cases (15) and (13).
The $b$-functions listed in that table of course agree with ours.\footnote{
   After making the translation of conventions $s \to s+1$:
   see footnote~\ref{footnote.shifted} above.
   See also F.~Sato and Sugiyama \cite[Section~3.1 and Lemma~4.2]{Sato_06}
   for an alternate approach to \reff{eq.cayley.1}, \reff{eq.antisymcayley.1}
   and \reff{eq.rectcayley.1}.
}
(The final ``generic'' case (30)
involves matrices $X$ of order $3 \times 2n$
with $P(X) = \tr (XJX^{\rm T})^2$,
which falls outside our methods
since it is neither a determinant nor a pfaffian.)
Finally, the most difficult result obtained in the present paper
--- namely, the multi-matrix rectangular Cayley identity identity
\reff{eq.multirectcayley.1} ---
has very recently been proven independently
by Sugiyama \cite[Theorem~0.1]{Sugiyama_11}
in the context of prehomogeneous vector spaces
associated to equioriented quivers of type~${\sf A}$;
his proof uses a decomposition formula found earlier
by himself and F.~Sato \cite{Sato_06}.
Indeed, Sugiyama proved an even more general result
\cite[Theorem~3.4]{Sugiyama_11},
applying to quivers of type~${\sf A}$ with arbitrary orientation.\footnote{
   We are grateful to Nero Budur for explaining to us
   the connection between our results and
   the theory of prehomogeneous vector spaces,
   and for drawing our attention to the work of Sugiyama \cite{Sugiyama_11}.
}

It is worth stressing that the Cayley identity \reff{eq.intro.1}
--- though not, as far as we can tell, the all-minors version \reff{eq.intro.2} ---
is an immediate consequence of a deeper identity due to
Capelli \cite{Capelli_1887,Capelli_1888,Capelli_1890},
in which the operator $H = (\det X)(\det \partial)$
is represented as a noncommutative determinant
involving the $\mathfrak{gl}(n)$ generators $X^{\rm T} \partial$:
see e.g.\ \cite[p.~53]{Umeda_98}, \cite[pp.~569--570]{Howe_91}
or \cite[Appendix]{CSS_Capelli_08}
for the easy deduction of Cayley from Capelli.
Likewise, the symmetric Cayley identity \reff{eq.symcayley.1}
follows from a symmetric Capelli-type identity due to
Turnbull \cite{Turnbull_48} (see also \cite{Wallace_53}),
and the antisymmetric Cayley identity
\reff{eq.antisymcayley.1}
follows from an antisymmetric Capelli-type identity due
independently to Kostant and Sahi \cite{Kostant_91}
and to Howe and Umeda \cite{Howe_91}
(see also \cite{Kinoshita_02}).
  
The one-matrix rectangular symmetric Cayley identity \reff{eq.rectcayley.1}
for $m=1$ follows from a Capelli-type identity
given in \cite[pp.~291--293]{Weyl_46}
and \cite[p.~61]{Umeda_98}.
For $m=n$ it of course follows from the ordinary Capelli identity.
For $2 \le m \le n-1$ we do not know any Capelli-type identity.

Proofs of the Capelli-type identities based on
group-representation theory have been given by
Howe and Umeda \cite{Howe_89,Howe_91}.
Combinatorial proofs of the Capelli and Turnbull identities
have been given by Foata and Zeilberger \cite{Foata_94}.
We have recently given very simple algebraic proofs
of these same identities as well as some generalizations \cite{CSS_Capelli_08}.
See also Weyl \cite[pp.~39--42]{Weyl_46} and
Fulton--Harris \cite[Appendix~F.3]{Fulton_91}
for more traditional proofs.
Further information on Capelli-type identities can be found in
\cite{Howe_91,Umeda_98}.

\subsection{Some algebraic preliminaries}  \label{subsec.alg}

%

A few words are needed about how the identities \reff{eq.cayley.1} ff.\ 
--- or more generally, Bernstein-type identities of the form
\be
   Q(s,x,\partial/\partial x) \, P(x)^s
   \;=\;
   b(s) \, P(x)^{s-1}
 \label{eq.sec2.bernstein}
\ee
where $x = (x_1,\ldots,x_n)$ --- are to be interpreted.
On the one hand, they can be interpreted as analytic identities for
functions of real or complex variables $x_1,\ldots,x_n$,
where $s$ is a real or complex number;
here $P(x)^s$ denotes any fixed branch on any open subset
of $\R^n$ or $\C^n$ where it is well-defined.
Alternatively, these formulae can be regarded as purely algebraic identities,
in several different ways:
\begin{itemize}
  \item[1)]
     For integer $s \ge 1$, as an identity in the ring $R[x_1,\ldots,x_n]$
     of polynomials in the indeterminates $x_1,\ldots,x_n$
     with coefficients in some commutative ring $R$
     (for instance, $R$ could be $\Z$, $\Q$, $\R$ or $\C$).\footnote{
        When $b(s)$ contains fractions
        [e.g.\ \reff{eq.symcayley.1}/\reff{eq.symcayley.2} and
         \reff{eq.para.sym.cayley1}/\reff{eq.para.sym.cayley2}]
        we should assume that the coefficient ring $R$ contains those fractions
        [i.e.\ in this case $\smhalf$].
}
  \item[2)]
     For any integer $s$ (positive or negative),
     as an identity in the field $K(x_1,\ldots,x_n)$
     of rational fractions in the indeterminates $x_1,\ldots,x_n$
     with coefficients in some field $K$
     (for instance, $K$ could be $\Q$, $\R$ or $\C$).\footnote{
        When $b(s)$ contains fractions
        [e.g.\ \reff{eq.symcayley.1}/\reff{eq.symcayley.2} and
         \reff{eq.para.sym.cayley1}/\reff{eq.para.sym.cayley2}]
        we should assume that the coefficient field $K$ contains those fractions
        [i.e.\ in this case $\smhalf$].
        Usually we will take $K$ to be a field of characteristic 0,
        so that $K$ contains the rationals $\Q$ as a subfield.
}
   \item[3)]
     For $s$ interpreted symbolically, as an identity in a module
     defined as follows
     \cite[pp.~93--94]{Coutinho_95} \cite[pp.~96~ff.]{Krause_00}:
     Let $K$ be a field of characteristic 0,
     let $x_1,\ldots,x_n$ and $s$ be indeterminates,
     and let $A_n(K)[s]$ be the $K$-algebra generated by
     $x_1,\ldots,x_n, \partial/\partial x_1, \ldots, \partial/\partial x_n$
     and $s$ with the usual commutation relations.
     (That is, it is the algebra of differential operators with respect to
     $x_1,\ldots,x_n$ in which the coefficients are polynomials in
     $x_1,\ldots,x_n,s$ with coefficients in $K$.)
     Now fix a nonzero polynomial $P \in K[x_1,\ldots,x_n]$,
     and let $K[x,s,P^{-1}]$ denote the ring of rational fractions
     in the indeterminates $x_1,\ldots,x_n,s$ whose denominators
     are powers of $P$.
     (It is a subring of the field $K(x_1,\ldots,x_n,s)$
     of all rational fractions in $x_1,\ldots,x_n,s$.)
     Let $K[x,s,P^{-1}] P^s$ be the free $K[x,s,P^{-1}]$-module
     consisting of objects of the form $f P^s$
     where $f \in K[x,s,P^{-1}]$;  here $P^s$ is treated as a formal symbol.
     We can define formal differentiation by
     \be
        {\partial \over \partial x_i} \, (f P^s)
        \;=\;
        \left( {\partial f \over \partial x_i} \,+\,
               s f \, {\partial P \over \partial x_i} \, P^{-1} \right) P^s
     \ee
     where $\partial f / \partial x_i$ is the standard formal derivative
     of a rational fraction.  This differentiation is easily extended
     to an action of $A_n(K)[s]$ on $K[x,s,P^{-1}] P^s$,
     making the latter into a left $A_n(K)[s]$-module.\footnote{
  More generally, we can proceed as follows:
   Let $R$ be an integral domain and let $A$ be an abelian group.
   We then define $R^A$ to be the commutative ring with identity
   generated by the symbols $x^a$ ($x \in R, a \in A$) subject to the relations
   $x^a x^b = x^{a+b}$,  $x^a y^a = (xy)^a$  and  $x^0 = 1$.
   In particular, if $A$ contains the integers as a subgroup,
   then we can consider $R$ as a subring of $R^A$
   by identifying $x \in R$ with $x^1 \in R^A$.
   
   Now suppose that $R$ is a polynomial ring $S[x_1,\ldots,x_n]$
   where $S$ is an integral domain of characteristic 0,
   and that $A$ is a subgroup of the additive group of $S$
   (where $1 \in \Z \subseteq A$ is identified with $1 \in S$).
   Then we can define an action of the differential operators
   $\partial/\partial x_i$ on $R^A$ by
   $$ {\partial \over \partial x_i} \, (P^a)  \;=\;
      a \, {\partial P \over \partial x_i} \, P^{a-1}
   $$
   together with the usual product rule.
   This makes $R^A$ into a left $A_n(S)$-module
   [where $A_n(S)$ is the Weyl algebra in $n$ variables over $S$].

   In order to handle Bernstein-type identities,
   we will introduce an indeterminate $s$
   and take $A = \Z + s\Z$ and $S = K[s]$ for some field $K$.
   Then we will work within the submodule of $R^A$ consisting of
   elements of the form $f P^{s+a}$ for $f \in K[x_1,\ldots,x_n,s]$,
   $a \in \Z$ and some \emph{fixed} $P \in K[x_1,\ldots,x_n]$.
   This submodule is isomorphic to $K[x,s,P^{-1}] P^s$.
}
\end{itemize}
Let us now show that all these interpretations are equivalent.

We begin by recalling some elementary facts.
Let $p(x_1,\ldots,x_n)$ be a polynomial with coefficients
in some commutative ring $R$,
and let $d_i$ be the degree of $p$ with respect to the variable $x_i$.
Suppose that there exist sets $X_1,\ldots,X_n \subseteq R$
with $|X_i| > d_i$ for all~$i$,
such that $p(x_1,\ldots,x_n) = 0$ whenever
$(x_1,\ldots,x_n) \in X_1 \times \ldots \times X_n$.
Then $p$ must be the zero polynomial,
i.e.\ all its coefficients are zero.
Note in particular that if the sets $X_1,\ldots,X_n$ are infinite,
then this reasoning applies to polynomials of arbitrary degree.
As a special case of this, if $R = \R$ or $\C$
and $p(x_1,\ldots,x_n) = 0$ for all $x = (x_1,\ldots,x_n)$
lying in some nonempty open set $U$ of $\R^n$ or $\C^n$,
then $p$ must be the zero polynomial.

Now let $K$ be a field of characteristic 0,
and let $P(x)$, $Q(s,x,\partial/\partial x)$ and $b(s)$
be polynomials with coefficients in $K$
[as always we use the shorthand $x = (x_1,\ldots,x_n)$].
Then elementary algebraic manipulations allow us to write
\be
   Q(s,x,\partial/\partial x) \, P(x)^s
   \,-\,
   b(s) \, P(x)^{s-1}
   \;=\;
   R(s,x) \, P(x)^{s-m}
\ee
for some polynomial $R(s,x)$ and some integer $m \ge 0$.
Combining this fact with the preceding observations,
we obtain immediately the following two propositions:

\begin{proposition}[Equivalence theorem for symbolic $s$]
   \label{prop.equiv.symbolic_s}
Let $K$ be a field of characteristic 0,
and let $P$, $Q$, $b$ and $R$ be as before.
Then the following are equivalent:
\begin{itemize}
   \item[(a)]  \reff{eq.sec2.bernstein} holds as an algebraic identity
        in the module $K[x,s,P^{-1}] P^s$, symbolically in $s$.
   \item[(b)]  $R(s,x) = 0$ in the polynomial ring $K[s,x]$.
\end{itemize}
Furthermore, if $K$ is infinite, then (a)--(b) are equivalent to:
\begin{itemize}
   \item[(c)]  For infinitely many $s \in K$,
        there exist infinite sets $X_1,\ldots,X_n \subseteq K$
        (possibly depending on $s$),
        such that $R(s,x_1,\ldots,x_n) = 0$
        for $(x_1,\ldots,x_n) \in X_1 \times \ldots \times X_n$
        and the given value of $s$.
\end{itemize}
Finally, if $K = \R$ or $\C$, then (a)--(c) are equivalent to:
\begin{itemize}
   \item[(d)]  For infinitely many $s \in K$,
        there exist a nonempty open set $U \subseteq K^n$
        (possibly depending on $s$)
        and a branch of $P(x)^s$ defined on $U$ such that
        \reff{eq.sec2.bernstein} holds for all $x \in U$.
   \item[(e)]  For every $s \in K$,
        every nonempty open set $U \subseteq K^n$
        and every branch of $P(x)^s$ defined on $U$,
        \reff{eq.sec2.bernstein} holds for all $x \in U$.
\end{itemize}
\end{proposition}

\begin{proposition}[Equivalence theorem for fixed $s$]
   \label{prop.equiv.fixed_s}
Let $K$ be a field of characteristic 0,
let $P$, $Q$, $b$ and $R$ be as before, and fix some element $s \in K$.
Then the following are equivalent:
\begin{itemize}
   \item[(b)]  $R(s,x) = 0$ in the polynomial ring $K[x]$
        (for the given value of $s$).
\end{itemize}
Furthermore, if $K$ is infinite, then (a)--(b) are equivalent to:
\begin{itemize}
   \item[(c)]  There exist infinite sets $X_1,\ldots,X_n \subseteq K$
        such that $R(s,x_1,\ldots,x_n) = 0$
        for $(x_1,\ldots,x_n) \in X_1 \times \ldots \times X_n$
        (for the given value of $s$).
\end{itemize}
Finally, if $K = \R$ or $\C$, then (a)--(c) are equivalent to:
\begin{itemize}
   \item[(d)]  There exist a nonempty open set $U \subseteq K^n$
        and a branch of $P(x)^s$ defined on $U$ such that
        \reff{eq.sec2.bernstein} holds for all $x \in U$
        (for the given value of $s$).
   \item[(e)]  For every nonempty open set $U \subseteq K^n$
        and every branch of $P(x)^s$ defined on $U$,
        \reff{eq.sec2.bernstein} holds for all $x \in U$
        (for the given value of $s$).
\end{itemize}
\end{proposition}


In particular, Proposition~\ref{prop.equiv.symbolic_s}
shows that it suffices to prove \reff{eq.sec2.bernstein}
for infinitely many positive or negative integers $s$,
using the elementary interpretation (1) or (2) above;
it then holds automatically for arbitrary $s$
as an identity in the module $K[x,s,P^{-1}] P^s$
and as an analytic identity.
Likewise, it suffices to specialize $x_1,\ldots,x_n$
to real or complex variables and to prove \reff{eq.sec2.bernstein}
for some nonempty open set in $\R^n$ or $\C^n$.
In what follows, we shall repeatedly take advantage of
these simplifications.
(Many previous authors --- especially the earlier ones ---
have done so as well, but without making
Proposition~\ref{prop.equiv.symbolic_s} explicit.)

\section{Elementary proofs of Cayley-type identities}  \label{sec.elementary}

In this section we give proofs of the three main
Cayley-type identities for square matrices
(Theorems~\ref{thm.cayley}--\ref{thm.antisymcayley})
that use nothing but elementary properties of determinants
(notably, Jacobi's identity for cofactors)
along with the elementary formulae for the derivative of a product or a power.

The general situation we have to handle
in all three cases is as follows:
Let $\Gamma$ be a finite index set,
let $(E_\gamma)_{\gamma \in \Gamma}$ be given $n \times n$ matrices
with elements in some field $K$,
and let $(x_\gamma)_{\gamma \in \Gamma}$ be indeterminates.
Now define the matrix $X = \sum_{\gamma \in \Gamma} x_\gamma E_\gamma$.
If $A = (A_\gamma)_{\gamma \in \Gamma}$ is a $K$-valued vector,
we write $\partial_A = \sum_{\gamma \in \Gamma} A_\gamma \,
                                                \partial/\partial x_\gamma$
and $E_A = \sum_{\gamma \in \Gamma} A_\gamma E_\gamma$,
so that $\partial_A X = E_A$.
We need a formula for successive derivatives of $(\det X)^s$:

\begin{lemma}
  \label{lemma.elem}
Let $(E_\gamma)_{\gamma \in \Gamma}$ be $n \times n$ matrices,
let $(x_\gamma)_{\gamma \in \Gamma}$ be indeterminates,
and define $X = \sum_{\gamma \in \Gamma} x_\gamma E_\gamma$.
For any sequence $A_1,\ldots,A_k$ of $K$-valued vectors, we have
\begin{subeqnarray}
\!\!\!\! & & \!\!\!\!\!
  \left( \prod\limits_{i=1}^k \partial_{A_i} \right)
  \,
  (\det X)^{s}
  \nonumber \\
& &
=\;
 (-1)^{k} \left( \det X \right)^{s}
\sum_{\tau \in \scrs_k} (-s)^{\#(\hbox{\scriptsize\rm cycles of } \tau)}
\!\!\!\!\!\!\!
 \prod_{\begin{scarray}
\hbox{\scriptsize\rm cycles of } \tau \\
C=(\alpha_{1},\ldots\alpha_{\ell})\end{scarray}}
\!\!\!\!\!\!\!
\tr \left( X^{-1}E_{A_{\alpha_{1}}}\ldots X^{-1}E_{A_{\alpha_{\ell}}} 
\right)
    \qquad \\[1mm]
& &
=\;
 \left( \det X \right)^{s}
\sum_{\tau \in \scrs_k} \sgn(\tau) \,
   s^{\#(\hbox{\scriptsize\rm cycles of } \tau)}
\!\!\!\!\!\!\!
 \prod_{\begin{scarray}
\hbox{\scriptsize\rm cycles of } \tau \\
C=(\alpha_{1},\ldots\alpha_{\ell})\end{scarray}}
\!\!\!\!\!\!\!
\tr \left( X^{-1}E_{A_{\alpha_{1}}}\ldots X^{-1}E_{A_{\alpha_{\ell}}} 
\right) \;.
\end{subeqnarray}
\end{lemma}

\proof
By induction on $k$.
The case $k=1$ follows from Cramer's rule or alternatively from the relation
\be
(\det X)^{s}   \;=\;  \exp  ( s \,\tr \log X )  \;.
\ee
For the inductive step we shall need the identity
\be
\partial_A (X^{-1})   \;=\;   - X^{-1} (\partial_A X) X^{-1}  \;,
  \label{eq.deriv.X-1}
\ee
which follows from $\partial_A (X^{-1} X) = 0$.
Assume now that the theorem is valid for $k$ and let us apply
$\partial_{A_{k+1}}$.
When this derivative hits $(\det X)^s$,
it creates a new cycle $(k+1)$,
with prefactor $s = (-1)(-s)$;
these terms correspond to permutations $\tau \in \scrs_{k+1}$
in which the element $k+1$ is fixed.
Alternatively, the derivative can hit one of the $X^{-1}$ factors
in one of the traces;
by \reff{eq.deriv.X-1} this inserts the element $k+1$
into one of the existing cycles at an arbitrary position,
and produces an extra factor $-1$;
these terms correspond to permutations $\tau \in \scrs_{k+1}$
in which the element $k+1$ is not fixed.
\qed

For $1 \le i,j \le n$, let $E^{ij}$ be the matrix with a 1 in position $ij$
and zeros elsewhere, i.e.
\be
   (E^{ij})_{i'j'}  \;=\;  \delta_{i,i'} \delta_{j,j'}  \;.
 \label{def.Eij}
\ee
We will express the matrices $E_\gamma$ in each of our three cases
in terms of the $E^{ij}$.

\subsection{Ordinary Cayley identity}  \label{sec.elementary.1}

\proofof{Theorem~\ref{thm.cayley}}
In this case the index set $\Gamma$ is simply $[n] \times [n]$,
and we write $X = \sum_{i,j=1}^n x_{ij} E^{ij}$.
Now let $I = \{i_1,\ldots,i_k\}$ with $i_1 < \ldots < i_k$
and $J = \{j_1,\ldots,j_k\}$ with $j_1 < \ldots < j_k$,
so that
\be
   \det(\partial_{IJ})  \;=\;
   \sum_{\sigma \in \scrs_k}  \sgn(\sigma) \prod_{r=1}^k
         {\bigpartial \over \bigpartial x_{i_r j_{\sigma(r)}}}
   \;.
\ee
For each fixed $\sigma \in \scrs_k$, we apply Lemma~\ref{lemma.elem}
with $\partial_{A_r} = \partial/\partial x_{i_r j_{\sigma(r)}}$.
In the traces we have
$E_{A_{\alpha_r}} = E^{i_{\alpha_r} j_{\sigma(\alpha_r)}}$
and hence
\begin{subeqnarray}
\tr \left( X^{-1} E^{i_{\alpha_{1}}j_{\sigma(\alpha_{1})}} \ldots 
X^{-1} E^{i_{\alpha_{\ell}}j_{\sigma(\alpha_{\ell})}} \right)  &= &
X^{-1}_{j_{\sigma(\alpha_{\ell})} i_{\alpha_{1}}} 
X^{-1}_{j_{\sigma(\alpha_{1})} i_{\alpha_{2}}} \ldots 
X^{-1}_{j_{\sigma(\alpha_{\ell-1})} i_{\alpha_{\ell}}} \\
& = & \prod_{r= 1}^{\ell} X^{-\rm T}_{i_{\tau(\alpha_r)} j_{\sigma(\alpha_r)}}
   \;,
\end{subeqnarray}
where $X^{-\rm T} \equiv (X^{-1})^{\rm T}$
and we have used the fact that, for $\tau$ as in Lemma~\ref{lemma.elem},
$\tau(\alpha_{i}) = \alpha_{i+1}$ for $i=1,\ldots,\ell-1$
and $\tau(\alpha_{\ell})=\alpha_{1}$.
We can now combine all the different traces into a single product.
We obtain
\be
   \det(\partial_{IJ}) \, (\det X)^{s}
   \;=\;
   (\det X)^{s}  \sum_{\sigma \in \scrs_k} \sgn(\sigma)
                 \sum_{\tau \in \scrs_k} \sgn(\tau)  \,
   s^{\#(\hbox{\scriptsize\rm cycles of } \tau)}
   \prod_{r=1}^k  (X^{-\rm T})_{i_{\tau(r)} j_{\sigma(r)}}
   \;.
\ee
Let us now define the permutation $\pi= \sigma \circ \tau^{-1}$
and change variables from $(\sigma,\tau)$ to $(\tau,\pi)$,
using $\sgn(\sigma) \, \sgn(\tau) = \sgn(\pi)$.
The product over $r$ can be written equivalently as a product over
$t = \tau(r)$.
We have
\be
  \sum_{\pi \in \scrs_k} \sgn(\pi)
  \prod_{t=1}^k  (X^{-\rm T})_{i_t j_{\pi(t)}}
  \;=\;
  \det( (X^{-\rm T})_{IJ} )
  \;=\;
  (\det X)^{-1} \, \epsilon(I,J) \, (\det X_{I^c J^c})
 \label{eq.sumoverpi}
\ee
by Jacobi's identity [Lemma~\ref{lemma.properties.det}(e)], while
\be
   \sum_{\tau \in \scrs_k}
   s^{\#(\hbox{\scriptsize\rm cycles of } \tau)}
   \;=\;
   s(s+1) \cdots (s+k-1)
 \label{eq.cycles.genfn}
\ee
(see e.g.\ \cite[p.~263, eq.~(6.11)]{Graham_94}
 or \cite[Proposition~1.3.4]{Stanley_97}
 for this well-known equality).
\qed

\subsection{Symmetric Cayley identity}

\proofof{Theorem~\ref{thm.symcayley}}
In this case the index set $\Gamma$ consists of
ordered pairs $(i,j) \in [n] \times [n]$ with $i \le j$,
and we write
$X^{\rm sym} = \sum_{i < j} x_{ij} (E^{ij}+E^{ji})
   + \sum_{i} x_{ii} E^{ii}$.
Then
\be
\partial_{ij}^{\rm sym} X^{\rm sym} \;=\; \smhalf (E^{ij}+E^{ji})
\ee
in all three cases ($i<j$, $i>j$ and $i=j$).
Now let us apply $\det (\partial_{IJ}^{\rm sym})$ to
$(\det X^{\rm sym})^{s}$ and compare to what we had
in the previous proof.
On the one hand we have a factor $2^{-k}$
coming from the $k$ derivatives.
On the other hand, each $E_A$ is now
a sum of two terms $E^{ij}$ and $E^{ji}$;
in each cycle $(\alpha_1,\ldots,\alpha_l)$ of $\tau$,
the argument of the trace becomes
\be
X^{-1}\left(  E^{i_{\alpha_{1}}j_{\sigma(\alpha_{1})}} +  
E^{j_{\sigma(\alpha_{1})}i_{\alpha_{1}}}\right)
X^{-1}\left(  E^{i_{\alpha_{2}}j_{\sigma(\alpha_{2})}} +  
E^{j_{\sigma(\alpha_{2})}i_{\alpha_{2}}}\right) \cdots
\ee
(to lighten the notation we have written $X$ instead of $X^{\rm sym}$).
We therefore need to sum over all the $2^l$ ways
of choosing $E^{ij}$ or $E^{ji}$ in each factor within the given trace
(hence $2^k$ choices overall).
Performing the trace, we will obtain terms of the form
\be
X^{-1}_{j_{\sigma(\alpha_{r})} i_{\alpha_{r+1}}} ,\;
X^{-1}_{i_{\alpha_{r}} i_{\alpha_{r+1}}} ,\;
X^{-1}_{j_{\sigma(\alpha_{r})}j_{\sigma(\alpha_{r+1})}} ,\;
X^{-1}_{i_{\alpha_{r}} j_{\sigma(\alpha_{r+1})} }   \;.
\ee
Let us now fix one of the $2^k$ choices and sum over
the permutation $\sigma$.
If one or more of the factors is of the form
$X^{-1}_{j_{\sigma(\alpha_{r})}j_{\sigma(\alpha_{r+1})}}$,
then the sum over $\sigma$ will vanish because
the exchange between $\sigma(\alpha_{r})$ and $\sigma(\alpha_{r+1})$
takes a $-1$ from $\sgn(\sigma)$.
Therefore, in each cycle of $\tau$ there are
only two nonvanishing contributions,
corresponding to the two ways of coherently orienting the cycle.
One of these has $X^{-\rm T}$ (as in the previous proof)
and the other has $X^{-1}$, but these are in fact equal
since $X$ is symmetric.
We thus have, compared to the previous proof,
an extra factor $2^{\#(\hbox{\scriptsize\rm cycles of } \tau)}$.
We now change variables, as before, in the sum over permutations.
The sum over $\pi$ gives \reff{eq.sumoverpi} exactly as before,
while the sum over $\tau$ now gives
\be
   2^{-k} \sum_{\tau \in \scrs_k}
   (2s)^{\#(\hbox{\scriptsize\rm cycles of } \tau)}
   \;=\;
   s(s+\smhalf) \cdots \left(s+ {k-1 \over 2} \right)
\ee
by \reff{eq.cycles.genfn}.
\qed

\subsection{Antisymmetric Cayley identity}

Let us recall the definition of the pfaffian
of a $2n \times 2n$ antisymmetric matrix:
\be
   \pf A  \;=\;
   {1 \over 2^n n!}
   \sum_{\sigma \in \scrs_{2n}}  \sgn(\sigma) \,
                  a_{\sigma(1) \sigma(2)} \cdots a_{\sigma(2n-1) \sigma(2n)}
   \;.
 \label{def.pfaffian}
\ee

\proofof{Theorem~\ref{thm.antisymcayley}}
In this case the index set $\Gamma$ consists of
ordered pairs $(i,j) \in [n] \times [n]$ with $i < j$,
and we write
$X^{\rm antisym} = \sum_{i < j} x_{ij} (E^{ij}-E^{ji})$. Then
\be
\partial_{ij}^{\rm antisym} X^{\rm antisym} =  E^{ij}-E^{ji}
  \;.
\ee
Now let $I = \{i_1,\ldots,i_{2k} \}$ with $i_1 < \ldots < i_{2k}$,
and let us apply $\pf (\partial_{II}^{\rm antisym})$
to $(\pf X^{\rm antisym})^{s} = (\det X^{\rm antisym})^{s/2}$:
using the representation \reff{def.pfaffian}
for $\pf (\partial_{II}^{\rm antisym})$, we obtain
\be
  \pf (\partial_{II}^{\rm antisym}) \, (\pf X^{\rm antisym})^{s}
  \;=\; 
  {1\over  2^{k}{k}!} \sum_{\sigma \in \scrs_{2k}}
   \sgn(\sigma) \,  \Bigl( \prod_{r=1}^k  {\partial \over \partial 
x_{i_{\sigma(2r-1)}i_{\sigma(2r)}}}  \Bigr)
  \, (\det X^{\rm antisym})^{s/2}
   \;.
\ee
Now apply Lemma~\ref{lemma.elem} as before.
In each cycle $(\alpha_1,\ldots,\alpha_\ell)$ of $\tau$,
the argument of the trace becomes
\be
   X^{-1} \left( E^{i_{\sigma(2 \alpha_{1}-1)}i_{\sigma(2 \alpha_{1})}}- 
  E^{i_{\sigma(2 \alpha_{1})}i_{\sigma(2 \alpha_{1}-1)}}\right)
  X^{-1} \left( E^{i_{\sigma(2 \alpha_{2}-1)}i_{\sigma(2 \alpha_{2})}} -
                 E^{i_{\sigma(2 \alpha_{2})}i_{\sigma(2 \alpha_{2}-1)}}\right)
  \cdots
  \;.
\ee
Once again we have $2^k$ choices in the $E$ factors;
but here these choices correspond simply to pre-multiplying $\sigma$
by one of the $2^k$ permutations that leave fixed the pairs
$\{1,2\}, \,\ldots, \, \{2k-1,2k\}$;
therefore, after summing over $\sigma$ we simply get a factor $2^k$.
Let us now introduce the permutation $\sigma^\tau$ defined by
\begin{subeqnarray}
  \sigma^\tau ( 2 r - 1) & = & \sigma(2 \tau(r)-1)  \\
  \sigma^\tau ( 2r) &= &\sigma(2r )
\end{subeqnarray}
For each $\tau \in \scrs_k$,
the map $\sigma \mapsto \sigma^\tau$ is an automorphism of $\scrs_{2k}$
and satisfies $\sgn(\sigma^\tau)=\sgn(\sigma)\sgn(\tau)$.
We have
\begin{eqnarray}
& &
   \pf(\partial_{II}^{\rm antisym}) \, (\pf X^{\rm antisym})^{s}
      \nonumber \\
& &
\qquad =\;
   2^k \, {1\over  2^k k!}
   \sum_{\sigma \in \scrs_{2k}}  \sgn(\sigma)
   \sum_{\tau \in \scrs_{k}}  \sgn(\tau) \,
   (s/2)^{\#(\hbox{\scriptsize\rm cycles of } \tau)}
   \prod_{r=1}^k (X^{-\rm T})_{i_{\sigma(2 \tau(r)-1)}i_{\sigma(2r)}}
\;.
   \nonumber \\
\end{eqnarray}
We now define $\pi = \sigma^\tau$
and change variables from $(\sigma,\tau)$ to $(\tau,\pi)$ as before.
The sum over $\pi$ gives
\be
   {1\over  2^k k!}
   \sum_{\pi \in \scrs_{2k}}  \sgn(\pi)
   \prod_{r=1}^k (X^{-\rm T})_{i_{\pi(2r-1)} i_{\pi(2r)}}
   \;=\;
   \pf((X^{-\rm T})_{II})
   \;=\;
   (\pf X)^{-1} \, \epsilon(I) \, (\pf X_{I^c I^c})
\ee
by the pfaffian version of Jacobi's identity
[cf.\ \reff{eq.app.pfaffian.jacobi}],
while the sum over $\tau$ gives
\be
   2^k \sum_{\tau \in \scrs_{k}}
   (s/2)^{\#(\hbox{\scriptsize\rm cycles of } \tau)}
   \;=\;
   s(s+2) \cdots (s+2k-2)
\ee
by \reff{eq.cycles.genfn}.
\qed

\section[Proofs of Cayley-type identities by Grassmann/Gaussian
         representation of $(\det X)^s$]
        {Proofs of Cayley-type identities by Grassmann/\hfill\break Gaussian
         representation of $(\det X)^s$}      \label{sec.grassmann.1}

Let us now give simple proofs of
Theorems~\ref{thm.cayley}--\ref{thm.antisymcayley},
based on representing $(\det X)^s$
as a fermionic or bosonic Gaussian integral.
A brief introduction to fermionic and bosonic Gaussian integration
can be found in Appendix~\ref{app.grassmann}.

\subsection{Ordinary Cayley identity}

\proofof{Theorem~\ref{thm.cayley}}
Assume that $s$ is a positive integer, and let us introduce
Grassmann variables $\psi_i^{(\alpha)}, \psibar_i^{(\alpha)}$
for $i=1,\ldots,n$ and $\alpha=1,\ldots,s$.  We can then write
\be
   (\det X)^s  \;=\;  \int \! \scrd(\psi,\psibar) \; e^{\psibar X \psi}
   \;,
 \label{eq.detXs}
\ee
where we have used the shorthand
\be
   \psibar X \psi  \;\equiv\;
   \sum_{\alpha=1}^s \sum_{i,j=1}^n
       \psibar_i^{(\alpha)} x_{ij} \psi_j^{(\alpha)}
   \;.
\ee
Now let $I = \{i_1,\ldots,i_k\}$ with $i_1 < \ldots < i_k$
and $J = \{j_1,\ldots,j_k\}$ with $j_1 < \ldots < j_k$,
so that
\be
   \det(\partial_{IJ})  \;=\;
   \sum_{\sigma \in \scrs_k}  \sgn(\sigma) \prod_{r=1}^k
         {\bigpartial \over \bigpartial x_{i_r j_{\sigma(r)}}}
   \;.
\ee
Applying this to \reff{eq.detXs}, we obtain
\be
   \det(\partial_{IJ}) \, (\det X)^s  \;=\;
   \int \! \scrd(\psi,\psibar) \;
   \sum_{\sigma \in \scrs_k}  \sgn(\sigma)
   \sum_{\alpha_1,\ldots,\alpha_k=1}^s
   \Biggl(
    \prod_{r=1}^k  \psibar_{i_r}^{(\alpha_r)} \psi_{j_{\sigma(r)}}^{(\alpha_r)} 
   \Biggr)
   \; e^{\psibar X \psi}
   \;.
\ee
When $X$ is an invertible real or complex matrix,
Wick's theorem for ``complex'' fermions
(Theorem~\ref{thm.wick.complexfermions}) gives
\be
   \int \! \scrd(\psi,\psibar) \;
   \Biggl(
    \prod_{r=1}^k  \psibar_{i_r}^{(\alpha_r)} \psi_{j_{\sigma(r)}}^{(\alpha_r)}
   \Biggr)
   \; e^{\psibar X \psi}
   \;=\;
   (\det X)^s
   \sum_{\tau \in \scrs_k}  \sgn(\tau)
   \prod_{r=1}^k  (X^{-\rm T})_{i_r j_{\sigma(\tau(r))}}
                  \delta_{\alpha_r \alpha_{\tau(r)}}
   \;.
\ee
Let us now define the permutation $\pi= \sigma \circ \tau$
and change variables from $(\sigma,\tau)$ to $(\tau,\pi)$,
using $\sgn(\sigma) \, \sgn(\tau) = \sgn(\pi)$.
We then have
\be
   \sum_{\pi \in \scrs_k} \sgn(\pi)
   \prod_{r=1}^k  (X^{-\rm T})_{i_r j_{\pi(r)}}
   \;=\;
   \det( (X^{-\rm T})_{IJ} )
   \;=\;
   (\det X)^{-1} \, \epsilon(I,J) \, (\det X_{I^c J^c})
\ee
by Jacobi's identity, while
\be
   \sum_{\tau \in \scrs_k} \sum_{\alpha_1,\ldots,\alpha_k=1}^s
     \prod_{r=1}^k \delta_{\alpha_r \alpha_{\tau(r)}}
   \;=\;
   \sum_{\tau \in \scrs_k} s^{\#(\hbox{\scriptsize\rm cycles of } \tau)}
   \;=\;
   s(s+1) \cdots (s+k-1)
\ee
by \reff{eq.cycles.genfn}.
This proves \reff{eq.cayley.2} for integer $s \ge 1$,
when $X$ is an invertible real or complex matrix.
By Proposition~\ref{prop.equiv.symbolic_s},
this is sufficient to prove the identity.
\qed

\alternateproofof{Theorem~\ref{thm.cayley}}
Instead of using ``complex'' fermions, we can use complex bosons.
So assume that $s = -m$ where $m$ is a positive integer,
and that $X$ is a complex matrix whose hermitian part is
positive-definite;
and let us introduce bosonic variables
$\varphi_i^{(\alpha)}, \varphibar_i^{(\alpha)}$
for $i=1,\ldots,n$ and $\alpha=1,\ldots,m$.  We then have
\be
   (\det X)^{-m}  \;=\;
   \int \! \scrd(\varphi,\varphibar) \; e^{-\varphibar X \varphi}
   \;.
 \label{eq.detXs.bis}
\ee
By the same method as before, we obtain
\be
   \det(\partial_{IJ}) \, (\det X)^{-m}  \;=\;
   (-1)^k \int \! \scrd(\varphi,\varphibar) \;
   \sum_{\sigma \in \scrs_k}  \sgn(\sigma)
   \sum_{\alpha_1,\ldots,\alpha_k=1}^m
   \Biggl(
       \prod_{r=1}^k
       \varphibar_{i_r}^{(\alpha_r)} \varphi_{j_{\sigma(r)}}^{(\alpha_r)} 
   \Biggr)
   \; e^{-\varphibar X \varphi}
   \;.
\ee
Wick's theorem for complex bosons
(Theorem~\ref{thm.wick.complexbosons}) now gives
\be
   \int \! \scrd(\varphi,\varphibar) \;
   \Biggl(
       \prod_{r=1}^k
       \varphibar_{i_r}^{(\alpha_r)} \varphi_{j_{\sigma(r)}}^{(\alpha_r)}
   \Biggr)
   \; e^{-\varphibar X \varphi}
   \;=\;
   (\det X)^{-m}
   \sum_{\tau \in \scrs_k}
   \prod_{r=1}^k  (X^{-\rm T})_{i_r j_{\sigma(\tau(r))}}
                  \delta_{\alpha_r \alpha_{\tau(r)}}
   \;.
\ee
Once again we define the permutation $\pi= \sigma \circ \tau$
and change variables,
using now $\sgn(\sigma) = \sgn(\tau) \, \sgn(\pi)$.
We again have
\be
   \sum_{\pi \in \scrs_k} \sgn(\pi)
   \prod_{r=1}^k  (X^{-\rm T})_{i_r j_{\pi(r)}}
   \;=\;
   \det( (X^{-\rm T})_{IJ} )
   \;=\;
   (\det X)^{-1} \, \epsilon(I,J) \, (\det X_{I^c J^c})
   \;,
\ee
while now
\begin{eqnarray}
   \sum_{\tau \in \scrs_k} \sgn(\tau) \sum_{\alpha_1,\ldots,\alpha_k=1}^m
     \prod_{r=1}^k \delta_{\alpha_r \alpha_{\tau(r)}}
   & = &
   \sum_{\tau \in \scrs_k} \sgn(\tau) \,
      m^{\#(\hbox{\scriptsize\rm cycles of } \tau)}
   \nonumber \\[1mm]
   & = &
   \sum_{\tau \in \scrs_k} (-1)^{\#(\hbox{\scriptsize\rm odd cycles of } \tau)}
      \, (-m)^{\#(\hbox{\scriptsize\rm cycles of } \tau)}
   \nonumber \\[1mm]
   & = &
   (-1)^k s(s+1) \cdots (s+k-1)
\end{eqnarray}
by \reff{eq.cycles.genfn} since $s=-m$.
This proves \reff{eq.cayley.2} for integer $s \le -1$,
when $X$ is a complex matrix whose hermitian part is
positive-definite;
we conclude by Proposition~\ref{prop.equiv.symbolic_s} as before.
\qed

\subsection{Symmetric Cayley identity}

\proofof{Theorem~\ref{thm.symcayley}}
In this case we use real bosons.
So assume that $s = -m/2$ where $m$ is a positive integer,
and that $X^{\rm sym}$ is a real symmetric positive-definite matrix;
and let us introduce bosonic variables $\varphi_i^{(\alpha)}$
for $i=1,\ldots,n$ and $\alpha=1,\ldots,m$.  We then have
\be
   (\det X^{\rm sym})^{-m/2}  \;=\;
   \int \! \scrd\varphi \; e^{-\smsmhalf\varphi X^{\rm sym} \varphi}
   \;.
 \label{eq.detXs.bis2}
\ee
The operator $\partial^{\rm sym}_{ij} =
  \smhalf (1 + \delta_{ij}) \partial/\partial x_{ij}$
is exactly what is needed to bring down a factor $-\smhalf \varphi_i \varphi_j$
when acting on $e^{-\smsmhalf\varphi X^{\rm sym} \varphi}$.
(To lighten the notation, let us henceforth write $X$ instead of
 $X^{\rm sym}$.)
We therefore have
\be
   \det(\partial^{\rm sym}_{IJ}) \, (\det X)^{-m/2}  \;=\;
   (-\smhalf)^k \int \! \scrd\varphi \;
   \sum_{\sigma \in \scrs_k}  \sgn(\sigma)
   \!\!
   \sum\limits_{\alpha_1,\ldots,\alpha_k=1}^m
   \Biggl(
       \prod_{r=1}^k
       \varphi_{i_r}^{(\alpha_r)} \varphi_{j_{\sigma(r)}}^{(\alpha_r)}
   \Biggr)
   \; e^{-\smsmhalf \varphi X \varphi}
   \;.
 \label{eq.int.3}
\ee
Wick's theorem for real bosons
(Theorem~\ref{thm.wick.realbosons})
applied to $\int \! \scrd\varphi \, \prod_{r=1}^k 
    \varphi_{i_r}^{(\alpha_r)} \varphi_{j_{\sigma(r)}}^{(\alpha_r)}$
now gives rise to two types of contractions:
those that only pair $i$'s with $j$'s,
and those that pair at least one $i$ with another $i$
(and hence also a $j$ with another $j$).
The pairings of the first class are given by a sum over permutations $\tau$,
and yield
\be
   (\det X)^{-m/2}
   \sum_{\tau \in \scrs_k}
   \prod_{r=1}^k  (X^{-1})_{i_r j_{\sigma(\tau(r))}}
                  \delta_{\alpha_r \alpha_{\tau(r)}}
   \;.
\ee
Changing variables as in the alternate proof of Theorem~\ref{thm.cayley},
the sum over $\pi$ again yields
$(\det X)^{-1} \, \epsilon(I,J) \, (\det X_{I^c J^c})$,
while the sum over $\tau$ yields
\be
   \sum_{\tau \in \scrs_k} (-1)^{\#(\hbox{\scriptsize\rm odd cycles of } \tau)}
      \, (-m)^{\#(\hbox{\scriptsize\rm cycles of } \tau)}
   \;=\;
   (-1)^k (2s)(2s+1) \cdots (2s+k-1)
   \qquad
 \label{eq.perms.3}
\ee
by \reff{eq.cycles.genfn} since $s=-m/2$.
Inserting these results into \reff{eq.int.3},
we obtain \reff{eq.symcayley.2}.

Let us now show that for each pairing of the second class,
the sum over $\sigma$ yields zero.
By hypothesis, at least one $j$ is paired with another $j$,
say $j_p$ with $j_q$.
Now let $\pi_{pq} \in \scrs_k$ be the permutation that interchanges
$p$ and $q$ while leaving all other elements fixed.
Then $\sigma \mapsto \sigma \circ \pi_{pq}$ is a sign-reversing involution
for the sum in question.

Once again we argue that because the equality \reff{eq.symcayley.2} holds
for infinitely many $s$ and for $X$ in a nonempty open set,
it must hold symbolically in $s$.
\qed

\subsection{Antisymmetric Cayley identity}

\proofof{Theorem~\ref{thm.antisymcayley}}
In this case we use ``real'' fermions.
So assume that $s$ is a positive integer, and let us introduce
``real'' Grassmann variables $\theta_i^{(\alpha)}$
for $i=1,\ldots,2m$ and $\alpha=1,\ldots,s$.  We can then write
\be
   (\pf X^{\rm antisym})^s  \;=\;
   \int \! \scrd\theta \; e^{\smsmhalf \theta X^{\rm antisym} \theta}
   \;.
 \label{eq.detXs.bis3}
\ee
Now let $I = \{i_1,\ldots,i_{2k}\}$ with $i_1 < \ldots < i_{2k}$, so that
\be
   \pf(\partial^{\rm antisym}_{II})  \;=\;
   {1 \over 2^k k!} 
   \sum_{\sigma \in \scrs_{2k}}  \sgn(\sigma) \prod_{r=1}^k
         {\bigpartial \over \bigpartial x_{i_{\sigma(2r-1)} i_{\sigma(2r)}}}
   \;.
\ee
Applying this to \reff{eq.detXs.bis3}, we obtain
\be
   \pf(\partial^{\rm antisym}_{II}) \, (\pf X^{\rm antisym})^s
   \;=\;
   \int \! \scrd\theta \;
   {1 \over 2^k k!}
   \sum_{\sigma \in \scrs_{2k}}  \sgn(\sigma)
   \sum_{\alpha_1,\ldots,\alpha_k=1}^s
   \Biggl(
    \prod_{r=1}^k  \theta_{i_{\sigma(2r-1)}}^{(\alpha_r)}
                   \theta_{i_{\sigma(2r)}}^{(\alpha_r)}
   \Biggr)
   \; e^{\smsmhalf \theta X^{\rm antisym} \theta}
   \;.
\ee
(To lighten the notation, let us henceforth write $X$ instead of
 $X^{\rm antisym}$.)
When $X$ is an invertible real or complex matrix,
Wick's theorem for ``real'' fermions
(Theorem~\ref{thm.wick.realfermions}) gives
\begin{eqnarray}
& &
   \int \! \scrd \theta 
   \Biggl(
       \prod_{r=1}^k
       \theta_{i_{\sigma(2r-1)}}^{(\alpha_r)}
       \theta_{i_{\sigma(2r)}}^{(\alpha_r)}
   \Biggr)
   \; e^{\smsmhalf \theta X \theta}
\nonumber \\
& &
\qquad =\;
   (\pf X)^s
   {1 \over 2^k k!}
   \sum_{\tau \in \scrs_{2k}}  \sgn(\tau)
   \prod_{r=1}^k  (X^{-\rm T})_{i_{\sigma(\tau(2r-1))} i_{\sigma(\tau(2r))}}
      \delta_{\alpha_{\lceil \tau(2r-1) \rceil}
              \alpha_{\lceil \tau(2r) \rceil}
             }
   \;.
\nonumber \\
\end{eqnarray}
Once again we define the permutation $\pi= \sigma \circ \tau$
and change variables from $(\sigma,\tau)$ to $(\tau,\pi)$,
using $\sgn(\sigma) \, \sgn(\tau) = \sgn(\pi)$.
The sum over $\pi$ gives
\be
   {1 \over 2^k k!}
   \sum_{\pi \in \scrs_{2k}}  \sgn(\pi)
   \prod_{r=1}^k  (X^{-\rm T})_{i_{\pi(2r-1)} i_{\pi(2r)}}
   \;=\;
   \pf( (X^{-\rm T})_{II} )
   \;=\;
   (\pf X)^{-1} \, \epsilon(I) \, (\pf X_{I^c I^c})
\ee
by the pfaffian version of Jacobi's identity
[cf.\ \reff{eq.app.pfaffian.jacobi}].
The remaining sums give
\be
   {1 \over 2^k k!}
   \sum_{\tau \in \scrs_{2k}}
   \sum_{\alpha_1,\ldots,\alpha_k=1}^s
   \prod_{r=1}^k
      \delta_{\alpha_{\lceil \tau(2r-1) \rceil}
              \alpha_{\lceil \tau(2r) \rceil}
             }
   \;=\;
   \sum_{M \in \scrm_{2k}}
   s^{\#(\hbox{\scriptsize\rm cycles of } M \cup M_0)}
  \;,
\ee
where the latter sum runs over all perfect matchings $M$ of $2k$ elements,
and $M_0$ is some fixed perfect matching
[in our case $(12)(34) \cdots (2k\!-\!1 \: 2k)$];
we observe that if $M$ and $M_0$ are thought of as edge sets on $[2k]$,
then $M \cup M_0$ is the edge set of a graph in which each vertex
has degree 2, and so is a union of cycles,
showing that the right-hand side is well-defined.
Let us now show that this latter sum equals $s(s+2) \cdots (s+2k-2)$:


Let $e_1,\ldots,e_k$ be the edges of $M_0$.
Then each matching $M \in \scrm_{2k}$ induces a decomposition of
the set $\{e_1,\ldots,e_k\}$ into cycles,
according to how those edges are traversed in $M \cup M_0$;
in other words, $M$ induces a permutation of $[k]$.
Moreover, for each cycle $C$ in this decomposition,
there are $2^{|C| - 1}$ ways of connecting up the vertices.
Thus each permutation $\pi$ of $[k]$ arises from
$\prod_{C \in \pi} 2^{|C|-1} = 2^{k-\#(\hbox{\scriptsize\rm cycles of } \pi)}$
different matchings $M$.  We therefore have
\begin{subeqnarray}
   \sum_{M \in \scrm_{2k}}
      s^{\#(\hbox{\scriptsize\rm cycles of } M \cup M_0)}
   & = &
   2^k \sum_{\pi \in \scrs_k}
      (s/2)^{\#(\hbox{\scriptsize\rm cycles of } \pi)}
   \\[1mm]
   & = &
   s(s+2) \cdots (s+2k-2)
   \;.
 \label{eq.M2k.Pk}
\end{subeqnarray}

Once again we invoke Proposition~\ref{prop.equiv.symbolic_s}
to conclude.\footnote{
   We thank Alex Scott for help in cleaning up
   our proof of \reff{eq.M2k.Pk}.
}
\qed

\section{Proofs of Cayley-type identities by Grassmann
   representation of $\det(\partial)$}      \label{sec.grassmann.2}

In this section we give an alternate Grassmann-based approach
to proving Cayley-type identities:
now it is the differential operator $\det(\partial)$
that is represented as a fermionic Gaussian integral.
This technique is in our opinion very powerful:
not only does it give the slickest proofs of the three main
Cayley-type identities for square matrices;
it also gives the only direct algebraic/combinatorial proofs (thus far) of the
(considerably more difficult) rectangular Cayley identities.

The basic fact we will need is that an operator
$\exp(a \cdot \partial)$ generates translation by $a$.
More precisely, let $R$ be a commutative ring containing the rationals;
then for any polynomial $P(z_1,\ldots,z_n)$ with coefficients in $R$
and any constants $a_1,\ldots,a_n \in R$, we have
the {\em translation formula}\/
\be
   \exp\left( \sum_i a_i {\partial \over \partial z_i} \right) \, P(z)
   \;=\;
   P(z+a)  \;.
 \label{eq.translation}
\ee
(Here $\exp$ is defined by Taylor series;
 note that all but finitely many terms will annihilate $P$.)
Indeed, the identity \reff{eq.translation} is nothing other than
Taylor's theorem for polynomials $P$.
In particular we will use \reff{eq.translation} when the
commutative ring $R$ consists of the even elements of
some Grassmann algebra.
Moreover, in our applications the elements $a_i$ will be nilpotent,
so that the Taylor series for $\exp(a \cdot \partial)$ is in fact finite.

We will also need a formula for the change of a determinant
under a low-rank perturbation:
see Appendix~\ref{app.det}.

Let us begin by explaining the general structure of all these proofs.
We introduce a Grassmann integral representation for 
the differential operator $\det(\partial)$
and let it act on $(\det X)^s$.
After a change of variables in the Grassmann integral,
we obtain the desired quantity $(\det X)^{s-1}$
[or its generalization for minors $I,J$]
multiplied by a purely combinatorial factor
that is independent of the matrix $X$.
We then proceed to calculate this combinatorial factor,
which turns out to be an explicit polynomial in $s$.

Unfortunately, the ``all-minors'' versions of these proofs
(i.e.\ those for $|I| = |J| = k < n$)
are slightly more complicated than the ``basic'' versions
(i.e.\ those for $I = J = [n]$).
We have therefore structured our presentation so as to
give the ``basic'' proof first,
and then indicate the modifications needed to handle the
``all-minors'' case.

\subsection{Ordinary Cayley identity}  \label{sec.grassmann.2.ordinary}

\proofof{Theorem~\ref{thm.cayley}}
We introduce Grassmann variables $\eta_i, \etabar_i$ ($1 \le i \le n$)
and use the representation
\be
   \det(\partial)  \;=\;
   \int \! \scrd_n(\eta,\etabar) \; e^{\etabar^{\rm T} \partial \eta}
   \;,
 \label{eq.detpartial}
\ee
where the subscript on $\scrd$ serves to remind us
of the length of the vectors in question,
and we have employed the shorthand notation
\be
   \etabar^{\rm T} \partial \eta
   \;\equiv\;
   \sum_{i,j=1}^n \etabar_i  \, {\partial \over \partial x_{ij}} \, \eta_j
   \;=\;
   \sum_{i,j=1}^n \etabar_i  \eta_j \, {\partial \over \partial x_{ij}}
   \;.
\ee
[Indeed, \reff{eq.detpartial} is simply a special case of
the fermionic Gaussian integral \reff{gaussint.fermionic.complex}
where the coefficient ring $R$ is the ring $\Q[\partial]$
of polynomials in the differential operators $\partial/\partial x_{ij}$.]
Applying \reff{eq.detpartial} to $(\det X)^s$
where $s$ is a positive integer
and using the translation formula \reff{eq.translation}, we obtain\footnote{
   We have assumed here that $s$ is a positive integer,
   because we have proven \reff{eq.translation} only for
   {\em polynomials}\/ $P$.
   Alternatively, we could avoid this assumption
   by proving \reff{eq.translation}
   also for more general functions (e.g.\ powers of polynomials)
   when all the $a_i$ are nilpotent.
}
\begin{subeqnarray}
   \det(\partial) \, (\det X)^s
   & = &
   \int \! \scrd_n(\eta,\etabar) \; \det(X + \etabar \eta^{\rm T})^s
      \\
   & = &
   (\det X)^s \int \! \scrd_n(\eta,\etabar) \;
          \det(I + X^{-1} \etabar \eta^{\rm T})^s
   \;.
\end{subeqnarray}
(We assume here that $X$ is an invertible real or complex matrix.)
Let us now change variables from $(\eta,\etabar)$ to
$(\eta',\etabar') \equiv (\eta, X^{-1} \etabar)$;
we pick up a Jacobian $\det(X^{-1}) = (\det X)^{-1}$
and thus have
\be
   \det(\partial) \, (\det X)^s
   \;=\;
   (\det X)^{s-1} \int \! \scrd_n(\eta',\etabar') \;
          \det(I + \etabar' \eta^{\prime\rm T})^s
   \;.
 \label{eq.grass1.star1}
\ee
Formula \reff{eq.grass1.star1} expresses $\det(\partial) \, (\det X)^s$
as the desired quantity $(\det X)^{s-1}$ times a purely combinatorial factor
\be
   P(s,n)  \;\equiv\;  \int \! \scrd_n(\eta,\etabar) \;
          \det(I + \etabar \eta^{\rm T})^s
  \;,
\ee
which we now proceed to calculate.
The matrix $I + \etabar \eta^{\rm T}$ is a rank-1 perturbation 
of the identity matrix;  by Lemma~\ref{lemma.lowrank} we have
\begin{subeqnarray}
   \det(I + \etabar \eta^{\rm T})^s
   & = &
   (1 - \etabar^{\rm T} \eta)^{-s}
      \\[1mm]
   & = & \sum_{\ell=0}^\infty
            (-1)^\ell \, {-s \choose \ell} \, (\etabar^{\rm T} \eta)^\ell
 \label{eq.grass1.star2}
\end{subeqnarray}
where
\be
   \etabar^{\rm T} \eta  \;\equiv\; \sum_{i=1}^n \etabar_i \eta_i
   \;.
\ee
Since
\be
   \int \! \scrd_n(\eta,\etabar) \; (\etabar^{\rm T} \eta)^\ell
   \;=\;
   n! \, \delta_{\ell,n} 
   \;,
 \label{eq.scalar.etaetabar}
\ee
it follows that
\begin{subeqnarray}
   P(s,n)
   & = &
   (-1)^n \, {-s \choose n} \, n!
      \\[1mm]
   & = &
   s(s+1) \cdots (s+n-1)
   \;.
 \label{eq.grass1.star3}
\end{subeqnarray}
This proves \reff{eq.cayley.1}
when $X$ is an invertible real or complex matrix
and $s$ is a positive integer;
the general validity of the identity then follows from
Proposition~\ref{prop.equiv.symbolic_s}.\footnote{
   Alternatively, one can observe
   (see the Remark after Theorem~\ref{thm.cayley})
   that the Cayley identity \reff{eq.cayley.1} is equivalent to
   $$
   \det(\partial) \det(X+A)^s
   \;=\;
   s(s+1)\cdots(s+n-1) \det(X+A)^{s-1}
   $$
   for any fixed matrix $A$.
   This latter identity can be proven exactly as above,
   with $X^{-1}$ replaced by $(X+A)^{-1}$,
   and it works whenever $X+A$ is invertible.
   Taking, for instance, $A=-cI$, we prove the identity
   for all real or complex matrices $X$ whose spectrum does not contain
   the point $c$.
   Putting together these identities for all $c \in \C$,
   we cover all matrices $X$.
}

Let us now indicate the modifications needed to prove \reff{eq.cayley.2}
for a minor $I,J$, where $I,J \subseteq [n]$ with $|I|=|J|=k$.
We begin by introducing a Grassmann representation for $\det(\partial_{IJ})$:
\be
\det(\partial_{IJ})
\;=\;
\epsilon(I,J)
\int \! \scrd_n(\eta, \etabar) \,
\left(\prod \etabar \eta \right)_{I^c,J^c} \;
e^{\etabar^{\rm T} \partial \eta}
\ee
where
\be
\left(\prod \etabar \eta \right)_{I^c,J^c}
 \;\equiv\;
\etabar_{i_1} \eta_{j_1} \cdots \etabar_{i_{n-k}} \eta_{j_{n-k}}
\ee
and $i_1,\ldots,i_{n-k}$ (resp.\ $j_1,\ldots,j_{n-k}$)
are the elements of $I^c$ (resp.\ $J^c$) in increasing order.
Applying the translation formula \reff{eq.translation} as before,
we obtain
\begin{subeqnarray}
   \det(\partial_{IJ}) \, (\det X)^s
   & = &
   \epsilon(I,J)\int \! \scrd_n(\eta,\etabar) \;
   \Big(\prod \etabar \eta \Big)_{I^c,J^c} \;
   \det(X + \etabar \eta^{\rm T})^s
       \slabel{eq.grassmann.ordinary.IJ.a}      \\
   & = &
   \epsilon(I,J) \, (\det X)^s \int \! \scrd_n(\eta,\etabar) \;
   \Big(\prod \etabar \eta \Big)_{I^c,J^c} \;
          \det(I + X^{-1} \etabar \eta^{\rm T})^s
   \;.
   \nonumber \\[-2mm]
       \label{eq.grassmann.ordinary.IJ}
\end{subeqnarray}
Once again we change variables from $(\eta,\etabar)$ to
$(\eta',\etabar') \equiv (\eta, X^{-1} \etabar)$, picking
up a Jacobian $\det(X^{-1}) = (\det X)^{-1}$;
and we use the identity \reff{eq.grass1.star2}, obtaining
\begin{subeqnarray}
(\det \partial_{IJ}) \, (\det X)^s
&=&
\epsilon(I,J) \, (\det X)^{s-1}
\int \! \scrd_n(\eta, \etabar)
\left(\prod (X\etabar) \eta \right)_{I^c,J^c} \,
\sum_{\ell=0}^{\infty}
(-1)^\ell \binom{-s}{\ell}
(\etabar^{\rm T} \eta)^{\ell}
    \nonumber \\ \\[-2mm]
 & = &
\epsilon(I,J) \, (\det X)^{s-1}
\sum_{r_1,\ldots,r_{n-k} \in [n]}
\left(\prod_{p=1}^{n-k}X_{i_p,r_p}\right)
\sum_{\ell=0}^{\infty} (-1)^\ell \binom{-s}{\ell} \times \nonumber \\
&& \qquad\qquad \int\!\scrd_n(\eta, \etabar)
\left(\prod_{p=1}^{n-k}\etabar_{r_p}\eta_{j_p}\right)
(\etabar^{\rm T} \eta)^{\ell}
   \;.
\end{subeqnarray}
The rules of Grassmann integration constrain the integral to be zero
unless $\ell = k$ and
$(r_1, \dots, r_{n-k}) = (j_{\sigma(1)}, \ldots, j_{\sigma(n-k)})$
for some permutation $\sigma \in \scrs_{n-k}$.
We therefore have
\begin{subeqnarray}
(\det \partial_{IJ}) \, (\det X)^s
&=&
\epsilon(I,J) \, (\det X)^{s-1} \!
\sum_{\sigma \in \scrs_{n-k}} \sgn(\sigma)
\left(\prod_{p=1}^{n-k}X_{i_p,j_{\sigma(p)}}\right)
  (-1)^k \binom{-s}{k} \times \nonumber\\
&& \qquad  \int \! \scrd_n(\eta, \etabar) \,
  \left(\prod_{p=1}^{n-k}\etabar_{j_p}\eta_{j_p}\right)
k!\left(\prod_{i\in J}\etabar_{i}\eta_{i}\right)   \\
&=&
\epsilon(I,J) \, (\det X)^{s-1}
 \, s(s+1)\cdots(s+k-1)
\sum_{\sigma \in \scrs_{n-k}}
\sgn(\sigma)
\prod_{p=1}^{n-k}X_{i_p,j_{\sigma(p)}}  \nonumber\\ \\
&=&
\epsilon(I,J) \, (\det X)^{s-1}
 \, s(s+1)\cdots(s+k-1) \, \det(X_{I^c J^c})
   \;.
       \label{eq.grassmann.ordinary.IJ.final}
\end{subeqnarray}
This proves \reff{eq.cayley.2}.
\qed

\subsection{Two useful lemmas for the all-minors case}
   \label{sec.grassmann.2.lemmas}

Let us now pause to abstract the type of reasoning that was just used
in proving the all-minors identity \reff{eq.cayley.2},
as similar reasoning will be needed to prove the all-minors
versions of the symmetric and rectangular Cayley identities
and the all-principal-minors versions of the antisymmetric identities.
The reader who is interested mainly in the case $I=J=[n]$
can skip this subsection on a first reading.

The key results of this subsection will be a pair of general formulae:
Lemma~\ref{lem.pbABp.Alan} for ``complex'' fermions
and Lemma~\ref{lem.pbCJCp.Alan} for ``real'' fermions.
Important special cases of these formulae
(which also have easier direct proofs) will be stated in
Corollaries~\ref{cor.pbABp.Marco} and \ref{cor.Atheta.Marco}, respectively.
These corollaries are, in fact, all we need to handle the
all-minors symmetric Cayley identity (Section~\ref{sec.grassmann.2.symmetric})
and the all-principal-minors antisymmetric Cayley identity
(Section~\ref{sec.grassmann.2.antisymm}).
However, the rectangular Cayley identities
(Sections~\ref{sec.grassmann.2.tmrect}--\ref{sec.grassmann.2.multi})
will need the full strength of Lemmas~\ref{lem.pbABp.Alan} and
\ref{lem.pbCJCp.Alan} to handle the all-minors case.

Let $I,J \subseteq [n]$ with $|I| = |J| = k$,
and let $I^c, J^c$ be the complementary subsets.
We denote by $i_1,\ldots,i_{n-k}$ (resp.\ $j_1,\ldots,j_{n-k}$)
the elements of $I^c$ (resp.\ $J^c$) in increasing order.
Let $\eta_i, \etabar_i$ ($1 \le i \le n$) be Grassmann variables as before.
Then, for any $n \times n$ matrices $A,B$,
we define
\be
\Big( \prod (A \etabar) (B\eta) \Big)_{I^c,J^c}
  \;=\;
(A \etabar)_{i_1} (B\eta)_{j_1} 
\cdots (A \etabar)_{i_{n-k}} (B\eta)_{j_{n-k}}
   \;.
\label{eq.prodAebeB.Marco}
\ee
In particular, $\bigl( (A \etabar) (B \eta) \bigr)_{\emptyset,\emptyset}=1$.

Now suppose that we have $N$ further sets of (real) Grassmann variables
$\theta_i^{(\alpha)}$ ($1 \le i \le n$, $1 \le \alpha \le N$)
--- the case $N=0$ is also allowed ---
and suppose that $f(\eta,\etabar,\theta)$
is a polynomial in the scalar products
$\etabar^{\rm T} \eta$,
$\etabar^{\rm T} \theta^{(\alpha)}$,
$\eta^{\rm T} \theta^{(\alpha)}$
and $\theta^{(\alpha) \rm T} \theta^{(\beta)}$.
Then a Grassmann integral of the form
\be
   \int \! \scrd_n(\eta,\etabar,\theta) \,
   \left( \prod_{i \in L} \etabar_i \eta_i \right)
   f(\eta,\etabar,\theta)
 \label{eq.pre-lemma}
\ee
obviously takes the same value for all sets $L \subseteq [n]$
of the same cardinality.

Let us now show how a Grassmann integral involving
$( \prod (A \etabar) (B\eta) )_{I^c,J^c}$ and $f(\eta,\etabar,\theta)$
can be written as a determinant containing $A$ and $B$
multiplied by a purely combinatorial factor:

\begin{lemma}
  \label{lem.pbABp.Alan}
Let $I,J \subseteq [n]$ with $|I| = |J| = k$,
let $A,B$ be $n \times n$ matrices,
and let $f(\eta,\etabar,\theta)$ be a polynomial in the scalar products
as specified above.
Suppose further that the number $N$ of additional fermion species
is {\em even}\/.
Then
\begin{eqnarray}
& &
\int \! \scrd_n(\eta,\etabar,\theta) \,
\Big( \prod (A \etabar) (B\eta) \Big)_{I^c,J^c} 
\, f(\eta,\etabar,\theta)
    \nonumber \\
& & \qquad =\;
  \det[(A B^{\rm T})_{I^c J^c}]
    \int \! \scrd_n(\eta,\etabar,\theta) \,
    \left( \prod_{i \in L} \etabar_i \eta_i \right)
    f(\eta,\etabar,\theta)
\end{eqnarray}
\label{eq.lemma7.1.Alan}
where $L \subseteq [n]$ is any set of cardinality $n-k$.
\end{lemma}

\proof
The expansion of the product \reff{eq.prodAebeB.Marco} produces
\be
   \Big( \prod (A \etabar) (B\eta) \Big)_{I^c,J^c}
   \;=\;
   \sum_{\begin{scarray}
            r_1,\ldots,r_{n-k} \in [n] \\
            s_1,\ldots,s_{n-k} \in [n]
         \end{scarray}}
   \prod_{p=1}^{n-k} A_{i_p,r_p} B_{j_p,s_p}  \etabar_{r_p} \eta_{s_p}
   \;.
 \label{proof.Alan.1}
\ee
Note that, by nilpotency of Grassmann variables,
the product on the right-hand side of \reff{proof.Alan.1}
is nonvanishing only if
the indices $r_1,\ldots,r_{n-k}$ are all distinct and
the indices $s_1,\ldots,s_{n-k}$ are also all distinct.
Let us now integrate \reff{proof.Alan.1} against
one of the monomials arising in the expansion of $f(\eta,\etabar,\theta)$.
For each ``site'' $i \in [n]$,
this monomial contains an even number of factors
$\eta_i$, $\etabar_i$ or $\theta_i^{(\alpha)}$;
and it must contain each one of the factors $\theta_i^{(\alpha)}$
($1 \le \alpha \le N$) precisely once if the Grassmann integral
over $\theta$ is to be nonvanishing.
Since $N$ is even, this means that at each site $i$
we must either have both of the factors $\eta_i$ and $\etabar_i$,
or neither.
In order to have a nonvanishing Grassmann integral over $\eta$ and $\etabar$,
the former situation must occur at $k$ sites and the latter at $n-k$
sites;  moreover, the latter $n-k$ sites must correspond precisely
to the factors $\etabar_{r_p} \eta_{s_p}$ in \reff{proof.Alan.1}.
We can therefore assume that
$(r_1,\ldots,r_{n-k}) = (s_{\sigma(1)},\ldots,s_{\sigma(n-k)})$
for some permutation $\sigma \in \scrs_{n-k}$;
the contributing terms in \reff{proof.Alan.1} are then
\be
   \sum_{\begin{scarray}
            r_1,\ldots,r_{n-k} \in [n] \\
            \hboxrm{all distinct}
         \end{scarray}}
   \sum_{\sigma \in \scrs_{n-k}}
   \prod_{p=1}^{n-k} A_{i_p,r_p} B_{j_p,r_{\sigma^{-1}(p)}}
         \etabar_{r_p} \eta_{r_{\sigma^{-1}(p)}}
   \;.
\ee
The factors $B_{j_p,r_{\sigma^{-1}(p)}}$ can be reordered freely;
reordering of the Grassmann factors $\eta_{r_{\sigma^{-1}(p)}}$
yields a prefactor $\sgn(\sigma)$.  We therefore obtain
\be
   \sum_{\begin{scarray}
            r_1,\ldots,r_{n-k} \in [n] \\
            \hboxrm{all distinct}
         \end{scarray}}
   \sum_{\sigma \in \scrs_{n-k}} \sgn(\sigma) \,
   \prod_{p=1}^{n-k} A_{i_p,r_p} B_{j_{\sigma(p)},r_p}
         \etabar_{r_p} \eta_{r_p}
   \;.
\ee
By the remarks around \reff{eq.pre-lemma}, the integral
\be
   \int \! \scrd_n(\eta,\etabar,\theta)
   \left( \prod_{p=1}^{n-k} \etabar_{r_p} \eta_{r_p} \right)
   f(\eta,\etabar,\theta)
\ee
is independent of the choice of $r_1,\ldots,r_{n-k}$
(provided that they are all distinct)
and hence can be pulled out.
We are left with the factor
\be
   \sum_{\begin{scarray}
            r_1,\ldots,r_{n-k} \in [n] \\
            \hboxrm{all distinct}
         \end{scarray}}
   \sum_{\sigma \in \scrs_{n-k}} \sgn(\sigma) \,
   \prod_{p=1}^{n-k} A_{i_p,r_p} B_{j_{\sigma(p)},r_p}
   \;.
 \label{proof.Alan.2}
\ee
We can now remove the restriction that the $r_1,\ldots,r_{n-k}$
be all distinct, because the terms with two or more $r_i$ equal
cancel out when we sum over permutations with the factor $\sgn(\sigma)$.
So \reff{proof.Alan.2} equals
\begin{subeqnarray}
   \sum_{r_1,\ldots,r_{n-k} \in [n]}
      \sum_{\sigma \in \scrs_{n-k}} \sgn(\sigma) \,
      \prod_{p=1}^{n-k} A_{i_p,r_p} B_{j_{\sigma(p)},r_p}
   & = &
   \sum_{\sigma \in \scrs_{n-k}} \sgn(\sigma) \,
      \prod_{p=1}^{n-k} (A B^{\rm T})_{i_p,j_{\sigma(p)}}
   \qquad\qquad \\[1mm]
   & = &
   \det[ (A B^{\rm T})_{I^c J^c} ]
   \;.
\end{subeqnarray}
\qed

\smallskip

{\bf Remark.}  If $N$ is odd, the situation is different:
at each site $i$ the monomial coming from $f$ must now provide
{\em exactly one}\/ of the factors $\eta_i$ and $\etabar_i$.
This means that we can get a nonzero contribution only when
$n$ is even and $k=n/2$.
For instance, suppose that $N=1$, $n=2$, $I^c = \{i\}$ and $J^c = \{j\}$,
and that
$f(\eta,\etabar,\theta) = (\eta^{\rm T} \theta) (\etabar^{\rm T} \theta)$.
Then
\be
\int \! \scrd_2(\eta,\etabar,\theta) \,
\Big( \prod (A \etabar) (B\eta) \Big)_{I^c,J^c}
\, (\eta^{\rm T} \theta) (\etabar^{\rm T} \theta)
  \;=\;
  - A_{i1} B_{j2} + A_{i2} B_{j1}
   \;.
\ee
We shall not consider this situation further,
as we shall not need it in the sequel.
\qed



The following identity is what was used in the proof of
the all-minors ordinary Cayley identity \reff{eq.cayley.2}
and will also be used in the proof of
the all-minors symmetric Cayley identity \reff{eq.symcayley.2}:

\begin{corollary}
\label{cor.pbABp.Marco}
Let $I,J \subseteq [n]$ with $|I| = |J| = k$,
let $A,B$ be $n \times n$ matrices,
let $M$ be an invertible $n \times n$ matrix,
and let $\ell$ be a nonnegative integer.
Then
\be
\int \! \scrd_n(\eta,\etabar) \,
\Big( \prod (A \etabar) (B\eta) \Big)_{I^c,J^c} 
(\etabar^{\rm T} M \eta)^\ell
  \;=\;
  k!\,\delta_{\ell,k} \, (\det M) \, \det[(A M^{-\rm T} B^{\rm T})_{I^c J^c}]
  \;.
   \label{eq.cor.pbABp.Marco}
\ee
\end{corollary}

\firstproof
The case $M=I$ is an easy consequence of Lemma~\ref{lem.pbABp.Alan}.
The case of a general invertible matrix $M$ can be reduced to
the case $M=I$ by the change of variables
$(\eta',\etabar') \equiv (M\eta, \etabar)$,
which picks up a Jacobian $\det M$ and replaces $B$ by $B M^{-1}$.
\qed

Here is a simple direct proof of Corollary~\ref{cor.pbABp.Marco}
that avoids the combinatorial complexity
of the full Lemma~\ref{lem.pbABp.Alan}
and instead relies on standard facts from the theory of
Grassmann--Berezin integration (Appendix~\ref{app.grassmann})\footnote{
   Said another way:
   The proof of Wick's theorem for ``complex'' fermions
   (Theorem~\ref{thm.wick.complexfermions})
   requires combinatorial work similar in difficulty to that
   occurring in the proof of Lemma~\ref{lem.pbABp.Alan}.
   But since Wick's theorem is a ``standard'' result,
   we can employ it without reproving it {\em ab initio}\/.
}:

\secondproof
It is easy to see that the Grassmann integral is nonvanishing
only when $\ell = k$.
So in this case we can replace $(\etabar^{\rm T} M \eta)^k$
in the integrand by $k! \, \exp(\etabar^{\rm T} M \eta)$.
The claim is then an immediate consequence of
Wick's theorem for ``complex'' fermions
(Theorem~\ref{thm.wick.complexfermions}).
\qed

\bigskip

In order to handle the all-principal-minors versions
of the antisymmetric identities,
we shall need formulae analogous to Lemma~\ref{lem.pbABp.Alan}
and Corollary~\ref{cor.pbABp.Marco},
but for ``real'' rather than ``complex'' fermions.
So let $I \subseteq [2m]$ with $|I|=2k$,
and let $I^c$ be the complementary subset.
We denote by $i_1,\ldots,i_{2m-2k}$
the elements of $I^c$ in increasing order.
Let $\theta_1, \ldots, \theta_{2m}$ be ``real'' Grassmann variables.
Then, for any $2m \times 2m$ matrix $C$, we define
\be
\left( \prod (C \theta)\right)_{I^c}
\;\equiv\;
(C \theta)_{i_1}
\cdots (C \theta)_{i_{2m-2k}}
   \label{eq.prodAtheta.Marco}
\ee
(note that the order of factors $\theta$ is crucial here
 because they anticommute).
In particular, $(\prod (C \theta))_{\emptyset}=1$.

Now suppose that we have $N$ further sets of (real) Grassmann
variables
$\chi_i^{(\alpha)}$ ($1 \le i \le 2m$, $1 \le \alpha \le N$)
--- the case $N=0$ is also allowed ---
and suppose that $f(\theta,\chi)$
is a polynomial in the symplectic products
$\smfrac{1}{2} \theta J \theta$,
$\theta J \chi^{(\alpha)}$,
$\smfrac{1}{2} \chi^{(\alpha)} J \chi^{(\alpha)}$
and $\chi^{(\alpha)} J \chi^{(\beta)}$,
where the $2m \times 2m$ matrix $J$ is defined by
\be
J \;=\;
\left(
\begin{array}{cc|cc|c}
 0 & 1 & \multicolumn{3}{c}{} \\
-1 & 0 & \multicolumn{3}{c}{} \\
\cline{1-4}
 &  & 0 & 1 & \\
 &  & -1 & 0 & \\
\cline{3-4}
\multicolumn{4}{c}{} & \ddots
\end{array}
\right)  \;.
  \label{def.J.lemma5.3}
\ee
Now consider a Grassmann integral of the form
\be
   \int \! \scrd_{2m} (\theta,\chi) \,
\left( \prod_{i \in L} \theta_i \right)
   f(\theta,\chi)
 \label{eq.pre-lemma.real}
\ee
for some set $L \subseteq [2m]$,
where the product $\prod_{i \in L} \theta_i$ is understood
to be written from left to right in increasing order of the indices.
Clearly $|L|$ must be even for the integral \reff{eq.pre-lemma.real}
to be nonvanishing.  But more is true:  $L$ must in fact be a
union of pairs $\{2j-1,2j\}$, otherwise \reff{eq.pre-lemma.real}
will again vanish.\footnote{
   This statement of course pertains to our specific choice of $J$,
   and would be modified in the obvious way if we had chosen a
   different convention for $J$.
}
To see this, it suffices to observe that all the variables
$\theta$ or $\chi^{(\alpha)}$ with index $2j-1$ or $2j$
appear {\em in pairs}\/ in $f(\theta,\chi)$;
therefore, in order to have an integrand that is {\em even}\/
in these variables, the set $L$ must contain either both of
$2j-1$ and $2j$ or neither.  Let us call a set $L \subseteq [2m]$
{\em well-paired}\/ if it is a union of pairs $\{2j-1,2j\}$.
Obviously the integral \reff{eq.pre-lemma.real} takes the
same value for all well-paired sets
$L \subseteq [2m]$ of the same cardinality.

Let us now show how a Grassmann integral involving
$( \prod (C \theta) )_{I^c}$ and $f(\theta,\chi)$
can be written as a pfaffian containing $C$
multiplied by a purely combinatorial factor:

\begin{lemma}
  \label{lem.pbCJCp.Alan}
Let $I \subseteq [2m]$ with $|I| = 2k$,
let $C$ be a $2m \times 2m$ matrix,
and let $f(\theta,\chi)$ be a polynomial in the symplectic products
as specified above.
Then
\begin{eqnarray}
& &
\int \! \scrd_{2m} (\theta,\chi) \,
\Big( \prod (C \theta) \Big)_{I^c} 
\, f(\theta,\chi)
    \nonumber \\
& & \qquad =\;
  \pf[(C J C^{\rm T})_{I^c I^c}]
    \int \! \scrd_{2m} (\theta,\chi) \,
    \left( \prod_{i \in L} \theta_i \right)
    f(\theta,\chi)
\label{eq.lemma7.1.Alan.now5.3}
\end{eqnarray}
where $L \subseteq [2m]$ is any well-paired set of cardinality $2m-2k$.
\end{lemma}

\proof
The expansion of the product \reff{eq.prodAtheta.Marco} produces
\be
   \Big( \prod (C \theta) \Big)_{I^c}
   \;=\;
   \sum_{r_1,\ldots,r_{2m-2k} \in [2m]}
      \left( \prod_{p=1}^{2m-2k}C_{i_p,r_p} \right)
      \theta_{r_1} \cdots \theta_{r_{2m-2k}}
   \;.
 \label{proof.Alan.1.65876}
\ee
By nilpotency of Grassmann variables,
the product on the right-hand side of \reff{proof.Alan.1.65876}
is nonvanishing only if the indices $r_1,\ldots,r_{2m-2k}$ are all distinct.
Moreover, as noted above, the set $R = \{r_1,\ldots,r_{2m-2k} \}$
must be well-paired if the integral is to be nonvanishing.
Let us denote by $r'_1,\ldots,r'_{2m-2k}$ the elements of $R$
in increasing order.  Then we can write
$(r_1,\ldots,r_{2m-2k}) = (r'_{\sigma(1)},\ldots,r'_{\sigma(2m-2k)})$
for some permutation $\sigma \in \scrs_{2m-2k}$.
The contributing terms in \reff{proof.Alan.1.65876} are then
\be
   \sum_{\begin{scarray}
            R \subseteq [2m] \\
            |R| = 2m-2k \\
            R \, \hboxrm{well-paired}
         \end{scarray}}
   \sum_{\sigma \in \scrs_{2m-2k}}
   \left( \prod_{p=1}^{2m-2k}C_{i_p,r'_{\sigma(p)}} \right)
      \theta_{r'_{\sigma(1)}} \cdots \theta_{r'_{\sigma(2m-2k)}}
   \;.
 \label{proof.Alan.1.65876.bis}
\ee
We can now reorder the factors $\theta$ into increasing order,
yielding a factor $\sgn(\sigma)$, so that \reff{proof.Alan.1.65876.bis} becomes
\be
   \sum_{\begin{scarray}
            R \subseteq [2m] \\
            |R| = 2m-2k \\
            R \, \hboxrm{\rm well-paired}
         \end{scarray}}
   \sum_{\sigma \in \scrs_{2m-2k}}
   \sgn(\sigma)
   \left( \prod_{p=1}^{2m-2k}C_{i_p,r'_{\sigma(p)}} \right)
   \left( \prod_{i \in R} \theta_i \right)
   \,.
\ee
We now multiply by $f(\theta,\chi)$ and integrate $ \scrd_{2m} (\theta,\chi)$.
As noted previously, the integral
\be
   \int \! \scrd_{2m} (\theta,\chi) \,
   \left( \prod_{i \in R} \theta_i \right)
   f(\theta,\chi)
\ee
is independent of the choice of $R$ 
(provided that it is well-paired and of cardinality $2m-2k$)
and hence can be pulled out.
In order to calculate the prefactor,
let us substitute the defining structure of well-paired sets $R$
by writing $r'_{2h-1} = 2j_h-1$ and $r'_{2h} = 2j_h$
for suitable indices $j_1 < \ldots < j_{m-k}$ in $[m]$.
We are left with
\be
   \sum_{\begin{scarray}
            1 \le j_1 < \ldots < j_{m-k} \le m \\
            r'_{2h-1} \equiv 2 j_h - 1 ,\; r'_{2h} \equiv 2 j_h 
         \end{scarray}}
   \sum_{\sigma \in \scrs_{2m-2k}} \sgn(\sigma) \,
   \prod_{p=1}^{2m-2k} C_{i_p,r'_{\sigma(p)}}
   \;.
 \label{proof.Alan.2.real}
\ee
Alternatively, we can sum over all {\em distinct}\/
$j_1, \ldots j_{m-k} \in [m]$
by inserting a factor $1/(m-k)!$.\footnote{
   Such a factor might not exist in a general commutative ring $R$.
   But we can argue as follows:
   First we prove the identity when $R = \R$;
   then we observe that since both sides of the identity are polynomials
   with integer coefficients in the matrix elements of $C$,
   the identity must hold {\em as polynomials}\/;
   therefore the identity holds when the matrix elements of $C$
   are specialized to arbitrary values in an arbitrary commutative ring $R$.
}
Remarkably, we can now remove the restriction that the
$j_1,\ldots,j_{m-k}$ be all distinct,
because the terms with two or more $j_\alpha$ equal cancel out
when we sum over permutations with the factor $\sgn(\sigma)$.\footnote{
   If $j_\alpha = j_\beta$ with $\alpha \neq \beta$,
   it suffices to use the involution
   exchanging $r'_{2\alpha-1}$ with $r'_{2\beta-1}$
   (or alternatively the involution exchanging
    $r'_{2\alpha}$ with $r'_{2\beta}$).
}
We can also reorder the factors $C$ freely.
Therefore \reff{proof.Alan.2.real} equals
\begin{eqnarray}
   &  &
   {1 \over (m-k)!}
   \!\!\!
   \sum_{\begin{scarray}
            j_1,\ldots,j_{m-k} \in [m] \\
            r'_{2h-1} \equiv 2 j_h - 1 ,\; r'_{2h} \equiv 2 j_h 
         \end{scarray}}
   \!
   \sum_{\sigma \in \scrs_{2m-2k}} \sgn(\sigma) \,
   \prod_{p=1}^{2m-2k}  C_{i_{\sigma^{-1}(p)}, r'_p}
       \nonumber \\
   &  &
   \quad =\;
   {1 \over (m-k)!}
   \!\!\!
   \sum_{\begin{scarray}
            j_1,\ldots,j_{m-k} \in [m] \\
            r'_{2h-1} \equiv 2 j_h - 1 ,\; r'_{2h} \equiv 2 j_h 
         \end{scarray}}
   \!
   \sum_{\sigma \in \scrs_{2m-2k}} \sgn(\sigma) \,
   \prod_{q=1}^{m-k}  C_{i_{\sigma^{-1}(2q-1)}, r'_{2q-1}} \,
                      C_{i_{\sigma^{-1}(2q)}, r'_{2q}}
       \nonumber \\
   &  &
   \qquad =\;
   {1 \over 2^{m-k} (m-k)!}
   \sum_{\begin{scarray}
             s_1,\ldots,s_{m-k} \in [2m] \\
             t_1,\ldots,t_{m-k} \in [2m]
         \end{scarray}}
   \sum_{\sigma \in \scrs_{2m-2k}} \sgn(\sigma) \,
   \prod_{q=1}^{m-k}  C_{i_{\sigma^{-1}(2q-1)}, s_q} \,
                      C_{i_{\sigma^{-1}(2q)}, t_q} \, J_{s_q t_q}
       \nonumber \\[1mm]
   &  &
   \qquad =\;
   \pf[ (C J C^{\rm T})_{I^c I^c} ]
   \;.
\end{eqnarray}
\qed

The following identity is what will be used in the proof of
the all-principal-minors antisymmetric Cayley identity
\reff{eq.antisymcayley.2}:

\begin{corollary}
  \label{cor.Atheta.Marco}
Let $I \subseteq [2m]$ with $|I|=2k$,
let $C$ be a $2m \times 2m$ matrix,
let $M$ be an invertible antisymmetric $2m \times 2m$ matrix,
and let $\ell$ be a nonnegative integer.
Then
\be
\label{eq.cor.Atheta.Marco}
\int \! \scrd_{2m}(\theta) \left(\prod (C\theta)\right)_{I^c}
   \, (\smhalf \theta^{\rm T} M \theta)^\ell
\;= \;
k!\,\delta_{\ell,k}\, (\pf M) \, \pf[(C M^{-\rm T} C^{\rm T} )_{I^c,I^c}]
   \;.
\ee
\end{corollary}

\firstproof
The case $M=J$ is an easy consequence of Lemma~\ref{lem.pbCJCp.Alan}.
For a general {\em real}\/ antisymmetric matrix $M$,
we can use the decomposition $M=AJA^{\rm T}$
(Lemma~\ref{lemma.decomp.1.FULL})
and the change of variables $\theta' = A^{\rm T} \theta$
(which picks up a Jacobian $\det A = \pf M$)
to reduce to the case $M=J$.

Finally, if $M$ is an antisymmetric matrix with coefficients
in an arbitrary commutative ring $R$,
we argue that both sides of \reff{eq.cor.Atheta.Marco}
are polynomials with integer coefficients in the matrix elements of $M$;
since they agree for all {\em real}\/ values of those matrix elements
(or even for real values in a nonempty open set),
they must agree {\em as polynomials}\/;
and this implies that they agree for arbitrary values of
the matrix elements in an arbitrary commutative ring $R$.
\qed

Here is a simple direct proof of Corollary~\ref{cor.Atheta.Marco}
that avoids the combinatorial complexity
of the full Lemma~\ref{lem.pbCJCp.Alan}
and instead relies on standard facts from the theory of
Grassmann--Berezin integration:

\secondproof
It is easy to see that the Grassmann integral is nonvanishing
only for $\ell = k$.
So in this case we can replace $(\smhalf \theta^{\rm T} M \theta)^k$
in the integrand by $k! \, \exp(\smhalf \theta^{\rm T} M \theta)$.
The claim is then an immediate consequence of
Wick's theorem for ``real'' fermions
(Theorem~\ref{thm.wick.realfermions}).
\qed

\subsection{Symmetric Cayley identity}   \label{sec.grassmann.2.symmetric}

With Corollary~\ref{cor.pbABp.Marco} in hand,
we can this time proceed directly to the proof of the all-minors identity.

\proofof{Theorem~\ref{thm.symcayley}}
%
Recall that the matrix $\partial^{\rm sym}$ is given by
\be
   (\partial^{\rm sym})_{ij}
    \;=\;  \cases{ \partial / \partial x_{ii}  & if $i=j$  \cr
                   \noalign{\vskip 4pt}
                   \smhalf \partial / \partial x_{ij}   & if $i < j$ \cr
                   \noalign{\vskip 4pt}
                   \smhalf \partial / \partial x_{ji}   & if $i > j$ \cr
                 }
\ee
As before, we introduce Grassmann variables
$\eta_i, \etabar_i$ ($1 \le i \le n$)
and use the representation
\begin{equation}
   \det(\partial^{\mathrm{sym}}_{IJ})
   \;=\;
   \epsilon(I,J)
   \int \! \scrd_n(\eta, \etabar) \,
     \Big(\prod \etabar \eta \Big)_{I^c,J^c} \;
     \exp\!\left[ \sum\limits_{i\leq j}
                  \smhalf (\etabar_i \eta_j + \etabar_j \eta_i)
                  \frac{\partial}{\partial x_{ij}}
           \right]
\ef.
\end{equation}
By the translation formula \reff{eq.translation}, we have
\begin{equation}
\det(\partial^{\mathrm{sym}}_{IJ}) \,
f(\{x_{ij}\}_{i \leq j})
\;=\;
\epsilon(I,J)
\int \! \scrd_n(\eta, \etabar) \,
\Big(\prod \etabar \eta \Big)_{I^c,J^c} \;
f ( \{ x_{ij}+ \smhalf (\etabar_i \eta_j - \eta_i \etabar_j) \}_{i \leq j}
  )
\end{equation}
for an arbitrary polynomial $f$.
We shall use this formula in the case
$f(\{x_{ij}\}_{i \leq j})=\det(X^{\rm sym})^s$
where $s$ is a positive integer, so that
\be
\det(\partial^{\mathrm{sym}}_{IJ}) \,
\det(X^{\rm sym})^s
\;=\;
 \epsilon(I,J)
\int \! \scrd_n(\eta, \etabar) \,
\Big( \prod \etabar \eta \Big)_{I^c,J^c} \;
\det\Bigl[
    X^{\rm sym}+ \smhalf (\etabar \eta^{\rm T}  - \eta \etabar^{\rm T} )
    \Bigr]^s
 \;.
   \label{eq.grassmann.sym.IJ}
\ee
It is convenient to introduce the shorthand
\be
(X^{\mathrm{sym}})^{\mathrm{trans}}
\;\equiv\;
X^{\mathrm{sym}} \,+\,  \smhalf (\etabar \eta^{\rm T}  - \eta \etabar^{\rm T} )
\ee
for the argument of $\det$.

Let us now change variables from $(\eta,\etabar)$ to
$(\eta', \etabar') \equiv (\eta, (X^{\rm sym})^{-1} \etabar)$,
with Jacobian $(\det X^{\rm sym})^{-1}$.
Dropping primes from the new variables, we observe that the expression for
the translated matrix can be written as
\begin{equation}
(X^{\mathrm{sym}})^{\mathrm{trans}}
\;=\;
X^{\mathrm{sym}}
\Big[I + \smhalf 
(\etabar \eta^{\rm T} - (X^{\mathrm{sym}})^{-1} \eta \etabar^{\rm T}
  X^{\mathrm{sym}} )\Big]
   \;,
\end{equation}
so that
\begin{equation}
\det (X^{\mathrm{sym}})^{\mathrm{trans}}
\;=\;
(\det X^{\mathrm{sym}}) \det\!\Big[I + \smhalf 
(\etabar \eta^{\rm T} - (X^{\mathrm{sym}})^{-1} \eta \etabar^{\rm T}
  X^{\mathrm{sym}} )\Big]
\ef.
\end{equation}
Applying Corollary~\ref{corol.lowranksym} to the rightmost determinant
yields
\be
\det\!\Big[I + \smhalf 
(\etabar \eta^{\rm T} - (X^{\mathrm{sym}})^{-1} \eta \etabar^{\rm T}
  X^{\mathrm{sym}} )\Big]
\;=\;
(1- \smfrac{1}{2} \etabar^{\rm T} \eta)^{-2}
\ee
so that we are left with the Grassmann-integral expression
\begin{eqnarray}
&&
\det(\partial^{\mathrm{sym}}_{IJ}) \,
(\det X^{\mathrm{sym}})^s
    \nonumber \\ 
& & \qquad = \;
(\det X^{\mathrm{sym}})^{s-1} \, \epsilon(I,J)
\int \! \scrd_n(\eta, \etabar) \;
\Big( \prod (X^{\mathrm{sym}} \etabar) (\eta) \Big)_{I^c,J^c} \,
(1- \smfrac{1}{2} \etabar^{\rm T} \eta)^{-2s}
\ef.
   \nonumber \\
  \label{eq.sym.111}
\end{eqnarray}
Now insert the expansion
\be
  (1-\smfrac{1}{2} \etabar^{\rm T} \eta)^{-2s}
  \;=\;
   \sum_{\ell=0}^\infty (-\smfrac{1}{2})^\ell \binom{-2s}{\ell}
        (\etabar^{\rm T} \eta)^\ell
 \label{eq.lowranksym2}
\ee
into \reff{eq.sym.111} and use Corollary~\ref{cor.pbABp.Marco}:
we obtain
\begin{equation}
\det(\partial^{\mathrm{sym}}_{IJ}) \,
(\det X^{\mathrm{sym}})^s
 \;=\;
(\det X^{\mathrm{sym}})^{s-1} \,
\epsilon(I,J) \, (\det X^{\mathrm{sym}}_{I^c J^c})
   \, k! \, (-\smfrac{1}{2})^k \binom{-2s}{k}
   \,.
 \label{eq.grassmann.sym.IJ.last}
\ee
This proves \reff{eq.symcayley.2}
when $X^{\rm sym}$ is an invertible real or complex
symmetric $n \times n$ matrix and $s$ is a positive integer;
the general validity of the identity then follows from
Proposition~\ref{prop.equiv.symbolic_s}.
\qed

{\bf Remark.}  A slight variant of this proof employs the
factorization $X^{\rm sym} = A A^{\rm T}$
when $X^{\rm sym}$ is a real symmetric positive-definite $n \times n$ matrix
(Lemma~\ref{lemma.cholesky});
we then use the change of variables
$(\eta', \etabar') \equiv (A^{-1} \eta, A^{-1} \etabar)$,
with Jacobian $(\det A)^{-2} = (\det X^{\rm sym})^{-1}$.
This slightly shortens the calculations,
but has the disadvantage that the proof is no longer purely algebraic/combinatorial,
because it invokes a decomposition
that is valid for {\em real symmetric positive-definite}\/ matrices
in order to define the needed change of variables
(which involves $A$ and not just $X^{\mathrm{sym}}$).
We therefore prefer to avoid matrix decompositions wherever we can
(which is unfortunately not always).
See also the Remark at the end of Section~\ref{sec.grassmann.2.antisymm},
and the discussion in Appendix~\ref{app.matrix.decomp}.

\subsection{Antisymmetric Cayley identity}  \label{sec.grassmann.2.antisymm}

With Corollary~\ref{cor.Atheta.Marco} in hand,
we can again proceed directly to the proof of the
all-(principal-)minors identity.

\proofof{Theorem~\ref{thm.antisymcayley}}
The $2m \times 2m$ matrix $\partial^{\rm antisym}$ is given by
\be
   (\partial^{\rm antisym})_{ij}
    \;=\;  \cases{ 0     & if $i=j$  \cr
                   \noalign{\vskip 4pt}
                   \partial / \partial x_{ij}   & if $i < j$ \cr
                   \noalign{\vskip 4pt}
                   -\, \partial / \partial x_{ji}   & if $i > j$ \cr
                 }
\ee
We introduce ``real'' Grassmann variables $\theta_i$ ($1 \le i \le 2m$)
and use the representation
\begin{equation}
\pf(\partial^{\mathrm{antisym}}_{II})
\;=\;
 \epsilon(I)
\int \! \scrd_{2m}(\theta) \,
\left( \prod \theta \right)_{I^c}  \,
\exp\!\left[ \sum\limits_{i < j} \theta_i \theta_j
                                 \frac{\partial}{\partial x_{ij}}
      \right]
\ef.
\end{equation}
By the translation formula \reff{eq.translation},
\begin{equation}
\pf(\partial^{\mathrm{antisym}}_{II}) \,
f(\{x_{ij}\}_{i < j})
\;=\;
 \epsilon(I)
\int \! \scrd_{2m}(\theta) \,
\left( \prod \theta \right)_{I^c}  \,
f ( \{ x_{ij} + \theta_i \theta_j \}_{i < j} )
\end{equation}
for an arbitrary polynomial $f$.
We shall use this formula in the case
$f(\{x_{ij}\}_{i < j})= (\pf X^{\rm antisym})^s$
where $s$ is a positive even integer.
It is convenient to introduce the shorthand
\begin{equation}
(X^{\mathrm{antisym}})^{\mathrm{trans}}
\;=\; X^{\mathrm{antisym}} + \theta \theta^{\rm T} 
\end{equation}
for the argument of $\pf$.

Suppose now that $X^{\rm antisym}$ is a $2m \times 2m$
real antisymmetric matrix of rank $2m$.
Then Lemma~\ref{lemma.decomp.1.FULL} guarantees that
we can find a matrix $A \in GL(2m)$
such that $X^{\mathrm{antisym}}=A J A^{\rm T}$, where
\be
J \;=\;
\left(
\begin{array}{cc|cc|c}
 0 & 1 & \multicolumn{3}{c}{} \\
-1 & 0 & \multicolumn{3}{c}{} \\
\cline{1-4}
 &  & 0 & 1 & \\
 &  & -1 & 0 & \\
\cline{3-4}
\multicolumn{4}{c}{} & \ddots
\end{array}
\right)
   \label{def.J}
\ee
is the standard $2m \times 2m$ symplectic form
(note that with this convention $\pf J = +1$ for any $m$).
%
%
We have
\begin{equation}
(X^{\mathrm{antisym}})^{\mathrm{trans}}
\;=\;
A (J + A^{-1} \theta \theta^{\rm T}  A^{-{\rm T}} ) A^{\rm T} 
\end{equation}
and
\begin{equation}
\pf (X^{\mathrm{antisym}})^{\mathrm{trans}}
\;=\;
(\pf X^{\mathrm{antisym}}) \pf ( J + A^{-1}
\theta \theta^{\rm T}  A^{-{\rm T}} )
\ef.
\end{equation}
Now change variables from $\theta$ to $\theta' \equiv A^{-1} \theta$,
with Jacobian $(\det A)^{-1}= (\pf X)^{-1}$.
Dropping primes, we are left with
\be
   \pf(\partial^{\rm antisym}_{II}) \, (\pf X^{\rm antisym})^s
   \;=\;
   \epsilon(I) \,
   (\pf X^{\rm antisym})^{s-1}
   \int \! \scrd_{2m}(\theta)
   \left( \prod (A\theta) \right)_{I^c}  \,
   \pf (J + \theta \theta^{\rm T})^s
\;.
\ee
We can now write
\be
   \pf (J + \theta \theta^{\rm T})^s
   \;=\;
   \det(J + \theta \theta^{\rm T})^{s/2}
   \;=\;
   \det(I  - J \theta \theta^{\rm T})^{s/2}
 \label{eq.pfs.detsover2}
\ee
since $\pf J = \det J = +1$ and $J^{-1} = -J$.\footnote{
   The equality
   $\pf (J + \theta \theta^{\rm T}) = \det(J + \theta \theta^{\rm T})^{1/2}$
   follows from the general fact $(\pf M)^2 = \det M$
   together with the observation that
   both $\pf (J + \theta \theta^{\rm T})$
   and $\det (J + \theta \theta^{\rm T})$
   are elements of the Grassmann algebra with constant term 1,
   and that an element of the Grassmann algebra with constant term 1
   has a {\em unique}\/ square root in the Grassmann algebra
   with this property.
}
The matrix $I  - J \theta \theta^{\rm T}$
is a rank-1 perturbation of the identity matrix;
applying Lemma~\ref{lemma.lowrank} with vectors
$u= -J \theta$ and $v = \theta$, we obtain
\begin{equation}
\det (I - J \theta \theta^{\rm T} )
\;=\;
(1- \theta^{\rm T} J \theta)^{-1}
\ef.
\end{equation}
Using Corollary~\ref{cor.Atheta.Marco}, we have
\begin{subeqnarray}
& &
   \int \! \scrd_{2m}(\theta)
      \left( \prod (A\theta) \right)_{I^c}  \,
      \pf (J + \theta \theta^{\rm T})^s
   \nonumber \\
& & \qquad\qquad =\;
   \int \! \scrd_{2m}(\theta)
      \left( \prod (A\theta) \right)_{I^c}  \,
      (1- \theta^{\rm T} J \theta)^{-s/2}
            \\
& & \qquad\qquad =\;
   (-2)^k k! \binom{-s/2}{k} (\pf J) \,
                             \pf[ (A J^{-\rm T} A^{\rm T})_{I^c I^c}]
            \nonumber \\
& & \qquad\qquad =\;
   s (s+2) \cdots (s+2m-2) \, \pf(X^{\mathrm{antisym}}_{I^c I^c})
 \label{eq.realfermion.scalartheta}
\end{subeqnarray}
since $\pf J = +1$
and $J^{-\rm T} = J$.
This proves \reff{eq.antisymcayley.2}
when $X^{\rm antisym}$ is a real antisymmetric matrix of rank $2m$
and $s$ is a positive even integer;
the general validity of the identity then follows from
Proposition~\ref{prop.equiv.symbolic_s}.
\qed

{\bf Remark.}  Although the statement of Theorem~\ref{thm.antisymcayley}
is purely algebraic/combinatorial, the foregoing proof is unfortunately
not {\em purely}\/ algebraic/combinatorial,
as it invokes the decomposition $X^{\mathrm{antisym}}=A J A^{\rm T}$
that is valid for {\em real}\/ full-rank antisymmetric matrices
in order to define the needed change of variables
(which involves $A$ and not just $X^{\mathrm{antisym}}$).
A similar decomposition will be invoked in the proofs
of the rectangular Cayley identities
(Sections~\ref{sec.grassmann.2.tmrect}, \ref{sec.grassmann.2.rect},
 \ref{sec.grassmann.2.antisymrect} and \ref{sec.grassmann.2.multi}).
This contrasts with the proofs of the ordinary and symmetric
Cayley identities
(Sections~\ref{sec.grassmann.2.ordinary} and \ref{sec.grassmann.2.symmetric}),
where we were able to define the needed change of variables
in terms of the original matrix $X$ or $X^{\mathrm{sym}}$.
The matrix factorization lemmas needed for the former proofs
are collected and proven in Appendix~\ref{app.matrix.decomp}.

\subsection{Two-matrix rectangular Cayley identity}
\label{sec.grassmann.2.tmrect}

In the remaining subsections of this section we shall prove the various
rectangular Cayley identities that were stated in
Section~\ref{sec.statement.rectangular}.
It is convenient to begin with the two-matrix rectangular Cayley identity
(Theorem~\ref{thm.tmrectcayley}),
whose proof is somewhat less intricate than that of the
corresponding one-matrix identities
(Theorems~\ref{thm.rectcayley} and \ref{thm.antisymrectcayley}).
Indeed, the one-matrix rectangular symmetric and antisymmetric Cayley identities
are related to the two-matrix identity
in roughly the same way as the symmetric and antisymmetric Cayley identities
are related to the ordinary one.

Here we will need the full strength of Lemma~\ref{lem.pbABp.Alan}
to handle the all-minors case.

\proofof{Theorem~\ref{thm.tmrectcayley}}
We begin once again by representing the differential operator
as a Grassmann integral: exploiting Corollary~\ref{cor.blockmatrixdet}, we have
\begin{equation}
\det (\partial_X \partial^{\rm T}_Y )   \;=\;
\det \left( \begin{array}{c|c}
0_m & \partial_X \\
\hline
-\partial^{\rm T}_Y  & I_n
\end{array} \right)
\;=\;
\int \! \scrd_m(\psi,\psibar) \, \scrd_n(\eta,\etabar)
\,
e^{\etabar^{\rm T} \! \eta
                  + \psibar^{\rm T} \! \partial_X \eta
                  + \psi^{\rm T} \! \partial_Y \etabar }
\ef.
\end{equation}
Here $\psi_i, \psibar_i$ ($1 \le i \le m$) and
$\eta_j, \etabar_j$ ($1 \le j \le n$) are Grassmann variables,
and the subscripts on $\scrd$ serve to remind us
of the length of each vector;
shorthand notations for index summations are understood, e.g.~%
$\psibar^{\rm T} \partial_X \eta \equiv \sum_{i=1}^m \sum_{j=1}^n
\psibar_i \eta_j \partial/\partial x_{ij}$.
For a general minor $I,J \subseteq [m]$ with $|I|=|J|=k$,
we have (writing $L = \{m+1,\ldots,m+n\}$)
\begin{subeqnarray}
\det[(\partial_X \partial^{\rm T}_Y )_{IJ}]
&=&
\det \left[ \left( \begin{array}{c|c}
0_m & \partial_X \\
\hline
-\partial_Y^{\rm T}  & I_n
\end{array} \right)_{I \cup L, J \cup L} \right]
\\ &=&
\epsilon(I,J) \int \! \scrd_m(\psi,\psibar) \, \scrd_n(\eta,\etabar)
\,
\Big( \prod \psibar \psi \Big)_{I^c,J^c}
e^{\etabar^{\rm T} \! \eta + \psibar^{\rm T} \! \partial_X \eta +
                          \psi^{\rm T} \! \partial_Y \etabar}
\ef.
   \nonumber \\
\end{subeqnarray}

Applying the translation formula \reff{eq.translation}
to the whole set of variables $\{x_{ij}, y_{ij} \}$ produces
\begin{equation}
\label{eq.76595.newNEW}
\det[(\partial_X \partial^{\rm T}_Y )_{IJ}]  \, f(X, Y)  \;=\;
\epsilon(I,J) \int \! \scrd_m(\psi,\psibar) \, \scrd_n(\eta,\etabar)
   \,
\Big( \prod \psibar \psi \Big)_{I^c,J^c}
e^{\etabar^{\rm T} \! \eta} \,
f(X + \psibar \eta^{\rm T}, Y + \psi \etabar^{\rm T})
\end{equation}
for an arbitrary polynomial $f$.
We shall use this formula in the case
$f(X, Y)=\det(X Y^{\rm T})^s$
where $s$ is a positive integer.
It is convenient to introduce the shorthands
\begin{subeqnarray}
X^{\rm trans} & \equiv & X + \psibar \eta^{\rm T}  \\
Y^{\rm trans} & \equiv & Y + \psi \etabar^{\rm T}
\end{subeqnarray}
for the arguments of $f$.

Suppose now that $X$ and $Y$ are real $m \times n$ matrices of rank $m$
that are sufficiently close to the matrix $\widehat{I}_{mn}$ defined by
\be
   (\widehat{I}_{mn})_{ij} 
   \;=\;
   \cases{  1  &  if $i=j$  \cr
            \noalign{\vskip 6pt}
            0  &  if $i \neq j$ \cr
         }
  \label{def.Ihat}
\ee
[Note that $\widehat{I}_{mn} \equiv (I_m,0_{m \times (n-m)})$ when $m \le n$,
 and $\widehat{I}_{mn} = \widehat{I}_{nm}^{\rm T}$ otherwise.
 Henceforth we shall drop the subscripts $mn$ on $\widehat{I}_{mn}$
 to lighten the notation.]
Then by Lemma~\ref{lemma.decomp.2}
we can find matrices $P, P' \in GL(m)$ and $Q \in GL(n)$
such that $X=P \widehat{I} Q$ and $Y=P' \widehat{I} Q^{-\rm T}$.
We have
\begin{equation}
\det(X Y^{\rm T} ) \;=\; \det
(P \widehat{I} Q Q^{-1}  \widehat{I}^{\rm T} {P'}^{\rm T} )
   \;=\; \det(P {P'}^{\rm T} ) \;=\; \det(P) \det(P')
\end{equation}
and
\begin{equation}
X^{\rm trans} (Y^{\rm trans})^{\rm T}
\;\equiv\;
(X + \psibar \eta^{\rm T})( Y^{\rm T} - \etabar \psi^{\rm T} )
  \;=\;
P [\widehat{I}+P^{-1}(\psibar \eta^{\rm T})Q^{-1} ] Q Q^{-1}
[\widehat{I}^{\rm T} - Q(\etabar \psi^{\rm T} ) {P'}^{-{\rm T}}]
{P'}^{\rm T}
\ef.
\end{equation}
Let us now change variables from $(\psi, \psibar, \eta, \etabar)$
to $(\psi', \psibar', \eta', \etabar') \equiv
    ({P'}^{-1} \psi, P^{-1} \psibar, Q^{-\rm T} \eta, Q \etabar)$,
with Jacobian $(\det P)^{-1} (\det P')^{-1} =\det(X Y^{\rm T} )^{-1}$.
In the new variables we have
(dropping now the primes from the notation)
\begin{equation}
X^{\rm trans} (Y^{\rm trans})^{\rm T}
  \;=\;
P (\widehat{I}+ \psibar \eta^{\rm T} )
(\widehat{I}^{\rm T} - \etabar \psi^{\rm T} )
{P'}^{\rm T}
\ef,
\end{equation}
and the translated determinant is given by
\begin{equation}
\det[(X^{\rm trans}) (Y^{\rm trans})^{\rm T} ] =
\det(XY^{\rm T} )
\det[(\widehat{I}+ \psibar \eta^{\rm T} )
(\widehat{I}^{\rm T} - \etabar \psi^{\rm T} )]
\ef,
\end{equation}
so that
\begin{eqnarray}
   & &
   \det[(\partial_X \partial^{\rm T}_Y )_{IJ}]  \,  \det(XY^{\rm T} )^s
       \nonumber \\
   & &
   \qquad =\;
   \epsilon(I,J) \, \det(XY^{\rm T} )^{s-1}
   \int \! \scrd_m(\psi,\psibar) \, \scrd_n(\eta,\etabar)
   \,
   \Big( \prod (P\psibar) (P' \psi) \Big)_{I^c,J^c}
   \, \times
        \nonumber \\
   & & \qquad\qquad\qquad
   e^{\etabar^{\rm T} \! \eta} \,
   \det[(\widehat{I}+ \psibar \eta^{\rm T} )
     (\widehat{I}^{\rm T} - \etabar \psi^{\rm T} )]^s
   \;.
\end{eqnarray}

Let us now split the vectors $\eta$ and $\etabar$ as
\begin{subeqnarray}
(\eta_1, \ldots, \eta_n) & = &
(\lambda_1, \ldots, \lambda_m , \chi_1, \ldots, \chi_{n-m})
\\
(\etabar_1, \ldots, \etabar_n)  & = &
(\lambdabar_1, \ldots, \lambdabar_m , \chibar_1, \ldots, \chibar_{n-m})
\end{subeqnarray}
so that
\begin{equation}
(\widehat{I}+ \psibar \eta^{\rm T} )
     (\widehat{I}^{\rm T} - \etabar \psi^{\rm T} )
 \;=\;
I_m +\psibar \lambda^{\rm T} - \lambdabar \psi^{\rm T}
+c \psibar \psi^{\rm T}
\end{equation}
with $c= \lambdabar^{\rm T} \lambda + \chibar^{\rm T} \chi$.
This matrix has the form of a low-rank perturbation
$I_m+\sum_{\alpha=1}^{2} u_{\alpha} v_{\alpha}^{\rm T} $,
with vectors $\{u_{\alpha}\}$, $\{v_{\alpha}\}$ given by
\begin{equation}
\begin{array}{c||c|c}
\alpha & u_{\alpha} & v_{\alpha}
\\
\hline
\rule{0pt}{4.5mm}
1 & \psibar & \lambda + c \psi \\
2 & \lambdabar & -\psi \\
\end{array}
   \label{eq.tmrectcayleyproof.rank2}
\end{equation}
By Lemma~\ref{lemma.lowrank} we can write the needed determinant
as the determinant of a $2 \times 2$ matrix:
\begin{subeqnarray}
& &
   \det[(\widehat{I}+ \psibar \eta^{\rm T} )
     (\widehat{I}^{\rm T} - \etabar \psi^{\rm T} )]
        \nonumber \\[2mm]
&  &  \qquad\qquad =\;
   \det{}^{-1}
   \left(
   \begin{array}{cc}
   \rule[-2.7mm]{0pt}{4.7mm}
  1 + (\lambda + c \psi)^{\rm T} \psibar &
                 (\lambda + c \psi)^{\rm T} \lambdabar \\
  - \psi^{\rm T} \psibar & 1 - \psi^{\rm T} \lambdabar
   \end{array}
   \right)
          \\[2mm]
&  &  \qquad\qquad =\;
   \det{}^{-1}
   \left(
   \begin{array}{cc}
   \rule[-2.7mm]{0pt}{4.7mm}
  1 + \lambda^{\rm T} \psibar & c + \lambda^{\rm T} \lambdabar \\
  - \psi^{\rm T} \psibar & 1 - \psi^{\rm T} \lambdabar
   \end{array}
   \right)
          \\[2mm]
&  &  \qquad\qquad =\;
 [ (1 - \psibar^{\rm T} \lambda) (1 + \lambdabar^{\rm T} \psi)
  - (\chibar^{\rm T} \chi) (\psibar^{\rm T} \psi) ]^{-1}
  \;,
\end{subeqnarray}
where the second equality is a row operation
(row 1 $\to$ row 1 + $c$ row 2).
We therefore have
\begin{eqnarray}
   & &
   \!\!\!\!\!\!
   \det[(\partial_X \partial^{\rm T}_Y )_{IJ}]  \,  \det(XY^{\rm T} )^s
   \;=\;
   \epsilon(I,J) \, \det(XY^{\rm T} )^{s-1}
   \int \! \scrd_{n-m}(\chi,\chibar) \; e^{\chibar^{\rm T} \! \chi}
      \,\times
      \nonumber \\
   & &
   \!\!
   \int \! \scrd_m(\psi,\psibar) \, \scrd_m(\lambda,\lambdabar)
   \;
   e^{\lambdabar^{\rm T} \! \lambda} \,
   \Big( \prod (P\psibar) (P' \psi) \Big)_{I^c,J^c}  \,
[(1 - \psibar^{\rm T} \lambda) (1 + \lambdabar^{\rm T} \psi)
- (\chibar^{\rm T} \chi) (\psibar^{\rm T} \psi)]^{-s}
   \;.
      \nonumber \\
\end{eqnarray}
Note in particular that the integrand
depends on $\psi, \psibar, \lambda, \lambdabar$
only via scalar products.
This allows us to apply Lemma~\ref{lem.pbABp.Alan}
to the integral over $\psi, \psibar, \lambda, \lambdabar$;
using also the fact that $P(P')^{\rm T} = XY^{\rm T} $, we obtain
\begin{eqnarray}
& & \!\!\!\!
\det[(\partial_X \partial^{\rm T}_Y )_{IJ}]  \, \det(XY^{\rm T} )^s
  \;=\;
  \epsilon(I,J) \, \det(XY^{\rm T} )^{s-1} \,
  \det[(XY^{\rm T})_{I^c J^c}]
    \,\times
    \nonumber \\
& &
\qquad
  \int \! \scrd_{n-m}(\chi,\chibar) \; e^{\chibar^{\rm T} \! \chi}
  \int \! \scrd_m(\psi,\psibar) \, \scrd_m(\lambda,\lambdabar) \;
     e^{\lambdabar^{\rm T} \! \lambda}
\left( \prod_{a=k+1}^{m} \psibar_a \psi_a \!\right)
    \,\times
    \nonumber \\
& &
\qquad\qquad
[(1 - \psibar^{\rm T} \lambda) (1 + \lambdabar^{\rm T} \psi)
- (\chibar^{\rm T} \chi) (\psibar^{\rm T} \psi)]^{-s}
  \;.
\label{eq.76595pippoNEW}
\end{eqnarray}
This formula expresses
$\det[(\partial_X \partial^{\rm T}_Y)_{IJ}]  \det(XY^{\rm T} )^s$
as the desired quantity
$\det(XY^{\rm T} )^{s-1} \times$ $\epsilon(I,J) \det[(XY^{\rm T} )_{I^c J^c}]$
multiplied by the purely combinatorial factor
\begin{eqnarray}
& &
P^{\textrm{\scriptsize tmrect}}(s,m,n,k)
 \;\equiv\;
\int \! \scrd_m(\psi,\psibar) \, \scrd_m(\lambda,\lambdabar) \,
  \scrd_{n-m}(\chi,\chibar)
  \,\times
     \nonumber \\
& & \qquad
 e^{\lambdabar^{\rm T} \! \lambda + \chibar^{\rm T} \! \chi}
\left( \prod_{a=k+1}^{m} \psibar_a \psi_a \right)
[(1 - \psibar^{\rm T} \lambda) (1 + \lambdabar^{\rm T} \psi)
- (\chibar^{\rm T} \chi) (\psibar^{\rm T} \psi)]^{-s}
    \;,  \qquad
\label{eq.834653865NEW}
\end{eqnarray}
which we now proceed to calculate.

First note that the factor $\prod_{a=k+1}^{m} \psibar_a \psi_a$
forces the Taylor expansion of the square bracket
to contain no
variables $\psi_a$, $\psibar_a$ with $k+1 \le a \le m$.
We can therefore drop the factor $\prod_{a=k+1}^{m} \psibar_a \psi_a$,
forget about the variables $\psi_a$, $\psibar_a$ with $k+1 \le a \le m$,
and consider $\psi, \psibar$ henceforth as vectors of length $k$.
Let us also rename the vectors $\lambda$, $\lambdabar$ as
\begin{subeqnarray}
(\lambda_1, \ldots, \lambda_m) & = &
(\lambda_1, \ldots, \lambda_{k}, \mu_1, \ldots, \mu_{m-k})
   \\
(\lambdabar_1, \ldots, \lambdabar_m)  & = &
(\lambdabar_1, \ldots, \lambdabar_{k}, \mubar_1, \ldots, \mubar_{m-k})
\end{subeqnarray}
Then the quantity \reff{eq.834653865NEW} is equivalent to
\begin{eqnarray}
\label{eq.834653865bisNEW}
& &
P^{\textrm{\scriptsize tmrect}}(s,m,n,k)
\;=\;
\int \! \scrd_k(\psi,\psibar) \, \scrd_k(\lambda,\lambdabar) \,
\scrd_{m-k}(\mu,\mubar) \,
  \scrd_{n-m}(\chi,\chibar) \,\times
      \nonumber \\
& & \qquad\qquad
  e^{\lambdabar^{\rm T} \! \lambda + \mubar^{\rm T} \! \mu
                                   + \chibar^{\rm T} \! \chi} \,
[(1 - \psibar^{\rm T} \lambda) (1 + \lambdabar^{\rm T} \psi)
- (\chibar^{\rm T} \chi) (\psibar^{\rm T} \psi)]^{-s}
\end{eqnarray}
where scalar products involving $\psi, \psibar, \lambda, \lambdabar$
are understood as referring only to the first $k$ variables.
Integration over the variables $\mu, \mubar$ is trivial
and produces just a factor 1. So we are left with
\begin{eqnarray}
& &
P^{\textrm{\scriptsize tmrect}}(s,m,n,k)
\;=\;
\int \! \scrd_k(\psi,\psibar) \, \scrd_k(\lambda,\lambdabar) \,
  \scrd_{n-m}(\chi,\chibar) \,\times
      \nonumber \\
& & \qquad\qquad
  e^{\lambdabar^{\rm T} \! \lambda + \chibar^{\rm T} \! \chi} \,
[(1 - \psibar^{\rm T} \lambda) (1 + \lambdabar^{\rm T} \psi)
- (\chibar^{\rm T} \chi) (\psibar^{\rm T} \psi)]^{-s}
   \;.
 \label{eq.834653865trisNEW}
\end{eqnarray}
Note that $P^{\textrm{\scriptsize tmrect}}(s,m,n,k)$ depends on $n$ and
$m$ only via the combination $n-m$.

The binomial expansion of the integrand in \reff{eq.834653865trisNEW} yields
\be
[(1 - \psibar^{\rm T} \lambda) (1 + \lambdabar^{\rm T} \psi)
- (\chibar^{\rm T} \chi) (\psibar^{\rm T} \psi)]^{-s}
\;=\;
\sum_{h=0}^\infty \!
\binom{-s}{h} \!
[(1 - \psibar^{\rm T} \lambda) (1 + \lambdabar^{\rm T} \psi)]^{-s-h} \,
[- (\chibar^{\rm T} \chi) (\psibar^{\rm T} \psi)]^{h}
  \;.
\ee
For each fixed value of $h$,
integration over variables $\chi$, $\chibar$ gives
\be
\int \! \scrd_{n-m}(\chi,\chibar) \;
  e^{\chibar^{\rm T} \! \chi} \,
(- \chibar^{\rm T} \chi)^{h}
\;=\; (-1)^h \binom{n-m}{h} h!
 \;.
\ee
Now the factor $(\psibar^{\rm T} \psi)^h$ will contain $h$ pairs of variables,
which can be chosen in $\binom{k}{h} h!$ ways (counting reorderings);
then in the term
$f(\psi,\psibar,\lambda,\lambdabar) =
 [(1 - \psibar^{\rm T} \lambda) (1 + \lambdabar^{\rm T} \psi)]^{-s-h}$
we will have to use only the other $k-h$ pairs of variables
for both $\psi,\psibar$ and $\lambda,\lambdabar$, 
and with a reasoning as in Lemma~\ref{lem.pbABp.Alan} we can choose them
to be the first $k-h$ indices,
so that
\begin{eqnarray}
& &
\int \! \scrd_k(\psi,\psibar) \, \scrd_k(\lambda,\lambdabar) \,
  e^{\lambdabar^{\rm T} \! \lambda} \,
(\psibar^{\rm T} \psi)^{h} f(\psi,\psibar,\lambda,\lambdabar)
     \nonumber \\
& & \qquad\qquad
=\;
\binom{k}{h} h!
\int \! \scrd_{k-h}(\psi,\psibar) \, \scrd_{k-h}(\lambda,\lambdabar) \,
  e^{\lambdabar^{\rm T} \! \lambda} \,
f(\psi,\psibar,\lambda,\lambdabar)
  \;,
\end{eqnarray}
and hence
\begin{eqnarray}
& &
P^{\textrm{\scriptsize tmrect}}(s,m,n,k)
\;=\;
\sum_{h=0}^\infty
(-1)^h
\binom{n-m}{h} \binom{k}{h} (h!)^2
\binom{-s}{h} \,\times
     \nonumber \\
& & \qquad\qquad
\int \! \scrd_{k-h}(\psi,\psibar) \, \scrd_{k-h}(\lambda,\lambdabar) \,
  e^{\lambdabar^{\rm T} \! \lambda} \,
[(1 - \psibar^{\rm T} \lambda) (1 + \lambdabar^{\rm T} \psi)]^{-s-h}
   \;.
   \qquad
\end{eqnarray}
Let us now make the change of variables $\psi_i \to \lambda_i$,
$\lambda_i \to - \psi_i$ (whose Jacobian is cancelled by the
reordering in the measure of integration);  we have
\begin{subeqnarray}
& &
\int \! \scrd_{k-h}(\psi,\psibar) \, \scrd_{k-h}(\lambda,\lambdabar) \,
  e^{\lambdabar^{\rm T} \! \lambda} \,
[(1 - \psibar^{\rm T} \lambda) (1 + \lambdabar^{\rm T} \psi)]^{-s-h}
       \\[2mm]
& & \qquad
  =\;
\int \! \scrd_{k-h}(\psi,\psibar) \, \scrd_{k-h}(\lambda,\lambdabar) \,
  e^{-\lambdabar^{\rm T} \! \psi} \,
[(1 + \psibar^{\rm T} \psi) (1 + \lambdabar^{\rm T} \lambda)]^{-s-h}
   \qquad
       \\[2mm]
& & \qquad
  =\;
\int \! \scrd_{k-h}(\psi,\psibar) \, \scrd_{k-h}(\lambda,\lambdabar) \,
[(1 + \psibar^{\rm T} \psi) (1 + \lambdabar^{\rm T} \lambda)]^{-s-h}
       \\[2mm]
& & \qquad
  =\;
\left[ \int \! \scrd_{k-h}(\psi,\psibar) \;
(1 + \psibar^{\rm T} \psi)^{-s-h} \right]^{\! 2}
       \\[2mm]
& & \qquad
  =\;
\left[ \binom{-s-h}{k-h} (k-h)! \right]^{\! 2}
   \;.
\end{subeqnarray}
Collecting all the factors, we have
\begin{subeqnarray}
& &
P^{\textrm{\scriptsize tmrect}}(s,m,n,k)
\;=\;
\sum_{h=0}^\infty
(-1)^h
\binom{n-m}{h} \binom{k}{h} (h!)^2
\binom{-s}{h}
\left[ \binom{-s-h}{k-h} (k-h)! \right]^2
    \nonumber \\  \\
& & \qquad
  =\;
\sum_{h=0}^\infty
(-1)^h
\binom{n-m}{h} (k!)^2 \binom{-s-h}{k-h}
\binom{-s}{k}
       \\[2mm]
& & \qquad
  =\;
\sum_{h=0}^\infty
(-1)^h
\binom{n-m}{h} (k!)^2 \, (-1)^{k-h} \binom{s+k-1}{k-h}
(-1)^k \binom{s+k-1}{k}
       \\[2mm]
& & \qquad
  =\;
(k!)^2 \binom{s+k-1}{k}
\sum_{h=0}^\infty
\binom{n-m}{h} \binom{s+k-1}{k-h}
       \\[2mm]
& & \qquad
  =\;
(k!)^2 \binom{s+k-1}{k} \binom{s+k+n-m-1}{k}
       \\[2mm]
& & \qquad
  =\;
  \prod_{j=0}^{k-1}  (s+j)(s+n-m+j)  \;,
\end{subeqnarray}
where the sum over $h$ was performed using the Chu--Vandermonde convolution
(Lemma~\ref{lemma.vandermonde}).

This proves Theorem~\ref{thm.tmrectcayley} when
$X$ and $Y$ are real $m \times n$ matrices of rank $m$
lying in a sufficiently small neighborhood of $\widehat{I}_{mn}$,
and $s$ is a positive integer.
The general validity of the identity then follows from
Proposition~\ref{prop.equiv.symbolic_s}.
\qed

\subsection{One-matrix rectangular symmetric Cayley identity}
    \label{sec.grassmann.2.rect}

The proof of the one-matrix rectangular symmetric Cayley identity
is extremely similar to that of the two-matrix identity,
but is slightly more complicated because it involves
a perturbation of rank 4 rather than rank 2
[compare \reff{eq.rectcayleyproof.rank4}
 with \reff{eq.tmrectcayleyproof.rank2}].
Luckily, the resulting $4 \times 4$ determinant turns out to be
the square of a quantity involving only a $2 \times 2$ determinant
[cf.\ \reff{eq.perfectsquare.1}/\reff{eq.perfectsquare.2}].
Once again, we will need the full strength of Lemma~\ref{lem.pbABp.Alan}
to handle the all-minors case.

\proofof{Theorem~\ref{thm.rectcayley}}
We begin once again by representing the differential operator
as a Grassmann integral: exploiting Corollary~\ref{cor.blockmatrixdet}, we have
\begin{equation}
\det (\partial \partial^{\rm T} )   \;=\;
\det \left( \begin{array}{c|c}
0_m & \partial \\
\hline
-\partial^{\rm T}  & I_n
\end{array} \right)
\;=\;
\int \! \scrd_m(\psi,\psibar) \, \scrd_n(\eta,\etabar)
\,
e^{\etabar^{\rm T} \! \eta + \psi^{\rm T} \! \partial \etabar +
                          \psibar^{\rm T} \! \partial \eta}
\ef.
\end{equation}
Here $\psi_i, \psibar_i$ ($1 \le i \le m$) and
$\eta_j, \etabar_j$ ($1 \le j \le n$) are Grassmann variables,
and we use the same conventions as in the preceding subsection.
For a general minor $I,J \subseteq [m]$ with $|I|=|J|=k$,
we have (writing $L = \{m+1,\ldots,m+n\}$)
\begin{subeqnarray}
\det[(\partial \partial^{\rm T} )_{IJ}]   
&=&
\det \left[ \left( \begin{array}{c|c}
0_m & \partial \\
\hline
-\partial^{\rm T}  & I_n
\end{array} \right)_{I \cup L, J \cup L} \right]
\\ &=&
\epsilon(I,J) \int \! \scrd_m(\psi,\psibar) \, \scrd_n(\eta,\etabar)
\,
\Big( \prod \psibar \psi \Big)_{I^c,J^c}
e^{\etabar^{\rm T} \! \eta + \psi^{\rm T} \! \partial \etabar +
                          \psibar^{\rm T} \! \partial \eta}
\ef.
   \nonumber \\
\end{subeqnarray}

Applying the translation formula \reff{eq.translation}
to the whole set of variables $\{x_{ij}\}$ produces
\begin{equation}
\label{eq.76595.new}
\det[(\partial \partial^{\rm T} )_{IJ}]  \, f(X)  \;=\;
\epsilon(I,J) \int \! \scrd_m(\psi,\psibar) \, \scrd_n(\eta,\etabar)
   \, 
\Big( \prod \psibar \psi \Big)_{I^c,J^c}
e^{\etabar^{\rm T} \! \eta} \,
f(X + \psibar \eta^{\rm T}  + \psi \etabar^{\rm T} )
\end{equation}
for an arbitrary polynomial $f$.
We shall use this formula in the case $f(X)=\det(XX^{\rm T} )^s$
where $s$ is a positive integer.
It is convenient to introduce the shorthand
\be
X^{\rm trans} \;\equiv\; X + \psibar \eta^{\rm T}  + \psi \etabar^{\rm T} 
\ee
for the argument of $f$.

Suppose now that $X$ is a real $m \times n$ matrix of rank $m$
that is sufficiently close to the matrix $\widehat{I}_{mn}$
defined in \reff{def.Ihat}.
Then by Lemma~\ref{lemma.decomp.3}
we can find matrices $P \in GL(m)$ and $Q \in O(n)$
such that $X=P \widehat{I} Q$
[we drop the subscripts $mn$ on $\widehat{I}_{mn}$ to lighten the notation].
We have
\begin{equation}
\det(X X^{\rm T} ) \;=\; \det (P \widehat{I} Q Q^{\rm T}  \widehat{I}^{\rm T}  P^{\rm T} )
   \;=\; \det(PP^{\rm T} ) \;=\; \det(P)^2
\end{equation}
and
\begin{equation}
X^{\rm trans} \;\equiv\; X + \psibar \eta^{\rm T}  + \psi \etabar^{\rm T} 
  \;=\;
P [\widehat{I}+P^{-1}(\psibar \eta^{\rm T}  + \psi \etabar^{\rm T} )Q^{\rm T} ] Q
\ef.
\end{equation}
Let us now change variables from $(\psi, \psibar, \eta, \etabar)$
to $(\psi', \psibar', \eta', \etabar') \equiv
    (P^{-1} \psi, P^{-1} \psibar, Q \eta, Q \etabar)$,
with Jacobian $(\det P)^{-2}=\det(XX^{\rm T} )^{-1}$.
In the new variables we have
(dropping now the primes from the notation)
\begin{equation}
X^{\rm trans} =
P (\widehat{I}+ \psibar \eta^{\rm T}  + \psi \etabar^{\rm T} ) Q
\ef,
\end{equation}
and the translated determinant is given by
\begin{equation}
\det[(X^{\rm trans}) (X^{\rm trans})^{\rm T} ] =
\det(XX^{\rm T} )
\det[(\widehat{I}+ \psibar \eta^{\rm T}  + \psi \etabar^{\rm T} )
     (\widehat{I}^{\rm T} -\eta \psibar^{\rm T}  -\etabar \psi^{\rm T} )]
\ef,
\end{equation}
so that
\begin{eqnarray}
   & &
   \det[(\partial \partial^{\rm T} )_{IJ}]  \,  \det(XX^{\rm T} )^s
       \nonumber \\
   & &
   \qquad =\;
   \epsilon(I,J) \, \det(XX^{\rm T} )^{s-1}
   \int \! \scrd_m(\psi,\psibar) \, \scrd_n(\eta,\etabar)
   \,
   \Big( \prod (P\psibar) (P\psi) \Big)_{I^c,J^c}
   \, \times
        \nonumber \\
   & & \qquad\qquad\qquad
   e^{\etabar^{\rm T} \! \eta} \,
   \det[(\widehat{I}+ \psibar \eta^{\rm T}  + \psi \etabar^{\rm T} )
     (\widehat{I}^{\rm T} -\eta \psibar^{\rm T}  -\etabar \psi^{\rm T} )]^s
   \;.
\end{eqnarray}    
   
Let us now split the vectors $\eta$ and $\etabar$ as
\begin{subeqnarray}
(\eta_1, \ldots, \eta_n) & = &
(\lambda_1, \ldots, \lambda_m , \chi_1, \ldots, \chi_{n-m})
\\
(\etabar_1, \ldots, \etabar_n)  & = &
(\lambdabar_1, \ldots, \lambdabar_m , \chibar_1, \ldots, \chibar_{n-m})
\end{subeqnarray}
so that
\begin{equation}
(\widehat{I}+ \psibar \eta^{\rm T}  + \psi \etabar^{\rm T} )
(\widehat{I}^{\rm T} -\eta \psibar^{\rm T}  -\etabar \psi^{\rm T} )
 \;=\;
I_m +\psibar \lambda^{\rm T}  + \psi \lambdabar^{\rm T} 
- \lambda \psibar^{\rm T}  - \lambdabar \psi^{\rm T} 
+c \psibar \psi^{\rm T}  -c \psi \psibar^{\rm T} 
\end{equation}
with $c= \lambdabar^{\rm T} \lambda + \chibar^{\rm T} \chi$.
This matrix has the form of a low-rank perturbation
$I_m +\sum_{\alpha=1}^{4} u_{\alpha} v_{\alpha}^{\rm T} $,
with vectors $\{u_{\alpha}\}$, $\{v_{\alpha}\}$ given by
\begin{equation}
\begin{array}{c||c|c}
\alpha & u_{\alpha} & v_{\alpha}
\\
\hline
\rule{0pt}{4.5mm}
1 & \psibar & \lambda + c \psi \\
2 & \psi & \lambdabar - c \psibar \\
3 & \lambdabar & -\psi \\
4 & \lambda & -\psibar
\end{array}
   \label{eq.rectcayleyproof.rank4}
\end{equation}
By Lemma~\ref{lemma.lowrank} we can write the needed determinant
as the determinant of a $4 \times 4$ matrix;
after a few row and column manipulations we can write
\begin{subeqnarray}
& &
   \det[(\widehat{I}+ \psibar \eta^{\rm T}  + \psi \etabar^{\rm T} )
        (\widehat{I}^{\rm T} -\eta \psibar^{\rm T}  -\etabar \psi^{\rm T} )]
        \nonumber \\[2mm]
&  &  \qquad\qquad =\;
   \det{}^{-1}
   \left(
   \begin{array}{c|c}
   A &
   \rule[-2.7mm]{0pt}{4.7mm}
   \begin{array}{cc}
   0 & \chibar^{\rm T} \chi \\
   -\chibar^{\rm T} \chi & 0
   \end{array}
   \\
   \hline
   \rule{0pt}{6.8mm}
   \begin{array}{cc}
   0 & -\psibar^{\rm T} \psi \\
   \psibar^{\rm T} \psi & 0
   \end{array}
   &
   A^{\rm T} 
   \end{array}
   \right)
          \\[2mm]
&  &  \qquad\qquad =\;
   [\det A - (\psibar^{\rm T} \psi) (\chibar^{\rm T} \chi)]^{-2}
 \label{eq.perfectsquare.1}
\end{subeqnarray}
where
\be
   A \;=\;
   \left(\!
   \begin{array}{cc}
      1+ \lambda^{\rm T} \psibar & \lambda^{\rm T} \psi \\[1mm]
      \lambdabar^{\rm T} \psibar & 1+ \lambdabar^{\rm T} \psi
   \end{array}
   \!\right)
 \label{eq.perfectsquare.2}
\ee
and hence
\be
\det A  \;=\;  1+ \lambda^{\rm T} \psibar + \lambdabar^{\rm T} \psi +
(\lambda^{\rm T} \psibar) (\lambdabar^{\rm T} \psi) +
(\lambdabar^{\rm T} \psibar) (\psi^{\rm T} \lambda)
\ef.
  \label{eq.detA}
\ee
We therefore have
\begin{eqnarray}
   & &
   \!\!\!\!
   \det[(\partial \partial^{\rm T} )_{IJ}]  \,  \det(XX^{\rm T} )^s
   \;=\;
   \epsilon(I,J) \, \det(XX^{\rm T} )^{s-1}
   \int \! \scrd_{n-m}(\chi,\chibar) \; e^{\chibar^{\rm T} \! \chi}
      \,\times
      \nonumber \\
   & &
   \;
   \int \! \scrd_m(\psi,\psibar) \, \scrd_m(\lambda,\lambdabar)
   \;
   e^{\lambdabar^{\rm T} \! \lambda} \,
   \Big( \prod (P\psibar) (P\psi) \Big)_{I^c,J^c}  \,
   [\det A - (\psibar^{\rm T} \psi) (\chibar^{\rm T} \chi)]^{-2s}
   \;.
   \qquad\quad
\end{eqnarray}    
Note in particular that $e^{\lambdabar^{\rm T} \! \lambda}$
and $\det A - (\psibar^{\rm T} \psi) (\chibar^{\rm T} \chi)$
depend on $\psi, \psibar, \lambda, \lambdabar$
only via scalar products.
This allows us to apply Lemma~\ref{lem.pbABp.Alan}
to the integral over $\psi, \psibar, \lambda, \lambdabar$;
using also the fact that $PP^{\rm T} = XX^{\rm T} $, we obtain
\begin{eqnarray}
& & \!\!\!\!
\det[(\partial \partial^{\rm T})_{IJ}]  \, \det(XX^{\rm T} )^s
  \;=\;
  \epsilon(I,J) \, \det(XX^{\rm T} )^{s-1} \,
  \det[(XX^{\rm T})_{I^c J^c}]
    \,\times
    \nonumber \\
& &
\;
  \int \! \scrd_{n-m}(\chi,\chibar) \; e^{\chibar^{\rm T} \! \chi}
  \int \! \scrd_m(\psi,\psibar) \, \scrd_m(\lambda,\lambdabar) \;
     e^{\lambdabar^{\rm T} \! \lambda}
\left( \prod_{a=k+1}^{m} \psibar_a \psi_a \!\right)
[\det A - (\psibar^{\rm T} \psi) (\chibar^{\rm T} \chi)]^{-2s}  \;.
   \nonumber \\
\label{eq.76595pippo}
\end{eqnarray}
This formula expresses 
$\det[(\partial \partial^{\rm T})_{IJ}]  \det(XX^{\rm T} )^s$
as the desired quantity
$\det(XX^{\rm T} )^{s-1} \times$ $\epsilon(I,J) \det[(XX^{\rm T} )_{I^c J^c}]$
multiplied by the purely combinatorial factor
\begin{eqnarray}
P^{\textrm{\scriptsize symrect}}(s,m,n,k)
& \equiv &
\int \! \scrd_m(\psi,\psibar) \, \scrd_m(\lambda,\lambdabar) \,
  \scrd_{n-m}(\chi,\chibar)
  \,\times
     \nonumber \\
& & \;
 e^{\lambdabar^{\rm T} \! \lambda + \chibar^{\rm T} \! \chi}
\left( \prod_{a=k+1}^{m} \psibar_a \psi_a \! \right)
[\det A - (\psibar^{\rm T} \psi) (\chibar^{\rm T} \chi)]^{-2s}
    \;,  \qquad
\label{eq.834653865}
\end{eqnarray}
which we now proceed to calculate.

First note that the factor $\prod_{a=k+1}^{m} \psibar_a \psi_a$
forces the Taylor expansion of 
$[\det A - (\psibar \psi) (\chibar \chi)]^{-2s}$ to contain no
variables $\psi_a$, $\psibar_a$ with $k+1 \le a \le m$.
We can therefore drop the factor $\prod_{a=k+1}^{m} \psibar_a \psi_a$,
forget about the variables $\psi_a$, $\psibar_a$ with $k+1 \le a \le m$,
and consider $\psi, \psibar$ henceforth as vectors of length $k$.
Let us also rename the vectors $\lambda$, $\lambdabar$ as
\begin{subeqnarray}
(\lambda_1, \ldots, \lambda_m) & = &
(\lambda_1, \ldots, \lambda_{k}, \mu_1, \ldots, \mu_{m-k})
   \\
(\lambdabar_1, \ldots, \lambdabar_m)  & = &
(\lambdabar_1, \ldots, \lambdabar_{k}, \mubar_1, \ldots, \mubar_{m-k})
\end{subeqnarray}
Then the quantity \reff{eq.834653865} is equivalent to
\begin{eqnarray}
\label{eq.834653865bis}
P^{\textrm{\scriptsize symrect}}(s,m,n,k)
&=&
\int \! \scrd_k(\psi,\psibar) \, \scrd_k(\lambda,\lambdabar) \,
\scrd_{m-k}(\mu,\mubar) \,
  \scrd_{n-m}(\chi,\chibar) \,\times
      \nonumber \\
& & \qquad
  e^{\lambdabar^{\rm T} \! \lambda + \mubar^{\rm T} \! \mu
                                   + \chibar^{\rm T} \! \chi} \,
[\det A - (\psibar^{\rm T} \psi) (\chibar^{\rm T} \chi)]^{-2s}
\end{eqnarray}
where $\det A$ has the same expression as in \reff{eq.detA}
but scalar products involving $\psi, \psibar, \lambda, \lambdabar$
are understood as referring only to the first $k$ variables.
Integration over the variables $\mu, \mubar$ is trivial
and produces just a factor 1. So we are left with
\begin{equation}
\label{eq.834653865tris}
P^{\textrm{\scriptsize symrect}}(s,m,n,k)
\;=\;
\int \! \scrd_k(\psi,\psibar) \, \scrd_k(\lambda,\lambdabar) \,
  \scrd_{n-m}(\chi,\chibar) \;
  e^{\lambdabar^{\rm T} \! \lambda + \chibar^{\rm T} \! \chi} \,
[\det A - (\psibar^{\rm T} \psi) (\chibar^{\rm T} \chi)]^{-2s}
   \;.
\end{equation}
Note that $P^{\textrm{\scriptsize symrect}}(s,m,n,k)$ depends on $n$ and
$m$ only via the combination $n-m$.

The multinomial expansion of the integrand in \reff{eq.834653865tris}
(except for $e^{\chibar^{\rm T} \! \chi}$) is
\begin{eqnarray}
& &
\sum_{\begin{scarray}
          r \geq 0 \\
          t_1, \ldots, t_5 \geq 0
      \end{scarray}}
\!\! 
\binom{-2s}{t_1 + \ldots + t_5}
\binom{t_1 + \ldots + t_5}{t_1 , \ldots , t_5}
e^{\chibar^{\rm T} \! \chi}
\frac{(-1)^{t_5}}{r!}
   \;\times
\nonumber \\
& & \qquad\qquad
(\lambdabar^{\rm T} \lambda)^{r}
(\chibar^{\rm T} \chi)^{t_5}
(\lambda^{\rm T} \psibar)^{t_1+t_3} (\lambdabar^{\rm T} \psi)^{t_2+t_3}
[(\lambdabar^{\rm T} \psibar) (\psi^{\rm T} \lambda)]^{t_4}
(\psibar^{\rm T} \psi)^{t_5}
\ef, \qquad
   \label{eq.372354263.new}
\end{eqnarray}
The integration
\be
\int \scrd_{n-m}(\chi,\chibar) \,
e^{\chibar^{\rm T} \! \chi}
(\chibar^{\rm T} \chi)^{\ell}
      \;=\;
\frac{(n-m)!}{(n-m-\ell)!}
\ee
is trivial. So we are left with two integrations over complex
Grassmann vectors of length $k$. This is performed through a
lemma that we shall prove at the end of this subsection:

\begin{lemma}
For all integers $a,a',b,b',c,c' \ge 0$, we have
\begin{eqnarray}
& &
\int \! \scrd_n(\psi,\psibar) \, \scrd_n(\lambda,\lambdabar) \;
   (\lambdabar^{\rm T} \lambda)^{a} \,
   (\psibar^{\rm T} \psi)^{a'} \,
   (\lambdabar^{\rm T} \psi)^{b} \,
   (\lambda^{\rm T} \psibar)^{b'} \, 
   (\lambdabar^{\rm T} \psibar)^{c} \,
   (\psi^{\rm T} \lambda)^{c'}
  \nonumber \\[2mm]
& & \qquad\qquad =\;\,
 \delta_{aa'} \, \delta_{bb'} \, \delta_{cc'} \, \delta_{a+b+c,n}
   \binom{n}{a,b,c} (a! b! c!)^2
   \;.
  \label{eq.lem.GVform1}
\end{eqnarray}
\label{lem.GVform1}
\end{lemma}


Combining formula \reff{eq.372354263.new} (after integration of
$\chibar$, $\chi$) and the statement of the lemma,
we are left with three independent summations;
we choose the summation indices to be $t_1$, $t_3$ and $t_5$,
which we relabel as $h$, $l$ and $j$, respectively.
The remaining indices are given by
$r = j$,
$t_2 = h$, $t_4 = m-j-h-l$.
The resulting expression is
\begin{equation}
P^{\textrm{\scriptsize symrect}}(s,m,n,k)
\;=\;
(k!)^2 \binom{-2s}{k}
\!\!
\sum_{\begin{scarray}
           j,h,l \geq 0 \\
           j+h+l \leq k
      \end{scarray}}
\!
(-1)^j
\binom{-2s-k}{h} \binom{n-m}{j} \binom{h+l}{h}
\ef.
\end{equation}
Renaming $-2s-k \equiv a$ and $n-m \equiv b$ highlights the fact
that there is no direct dependence of the summands on $k$.
{}From Lemma~\ref{lemma2} we have
\begin{equation}
\sum_{\begin{scarray}
          j,h,l \geq 0 \\
          j+h+l \leq k
      \end{scarray}}
(-1)^j
\binom{a}{h} \binom{b}{j} \binom{h+l}{h}
\;=\;
\binom{a-b+k+1}{k}
\ef,
\end{equation}
which finally gives
\begin{subeqnarray}
P^{\textrm{\scriptsize symrect}}(s,m,n,k)
& = &
(k!)^2 \binom{-2s}{k} \binom{-(2s+n-m-1)}{k}
\\
& = &
\prod_{j=0}^{k-1} (2s+j)(2s+n-m-1+j)
\ef.
\label{eq.polinPsmn}
\end{subeqnarray}
This proves \reff{eq.rectcayley.2}
when $X$ is a real $m \times n$ matrix of rank $m$
lying in a sufficiently small neighborhood of $\widehat{I}_{mn}$,
and $s$ is a positive integer;
the general validity of the identity then follows from
Proposition~\ref{prop.equiv.symbolic_s}.
\qed



\proofof{Lemma~\ref{lem.GVform1}}
%
Let us rewrite \reff{eq.lem.GVform1} by forming an exponential
generating function: that is, we multiply both sides by
$\alpha^a (\alpha')^{a'} \beta^b (\beta')^{b'} \gamma^c (\gamma')^{c'}/
 [a! a'! b! b'! c! c'!]$
and sum over $a,a',b,b',c,c' \ge 0$.
So \reff{eq.lem.GVform1} is equivalent to
\be
   \int \! \scrd_n(\psi,\psibar) \, \scrd_n(\lambda,\lambdabar) \;
     e^{\alpha \lambdabar^{\rm T} \lambda +
        \alpha' \psibar^{\rm T} \psi +
        \beta \lambdabar^{\rm T} \psi +
        \beta' \lambda^{\rm T} \psibar +
        \gamma \lambdabar^{\rm T} \psibar +
        \gamma' \psi^{\rm T} \lambda}
   \;=\;
   (\alpha\alpha' + \beta\beta' + \gamma\gamma')^n
   \;,
 \label{eq.lem.GVform1.genfnAA}
\ee
and this is the formula that we shall prove.
Note first that the measure can be rewritten as
\be
\scrd_n(\psi,\psibar) \, \scrd_n(\lambda,\lambdabar)
=
d \psi_n d \psibar_n d \lambda_n d \lambdabar_n
\cdots
d \psi_1 d \psibar_1 d \lambda_1 d \lambdabar_1
\ee
with no minus signs.
So let us assemble $\psi,\psibar,\lambda,\lambdabar$
into a single Grassmann vector $\theta$ of length $4n$,
with $\psi_k = \theta_{4k}$, $\psibar_k = \theta_{4k-1}$,
$\lambda_k = \theta_{4k-2}$ and $\lambdabar_k = \theta_{4k-3}$;
then the measure becomes $\scrd_{4n}(\theta)$.
In the exponential we have an antisymmetric bilinear form 
$\frac{1}{2} \theta Q \theta$,
in which $Q$ is block-diagonal with $n$ identical $4 \times 4$ blocks
that we call $Q'$:
\be
Q' \;=\; 
\left(\!
\begin{array}{cccc}
    0   & \alpha  &  \gamma  &  \beta   \\
-\alpha &    0    &  \beta'  & -\gamma' \\
-\gamma & -\beta' &    0     & \alpha'  \\
-\beta  & \gamma' & -\alpha' &    0
\end{array}
\!\right)
  \;.
\ee
So the result of the integration is $\pf Q = (\pf Q')^n$.
And indeed, $\pf Q' = \alpha\alpha' + \beta\beta' + \gamma\gamma'$,
as was to be proven.
\qed

\subsection{One-matrix rectangular antisymmetric Cayley identity}
    \label{sec.grassmann.2.antisymrect}

Let us now prove the one-matrix rectangular antisymmetric Cayley identity
(Theorem~\ref{thm.antisymrectcayley}).
Again, this proof is extremely similar to that of the two-matrix and
the symmetric one-matrix rectangular identities.
As always with pfaffians, we will deal only with principal minors,
for which we will need the full strength of Lemma~\ref{lem.pbCJCp.Alan}.

\proofof{Theorem~\ref{thm.antisymrectcayley}}
We begin once again by representing the differential operator
as a Grassmann integral: exploiting Corollary~\ref{cor.blockmatrixpf},
we have
\begin{subeqnarray}
\pf (\partial J \partial^{\rm T} )
&=&
(-1)^m
\pf \left( \begin{array}{c|c}
0_{2m} & \partial \\
\hline
-\partial^{\rm T}  & J_{2n}
\end{array} \right)
\\[1mm]
&=&
(-1)^m
\int \! \scrd_{2m}(\psi) \, \scrd_{2n}(\eta)
\,
e^{\frac{1}{2} \eta^{\rm T} J \eta + \psi^{\rm T} \! \partial \eta}
   \;,
\end{subeqnarray}
where $\psi_i$ ($1 \le i \le 2m$) and $\eta_j$ ($1 \le j \le 2n$)
are ``real'' Grassmann variables,
and we use the same conventions as in the preceding subsection.
For a general even-dimensional principal minor $I \subseteq [2m]$ with $|I|=2k$,
we similarly have (writing $L = \{m+1,\ldots,m+n\}$)
\begin{subeqnarray}
\pf[(\partial J \partial^{\rm T} )_{I I}]
&=&
(-1)^k
\pf \left[ \left( \begin{array}{c|c}
0_{2m} & \partial \\
\hline
-\partial^{\rm T}  & J_{2n}
\end{array} \right)_{\!\! I \cup L, I \cup L} \right]
   \\[1mm]
 &=&
(-1)^k \epsilon(I)
\int \! \scrd_{2m}(\psi) \, \scrd_{2n}(\eta)
\,
\Big( \prod \psi \Big)_{\! I^c} \:
e^{\frac{1}{2} \eta^{\rm T} J \eta + 
\psi^{\rm T} \! \partial \eta}
  \;. \qquad
\end{subeqnarray}

Applying the translation formula \reff{eq.translation}
to the whole set of variables $\{x_{ij}\}$ produces
\begin{equation}
\label{eq.76595.new.skew}
\pf[(\partial J \partial^{\rm T} )_{II}]  \, f(X)  \;=\;
(-1)^k \epsilon(I)
\int \! \scrd_{2m}(\psi) \, \scrd_{2n}(\eta)
   \, 
\Big( \prod \psi \Big)_{\! I^c}
e^{\frac{1}{2} \eta^{\rm T} J \eta} \,
f(X + \psi \eta^{\rm T})
\end{equation}
for an arbitrary polynomial $f$.
We shall use this formula in the case $f(X)=\pf(XJX^{\rm T} )^s$
where $s$ is a positive integer.
It is convenient to introduce the shorthand
\be
X^{\rm trans} \;\equiv\; X + \psi \eta^{\rm T}
\ee
for the argument of $f$.

Suppose now that $X$ is a real $2m \times 2n$ matrix of rank $2m$
that is sufficiently close to the matrix $\widehat{I}_{2m,2n}$
defined in \reff{def.Ihat}.
Then by Lemma~\ref{lemma.decomp.4}
we can find matrices $P \in GL(2m)$ and $Q \in Sp(2n)$
such that $X=P \widehat{I}_{2m,2n} Q$.
We recall the defining property of $Sp(2n)$,
namely that $Q J_{2n} Q^{\rm T}= J_{2n}$.
We have
\begin{equation}
X J_{2n} X^{\rm T} \;=\; 
P \widehat{I}_{2m,2n} Q J_{2n} Q^{\rm T}  \widehat{I}^{\rm T}_{2m,2n}
  P^{\rm T} 
   \;=\;
P J_{2m} P^{\rm T}
   \;,
\label{eq.89857687}
\end{equation}
so that
\begin{equation}
\pf(X J_{2n} X^{\rm T} ) \;=\; 
\pf(P J_{2m} P^{\rm T} )  \;=\;
\det P
\end{equation}
and
\begin{equation}
X^{\rm trans} \;\equiv\; X + \psi \eta^{\rm T}
  \;=\;
P [\widehat{I}+P^{-1} \psi \eta^{\rm T} Q^{-1} ] Q
\ef.
\end{equation}
(we drop the subscripts on $\widehat{I}_{2m,2n}$ to lighten the notation).
Let us now change variables from $(\psi, \eta)$
to $(\psi', \eta') \equiv
    (P^{-1} \psi, Q^{\rm -T} \eta)$,
with Jacobian $(\det P)^{-1}=\pf(X J X^{\rm T} )^{-1}$.
In the new variables we have
(dropping now the primes from the notation)
\begin{equation}
X^{\rm trans} =
P (\widehat{I}+ \psi \eta^{\rm T} ) Q
\ef,
\end{equation}
and the translated pfaffian is given by
\begin{eqnarray}
\pf[(X^{\rm trans}) J (X^{\rm trans})^{\rm T} ] 
&=&
\pf[P (\widehat{I}+ \psi \eta^{\rm T} ) J
     (\widehat{I}^{\rm T} -\eta \psi^{\rm T} )
P^{\rm T}]
   \nonumber \\
&=&
\det (P) \, 
\pf[(\widehat{I}+ \psi \eta^{\rm T} ) J
     (\widehat{I}^{\rm T} -\eta \psi^{\rm T} )]
   \nonumber \\
&=&
\pf(XJX^{\rm T} )
\pf[(\widehat{I}+ \psi \eta^{\rm T} ) J
     (\widehat{I}^{\rm T} -\eta \psi^{\rm T} )]
\ef,
\end{eqnarray}
so that
\begin{eqnarray}
   & &
   \pf[(\partial J \partial^{\rm T} )_{II}]  \,  \pf(X J X^{\rm T} )^s
       \nonumber \\
   & &
   \qquad =\;
\pf(X J X^{\rm T} )^{s-1}
(-1)^k \epsilon(I)
   \int \! \scrd_{2m}(\psi) \, \scrd_{2n}(\eta)
   \,
   \Big( \prod (P\psi) \Big)_{\! I^c}
   \, \times
        \nonumber \\
   & & \qquad\qquad\qquad
   e^{ \frac{1}{2} \eta^{\rm T} J \eta} \,
   \pf[(\widehat{I}+ \psi \eta^{\rm T} ) J
     (\widehat{I}^{\rm T} -\eta \psi^{\rm T} )]^s
   \;.
\end{eqnarray}    
   
Let us now split the vector $\eta$ as
\be
   (\eta_1, \ldots, \eta_{2n})
   \;=\;
   (\lambda_1, \ldots, \lambda_{2m} , \chi_1, \ldots, \chi_{2(n-m)})
   \;,
\ee
so that
\begin{equation}
(\widehat{I}+ \psi \eta^{\rm T} ) J_{2n}
(\widehat{I}^{\rm T} -\eta \psi^{\rm T} )
 \;=\;
J_{2m} - \psi \lambda^{\rm T} J_{2m}^{\rm T} 
- J_{2m} \lambda \psi^{\rm T}
-c \psi \psi^{\rm T}
  \label{eq.antisym.atsight}
\end{equation}
with $c= \eta^{\rm T} J_{2n} \eta =
\lambda^{\rm T} J_{2m} \lambda + \chi^{\rm T} J_{2(n-m)} \chi$.
The matrix \reff{eq.antisym.atsight} is manifestly antisymmetric,
and has the form of a low-rank perturbation of $J_{2m}$.
Since $\pf J_{2m} = 1$ and the perturbation is purely Grassmannian,
we can write
\begin{equation}
\pf (J - \psi \lambda^{\rm T} J^{\rm T}
- J \lambda \psi^{\rm T}
-c \psi \psi^{\rm T} )
  \;=\;
{\det}^{1/2} (I + J \psi \lambda^{\rm T} J^{\rm T}
- \lambda \psi^{\rm T}
+ c J \psi \psi^{\rm T} )
   \;.
\end{equation}
Here the argument of $\det$ is a matrix of the form
$I_{2m} +\sum_{\alpha=1}^{2} u_{\alpha} v_{\alpha}^{\rm T}$,
with vectors $\{u_{\alpha}\}$, $\{v_{\alpha}\}$ given by
\begin{equation}
\begin{array}{c||c|c}
\alpha & u_{\alpha} & v_{\alpha}
\\
\hline
\rule{0pt}{4.5mm}
1 & J \psi & J \lambda + c \psi \\
2 & -\lambda & \psi \\
\end{array}
   \label{eq.rectcayleyproof.rank2skew}
\end{equation}
By Lemma~\ref{lemma.lowrank} we can write the needed determinant
as the determinant of a $2 \times 2$ matrix, which, after a slight
row manipulation, can be written as
\be
\det (I + J \psi \lambda^{\rm T} J^{\rm T}
- \lambda \psi^{\rm T}
+ c J \psi \psi^{\rm T} )
\;=\;
   \det{}^{-1}
   \left(
   \begin{array}{cc}
1 - \psi^{\rm T} \lambda & - \chi^{\rm T} J \chi \\
\psi^{\rm T} J \psi      & 1 - \psi^{\rm T} \lambda
   \end{array}
   \right)
   \,.
\ee
We therefore have
\begin{eqnarray}
   & &
   \!\!
   \pf[(\partial J \partial^{\rm T} )_{II}]  \,  
\pf(X J X^{\rm T} )^s
   \;=\;
\pf(X J X^{\rm T} )^{s-1}
(-1)^k \epsilon(I)
   \int \! \scrd_{2(n-m)}(\chi) \; e^{\frac{1}{2} \chi^{\rm T} J \chi}
      \,\times
      \nonumber \\
   & &
   \quad
   \int \! \scrd_{2m}(\psi) \, \scrd_{2m}(\lambda)
   \;
   e^{\frac{1}{2} \lambda^{\rm T} J \lambda} \,
   \Big( \prod (P\psi) \Big)_{\! I^c}  \,
[(1 - \psi^{\rm T} \lambda)^2 +
  (\psi^{\rm T} J \psi) (\chi^{\rm T} J \chi)]^{-s/2}
   \;.
   \qquad
      \nonumber \\
\end{eqnarray}
In this expression we have both an ordinary scalar product
($\psi^{\rm T} \lambda$)
and symplectic scalar products
($\lambda^{\rm T} J \lambda$, $\psi^{\rm T} J \psi$ and $\chi^{\rm T} J \chi$).
By a further change of variables
$\lambda \to \lambda' = -J \lambda$,
we can reduce to symplectic products only (the Jacobian is 1).
Dropping the primes, we have (also using $J^{\rm T} J J = J$)
\begin{eqnarray}
   & &
   \!\!
   \pf[(\partial J \partial^{\rm T} )_{II}]  \,  
\pf(X J X^{\rm T} )^s
   \;=\;
\pf(X J X^{\rm T} )^{s-1}
(-1)^k \epsilon(I)
   \int \! \scrd_{2(n-m)}(\chi) \; e^{\frac{1}{2} \chi^{\rm T} J \chi}
      \,\times
      \nonumber \\
   & &
   \quad
   \int \! \scrd_{2m}(\psi) \, \scrd_{2m}(\lambda)
   \;
   e^{\frac{1}{2} \lambda^{\rm T} J \lambda} \,
   \Big( \prod (P\psi) \Big)_{\! I^c}  \,
[(1 - \psi^{\rm T} J \lambda)^2 +
   (\psi^{\rm T} J \psi) (\chi^{\rm T} J \chi)]^{-s/2}
   \;.
   \qquad
     \nonumber \\
\end{eqnarray}
Now the integral on the $\psi$ and $\lambda$ fields is of the form
described in Lemma~\ref{lem.pbCJCp.Alan}. Applying this lemma, 
and making use of (\ref{eq.89857687}), we have
\begin{eqnarray}
   & &
   \pf[(\partial J \partial^{\rm T} )_{II}]  \,  
\pf(X J X^{\rm T} )^s
   \;=\;
\epsilon(I) 
\pf(X J X^{\rm T} )^{s-1}
\pf[(X J X^{\rm T} )_{I^c I^c}]
      \nonumber \\
   & &
   \qquad \times\;
(-1)^k 
   \int \! \scrd_{2(n-m)}(\chi) \; e^{\frac{1}{2} \chi^{\rm T} J \chi}
   \int \! \scrd_{2m}(\psi) \, \scrd_{2m}(\lambda)
   \;
   e^{\frac{1}{2} \lambda^{\rm T} J \lambda} \,
      \nonumber \\
   & &
   \qquad \times\;
   \Big( \prod \psi \Big)_{ \{2k+1,\ldots,2m\} }  \,
[(1 - \psi^{\rm T} J \lambda)^2 +
   (\psi^{\rm T} J \psi) (\chi^{\rm T} J \chi)]^{-s/2}
   \;.
   \qquad\quad
\end{eqnarray}
This formula expresses 
$\pf[(\partial J \partial^{\rm T})_{II}]  \pf(X J X^{\rm T} )^s$
as the desired quantity
$\pf(X J X^{\rm T} )^{s-1} \times \epsilon(I) \pf [(X J X^{\rm T} )_{I^c I^c}]$
multiplied by the purely combinatorial factor
\begin{eqnarray}
   & &
P^{\textrm{\scriptsize asrect}}(s,m,n,k)
   \;\equiv\;
(-1)^k 
   \int \! \scrd_{2(n-m)}(\chi) \; e^{\frac{1}{2} \chi^{\rm T} J \chi}
   \int \! \scrd_{2m}(\psi) \, \scrd_{2m}(\lambda)
      \, e^{\frac{1}{2} \lambda^{\rm T} J \lambda}
      \nonumber \\
   & &
   \qquad \times\;
   \Big( \prod \psi \Big)_{ \{2k+1,\ldots,2m\} }  \,
[(1 - \psi^{\rm T} J \lambda)^2 +
   (\psi^{\rm T} J \psi) (\chi^{\rm T} J \chi)]^{-s/2}
    \;,  \qquad
\label{eq.834653865.ANTISYM}
\end{eqnarray}
which we now proceed to calculate.

Note, first of all, that the factor $\prod_{a=2k+1}^{2m} \psi_a$
forces the Taylor expansion of 
$[(1 - \psi^{\rm T} J \lambda)^2 +
   (\psi^{\rm T} J \psi) (\chi^{\rm T} J \chi)]^{-s/2}$
to contain no variables $\psi_a$ with $2k+1 \le a \le 2m$,
and thus also no variables $\lambda_a$ with $2k+1 \le a \le 2m$. The
latter must all come from the 
trivial diagonal exponential
$e^{\frac{1}{2} \lambda^{\rm T} J \lambda}$, and their integration can
then be performed easily (it gives 1).
We can therefore drop the factor $\prod_{a=2k+1}^{2m} \psi_a$,
forget about the variables $\psi_a$ and $\lambda_a$ with $2k+1 \le a \le 2m$,
and consider $\psi$ and $\lambda$ henceforth as vectors of length $2k$.
We have
\begin{eqnarray}
P^{\textrm{\scriptsize asrect}}(s,m,n,k)
& \equiv &
(-1)^k 
   \int \! \scrd_{2(n-m)}(\chi) \; e^{\frac{1}{2} \chi^{\rm T} J \chi}
   \int \! \scrd_{2k}(\psi) \, \scrd_{2k}(\lambda)
      \,\times
      \nonumber \\
   & &
  \qquad
   \;
   e^{\frac{1}{2} \lambda^{\rm T} J \lambda} \,
[(1 - \psi^{\rm T} J \lambda)^2 +
   (\psi^{\rm T} J \psi) (\chi^{\rm T} J \chi)]^{-s/2}
    \;.  \qquad
\label{eq.834653865.ANTISYM.2}
\end{eqnarray}
We see that $P^{\textrm{\scriptsize asrect}}(s,m,n,k)$ depends on $n$ and
$m$ only via the combination $n-m$.

In the expression \reff{eq.834653865.ANTISYM.2}
we have symplectic scalar products of ``real'' Grassmann fields.
However, our summation lemmas (such as Lemma~\ref{lem.GVform1})
have heretofore been developed only for ordinary scalar products
of ``complex'' Grassmann fields.
However, through a relabeling of the fields,
we can rewrite \reff{eq.834653865.ANTISYM.2}
in terms of ordinary scalar products of ``complex'' fields
of half the dimensions.
More precisely, we relabel
$(\psi_1, \ldots, \psi_{2k}) \to 
 (\psibar_1, \psi_1, \ldots, \psibar_k, \psi_k)$,
$(\lambda_1, \ldots, \lambda_{2k}) \to 
 (\lambdabar_1, \lambda_1, \ldots, \lambdabar_k, \lambda_k)$
and
$(\chi_1, \ldots, \chi_{2(n-m)}) \to 
 (\chibar_1, \chi_1, \ldots, \chibar_{n-m}, \chi_{n-m})$. 
No signs arise in the measure of integration, and we have the correspondences
\begin{subeqnarray}
\smfrac{1}{2} \psi^{\rm T} J \psi
   &\longrightarrow & \psibar^{\rm T} \psi
\\
\smfrac{1}{2} \lambda^{\rm T} J \lambda
   &\longrightarrow & \lambdabar^{\rm T} \lambda
\\
\psi^{\rm T} J \lambda
   &\longrightarrow & \psibar^{\rm T} \lambda + \lambdabar^{\rm T} \psi
\\
\smfrac{1}{2} \chi^{\rm T} J \chi
   &\longrightarrow & \chibar^{\rm T} \chi
\end{subeqnarray}
which rewrites (\ref{eq.834653865.ANTISYM.2}) as
\begin{eqnarray}
   & &
   P^{\textrm{\scriptsize asrect}}(s,m,n,k)
   \;=\;
   (-1)^k
   \int \! \scrd_k(\psi,\psibar) \, \scrd_k(\lambda,\lambdabar) \,
       \scrd_{n-m}(\chi,\chibar) \;
     \nonumber \\
   & &
   \qquad \times\;
  e^{\lambdabar^{\rm T} \! \lambda + \chibar^{\rm T} \! \chi} \,
[(1 - \psibar^{\rm T} \lambda - \lambdabar^{\rm T} \psi )^2 
+ 4 (\psibar^{\rm T} \psi) (\chibar^{\rm T} \chi)]^{-s/2}
   \;.  \qquad
  \label{eq.834653865tris.ANTISYM}
\end{eqnarray}

We begin by getting rid of variables $\chi$ and $\chibar$,
writing
\begin{eqnarray}
& &
\int \!
\scrd_{n-m}(\chi,\chibar) \;
  e^{\chibar^{\rm T} \! \chi} \,
[(1 - \psibar^{\rm T} \lambda - \lambdabar^{\rm T} \psi )^2 
+ 4 (\psibar^{\rm T} \psi) (\chibar^{\rm T} \chi)]^{-s/2}
         \nonumber \\
& &
\qquad =\;
\sum_{a \geq 0}
\binom{-\frac{s}{2}}{a} (4 \psibar^{\rm T} \psi)^a 
(1 - \psibar^{\rm T} \lambda - \lambdabar^{\rm T} \psi)^{-s-2a}
\binom{n-m}{a} a!
   \;.  \qquad
  \label{eq.6754865}
\end{eqnarray}
So we have
\begin{eqnarray}
   & &
P^{\textrm{\scriptsize asrect}}(s,m,n,k)
   \;=\;
(-1)^k
\int \! \scrd_k(\psi,\psibar) \, \scrd_k(\lambda,\lambdabar) \,
  e^{\lambdabar^{\rm T} \! \lambda}
      \;\times \nonumber \\
& &
\qquad\qquad
\sum_{a \geq 0}
\binom{-\frac{s}{2}}{a} (4 \psibar^{\rm T} \psi)^a 
(1 - \psibar^{\rm T} \lambda - \lambdabar^{\rm T} \psi)^{-s-2a}
\binom{n-m}{a} a!
   \;.  \qquad
  \label{eq.834653865trisX}
\end{eqnarray}
Next we expand fully the integrand in \reff{eq.834653865trisX}, yielding
\begin{eqnarray}
   & &
P^{\textrm{\scriptsize asrect}}(s,m,n,k)
   \;= \!
\sum_{\begin{scarray}
    a, a' \geq 0 \\
    b, b' \geq 0
      \end{scarray}}
\!\! 
4^a 
(-1)^{k-b}
\binom{-\frac{s}{2}}{a} 
\binom{n-m}{a} a!
\binom{-s-2a}{b+b'}
\binom{b+b'}{b}
\frac{1}{a'!}
   \;\times
\nonumber \\
& & \qquad\qquad
\int \! \scrd_k(\psi,\psibar) \, \scrd_k(\lambda,\lambdabar) \,
(\psibar^{\rm T} \psi)^{a}
(\lambdabar^{\rm T} \lambda)^{a'}
(\psibar^{\rm T} \lambda)^{b}
(\psi^{\rm T} \lambdabar)^{b'}
\ef. \qquad
\end{eqnarray}
The fermionic integration here is a special case of
Lemma~\ref{lem.GVform1}, which thus gives
\begin{eqnarray}
  P^{\textrm{\scriptsize asrect}}(s,m,n,k)
  & = &  \!\!\!\!
\sum_{\begin{scarray}
    a, a' \geq 0 \\
    b, b' \geq 0
      \end{scarray}}
\!\! 
4^a 
(-1)^{k-b}
\binom{-\frac{s}{2}}{a} 
\binom{n-m}{a} a!
\binom{-s-2a}{b+b'}
\binom{b+b'}{b}
\frac{1}{a'!}
\nonumber \\
& & \qquad\qquad \times\;
 \delta_{aa'} \, \delta_{bb'} \, \delta_{a+b,k}
   \binom{k}{a} (a! b!)^2
\nonumber \\
   & = &
\sum_{a=0}^k (-4)^a
   (-\smfrac{s}{2})^{\underline{a}} (-s-2a)^{\underline{2k-2a}}
\binom{n-m}{a} 
\binom{k}{a}
a!
   \label{eq.372354263.new.1}
\end{eqnarray}
where $x^{\underline{k}}=x(x-1)\cdots(x-k+1)$.
Notice now that
\begin{equation}
(-\smfrac{s}{2})^{\underline{a}}  (-s-2a)^{\underline{2k-2a}}
 \;=\;
 (-\smfrac{1}{2})^a
s(s+2)\cdots(s+2k-2)
\,
(s+2a+1)(s+2a+3)\cdots(s+2k-1)
  \;.
\end{equation}
The factors of the form $s+2j$, namely $s(s+2)\cdots(s+2k-2)$,
are independent of the summation variable $a$.
What remains is
\begin{eqnarray}
  & &
\sum_{a=0}^k
2^a 
(s+2a+1)(s+2a+3)\cdots(s+2k-1)
\binom{n-m}{a} 
\binom{k}{a}
a!
    \nonumber \\
  & & \qquad =\;
\sum_{a=0}^k
2^k
\binom{\frac{s-1}{2}+k}{k-a}
(k-a)!
\binom{n-m}{a} 
\binom{k}{a}
a!
    \nonumber \\
  & & \qquad =\;
2^k k!
\sum_{a=0}^k
\binom{\frac{s-1}{2}+k}{k-a}
\binom{n-m}{a}
    \nonumber \\
  & & \qquad =\;
2^k
k!
\binom{\frac{s-1}{2}+n-m+k}{k}
    \nonumber \\
  & & \qquad =\;
\prod_{j=0}^{k-1} (s+1+2n-2m+2j)
\end{eqnarray}
where the sum over $a$ was a Chu--Vandermonde convolution
(Lemma~\ref{lemma.vandermonde}).
Restoring the factors of the form $s+2j$,
we conclude that
\be
   P^{\textrm{\scriptsize asrect}}(s,m,n,k)
   \;=\;
   \prod_{j=0}^{k-1} (s+2j) (s+1+2n-2m+2j)
   \;.
\ee

This proves \reff{eq.antisymrectcayley.2}
when $X$ is a real $2m \times 2n$ matrix of rank $2m$
lying in a sufficiently small neighborhood of $\widehat{I}_{2m,2n}$,
and $s$ is a positive integer;
the general validity of the identity then follows from
Proposition~\ref{prop.equiv.symbolic_s}.
\qed

\subsection{A lemma on Grassmann integrals of scalar products}
\label{sec.grassmann.lemmascalar}

In the preceding proofs we have frequently had to evaluate
Grassmann integrals over one or more sets of Grassmann variables,
in which the integrand depends only on scalar products among those sets.
A simple version of such a formula occurred already in
\reff{eq.scalar.etaetabar},
while more complicated versions arose in
\reff{eq.834653865trisNEW} ff.\ and \reff{eq.lem.GVform1}.
[Also, a version for real fermions arose in
 \reff{eq.realfermion.scalartheta}.]
So far we have simply treated such integrals ``by hand''.
But since a much more complicated such integral will arise
in the proof of the multi-matrix rectangular Cayley identity
(Theorem~\ref{thm.multirectcayley}) in the next subsection,
it is worth stating now a general lemma that allows us
to systematize such calculations.
Indeed, we think that this lemma
(Proposition~\ref{prop.intscalar} below) is of some independent interest.

Let us start by rephrasing the simplest case \reff{eq.scalar.etaetabar}
in a suggestive way.
Let $n$ be a positive integer, and introduce Grassmann variables
$\psi_i,\psibar_i$ ($1 \le i \le n$).
Let $f(x) = \sum_{k=0}^\infty a_k x^k$ be a formal power series
in one indeterminate.  We then have by \reff{eq.scalar.etaetabar}
\be
   \int \! \scrd_n(\psi,\psibar) \: f(\psibar^{\rm T} \psi)
   \;=\;
   n! \, a_n
   \;.
\ee
On the other hand, we also have trivially
\be
   \left. {d^n \over dx^n} f(x) \right| _{x=0}
   \;=\;
   n! \, a_n
\ee
and
\be
   \left. f\Bigl( {d \over dx} \Bigr) \, x^n \right| _{x=0}
   \;=\;
   n! \, a_n
   \;.
\ee
Hence
\be
   \int \! \scrd_n(\psi,\psibar) \: f(\psibar^{\rm T} \psi)
   \;=\;
   \left. {d^n \over dx^n} f(x) \right| _{x=0}
   \;=\;
   \left. f\Bigl( {d \over dx} \Bigr) \, x^n \right| _{x=0}
   \;.
\ee

Surprisingly enough, a similar formula exists for Grassmann integrals
involving multiple sets of fermionic variables:

\begin{proposition}
  \label{prop.intscalar}
Let $\ell$ and $n$ be positive integers, and introduce Grassmann variables
$\psi^\alpha_i,\psibar^\alpha_i$ ($1 \le \alpha \le \ell$, $1 \le i \le n$).
Let $f(X)$ be a formal power series in commuting indeterminates
$X = (x_{\alpha\beta})_{\alpha,\beta=1}^\ell$,
and write $\partial = (\partial/\partial x_{\alpha\beta})$.
We then have
\begin{subeqnarray}
   \int \! \scrd_n(\psi^1,\psibar^1) \cdots \scrd_n(\psi^\ell,\psibar^\ell)
      \; f( \{ \psibar^{\alpha \rm T} \psi^\beta \} )
   & = &
   \left. \det(\partial)^n \, f(X) \right|_{X=0}
 \slabel{eq.prop.intscalar.a}
       \\[2mm]
   & = &
   \left. f(\partial) \, (\det X)^n \right|_{X=0}
   \;. \qquad
 \slabel{eq.prop.intscalar.b}
 \label{eq.prop.intscalar}
\end{subeqnarray}
\end{proposition}

\proof
We write
\be
   f(X) \;=\; \sum_{N \in \mathbb{N}^{\ell \times \ell}}
       f_N \prod_{\alpha, \beta=1}^{\ell} x_{\alpha \beta}^{n_{\alpha \beta}}
\ee
where the sum runs over matrices
$N=( n_{\alpha \beta} )_{\alpha, \beta=1}^{\ell}$
of nonnegative integers.
Note first that it suffices to prove \reff{eq.prop.intscalar}
for {\em polynomials}\/ $f$, since all three expressions in
\reff{eq.prop.intscalar} have the property that only finitely many
coefficients $N$ contribute (namely, those matrices $N$ with
all row and column sums equal to $n$).
So it suffices to prove \reff{eq.prop.intscalar} for all monomials
$X^N \bydef \prod_{\alpha, \beta=1}^{\ell} x_{\alpha \beta}^{n_{\alpha \beta}}$.
But then it suffices to prove \reff{eq.prop.intscalar}
for the exponential generating function
\be
   \phi_{\Omega}(X)
   \;=\;
   \exp\tr(\Omega^{\rm T} X)
   \;=\;
   \sum_{N \in \mathbb{N}^{\ell \times \ell}} \, \prod_{\alpha, \beta=1}^{\ell} 
          \frac{(\omega_{\alpha \beta} x_{\alpha \beta})^{n_{\alpha \beta}}
               }{ n_{\alpha \beta} !}
\ee
where $\Omega=( \omega_{\alpha \beta} )_{\alpha, \beta=1}^{\ell}$
are indeterminates, since the value at the monomial $X^N$
can be obtained by extracting the coefficient $[\Omega^N]$.
Our goal is to prove that all the three expressions in
\reff{eq.prop.intscalar}, specialized to $f(X)=\phi_{\Omega}(X)$,
equal $(\det \Omega)^n$.

For the Grassmann-integral expression we have
\begin{subeqnarray}
&&
   \int \! \scrd_n(\psi^1,\psibar^1) \cdots \scrd_n(\psi^\ell,\psibar^\ell)
      \; \phi_{\Omega}( \{ \psibar^{\alpha \rm T} \psi^\beta \} )
   \nonumber \\
&&\qquad
\;=\;
   \int \! \scrd_n(\psi^1,\psibar^1) \cdots \scrd_n(\psi^\ell,\psibar^\ell)
      \; 
\exp \left(
\sum_{i=1}^n \sum_{\alpha, \beta=1}^{\ell} 
\psibar^{\alpha}_i \omega_{\alpha \beta} \psi^\beta_i
\right)
    \qquad \\
&&\qquad
\;=\;
(\det \Omega)^n
\ef.
\end{subeqnarray}
For the power series in derivative operators applied to a power of a
determinant, a simple application of the translation formula
\reff{eq.translation} yields
\begin{subeqnarray}
   \left.
\phi_{\Omega}(\partial) \, (\det X)^n \right|_{X=0}
&=&
   \left. 
\exp \left(
\sum_{\alpha, \beta=1}^{\ell} 
\omega_{\alpha \beta} \frac{\partial}{\partial x_{\alpha \beta}}
\right)
\, (\det X)^n \right|_{X=0}
\\[1mm]
&=&
   \left. 
(\det (X+\Omega))^n
\right|_{X=0}
\\[1mm]
&=&
(\det \Omega)^n
\ef.
 \label{eq.proof.prop.intscalar}
\end{subeqnarray}
For the power of a determinant in derivative operators applied to a
power series, it suffices to cite the proof \reff{eq.proof.prop.intscalar}
and invoke transposition duality $X \leftrightarrow \partial$
in the Weyl algebra;  but for completeness let us give a direct proof.
Use a fermionic representation of the differential operator
\be
\det(\partial)
\;=\;
   \int \! \scrd_{\ell}(\eta, \etabar)
\exp \left(
\sum_{\alpha, \beta=1}^{\ell} 
\etabar_{\alpha}
\frac{\partial}{\partial x_{\alpha \beta}}
\eta_{\beta}
\right)
\ee
and apply again the translation formula \reff{eq.translation}:
\begin{subeqnarray}
 \!\!\!\!\!
\det(\partial)  \:
[\exp \tr (\Omega^{\rm T} X)]
&=& \!
   \int \! \scrd_{\ell}(\eta, \etabar)
\exp\!\left(
\sum_{\alpha, \beta=1}^{\ell} 
\etabar_{\alpha}
\frac{\partial}{\partial x_{\alpha \beta}}
\eta_{\beta}
\right)
[\exp \tr (\Omega^{\rm T} X)]
  \qquad
  \\
&=& \!
   \int \! \scrd_{\ell}(\eta, \etabar)
\exp \tr(\Omega^{\rm T} (X+ \etabar \eta^{\rm T} ))
  \\
& = &
   [\exp \tr (\Omega^{\rm T} X)]
   \int \! \scrd_{\ell}(\eta, \etabar)
\exp\!\left(
\sum_{\alpha, \beta=1}^{\ell} 
\etabar_{\alpha} \omega_{\alpha\beta} \eta_{\beta}
\right)
  \\
&=&
(\det \Omega) \,
\exp \tr(\Omega^{\rm T} X )
\ef.
\end{subeqnarray}
Iterating $n$ times, we get
\be
   \left. 
\det(\partial)^n \, \phi_{\Omega}(X) \right|_{X=0}
\;=\;
   \left. 
(\det \Omega)^n
\exp \tr(\Omega^{\rm T} X )
\right|_{X=0}
\;=\;
(\det \Omega)^n
\ef.
\ee
\qed


\subsection{Multi-matrix rectangular Cayley identity}
\label{sec.grassmann.2.multi}


In this subsection we shall prove
the multi-matrix rectangular Cayley identity
(Theorem~\ref{thm.multirectcayley}),
which is the most difficult of the identities proven in this paper.
The difficulty arises principally from the fact that the number of matrices
appearing in the identity (which we call $\ell$) can be arbitrarily large,
and we are required to provide a proof valid for all $\ell$.
The proof nevertheless follows the basic pattern of the proofs
of the other identities (notably the two-matrix rectangular identity),
and divides naturally into two parts:
\begin{itemize}
   \item[(i)]  We represent the differential operator as a Grassmann
       integral, and after several manipulations we are able to express
       $\det[(\partial^{(1)} \cdots \partial^{(\ell)})_{IJ}] \,
        \det(X^{(1)} \cdots X^{(\ell)})^s$
       as the desired quantity
       $\det (X^{(1)} \cdots X^{(\ell)})^{s-1}
        \epsilon(I,J) \det[(X^{(1)} \cdots X^{(\ell)})_{I^c J^c}]$
       multiplied by a purely combinatorial factor $b(s)$
       that is given as a Grassmann integral.
   \item[(ii)]  We evaluate this Grassmann integral
       and prove that
       $b(s) = \prod\limits_{\alpha=1}^\ell
               \prod\limits_{j=1}^{k} (s+n_\alpha-j)$.
\end{itemize}
The major new complications, as compared to the preceding proofs,
come from the fact that our Grassmann integral
involves $\ell$ sets of fermionic fields,
while in the preceding proofs we only had one or two,
and the matrix arising from the application of the
``low-rank perturbation lemma'' (Lemma~\ref{lemma.lowrank})
is of size $\ell \times \ell$,
while in the preceding proofs we never had a matrix larger than $4 \times 4$.
Despite these complications,
step (i) follows closely the model established in the preceding proofs
and is not much more difficult than them:
the main novelty is that we need to use a variant of
the low-rank perturbation lemma (Corollary~\ref{cor.lowrankrect})
that is specially adapted to perturbations of the product form
\reff{def.Mn1nl} encountered here.
The big trouble arises in step (ii):
to evaluate the Grassmann integral involving scalar products
among $\ell$ sets of fermionic fields,
we shall first rewrite it as a differential operator
acting on a determinant (using Proposition~\ref{prop.intscalar})
and then exploit the fact that the matrix arising in this
differential operator is upper Hessenberg
(i.e.\ has zero entries below the first subdiagonal);
the latter leads to some special combinatorial/algebraic computations
(Theorem~\ref{thm.quasitriang} and Corollary~\ref{cor.solvingmulti}).

It is convenient to formally divide the statement of
Theorem~\ref{thm.multirectcayley} into two parts,
corresponding to steps (i) and (ii) above:

\bigskip
\noindent
{\bf Theorem~\ref{thm.multirectcayley}, part (i)}\quad 
{\em
Fix integers $\ell \ge 1$ and $n_1,\ldots,n_{\ell} \ge 0$
and write $n_{\ell+1} = n_1$.
For $1 \le \alpha \le \ell$,
let $X^{(\alpha)}$ be an $n_\alpha \times n_{\alpha+1}$
matrix of indeterminates, and let $\partial^{(\alpha)}$
be the corresponding matrix of partial derivatives.
If $I,J \subseteq [n_1]$ with $|I| = |J| = k$, then
\begin{eqnarray}
   & & \!\!\!
   \det[(\partial^{(1)} \cdots \partial^{(\ell)})_{IJ}] \,
      \det(X^{(1)} \cdots X^{(\ell)})^s
      \nonumber \\
   & & \quad =\;
   b_{n_1,\ldots,n_\ell;k}(s)
      \, \det (X^{(1)} \cdots X^{(\ell)})^{s-1}
   \, \epsilon(I,J) \, \det[(X^{(1)} \cdots X^{(\ell)})_{I^c J^c}]
      \qquad
 \label{eq.multirectcayley.2.part_i}
\end{eqnarray}
where
\begin{eqnarray}
   &  & \!\!\!\!
   b_{n_1,\ldots,n_\ell;k}(s)  \;=\;
(-1)^{k\ell}
\int \! 
\scrd_{k}(\psi^1,\psibar^1) \cdots
\scrd_{k}(\psi^{\ell},\psibar^{\ell})
\,
\scrd_{n_2-n_1}(\eta^2,\etabar^2) \cdots
\scrd_{n_{\ell}-n_1}(\eta^{\ell},\etabar^{\ell})
  \nonumber \\
  &  & \qquad\qquad\qquad\qquad\qquad \times\:
\exp \!\left[ \sum_{\alpha=2}^{\ell} 
   (\psibar^{\alpha \rm T} \psi^{\alpha-1} 
    + \etabar^{\alpha \rm T} \eta^{\alpha} )
\right]  \,
\det(I_\ell + M)^{-s}
 \label{eq.def.bn1nk}
\end{eqnarray}
and
\be
   M_{\alpha \beta}
   \;=\;  \cases{ \psibar^{\alpha \rm T} \psi^{\beta}
                         & \hbox{if } $\alpha \le \beta$ \cr
                  \noalign{\vskip 6pt}
                  - \etabar^{\alpha \rm T} \eta^{\alpha}
                         & \hbox{if } $\alpha = \beta + 1$ \cr
                  \noalign{\vskip 6pt}
                  0      & \hbox{otherwise} \cr
                }
 \label{eq.def.bn1nk.2}
\ee
}

\bigskip
\noindent
{\bf Theorem~\ref{thm.multirectcayley}, part (ii)}\quad 
{\em
We have
\be
   b_{n_1,\ldots,n_\ell;k}(s)  \;=\;
    \prod\limits_{\alpha=1}^\ell
          \prod\limits_{j=1}^{k} (s+n_\alpha-j)
   \;.
 \label{eq.multirectcayley.2.part_ii}
\ee
}

Let us now begin the proof of part (i) of Theorem~\ref{thm.multirectcayley}.
As before, we will need the full strength of Lemma~\ref{lem.pbABp.Alan}
to handle the all-minors case.

\proofof{Theorem~\ref{thm.multirectcayley}, part (i)}
We can assume that $n_{\alpha} \geq n_1$ for $2 \leq \alpha \leq \ell$,
since otherwise $\det(X^{(1)} \cdots X^{(\ell)})$ is the zero polynomial.

We begin, as usual, by representing the differential operator
as a Grassmann integral, this time exploiting 
Lemma~\ref{lem.blockmatrixdet.multi}
(rather than just Corollary~\ref{cor.blockmatrixdet}). We have
\begin{eqnarray}
\det(\partial^{(1)} \cdots \partial^{(\ell)})
&=&
\int \! \scrd_{n_1}(\psi^1,\psibar^1) \cdots
\scrd_{n_{\ell}}(\psi^{\ell},\psibar^{\ell})
\nonumber
\\
&&\quad \times\;
\exp \!\left[
\sum_{\alpha=2}^{\ell}
( \psibar^{\alpha \rm T} \psi^{\alpha} 
- \psibar^{\alpha-1\,\rm T} \partial^{(\alpha-1)} \psi^{\alpha} )
+ \psibar^{\ell \rm T} \partial^{(\ell)} \psi^{1}
\right]  . \qquad
\label{eq.3886594}
\end{eqnarray}
Here $\psi_i^{\alpha}, \psibar_i^{\alpha}$ 
($1 \leq \alpha \leq \ell$, $1 \le i \le n_{\alpha}$)
are Grassmann variables,
and the subscripts on $\scrd$ serve to remind us
of the length of each vector;
shorthand notations for index summations are understood,
e.g.\ $\psibar^{a \rm T} \partial^{(a)} \psi^{a+1} \equiv 
\sum_{i=1}^{n_a} \sum_{j=1}^{n_{a+1}}
\psibar^{a}_i \psi^{a+1}_j \partial/\partial x^{(a)}_{ij}$.
For a general minor $I,J \subseteq [n_1]$ with $|I|=|J|=k$,
the quantity
$\det[(\partial^{(1)} \cdots \partial^{(\ell)})_{IJ}]$ has a
representation like \reff{eq.3886594} but with an extra factor
$\epsilon(I,J) 
\Big( \prod \psibar^1 \psi^1 \Big)_{I^c,J^c}$ in the integrand.

It is convenient to make the change of variables $\psi^1 \to -\psi^1$,
as this makes the summand
$\psibar^{\ell} \partial^{(\ell)} \psi^{1}$
analogous to the
$\psibar^{\alpha-1} \partial^{(\alpha-1)} \psi^{\alpha}$
arising for $2 \le \alpha \leq \ell$.
We shall exploit this structure by writing
$\psi^{\ell+1}$ as a synonym for $\psi^1$.
This change of variables
introduces an overall factor $(-1)^{n_1} (-1)^{n_1-k} = (-1)^k$.

Applying the translation formula \reff{eq.translation}
to the whole set of variables $\{x^{(\alpha)}_{ij} \}$ produces
\begin{eqnarray}
&&
\det[(\partial^{(1)} \cdots \partial^{(\ell)} )_{IJ}]  
\: f(X^{(1)}, \ldots, X^{(\ell)})
\;=\;
\epsilon(I,J) 
\int \! \scrd_{n_1}(\psi^1,\psibar^1) \cdots
\scrd_{n_{\ell}}(\psi^{\ell},\psibar^{\ell})
\nonumber
\\
&&\qquad \times\;
(-1)^k
\exp \!\left[
\sum_{\alpha=2}^{\ell}
\psibar^{\alpha} \psi^{\alpha} 
\right]
\Big( \prod \psibar^1 \psi^1 \Big)_{I^c,J^c}
 \: f(X^{(1)}_{\rm trans}, \ldots, X^{(\ell)}_{\rm trans})
 \qquad
\label{eq.76595.newNEWmulti}
\end{eqnarray}
for an arbitrary polynomial $f$, where we have introduced the shorthand
\be
X^{(\alpha)}_{\rm trans}   \;=\;
X^{(\alpha)} - \psibar^{\alpha} (\psi^{\alpha+1})^{\rm T}
\ee
for the arguments of $f$.
We shall use the formula \reff{eq.76595.newNEWmulti} in the case
$f(X^{(1)}, \cdots, X^{(\ell)}) = \det(X^{(1)} \cdots X^{(\ell)})^s$
where $s$ is a positive integer.

Suppose now that the $X^{(\alpha)}$ are real matrices
of rank $\min(n_\alpha, n_{\alpha+1})$
that are sufficiently close to the matrix
$\widehat{I}_{n_\alpha n_{\alpha+1}}$
defined in \reff{def.Ihat}.
Then by Lemma~\ref{lemma.decomp.5}
we can find matrices $P_{\alpha} \in GL(n_{\alpha})$
for $1 \le \alpha \le \ell+1$ such that
$X^{(\alpha)} =
 P_{\alpha} \, \widehat{I}_{n_{\alpha} n_{\alpha+1}} \, P_{\alpha+1}^{-1}$.
We have $X^{(1)} \cdots X^{(\ell)} = P_{1} P_{\ell+1}^{-1}$
since $\widehat{I}_{n_1 n_2} \widehat{I}_{n_2 n_3} \cdots
       \widehat{I}_{n_\ell n_1} = I_{n_1}$
as a consequence of the fact that $n_{\alpha} \geq n_1$ for 
$2 \leq \alpha \leq \ell$.
Therefore
\begin{equation}
\det(X^{(1)} \cdots X^{(\ell)})
\;=\; 
\det(P_{1}) \det(P_{\ell+1})^{-1}
  \;.
  \label{eq.95876739}
\end{equation}
We also have
\begin{equation}
X_{\rm trans}^{(1)} \cdots X_{\rm trans}^{(\ell)}
\;=\; 
\prod_{\alpha=1}^{\ell}
P_{\alpha} 
\, [\widehat{I}_{n_{\alpha} n_{\alpha+1}} 
- 
P_{\alpha}^{-1}
\psibar^{\alpha} (\psi^{\alpha+1})^{\rm T}
P_{\alpha+1}
] \,
P_{\alpha+1}^{-1}
\end{equation}
where the product is taken from left ($\alpha=1$) to right ($\alpha=\ell$).

Let us now change variables from $(\psi^{\alpha}, \psibar^{\alpha})$
to $({\psi'}^{\alpha}, \psibar'{}^{\alpha})$ defined by
\begin{eqnarray}
   {\psi'}^{\alpha}
   & = &
   \cases{ P_{\alpha}^{\rm T} \psi^{\alpha}
                       & for $2 \le \alpha \le \ell$  \cr
           \noalign{\vskip 6pt}
           P_{\ell+1}^{\rm T} \psi^1
                       & for $\alpha = 1$  \cr
         }
    \\ [2mm]
   \psibar'{}^{\alpha}  
   & = &
   P_{\alpha}^{-1} \psibar^{\alpha}
\end{eqnarray}
The Jacobian is
$(\det P_1^{-1}) (\det P_{\ell+1}) = \det(X^{(1)} \cdots X^{(\ell)})^{-1}$
[using \reff{eq.95876739}].
In the new variables we have
(dropping now the primes from the notation)
\begin{equation}
X_{\rm trans}^{(1)} \cdots X_{\rm trans}^{(\ell)}
\;=\; 
P_{1}
\left(
\prod_{\alpha=1}^{\ell}
[\widehat{I}_{n_{\alpha} n_{\alpha+1}} 
- \psibar^{\alpha} (\psi^{\alpha+1})^{\rm T} ]
\right)
P_{\ell+1}^{-1}
  \;,
\end{equation}
so again using \reff{eq.95876739} we see that
the translated determinant is given by
\begin{equation}
\det(
X_{\rm trans}^{(1)} \cdots X_{\rm trans}^{(\ell)}
)  \;=\;
\det( X^{(1)} \cdots X^{(\ell)} ) \:
(\det M_{n_1, \ldots, n_{\ell}})
\end{equation}
where the matrix $M_{n_1, \ldots, n_{\ell}}$
depends only on the Grassmann variables:
\be
M_{n_1, \ldots, n_{\ell}}
\;=\;
\prod_{\alpha=1}^{\ell}
[\widehat{I}_{n_{\alpha} n_{\alpha+1}} 
- \psibar^{\alpha} (\psi^{\alpha+1})^{\rm T} ]
   \;.
 \label{def.Mn1nl}
\ee
Therefore
\begin{eqnarray}
   & & \!\!\!
\det[(\partial^{(1)} \cdots \partial^{(\ell)} )_{IJ}]
\,
\det( X^{(1)} \cdots X^{(\ell)} )^s
       \nonumber \\
   & &
   \quad =\;
   \epsilon(I,J) \, 
\det( X^{(1)} \cdots X^{(\ell)} )^{s-1} \,
(-1)^k
\int \! \scrd_{n_1}(\psi^1,\psibar^1) \cdots
\scrd_{n_{\ell}}(\psi^{\ell},\psibar^{\ell})
\nonumber
\\
&& \qquad \times\;
\Big( \prod (P_{1} \psibar^1) (P_{\ell+1}^{\rm -T} \psi^1) \Big)_{I^c,J^c}
\:
\exp \!\left[ \sum_{\alpha=2}^{\ell} \psibar^{\alpha} \psi^{\alpha} \right]
\det(M_{n_1, \ldots, n_{\ell}})^s
   \;.  \qquad
\end{eqnarray}

To evaluate $\det(M_{n_1, \ldots, n_{\ell}})$,
we shall use Corollary~\ref{cor.lowrankrect},
which is a variant of the low-rank perturbation lemma
that is specially adapted to matrices of the form \reff{def.Mn1nl}.
Before doing so, it is convenient to split the vectors $\psi^{\alpha}$ 
and $\psibar^{\alpha}$ as
\begin{subeqnarray}
(\psi^{\alpha}_1, \ldots, \psi^{\alpha}_{n_{\alpha}}) & = &
(\lambda^{\alpha}_1, \ldots, \lambda^{\alpha}_{n_1} , 
\zeta^{\alpha}_1, \ldots, \zeta^{\alpha}_{m_{\alpha}})
  \\[1mm]
(\psibar^{\alpha}_1, \ldots, \psibar^{\alpha}_{n_{\alpha}}) & = &
(\lambdabar^{\alpha}_1, \ldots, \lambdabar^{\alpha}_{n_1} , 
\zetabar^{\alpha}_1, \ldots, \zetabar^{\alpha}_{m_{\alpha}})
   \label{eq.split.psi}
\end{subeqnarray}
where $m_\alpha \bydef n_\alpha - n_1$
(note in particular that $\psi^1 = \lambda^1$
 and $\psibar^1 = \lambdabar^1$).
We now apply Corollary~\ref{cor.lowrankrect}
with $x_\alpha = \psibar^\alpha$, $y_\alpha = \psi^{\alpha+1}$
and $\epsilon = -1$ to obtain
\be
   \det(M_{n_1, \ldots, n_{\ell}})  \;=\;  (\det N)^{-1}
\ee
where the $\ell\times\ell$ matrix $N$ is defined by
\be
   N_{\alpha \beta}
   \;=\;
   \cases{ \sum\limits_{i=1}^{m_{\alpha+1,\beta}}
              \zeta^{\alpha+1}_i \zetabar^\beta_i
                   & \hbox{\rm if } $\alpha < \beta$  \cr
           \noalign{\vskip 6pt}
           \delta_{\alpha\beta} -
              \lambda^{\alpha+1 \, \rm T} \lambdabar^\beta
                   & \hbox{\rm if } $\alpha \ge \beta$  \cr
         }
\label{apply.cor.lowrankrect.defN}
\ee
and $m_{\alpha,\beta} \bydef
     \min\limits_{\alpha \le \gamma \le \beta} m_\gamma$.
We now fall into the conditions for the application of 
Lemma~\ref{lem.pbABp.Alan} with
$\eta = \psi^1 = \lambda^1$, $\etabar = \psibar^1 = \lambdabar^1$
and $\theta = \{\lambda^2,\ldots,\lambda^\ell,
                \lambdabar^2,\ldots,\lambdabar^\ell\}$
(here the variables $\zeta$ and $\zetabar$ just go for the ride), yielding
\begin{eqnarray}
   & &
\det[(\partial^{(1)} \cdots \partial^{(\ell)} )_{IJ}]
\,
\det[ X^{(1)} \cdots X^{(\ell)} ]^s
       \nonumber \\[2mm]
   & &
   \qquad =\;
   \epsilon(I,J) \, 
\det[ X^{(1)} \cdots X^{(\ell)} ]^{s-1}
\det[ (X^{(1)} \cdots X^{(\ell)})_{I^c,J^c} ]
\, (-1)^k
\nonumber
\\
&&\qquad \qquad \times\;
\int \! \scrd_{n_1}(\lambda^1,\lambdabar^1) \cdots
           \scrd_{n_1}(\lambda^{\ell},\lambdabar^{\ell}) \,
        \scrd_{m_2}(\zeta^2,\zetabar^2) \cdots
           \scrd_{m_\ell}(\zeta^{\ell},\zetabar^{\ell}) \,
\nonumber
\\
&&\qquad \qquad \times
\left( \prod_{j=k+1}^{n_1} \lambdabar^1_j \lambda^1_j \right)
\exp \!\left[ \sum_{\alpha=2}^{\ell}
   (\lambdabar^{\alpha \rm T} \lambda^{\alpha} +
   \zetabar^{\alpha \rm T} \zeta^{\alpha})
 \right]
   \det(N)^{-s}
   \;.
\end{eqnarray}
We have thus represented
$\det[(\partial^{(1)} \cdots \partial^{(\ell)})_{IJ}] \,
        \det(X^{(1)} \cdots X^{(\ell)})^s$
as the desired quantity
$\det (X^{(1)} \cdots X^{(\ell)})^{s-1}
        \epsilon(I,J) \det[(X^{(1)} \cdots X^{(\ell)})_{I^c J^c}]$
multiplied by a purely combinatorial factor $b_{n_1,\ldots,n_\ell;k}(s)$
that is given as a Grassmann integral:
\begin{eqnarray}
b_{n_1, \ldots, n_{\ell};k}(s)
&:=&
(-1)^k
\int \! \scrd_{n_1}(\lambda^1,\lambdabar^1) \cdots
           \scrd_{n_1}(\lambda^{\ell},\lambdabar^{\ell}) \,
        \scrd_{m_2}(\zeta^2,\zetabar^2) \cdots
           \scrd_{m_\ell}(\zeta^{\ell},\zetabar^{\ell})
\nonumber
\\
&&\quad \times
\left( \prod_{j=k+1}^{n_1} \lambdabar^1_j \lambda^1_j \right)
\exp \!\left[ \sum_{\alpha=2}^{\ell}
   (\lambdabar^{\alpha \rm T} \lambda^{\alpha} +
   \zetabar^{\alpha \rm T} \zeta^{\alpha})
 \right]
   \det(N)^{-s}
   \;. \qquad
 \label{eq.star5180}
\end{eqnarray}

In order to handle the factor
$\prod_{j=k+1}^{n_1} \lambdabar^1_j \lambda^1_j$,
it is convenient to further split the fields $\lambda$ and $\lambdabar$ as
\begin{subeqnarray}
(\lambda^{\alpha}_1, \ldots, \lambda^{\alpha}_{n_1}) & = &
(\lambda^{\alpha}_1, \ldots, \lambda^{\alpha}_{k} , 
\chi^{\alpha}_1, \ldots, \chi^{\alpha}_{n_1-k})
  \\[1mm]
(\lambdabar^{\alpha}_1, \ldots, \lambdabar^{\alpha}_{n_1}) & = &
(\lambdabar^{\alpha}_1, \ldots, \lambdabar^{\alpha}_{k} , 
\chibar^{\alpha}_1, \ldots, \chibar^{\alpha}_{n_1-k})
\end{subeqnarray}
Notice now that the overall factor 
$\prod_{j=k+1}^{n_1} \lambdabar^1_j \lambda^1_j
 = \prod_{i=1}^{n_1-k} \chibar^1_i \chi^1_i$
in the integrand kills all monomials in the expansion of the
rest of the integrand that contain
any field $\chi^1$ or $\chibar^1$.
Now, the factor $(\det N)^{-s}$ in \reff{eq.star5180} depends on the fields
$\lambda,\lambdabar,\chi,\chibar,\zeta,\zetabar$ only through products of the forms
$\{\lambda^\alpha_i \lambdabar^\beta_i\}_{1 \le \beta < \alpha \le \ell}$,
$\{\lambda^1_i \lambdabar^\beta_i\}_{1 \le \beta \le \ell}$,
$\{\chi^\alpha_i \chibar^\beta_i\}_{1 \le \beta < \alpha \le \ell}$,
$\{\chi^1_i \chibar^\beta_i\}_{1 \le \beta \le \ell}$
and $\{\zeta^\alpha_i \zetabar^\beta_i\}_{2 \le \alpha \le \beta \le \ell}$,
while the exponential depends only on combinations
$\lambdabar^\alpha \lambda^\alpha$, $\chibar^\alpha \chi^\alpha$
and $\zetabar^\alpha \zeta^\alpha$ that are ``charge-neutral''
in each field separately.
Therefore, the only monomials in the expansion of $(\det N)^{-s}$
that can contribute to the integral must also be
charge-neutral in each field separately.
But since no monomial containing any $\chi^1$ can arise,
it is impossible to make such a charge-neutral combination
using any other $\chi$ or $\chibar$
since all such terms
are of the form $\chi^\alpha \chibar^\beta$ with $\beta < \alpha$.
In a similar way, the terms
$\zeta^\alpha \zetabar^\beta$ with $2 \le \alpha < \beta \le \ell$
cannot contribute.
(The combinations $\zeta^\alpha \zetabar^\alpha$ do survive.)
The integral \reff{eq.star5180} will therefore be unchanged
if we replace $N$ by a new matrix $N'$ in which all these
``forbidden combinations'' are set to zero:
\be
   N'_{\alpha \beta}
   \;=\;
   \cases{ \zeta^{\beta \rm T} \zetabar^\beta
                   & \hbox{\rm if } $\alpha = \beta-1$  \cr
           \noalign{\vskip 6pt}
           0      & \hbox{\rm if } $\alpha < \beta-1$  \cr
           \noalign{\vskip 6pt}
           \delta_{\alpha\beta} -
              \lambda^{\alpha+1 \, \rm T} \lambdabar^\beta
                   & \hbox{\rm if } $\alpha \ge \beta$  \cr
         }
\label{apply.cor.lowrankrect.defNprime}
\ee
Note that the matrix $N'$ is lower Hessenberg
(i.e.\ has zero entries above the first superdiagonal).
Summarizing, we have
\begin{eqnarray}
& &
b_{n_1, \ldots, n_{\ell};k}(s)
   \;=\;
(-1)^k
\int \! \scrd_{k}(\lambda^1,\lambdabar^1) \cdots
           \scrd_{k}(\lambda^{\ell},\lambdabar^{\ell}) \,
   \nonumber \\
& & \qquad \times \;
        \scrd_{n_1-k}(\chi^1,\chibar^1) \cdots
           \scrd_{n_1-k}(\chi^{\ell},\chibar^{\ell}) \,
        \scrd_{m_2}(\zeta^2,\zetabar^2) \cdots
           \scrd_{m_\ell}(\zeta^{\ell},\zetabar^{\ell})
   \nonumber \\
& & \qquad \times \;
\left( \prod_{i=1}^{n_1-k} \chibar^1_i \chi^1_i \right)
\exp \!\left[ \sum_{\alpha=2}^{\ell}
   (\lambdabar^{\alpha \rm T} \lambda^{\alpha} +
    \chibar^{\alpha \rm T} \chi^{\alpha} +
    \zetabar^{\alpha \rm T} \zeta^{\alpha})
 \right]
   \det(N')^{-s}
   \;. \qquad
 \label{eq.star5180.bis}
\end{eqnarray}

Since $N'$ does not contain $\chi$ or $\chibar$,
we can immediately perform the integrations
over these fields, yielding 1.

To make the indices in the matrix $N'$ look nicer,
we perform the change of variables from $\lambda$ to $\lambda'$ defined by
$(\lambda')^\alpha = \lambda^{\alpha+1}$ for $1 \le \alpha \le \ell$
(and recalling that $\lambda^{\ell+1}$ is a shorthand for $\lambda^1$);
the variables $\lambdabar$ are left as is.
The Jacobian is $(-1)^{k(\ell-1)}$.
So, dropping primes,
we have
\begin{eqnarray}
& &
b_{n_1, \ldots, n_{\ell};k}(s)
   \nonumber \\[1mm]
& & \qquad = \;
(-1)^k
\int \! \scrd_{k}(\lambda^1,\lambdabar^1) \cdots
           \scrd_{k}(\lambda^{\ell},\lambdabar^{\ell}) \,
        \scrd_{m_2}(\zeta^2,\zetabar^2) \cdots
           \scrd_{m_\ell}(\zeta^{\ell},\zetabar^{\ell})
\nonumber
\\
&&\qquad\qquad\qquad\qquad \;\times\;
\exp \!\left[ \sum_{\alpha=2}^{\ell}
   (\lambdabar^{\alpha \rm T} \lambda^{\alpha-1} +
    \zetabar^{\alpha \rm T} \zeta^{\alpha})
 \right]
   \det(N'')^{-s}
   \qquad
 \label{eq.star5180.bisbis}
\end{eqnarray}
where
\be
   N''_{\alpha \beta}
   \;=\;
   \cases{ \zeta^{\beta \rm T} \zetabar^\beta
                   & \hbox{\rm if } $\alpha = \beta-1$  \cr
           \noalign{\vskip 6pt}
           0      & \hbox{\rm if } $\alpha < \beta-1$  \cr
           \noalign{\vskip 6pt}
           \delta_{\alpha\beta} -
              \lambda^{\alpha \, \rm T} \lambdabar^\beta
                   & \hbox{\rm if } $\alpha \ge \beta$  \cr
         }
\label{apply.cor.lowrankrect.defNdoubleprime}
\ee

Finally, in order to put our scalar products of complex fermions
in the standard forms $\zetabar^{\rm T} \zeta$ and $\lambdabar^{\rm T} \lambda$,
we anticommute all the bilinears (obtaining a minus sign);
and in order to keep the indices in a natural notation,
we replace $N''$ by its transpose.
After a renaming $\lambda \to \psi$ and $\zeta \to \eta$,
the result is \reff{eq.def.bn1nk}/\reff{eq.def.bn1nk.2},
where $(N'')^{\rm T} = I_\ell + M$.

This proves part (i) of Theorem~\ref{thm.multirectcayley}
when the $X^{(\alpha)}$ are real matrices
of rank $\min(n_\alpha, n_{\alpha+1})$
lying in a sufficiently small neighborhood of
$\widehat{I}_{n_\alpha n_{\alpha+1}}$,
and $s$ is a positive integer;
the general validity of the identity then follows from
Proposition~\ref{prop.equiv.symbolic_s}.
\qed

We now turn to the evaluation of the Grassmann integral
\reff{eq.def.bn1nk}/\reff{eq.def.bn1nk.2} for
$b_{n_1, \ldots, n_{\ell};k}(s)$.  
The integrand depends on the fields only through the scalar products
$\{ \psibar^{\alpha\rm T} \psi^\beta \}_{\alpha,\beta=1}^\ell$
and $\{ \etabar^{\alpha\rm T} \eta^\alpha \}_{\alpha=2}^\ell$.
We can therefore apply Proposition~\ref{prop.intscalar}
once to the entire set of variables $\psi,\psibar$
and separately to the variables $\eta^\alpha,\etabar^\alpha$
for each $\alpha$ ($2 \le \alpha \le \ell$).
We therefore introduce indeterminates
$X = (x_{\alpha\beta})_{\alpha,\beta=1}^\ell$
and $y = (y_\alpha)_{\alpha=2}^{\ell}$,
along with the corresponding differential operators
$\partial/\partial x_{\alpha\beta}$ and $\partial/\partial y_\alpha$.
Using Proposition~\ref{prop.intscalar}
in the form \reff{eq.prop.intscalar.b}, we obtain
(after renaming $y_2,\ldots,y_\ell$ as $y_1,\ldots,y_{\ell-1}$)
\be
b_{n_1, \ldots, n_{\ell};k}(s)
 \;=\;
(-1)^{k \ell}
(\det \widehat{M})^{-s}
\exp \!\left[ \sum_{\alpha=1}^{\ell-1} 
\left(
\frac{\partial}
{\partial x_{\alpha+1 ,\alpha}}
+
\frac{\partial}
{\partial y_{\alpha}}
\right)
\right]
\left.
\det(X)^k
\prod_{\alpha=1}^{\ell-1} 
y_{\alpha}^{m_{\alpha}}
\right|_{X=y=0}
   ,
\label{eq.476875686.tris}
\ee
where now $\widehat{M}$ reads
\be
\widehat{M}
\;=\;
\left(
\begin{array}{ccccc}
1+\frac{\partial}
{\partial x_{11}}  &  \frac{\partial}
{\partial x_{12}}   &  \frac{\partial}
{\partial x_{13}}  & \cdots  & \frac{\partial}
{\partial x_{1\ell}} \\[1mm]
-\frac{\partial}
{\partial y_{1}}  & 1+\frac{\partial}
{\partial x_{22}}  &  \frac{\partial}
{\partial x_{23}}  & \cdots  & \frac{\partial}
{\partial x_{2\ell}} \\[1mm]
0       & -\frac{\partial}
{\partial y_{2}}  & 1+\frac{\partial}
{\partial x_{33}} & \cdots  & \frac{\partial}
{\partial x_{3\ell}} \\[1mm]
\vdots  & \ddots  & \ddots & \ddots  &      \\[1mm]
0       & \cdots  & 0      & -\frac{\partial}
{\partial y_{\ell-1}} & 1+\frac{\partial}
{\partial x_{\ell\ell}}
\end{array}
\right)
  \,.
\ee
We can now apply the translation formula \reff{eq.translation}
to the exponential of the differential operator:
this transforms
$\prod_{\alpha=1}^{\ell-1} y_{\alpha}^{m_{\alpha}}$
into
$\prod_{\alpha=1}^{\ell-1} (1+y_{\alpha})^{m_{\alpha}}$
and $x_{\alpha\beta}$ into
$x'_{\alpha\beta}=x_{\alpha\beta}+\delta_{\alpha,\beta+1}$.
Thus
\be
b_{n_1, \ldots, n_{\ell};k}(s)
 \;=\;
(-1)^{k \ell}
(\det \widehat{M})^{-s}
\,
\left.
\det(X')^k
\prod_{\alpha=1}^{\ell-1} 
(1+y_{\alpha})^{m_{\alpha}}
\right|_{X=y=0}
   \;.
\label{eq.476875686.quadris}
\ee
Finally, we observe that $\widehat{M}$ does not contain
the differential operators $\partial/\partial x_{\alpha\beta}$
with $\alpha>\beta$, so in $X'$ we can set those variables $x_{\alpha\beta}$
to zero immediately.
We thus have \reff{eq.476875686.quadris} where now
\be
X'
\;=\;
\left(
\begin{array}{ccccc}
x_{11}  & x_{12}  & x_{13}  & \cdots & x_{1\ell} \\
1       & x_{22}  & x_{23}  & \cdots & x_{2\ell} \\
0       &    1    & x_{33}  & \cdots & x_{3\ell} \\
\vdots  & \ddots  & \ddots  & \ddots &        \\
0       & \cdots  & 0       & 1      & x_{\ell\ell}
\end{array}
\right)
  \,.
\ee


Since the evaluation of \reff{eq.476875686.quadris}
will involve a recursion in $\ell$,
it is convenient to introduce an infinite set of indeterminates
$\{ x_{\alpha \beta} \}_{1 \leq \alpha \leq \beta < \infty}$,
along with the corresponding set of differential operators
$\de_{\alpha \beta}=\de/\de x_{\alpha \beta}$,
as well as another infinite set of indeterminates
$\{a_{\alpha} \}_{1 \leq \alpha \le n-1}$
(only finitely many of these will play any role at any given stage).
Then define, for each $\ell \ge 1$, the quantities
\be
D_\ell(a) 
\;=\;
\det
\left(
\begin{array}{ccccc}
1+\de_{11} & \de_{12}  & \de_{13}  & \cdots & \de_{1\ell} \\
a_1       & 1+\de_{22} & \de_{23}  & \cdots & \de_{2\ell} \\
0         &    a_2    & 1+\de_{33} & \cdots & \de_{3\ell} \\
\vdots    & \ddots    & \ddots     & \ddots & \\
0         & \cdots    & 0          & a_{\ell-1} & 1+\de_{\ell\ell}
\end{array}
\right)
\ee
and
\be
X_\ell
\;=\;
\det
\left(
\begin{array}{ccccc}
x_{11}  & x_{12}  & x_{13}  & \cdots & x_{1\ell} \\
1       & x_{22}  & x_{23}  & \cdots & x_{2\ell} \\
0       &    1    & x_{33}  & \cdots & x_{3\ell} \\
\vdots  & \ddots  & \ddots  & \ddots &        \\
0       & \cdots  & 0       & 1      & x_{\ell\ell}
\end{array}
\right)
\ee
where we also set $D_0(a)=1$ and $X_0=1$.
Note also that $D_\ell(a)$ [resp.\ $X_\ell$]
involves only those $\de_{\alpha\beta}$ [resp.\ $x_{\alpha\beta}$]
with $\alpha \le \beta \le \ell$.

Given a formal indeterminate $s$ and a nonnegative integer $k$,
our goal in the remainder of this subsection is to compute the expression
\be
\left.
D_\ell(a)^{-s} X_\ell^k \right|_{X=0}
  \;,
  \label{eq.2876597649}
\ee
which abstracts the relevant features of \reff{eq.476875686.quadris}
[as $s$ is a formal variable, the choice between $s$ and $-s$ has no
special role and is made here just for convenience].
Since $D_\ell(a)$ is a polynomial
in the quantities $\{ \de_{\alpha\beta} \}$ and $\{a_\alpha\}$
with constant term 1,
$D_\ell(a)^{-s}$ is here to be understood as the series
\be
\label{eq.769872784}
D_\ell(a)^{-s} 
\;=\; 
\sum_{h=0}^\infty \binom{-s}{h} [D_\ell(a)-1]^h
\ee
which, when applied to $X_\ell^k$ as in \reff{eq.2876597649},
can be truncated to $h \le k\ell$.
We shall prove the following:

\begin{theorem}
\label{thm.quasitriang}
With the definitions above, we have
\be
  \left.  D_\ell(a)^{-s} X_\ell^k \right|_{X=0}
  \;=\;
  k! \, \binom{-s}{k}
  \prod_{\alpha=1}^{\ell-1} \, \sum_{b=0}^k
        k! \binom{-s-b}{k-b} \frac{a_{\alpha}^{b}}{b!}
   \;.
 \label{eq.8765867}
\ee
\end{theorem}


We remark that
\be
  \sum_{b=0}^k \binom{-s-b}{k-b} \frac{z^b}{b!}
   \;=\;
  \binom{-s}{k} \ofo(-k;s;z)
\ee
although we will not use this.

What we shall actually need is a specific corollary of
Theorem~\ref{thm.quasitriang}.
Let us introduce the variables $\{y_{\alpha} \}_{1 \leq \alpha < \infty}$
and the associated derivatives $\widehat{\de}_{\alpha} = \de/\de y_{\alpha}$,
as well as the further indeterminates (or nonnegative integers)
$\{m_{\alpha} \}_{1 \leq \alpha < \infty}$.

\begin{corollary}
\label{cor.solvingmulti}
With the definitions above, we have
\be
\left.
D_\ell(-\widehat{\de})^{-s} X_\ell^k 
\prod_{\alpha=1}^{\ell-1} (1+y_{\alpha})^{m_{\alpha}}
\right|_{X=y=0}
  \;=\;
(-1)^{k\ell}
\prod_{\alpha=0}^{\ell-1} \,
\prod_{i=0}^{k-1}
(s+m_{\alpha}+i)
\label{eq.cor.solvingmulti}
\ee
with the convention $m_0 = 0$.
\end{corollary}

Given Corollary~\ref{cor.solvingmulti}, the proof of
Theorem~\ref{thm.multirectcayley}, part (ii) is a triviality:

\bigskip
\noindent
{\sc Proof of  Theorem~\ref{thm.multirectcayley}, part (ii),
   given Corollary~\ref{cor.solvingmulti}.\ }
It suffices to recognize $\det[\widehat{M}]$ and $\det(X')$
as the operators $D_{\ell}(-\widehat{\de})$ and $X_{\ell}$
in Corollary~\ref{cor.solvingmulti}.
The sign $(-1)^{k \ell}$ in \reff{eq.476875686.quadris}
combines with the one in (\ref{eq.cor.solvingmulti}),
leaving exactly the prefactor claimed in Theorem~\ref{thm.multirectcayley}.
\qed

Let us next show how to deduce Corollary~\ref{cor.solvingmulti}
from Theorem~\ref{thm.quasitriang}:

\bigskip
\noindent
{\sc Proof of Corollary~\ref{cor.solvingmulti},
   given Theorem~\ref{thm.quasitriang}.\ }
We evaluate the left-hand side of \reff{eq.cor.solvingmulti}
by using \reff{eq.8765867} with $a_\alpha$
replaced by $-\widehat{\de}_\alpha$.
Since
\be
   \left. (-\widehat{\de}_\alpha)^b (1+y_\alpha)^{m_\alpha}
       \right|_{y_\alpha=0} 
   \;=\;
   (-1)^b  b! \binom{m_{\alpha}}{b}
\ee
and everything is factorized over $\alpha$, we obtain
\be
\left.
D_\ell(-\widehat{\de})^{-s} X_\ell^k 
\prod_{\alpha=1}^{\ell-1} (1+y_{\alpha})^{m_{\alpha}}
\right|_{X=y=0}
  \;=\;
  k! \, \binom{-s}{k}
  \prod_{\alpha=1}^{\ell-1} \, \sum_{b=0}^k
        k! \binom{-s-b}{k-b} 
   (-1)^b  \binom{m_{\alpha}}{b}
  \;.
\ee
Then for each $\alpha$ we have
\begin{eqnarray}
k!
\sum_{b=0}^k
\binom{-s-b}{k-b}
(-1)^b
\binom{m_{\alpha}}{b}
& = &
k!
\sum_{b=0}^k
(-1)^{k-b}
\binom{s+k-1}{k-b}
(-1)^b
\binom{m_{\alpha}}{b}
\nonumber  \\
  & = &
(-1)^k
k!
\sum_{b=0}^k
\binom{s+k-1}{k-b}
\binom{m_{\alpha}}{b}
\nonumber  \\
  & = &
(-1)^k
k!
\binom{s+k+m_{\alpha}-1}{k}
\nonumber  \\
  & = &
(-1)^k
   \prod_{i=0}^{k-1} (s+m_\alpha + i)
\end{eqnarray}
where we used the Chu--Vandermonde convolution (Lemma~\ref{lemma.vandermonde})
in going from the second to the third line.
The prefactor $k! \binom{-s}{k}$
gives an analogous contribution with $m_0 = 0$.
\qed

Finally, we turn to the proof of Theorem~\ref{thm.quasitriang}.
We shall need two main lemmas:
one that essentially provides an inductive step, and another
dealing with a special sum of multinomial coefficients.
Henceforth we shall write $D_\ell$ as a shorthand for $D_\ell(a)$.

We start with a pair of easy recursive formulae,
obtained by expansion of the determinant on the last column:

\begin{lemma}
   \label{lemma.expDetDn}
\be
D_\ell  \;=\; D_{\ell-1}  \,+\, \sum_{\alpha=1}^\ell (-1)^{\ell-\alpha}
a_{\alpha} \cdots a_{\ell-1} D_{\alpha-1} \de_{\alpha \ell}
\ee
(the empty product $a_{\alpha} \cdots a_{\ell-1}$ for $\alpha=\ell$ should
be understood as 1)
and
\be
X_\ell \;=\; \sum_{\alpha=1}^\ell (-1)^{\ell-\alpha}
X_{\alpha-1} 
x_{\alpha \ell}
   \;.
 \label{eq.det.Hessenberg}
\ee
\end{lemma}

\noindent
Related formulae for the determinant of a Hessenberg matrix
can be found in \cite{Tamm_1,Tamm_2}.


The induction lemma is the following:

\begin{lemma}
  \label{lemma.quatri.1}
Let $m \ge 1$ and $c_1,\ldots,c_m \ge 0$ be integers,
and let $t$ be an indeterminate.  Then
\begin{eqnarray}
D_m^t
\left.
\prod_{\alpha=1}^m 
D_{\alpha}^{c_{\alpha}}
\prod_{\alpha=1}^m 
X_{\alpha}^{c_{\alpha}}
\right|_{x_{\alpha m}=0}
&=&
\!\!\!\!\!\!
\sum_{
\begin{scarray}
b_1,\ldots,b_m \ge 0
\\
b_1+\cdots+b_m=c_m
\end{scarray}
}
\!\!\!\!\!\!
\binom{t+c_m}{c_m}c_m!
\binom{c_m}{b_1,\ldots, b_m}
\prod_{\alpha=1}^{m-1} 
a_{\alpha}^{\sum_{\beta=1}^{\alpha} b_{\beta}}
\nonumber
\\
&&
\quad
\times\;
D_{m-1}^t
\prod_{\alpha=1}^{m-1} 
D_{\alpha}^{c_{\alpha}+b_{\alpha+1}}
\prod_{\alpha=1}^{m-1} 
X_{\alpha}^{c_{\alpha}+b_{\alpha+1}}
  \label{eq.lemma.quatri.1}
\end{eqnarray}
\end{lemma}

\proof
Notice, first of all, that the factors $D_{\alpha}$ and $X_{\alpha}$
for $\alpha<m$ do not play any role, i.e.\ we can rewrite
the left-hand side of \reff{eq.lemma.quatri.1} as
\be
\label{eq.74698680870}
\left(
\prod_{\alpha=1}^{m-1} 
D_{\alpha}^{c_{\alpha}}
\right)
\left(
\left.
D_{m}^{c_{m}+ t}
X_m^{c_m}
\right|_{x_{\alpha m}=0}
\right)
\left(
\prod_{\alpha=1}^{m-1} 
X_{\alpha}^{c_{\alpha}}
\right)
\ee
and concentrate on the central factor alone.
To compute
$\displaystyle{
   \left.  D_{m}^{c_{m}+ t} X_m^{c_m} \right|_{x_{\alpha m}=0}
  }$,
we expand $D_m^{c_m+t}$ and $X_m^{c_m}$
using Lemma~\ref{lemma.expDetDn}:
\begin{eqnarray}
D_m^{t+c_m}
&=&
\sum_{b_1,\ldots,b_m \ge 0}
\binom{t+c_m}{b_1+ \cdots +b_m}
\binom{b_1+ \cdots +b_m}{b_1,\ldots,b_m}
D_{m-1}^{t+c_m-(b_1+ \cdots +b_m)}
\nonumber \\
&& \qquad \times\;
\prod_{\alpha=1}^m 
\left[
(-1)^{m-\alpha}
(a_{\alpha} \cdots a_{m-1})
D_{\alpha-1} 
\de_{\alpha m}
\right]^{b_{\alpha}}
\\[3mm]
X_m^{c_m}
&=&
\sum_{\begin{scarray}
b_1,\ldots,b_m \ge 0 \\
b_1+\cdots+b_m=c_m
\end{scarray}
}
\binom{c_m}{b_1,\ldots,b_m}
\prod_{\alpha=1}^m 
\left[
(-1)^{m-\alpha}
x_{\alpha m}
X_{\alpha-1} 
\right]^{b_{\alpha}}
\end{eqnarray}
Since $\de_{\alpha m}^b x_{\alpha m}^{b'}|_{x_{\alpha m}=0}
 = \delta_{b,b'}\, b!$, 
the two sets of summation indices,
when combined inside
$\displaystyle{
   \left.  D_{m}^{c_{m}+ t} X_m^{c_m} \right|_{x_{\alpha m}=0}
  }$,
must coincide, and we get
\begin{eqnarray}
\left.
D_{m}^{c_{m}+ t}
X_m^{c_m}
\right|_{x_{\alpha m}=0}
&=&
\sum_{\begin{scarray}
b_1,\ldots,b_m \ge 0 \\
b_1+\cdots+b_m=c_m
\end{scarray}}
\binom{t+c_m}{c_m}
c_m!
\binom{c_m}{b_1,\ldots,b_m}
D_{m-1}^{t}
\nonumber
\\
&& \qquad\times\;
\prod_{\alpha=1}^m 
\left[
(a_{\alpha} \cdots a_{m-1})
D_{\alpha-1} 
\right]^{b_{\alpha}}
\prod_{\alpha=1}^m 
\left(
X_{\alpha-1} 
\right)^{b_{\alpha}}
   \;.
\end{eqnarray}
In the two final products, we can drop the factors $D_0^{b_1}$ and
$X_0^{b_1}$ since $D_0=X_0=1$.
Reintroducing the missing factors from (\ref{eq.74698680870}), we
obtain \reff{eq.lemma.quatri.1}.
\qed

We will now apply Lemma~\ref{lemma.quatri.1} for $\ell$ ``rounds'',
starting with the initial conditions $m=\ell$, $t=-s-k$,
$c_1 = \ldots = c_{\ell-1} = 0$ and $c_\ell = k$.
At round $i$ we have $m = \ell+1-i$.
Let us denote by $c^i_1,\ldots,c^i_{\ell+1-i}$ the parameters $\{c_j\}$
immediately before entering round $i$,
and let us denote by $b^i_1,\ldots,b^i_{\ell+1-i}$
the summation indices in round $i$.
We therefore have the initial conditions
\be
   c^1_1 = \ldots = c^1_{\ell-1} = 0, \quad c^1_\ell = k
 \label{eq.initialcond}
\ee
and the recursion
\be
   c^{i+1}_j  \;=\;  c^i_j \,+\, b^i_{j+1}
      \qquad\hbox{for } 1 \le j \le \ell-i
 \label{eq.bc.recursion}
\ee
[cf.\ \reff{eq.lemma.quatri.1}].
The summation indices $b^i_1,\ldots,b^i_{\ell+1-i}$ obey the constraint
\be
   \sum_{j=1}^{\ell+1-i}  b^i_j  \;=\;  c^i_{\ell+1-i}  \;.
 \label{eq.bc.constraint}
\ee
Using \reff{eq.initialcond} and \reff{eq.bc.recursion}
we prove by induction that
\be
   c^i_{\ell+1-i}  \;=\;  \sum_{h=1}^{i-1} b^h_{\ell+2-i}
     \qquad\hbox{for } 2 \le i \le \ell  \;.
\ee
It is convenient to arrange the summation indices
$\{b^i_j\}_{i+j \le \ell+1}$ into a matrix
\be
B  \;=\;
\left(
\begin{array}{ccccc}
b^1_1 & b^1_2 & b^1_3 & \cdots & \multicolumn{1}{c|}{b^1_\ell} \\[1mm]
\cline{5-5}
b^2_1 & b^2_2 & \cdots & \multicolumn{1}{c|}{b^2_{\ell-1}} & \\[1mm]
\cline{4-4}
b^3_1 & \vdots & \multicolumn{1}{c|}{\ddotsinverse} && \\[1mm]
\cline{3-3}
\vdots & \multicolumn{1}{c|}{b^{\ell-1}_2} &&& \\[1mm]
\cline{2-2}
\multicolumn{1}{c|}{b^\ell_1} &&&& \\[1mm]
\cline{1-1}
\end{array}
\right)
   \label{eq.def.Bmatrix}
\ee
in which the first row sums to $k$
and, for $2 \le i \le \ell$,
row $i$ and column $\ell+2-i$ have the same sum.
Such matrices can be characterized as follows:

\begin{lemma}
   \label{lemma.Bmatrix}
Fix integers $\ell \ge 1$ and $k \ge 0$.
For a matrix $B = (b^i_j)_{i+j \le \ell+1}$
of nonnegative integers as in \reff{eq.def.Bmatrix},
the following conditions are equivalent:
\begin{itemize}
   \item[(a)] The first row sums to $k$ and, for $2 \le i \le \ell$,
row $i$ and column $\ell+2-i$ have the same sum.
   \item[(b)] For $1 \le h \le \ell$ we have
$\displaystyle{ \!\!\!\sum\limits_{\begin{scarray}
                                1 \le i \le h \\
                                1 \le j \le \ell+1-h
                             \end{scarray}}
                \!\! b^i_j \:=\: k}$.
   \item[(c)] There exist nonnegative integers
$b^2_\ell, b^3_{\ell-1},\ldots,b^\ell_2$ completing the matrix $B$
to one in which all the row and column sums are equal to $k$.
(Such numbers are obviously unique if they exist,
 and must lie in the interval $[0,k]$.)
\end{itemize}
\end{lemma}

\noindent
Note in particular that statement (b) with $h=\ell$
tells us that the first column sums to $k$.

\proof
(a) $\implies$ (b):  By induction on $h$.
By hypothesis the equality holds for $h=1$.
Then for $h \ge 2$ we have
\be
   \sum\limits_{\begin{scarray}
                    1 \le i \le h \\
                    1 \le j \le \ell+1-h
                \end{scarray}}
         \!\! b^i_j
   \:-\!\!
   \sum\limits_{\begin{scarray}
                    1 \le i \le h-1 \\
                    1 \le j \le \ell+2-h
                \end{scarray}}
         \!\! b^i_j
   \;=\;
   \sum\limits_{1 \le j \le \ell-h}  \!\! b^h_j
   \:-\!
   \sum\limits_{1 \le i \le h-1}  \!\! b^i_{\ell+2-h}
   \;=\;
   0 \quad
\ee
by hypothesis.

(b) $\implies$ (c):  It is easily checked that the definition
\be
   b^h_{\ell+2-h}
   \;=\;
   \sum\limits_{\begin{scarray}
                    1 \le i \le h-1 \\
                    1 \le j \le \ell+1-h
                \end{scarray}}
   \!\! b^i_j
   \qquad\hbox{for } 2 \le h \le \ell
\ee
does what is required.

(c) $\implies$ (a) is obvious.
\qed

{\bf Remark.}  An analogous equivalence holds, with the same proof
({\em mutatis mutandis}\/),
for matrices $B$ of nonnegative {\em real}\/ numbers
where $k$ is a fixed nonnegative {\em real}\/ number.
\qed

The $\ell$-fold application of Lemma~\ref{lemma.quatri.1}
with the initial conditions \reff{eq.initialcond}
gives rise to a sum over matrices $B$
satisfying the equivalent conditions (a)--(c)
of Lemma~\ref{lemma.Bmatrix}.
To each such matrix there corresponds an $\ell$-tuple
$\bar{c} = (\bar{c}^1,\ldots,\bar{c}^\ell)$
of integers in the range $[0,k]$,
where $\bar{c}^i = \sum_{j=1}^{\ell+1-i} b^i_j$
is the sum of row $i$ (it also equals $c^i_{\ell+1-i}$)
and of course $\bar{c}^1 = k$.
It will be convenient to partition the sum over matrices $B$
according to the vector $\bar{c}$,
so let us denote by $\scrb(\bar{c})$ be the set of matrices satisfying the
given conditions with the row sums $\bar{c}$:
\begin{eqnarray}
   & &
   \scrb(\bar{c})  \;=\;
   \biggl\{ 
       B = (b^i_j)_{i+j \le \ell+1} \in \N^{\ell(\ell+1)/2} \colon\;
      \sum_{j=1}^{\ell+1-i}  b^i_j = \bar{c}^i
      \hbox{ for } 1 \le i \le \ell
     \nonumber \\
   & &
   \qquad\qquad\qquad\qquad\qquad\qquad\qquad
   \hbox{ and }
      \sum_{h=1}^{i-1} b^h_{\ell+2-i} = \bar{c}^i
      \hbox{ for } 2 \le i \le \ell
   \biggr\}
   \;. \qquad
 \label{def.scrbbarc}
\end{eqnarray}
The summand is then
\be
   \left( \prod_{i=1}^\ell  \binom{t+\bar{c}^i}{\bar{c}^i}
                \, \bar{c}^i ! \,
                \binom{\bar{c}^i}{b^i_1,\ldots,b^i_{\ell+1-i}}
   \right)
   \left( \prod_{\alpha=1}^{\ell-1}
            a_\alpha^{\sum\limits_{i=1}^{\ell-\alpha}
                      \sum\limits_{j=1}^\alpha  b^i_j}
   \right)
   \;.
  \label{eq.bc.star4}
\ee
Furthermore, it follows from Lemma~\ref{lemma.Bmatrix}(b) that
\be
   \sum\limits_{i=1}^{\ell-\alpha} \sum\limits_{j=1}^\alpha  b^i_j
   \;=\;
   k \,-\, \bar{c}^{\ell+1-\alpha}
   \;.
  \label{eq.bc.star5}
\ee
Let us now show how to perform the sum over matrices $B$
with a given vector $\bar{c}$:

\begin{lemma}
  \label{lem.quatri.2}
Let $\ell \ge 1$, $k \ge 0$ and $\bar{c}^1,\ldots,\bar{c}^\ell \ge 0$
be integers, with $\bar{c}^1 = k$.  Then
\be
   \sum\limits_{B = (b^i_j)_{i+j \le \ell+1} \in \scrb(\bar{c})} \:
   \prod_{i=1}^\ell  \binom{\bar{c}^i}{b^i_1,\ldots,b^i_{\ell+1-i}}
   \;=\;
   \prod_{i=1}^\ell  \binom{k}{\bar{c}^i}
   \;.
  \label{eq.bc.star6}
\ee
\end{lemma}

Before proving Lemma~\ref{lem.quatri.2},
let us show how it can be used to complete the proof of
Theorem~\ref{thm.quasitriang}.

\bigskip
\noindent
{\sc Proof of Theorem~\ref{thm.quasitriang},
   given Lemma~\ref{lem.quatri.2}.\ }
Summing \reff{eq.bc.star4} over $B \in \scrb(\bar{c})$
and using \reff{eq.bc.star5} and \reff{eq.bc.star6}
along with $t=-s-k$ and $\bar{c}^1 = k$, we obtain
\be
   \left( \prod_{i=1}^\ell  \binom{-s-k+\bar{c}^i}{\bar{c}^i}
                \, \bar{c}^i ! \,
                \binom{k}{\bar{c}^i}
   \right)
   \left( \prod_{\alpha=1}^{\ell-1} a_\alpha^{k-\bar{c}^{\ell+1-\alpha}}
   \right)
   \;=\;
   k! \, \binom{-s}{k}
   \prod_{\alpha=1}^{\ell-1}
      \binom{-s-\widehat{c}^\alpha}{k-\widehat{c}^\alpha}
      \, k! \, {a_\alpha^{\widehat{c}^\alpha}  \over \widehat{c}^\alpha !}
\ee
where $\widehat{c}^\alpha = k-\bar{c}^{\ell+1-\alpha}$.
The sum over $\bar{c}^2,\ldots,\bar{c}^\ell$
--- or equivalently over $\widehat{c}^1,\ldots,\widehat{c}^{\ell-1}$ ---
now factorizes and gives precisely \reff{eq.8765867}.
\qed

\proofof{Lemma~\ref{lem.quatri.2}}
The sum on the left-hand side of \reff{eq.bc.star6} is nontrivial
because the ``row'' and ``column'' constraints \reff{def.scrbbarc}
are entangled.
We shall replace one of the two sets of constraints
(say, the ``column'' one) by a generating function:
that is, we introduce indeterminates $\xi_j$ ($1 \le j \le \ell$)
and consider
\be
   \Phi_{\ell,k,\bar{c}}(\xi)
   \;=\;
   \sum_{\begin{scarray}
           \{ b^i_j \}_{i+j \leq \ell+1} \\[2mm]
           \sum_{j=1}^{\ell+1-i} b^i_j = \bar{c}^i
         \end{scarray}}
   \prod_{j=1}^{\ell} \xi_j^{\sum_{i=1}^{\ell+1-j} b^i_j}
   \prod_{i=1}^\ell  \binom{\bar{c}^i}{b^i_1,\ldots,b^i_{\ell+1-i}}
  \;.
\ee
The sum now factorizes over rows:  we have
\begin{subeqnarray}
   \Phi_{\ell,k,\bar{c}}(\xi)
   & = &
   \prod_{i=1}^\ell
   \sum_{\begin{scarray}
            b^i_1,\ldots,b^i_{\ell+1-i} \ge 0 \\[1mm]
            b^i_1 + \ldots + b^i_{\ell+1-i} = \bar{c}^i
         \end{scarray}}
   \left( \prod_{j=1}^{\ell+1-i} \xi_j^{b^i_j} \right)
   \binom{\bar{c}^i}{b^i_1,\ldots,b^i_{\ell+1-i}}
         \\[1mm]
   & = &
   \prod_{i=1}^\ell
   \left( \sum_{j=1}^{\ell+1-i} \xi_j \right) ^{\! \bar{c}^i}
         \\[1mm]
   & = &
   (\xi_1 + \ldots + \xi_\ell)^k
   (\xi_1 + \ldots + \xi_{\ell-1})^{\bar{c}^2}
   (\xi_1 + \ldots + \xi_{\ell-2})^{\bar{c}^3}
      \,\cdots\,
   \xi_1^{\bar{c}^\ell}
   \,.
   \qquad
 \slabel{eq.prodxi}
\end{subeqnarray}
We must now extract the coefficient of the monomial
\be
   \xi_1^k \prod_{j=2}^\ell \xi_j^{\bar{c}^{\ell+2-j}}
   \;=\;
   \xi_1^k \, \xi_2^{\bar{c}^\ell} \, \xi_3^{\bar{c}^{\ell-1}}
     \,\cdots\, \xi_\ell^{\bar{c}^2}
\ee
in $\Phi_{\ell,k,\bar{c}}(\xi)$.
We first extract $[\xi_\ell^{\bar{c}^2}]$ from \reff{eq.prodxi}:
here $\xi_\ell$ occurs only in the first factor, so we get
$\binom{k}{\bar{c}^2}  (\xi_1 + \ldots + \xi_{\ell-1})^{k-\bar{c}^2}$
times the remaining factors, i.e.
\be
   \binom{k}{\bar{c}^2}
   (\xi_1 + \ldots + \xi_{\ell-1})^{k}
   (\xi_1 + \ldots + \xi_{\ell-2})^{\bar{c}^3}
      \,\cdots\,
   \xi_1^{\bar{c}^\ell}
   \;.
\ee
We can then extract $[\xi_{\ell-1}^{\bar{c}^3}]$ in the same way,
and so forth until the end,
yielding right-hand side of \reff{eq.bc.star6}.
\qed

Let us conclude by apologizing for the combinatorial complexity
involved in the proof of Theorem~\ref{thm.quasitriang}.
The simplicity of the final formula \reff{eq.8765867},
together with the simplicity of \reff{eq.bc.star6}
and the miraculous simplifications observed in its proof,
suggest to us that there ought to exist a much simpler and shorter
proof of Theorem~\ref{thm.quasitriang}.
But we have thus far been unable to find one.

\section{Proofs of diagonal-parametrized Cayley identities}  \label{sec.param}

In this section we prove the diagonal-parametrized Cayley identities
(Theorems~\ref{thm.para.cayley} and \ref{thm.para.sym.cayley}).
We give two proofs of each result:
the first deduces the diagonal-parametrized Cayley identity
from the corresponding Cayley identity
by a change of variables;
the second is a direct proof using a Grassmann representation
of the differential operator.

\subsection{Diagonal-parametrized ordinary Cayley identity} \label{sec.param.1}

\proofof{Theorem~\ref{thm.para.cayley}}
We change variables from $(x_{ij})_{i,j=1}^n$ to
new variables $(t_i)_{i=1}^n$ and $(y_{ij})_{1 \le i \neq j \le n}$
defined by
\begin{subeqnarray}
   t_i     & = &  x_{ii}   \\
   y_{ij}  & = &  x_{ii}^{-\alpha_i} x_{jj}^{-(1-\alpha_j)} x_{ij}
\end{subeqnarray}
We also set $y_{ii} = 1$ for all $i$
and define the matrices $T_\balpha = \diag(t_i^{\alpha_i})$,
$T_{1-\balpha} = \diag(t_i^{1-\alpha_i})$ and $Y = (y_{ij})$,
so that
\be
    X  \;=\;  T_\balpha Y T_{1-\balpha}
\ee
and hence
\be
   (\det X)^s  \;=\;  \left( \prod_{i=1}^n t_i^s \right) (\det Y)^s  \;.
 \label{eq.detTY}
\ee
A straightforward computation shows that the differential operators
(vector fields) $\partial/\partial x_{ij}$ can be rewritten in the
new variables as
\be
   {\partial \over \partial x_{ij}}  \;=\;
   \cases{  t_i^{-1} \left[ t_i \,
                        {\bigpartial \over \bigpartial t_i}
               - \alpha_i
                 \sum\limits_{k \neq i} y_{ik} \,
                     {\bigpartial \over \bigpartial y_{ik}}
               - (1-\alpha_i)
                 \sum\limits_{l \neq i} y_{li} \,
                     {\bigpartial \over \bigpartial y_{li}}
               \right]
                    &  if $i = j$  \cr
            \noalign{\vskip 4pt}
            t_i^{-\alpha_i} t_j^{-(1-\alpha_j)} \,
                     {\bigpartial \over \bigpartial y_{ij}}
                    &  if $i \neq j$  \cr
          }
 \label{def.Delta.alpha}
\ee
Let us denote by $\Delta_\balpha$ the matrix of
differential operators whose elements are given by the
right-hand side of \reff{def.Delta.alpha};
please note that each element commutes with
each other element not in the same row or column
[that is, $(\Delta_\balpha)_{ij}$ commutes with $(\Delta_\balpha)_{i'j'}$
 whenever $i \neq i'$ and $j \neq j'$].

Now \reff{eq.detTY}, considered as a function of $t = (t_i)$,
   is of the form ${\rm const} \times \prod_i t_i^s$;
   therefore, acting on \reff{eq.detTY},
   each operator $t_i \partial/\partial t_i$
   is equivalent to multiplication by $s$.
   It follows that the action of $\Delta_\balpha$ on \reff{eq.detTY}
   is identical to that of $\Delta_{\balpha,s}$ defined by
\be
   (\Delta_{\balpha,s})_{ij}  \;=\;
   \cases{  t_i^{-1} \left[ s
               - \alpha_i
                 \sum\limits_{k \neq i} y_{ik} \,
                     {\bigpartial \over \bigpartial y_{ik}}
               - (1-\alpha_i)
                 \sum\limits_{l \neq i} y_{li} \,
                     {\bigpartial \over \bigpartial y_{li}}
               \right]
                    &  if $i = j$  \cr
            \noalign{\vskip 4pt}
            t_i^{-\alpha_i} t_j^{-(1-\alpha_j)} \,
                     {\bigpartial \over \bigpartial y_{ij}}
                    &  if $i \neq j$  \cr
          }
 \label{def.Delta.balpha.s}
\ee
Furthermore, 
$(\Delta_\balpha)_{ij}$ and $(\Delta_{\balpha,s})_{ij}$
both commute with $(\Delta_{\balpha,s})_{i'j'}$
whenever $i \neq i'$ and $j \neq j'$.

Now suppose we have a product
$(\Delta_\balpha)_{i_1 j_1} \cdots (\Delta_\balpha)_{i_\ell j_\ell}$
acting on \reff{eq.detTY},
in which $i_1, \ldots, i_\ell$ are all distinct
and also $j_1, \ldots, j_\ell$ are all distinct
(this is the case that will arise when we take the determinant
    of any submatrix).
Then the rightmost factor $(\Delta_\balpha)_{i_\ell j_\ell}$
can be replaced by $(\Delta_{\balpha,s})_{i_\ell j_\ell}$;
and the same can be done for the other factors
by commuting them to the right, changing them from
$\Delta_\balpha$ to $\Delta_{\balpha,s}$, and commuting them back.
It follows that, under the given condition on the indices,
$(\Delta_\balpha)_{i_1 j_1} \cdots (\Delta_\balpha)_{i_\ell j_\ell} (\det X)^s
 = 
 (\Delta_{\balpha,s})_{i_1 j_1} \cdots (\Delta_{\balpha,s})_{i_\ell j_\ell}
 (\det X)^s$.

Translating back to the original variables $(x_{ij})_{i,j=1}^n$,
$\Delta_{\balpha,s}$ equals $D_{\balpha,s}$ defined by
\be
   (D_{\balpha,s})_{ij}  \;=\;
   \cases{  x_{ii}^{-1} \left[ s
               - \alpha_i
                 \sum\limits_{k \neq i} x_{ik} \,
                     {\bigpartial \over \bigpartial x_{ik}}
               - (1-\alpha_i)
                 \sum\limits_{l \neq i} x_{li} \,
                     {\bigpartial \over \bigpartial x_{li}}
               \right]
                    &  if $i = j$  \cr
            \noalign{\vskip 4pt}
            {\bigpartial \over \bigpartial x_{ij}}
                    &  if $i \neq j$  \cr
          }
 \label{def.D.alpha.s}
\ee
and the Cayley formula \reff{eq.cayley.2} tells us that
\be
   \det((D_{\balpha,s})_{IJ}) \, (\det X)^s  \;=\;
      s(s+1) \cdots (s+k-1) \, (\det X)^{s-1} \,
      \epsilon(I,J) \, (\det X_{I^c J^c})
      \;.
 \label{eq.para.cayley.proof1}
\ee
On the other hand,
$D_{\balpha,\bbeta,s} = \widehat{X}_\bbeta D_{\balpha,s} \widehat{X}_{1-\bbeta}$
where $\widehat{X}_\bbeta = \diag(x_{ii}^{\beta_i})$
and $\widehat{X}_{1-\bbeta} = \diag(x_{ii}^{1-\beta_i})$,
so that
\be
   \det((D_{\balpha,\bbeta,s})_{IJ})  \;=\;
      \Biggl( \prod\limits_{i \in I} x_{ii}^{\beta_i} \Biggr) \,
      \Biggl( \prod\limits_{j \in J} x_{jj}^{1-\beta_j} \Biggr) \,
      \det((D_{\balpha,s})_{IJ})
   \;.
 \label{eq.para.cayley.proof2}
\ee
Combining \reff{eq.para.cayley.proof1} and \reff{eq.para.cayley.proof2},
we obtain \reff{eq.para.cayley2}.
\qed

Now let us show an {\em ab initio}\/ proof based on
Grassmann representation of the differential operator
$\det((D_{\balpha,s})_{IJ})$.
Since $D_{\balpha,s}$ contains not only terms $\partial/\partial x_{ij}$
but also terms $x_{ij} \, \partial/\partial x_{ij}$,
the role played in Section~\ref{sec.grassmann.2}
by the translation formula \reff{eq.translation}
will here be played by the dilation-translation formula \reff{eq.transl2}:
see Section~\ref{app.gen.transl} for discussion.

\alternateproofof{Theorem~\ref{thm.para.cayley}}
For notational simplicity let us assume that the diagonal elements $x_{ii}$
are all equal to 1;  the general case can be recovered by a simple scaling.

Consider the matrix of differential operators $D_{\balpha,s}$
defined by \reff{def.D.alpha.s}.
In terms of the matrices $E^{ij}$ defined by \reff{def.Eij},
we can write
\begin{equation}
   D_{\balpha,s}
   \;=\;
   s I  \,+\,
   \sum_{i \neq j} \Big( E^{ij} - x_{ij}[\alpha_i E^{ii} + (1-\alpha_j) E^{jj}]
                   \Big)
              \frac{\partial}{\partial x_{ij}}
   \;\,.
\end{equation}
As before, we introduce Grassmann variables $\eta_i, \etabar_i$
($1 \le i \le n$) and use the representation
\begin{eqnarray}
   & &
   \det((D_{\balpha,s})_{IJ})
   \;=\;
   \epsilon(I,J)
   \int \scrd_n(\eta, \etabar)  \,
        \Big( \prod \etabar \eta \Big)_{I^c,J^c}
        \;\times
   \nonumber \\[2mm]
   & &
   \qquad\qquad
   \exp\!\left[
        s \sum\limits_i \etabar_i \eta_i
           \,+\, \sum\limits_{i \neq j} [\etabar_i \eta_j - x_{ij}
          (\alpha_i \etabar_i \eta_i + (1-\alpha_j) \etabar_j \eta_j )]
          \frac{\partial}{\partial x_{ij}}
       \right]
   \;. \qquad
\end{eqnarray}
Let us now apply this operator to a generic polynomial $f(X)$:
using the dilation-translation formula \reff{eq.transl2} on all the
variables $x_{ij}$ ($i \neq j$), we obtain
\begin{eqnarray}
   & &
   \det((D_{\balpha,s})_{IJ}) \, f(X)
   \;=\;
   \epsilon(I,J)
   \int \scrd_n(\eta, \etabar) \,
        \Big( \prod \etabar \eta \Big)_{I^c,J^c} \,
        \exp\Biggl[ s \sum\limits_i \etabar_i \eta_i \Biggr] \;
        \times
   \nonumber \\[2mm]
   &  &  \qquad\qquad\qquad\qquad
   f \Biggl( \biggl\{
(1- \alpha_i \etabar_i \eta_i) x_{ij} (1- (1-\alpha_j) \etabar_j \eta_j )
+\etabar_i \eta_j \biggl\}_{i \neq j} \Biggl)
\ef.
\end{eqnarray}
Note that the diagonal terms remain unchanged
at their original value $x_{ii}=1$.
Defining the diagonal matrices
$M_\balpha = \diag(1- \alpha_i \etabar_i \eta_i)$
and $M_{1-\balpha} = \diag[1- (1-\alpha_i) \etabar_i \eta_i]$,
we see that the argument of $f$ is
\begin{equation}
   X' \;=\; M_\balpha X M_{1-\balpha}  \,+\, \etabar \eta^{\rm T} 
\end{equation}
(note that also the diagonal elements $i=j$ come out right).
We are interested in $f(X)=\det(X)^s$, and we have
\begin{subeqnarray}
   \det X'
   & = &
   (\det M_\balpha) (\det M_{1-\balpha})
      \det (X + M_\balpha^{-1} \etabar \eta^{\rm T}  M_{1-\balpha}^{-1})
      \\[2mm]
   & = &
   (\det M_\balpha) (\det M_{1-\balpha})
      \det (X + \etabar \eta^{\rm T})
\end{subeqnarray}
where $M_\balpha^{-1} \etabar = \etabar$
and $\eta^{\rm T}  M_{1-\balpha}^{-1} = \eta^{\rm T}$ by nilpotency.
The factor $(\det M_\balpha)^{s} (\det M_{1-\balpha})^s$
exactly cancels the factor $\exp[s \sum_i \etabar_i \eta_i]$
in the integrand, and we are left with
\begin{equation}
  \det((D_{\balpha,s})_{IJ}) \det(X)^s
  \;=\;
  \epsilon(I,J)
\int \scrd_n(\eta, \etabar)
\Big( \prod \etabar \eta \Big)_{I^c,J^c}
\det (X + \etabar \eta^{\rm T} )^s
\ef,
\end{equation}
which coincides with \reff{eq.grassmann.ordinary.IJ.a}.
Arguing exactly as in
\reff{eq.grassmann.ordinary.IJ}--\reff{eq.grassmann.ordinary.IJ.final},
we therefore obtain
\be
   \det((D_{\balpha,s})_{IJ}) \, (\det X)^s  \;=\;
      s(s+1) \cdots (s+k-1) \, (\det X)^{s-1} \,
      \epsilon(I,J) \, (\det X_{I^c J^c})
      \;.
\ee
\qed

\subsection{Diagonal-parametrized symmetric Cayley identity}
   \label{sec.param.2}

An analogous proof gives the symmetric analogue:

\proofof{Theorem~\ref{thm.para.sym.cayley}}
We change variables from $(x_{ij})_{1 \le i \leq j \le n}$ to
new variables $(t_i)_{i=1}^n$ and $(y_{ij})_{1 \le i < j \le n}$
defined by
\begin{subeqnarray}
   t_i     & = &  x_{ii}   \\
   y_{ij}  & = &  x_{ii}^{-\frac{1}{2}} x_{jj}^{-\frac{1}{2}} x_{ij}
\end{subeqnarray}
We also set $y_{ii} = 1$ for all $i$, use the synonymous $y_{ji} =
y_{ij}$ for $j > i$,
and define the matrices $T = \diag(t_i^{\frac{1}{2}})$
and $Y = (y_{ij})$,
so that
\be
    X  \;=\;  T Y T
\ee
and hence
\be
   (\det X)^s  \;=\;  \left( \prod_{i=1}^n t_i^s \right) (\det Y)^s  \;.
 \label{eq.detTYs}
\ee
A straightforward computation shows that the differential operators
(vector fields) $\partial/\partial x_{ij}$ can be rewritten in the
new variables as
\be
   {\partial \over \partial x_{ij}}  \;=\;
   \cases{  t_i^{-1} \left[ t_i \,
                        {\bigpartial \over \bigpartial t_i}
               - \frac{1}{2}
                 \sum\limits_{k \neq i} y_{ik} \,
                     {\bigpartial \over \bigpartial y_{ik}}
               - \frac{1}{2}
                 \sum\limits_{l \neq i} y_{li} \,
                     {\bigpartial \over \bigpartial y_{li}}
               \right]
                    &  if $i = j$  \cr
            \noalign{\vskip 4pt}
            t_i^{-\frac{1}{2}} t_j^{-\frac{1}{2}} \,
                     {\bigpartial \over \bigpartial y_{ij}}
                    &  if $i \neq j$  \cr
          }
 \label{def.Delta.alphas}
\ee
Let us denote by $\Delta$ the matrix of
differential operators whose elements are given by the
right-hand side of \reff{def.Delta.alphas};
please note that each element commutes with
each other element not in the same row or column
[that is, $\Delta_{ij}$ commutes with $\Delta_{i'j'}$
 whenever $i \neq i'$ and $j \neq j'$].

Now \reff{eq.detTYs}, considered as a function of $t = (t_i)$,
   is of the form ${\rm const} \times \prod_i t_i^s$;
   therefore, acting on \reff{eq.detTYs},
   each operator $t_i \partial/\partial t_i$
   is equivalent to multiplication by $s$.
   It follows that the action of $\Delta$ is identical to that of
   $\Delta_{s}$ defined by
\be
   (\Delta_{s})_{ij}  \;=\;
   \cases{  t_i^{-1} \left[ s
               - \frac{1}{2}
                 \sum\limits_{k \neq i} y_{ik} \,
                     {\bigpartial \over \bigpartial y_{ik}}
               - \frac{1}{2}
                 \sum\limits_{l \neq i} y_{li} \,
                     {\bigpartial \over \bigpartial y_{li}}
               \right]
                    &  if $i = j$  \cr
            \noalign{\vskip 4pt}
            t_i^{-\frac{1}{2}} t_j^{-\frac{1}{2}} \,
                     {\bigpartial \over \bigpartial y_{ij}}
                    &  if $i \neq j$  \cr
          }
 \label{def.Delta.s}
\ee
Furthermore,
$\Delta_{ij}$ and $(\Delta_{s})_{ij}$
both commute with $(\Delta_{s})_{i'j'}$
whenever $i \neq i'$ and $j \neq j'$.

Now suppose we have a product
$\Delta_{i_1 j_1} \cdots \Delta_{i_\ell j_\ell}$
acting on \reff{eq.detTYs},
in which $i_1, \ldots, i_\ell$ are all distinct
and also $j_1, \ldots, j_\ell$ are all distinct.
Then the rightmost factor $\Delta_{i_\ell j_\ell}$
can be replaced by $(\Delta_{s})_{i_\ell j_\ell}$;
and the same can be done for the other factors
by commuting them to the right, changing them from
$\Delta$ to $\Delta_{s}$, and commuting them back.
It follows that, under the given condition on the indices,
$\Delta_{i_1 j_1} \cdots \Delta_{i_\ell j_\ell} (\det X)^s
 = 
 (\Delta_{s})_{i_1 j_1} \cdots (\Delta_{s})_{i_\ell j_\ell} (\det X)^s$.

Translating back to the original variables $(x_{ij})_{1 \le i \leq j \le n}$,
$\Delta_{s}$ equals $D^{\rm sym}_{s}$ defined by
\be
   (D^{\rm sym}_{s})_{ij}  \;=\;
   \cases{  x_{ii}^{-1} \left[ s
               - \frac{1}{2}
                 \sum\limits_{k \neq i} x_{ik} \,
                     {\bigpartial \over \bigpartial x_{ik}}
               - \frac{1}{2}
                 \sum\limits_{l \neq i} x_{li} \,
                     {\bigpartial \over \bigpartial x_{li}}
               \right]
                    &  if $i = j$  \cr
            \noalign{\vskip 4pt}
            {\bigpartial \over \bigpartial x_{ij}}
                    &  if $i \neq j$  \cr
          }
 \label{def.D.alpha.ss}
\ee
and the Cayley formula \reff{eq.cayley.2} tells us that
\be
   \det((D^{\rm sym}_{s})_{IJ}) \, (\det X^{\rm sym})^s  \;=\;
      s \left( s+\smfrac{1}{2} \right) \cdots \left( s+\frac{k-1}{2} \right) 
     \, (\det X^{\rm sym})^{s-1} \,
      \epsilon(I,J) \, (\det X^{\rm sym}_{I^c J^c})
      \;.
 \label{eq.para.cayley.proof1s}
\ee
On the other hand,
$D^{\rm sym}_{\bbeta,s} = \widehat{X}_\bbeta D^{\rm sym}_{s} \widehat{X}_{1-\bbeta}$
where $\widehat{X}_\bbeta = \diag(x_{ii}^{\beta_i})$
and $\widehat{X}_{1-\bbeta} = \diag(x_{ii}^{1-\beta_i})$,
so that
\be
   \det((D^{\rm sym}_{\bbeta,s})_{IJ})  \;=\;
      \Biggl( \prod\limits_{i \in I} x_{ii}^{\beta_i} \Biggr) \,
      \Biggl( \prod\limits_{j \in J} x_{jj}^{1-\beta_j} \Biggr) \,
      \det((D^{\rm sym}_{s})_{IJ})
   \;.
 \label{eq.para.cayley.proof2s}
\ee
Combining \reff{eq.para.cayley.proof1s} and \reff{eq.para.cayley.proof2s},
we obtain \reff{eq.para.sym.cayley2}.
\qed

Now let us show the proof based on
Grassmann representation of the differential operator:

\alternateproofof{Theorem~\ref{thm.para.sym.cayley}}
For notational simplicity let us assume that the diagonal elements $x_{ii}$
are all equal to 1;  the general case can be recovered by a simple scaling.
We are therefore using the matrix $X^{\rm sym}$ defined by
\begin{equation}
(X^{\rm sym})_{ij}  \;=\;
\cases{ 1      & for $i=j$  \cr
        x_{ij} & for $i<j$  \cr
        x_{ji} & for $i>j$  \cr
      }
\end{equation}

Consider the matrix of differential operators $D^{\rm sym}_{\bbeta,s}$
defined by \reff{def.para.sym.cayley}.
In terms of the matrices $E^{ij}$ defined by \reff{def.Eij},
we can write
\begin{equation}
   D^{\rm sym}_{\bbeta,s}
   \;=\;
   s I  \,+\,
   \smhalf
   \sum_{i<j} \Big[ E^{ij} + E^{ji} - x_{ij} (E^{ii} + E^{jj}) \Big]
              \frac{\partial}{\partial x_{ij}}
   \;\,.
\end{equation}
We introduce Grassmann variables $\eta_i, \etabar_i$
($1 \le i \le n$) and use the representation
\begin{eqnarray}
   & &
   \det((D^{\rm sym}_{\bbeta,s})_{IJ})
   \;=\;
   \epsilon(I,J)
   \int \scrd_n(\eta, \etabar) \,
   \Big( \prod \etabar \eta \Big)_{I^c,J^c}
   \;\times
       \nonumber \\[2mm]
   & &
   \qquad
   \exp\!\left[
         s \sum\limits_i \etabar_i \eta_i
         \,+\, \smhalf \sum\limits_{i<j}
             [\etabar_i \eta_j - \eta_i \etabar_j
              - x_{ij} (\etabar_i \eta_i +\etabar_j \eta_j )]
             \frac{\partial}{\partial x_{ij}}
       \right]
   \;. \qquad
\end{eqnarray}
Applying the dilation-translation formula \reff{eq.transl2} on all the
variables $x_{ij}$ ($i < j$), we obtain
\begin{eqnarray}
   & &
   \det((D^{\rm sym}_{\bbeta,s})_{IJ})  \, f(\{x_{ij}\}_{i<j})
   \;=\;
   \epsilon(I,J)
   \int \scrd_n(\eta, \etabar) \,
      \Big( \prod \etabar \eta \Big)_{I^c,J^c} \,
      \exp\Biggl[ s \sum\limits_i \etabar_i \eta_i \Biggr]
      \,\times
   \nonumber \\[2mm]
   & & \qquad\qquad
   f \Biggl( \biggl\{
(1- \smhalf \etabar_i \eta_i) x_{ij} (1- \smhalf \etabar_j \eta_j )
+\smhalf (\etabar_i \eta_j - \eta_i \etabar_j) \biggr\}_{i<j} \Biggr)
\ef.
\end{eqnarray}
Note that the diagonal terms remain unchanged
at their original value $x_{ii}=1$.
Defining the diagonal matrix
$M = \diag(1- \smhalf \etabar_i \eta_i)$,
we see that the argument of $f$ is
\begin{equation}
  {X^{\rm sym}}' \;=\; MX^{\rm sym} M \,+\,
       \smhalf (\etabar \eta^{\rm T}  - \eta \etabar^{\rm T} )
\end{equation}
(note that also the diagonal elements $i=j$ come out right).
We are interested in $f(X^{\rm sym})=\det(X^{\rm sym})^s$, and we have
\begin{subeqnarray}
   \det {X^{\rm sym}}'
   & = &
   (\det M)^2 \det [X + \smhalf M^{-1}( \etabar \eta^{\rm T}
                                        - \eta \etabar^{\rm T} )M^{-1}]
      \\[2mm]
   & = &
   (\det M)^2 \det [X + \smhalf ( \etabar \eta^{\rm T}  - \eta \etabar^{\rm T} )]
\end{subeqnarray}
where the last equality again follows by nilpotency.
The factor $(\det M)^{2s}$ exactly cancels
the factor $\exp[s \sum_i \etabar_i \eta_i]$ in the integrand,
and we are left with
\begin{equation}
   \det((D^{\rm sym}_{\bbeta,s})_{IJ}) \, \det(X^{\rm sym})^s
   \;=\;
   \epsilon(I,J)
   \int \scrd_n(\eta, \etabar) \,
   \Big( \prod \etabar \eta \Big)_{I^c,J^c} \,
   \det [X^{\rm sym} + \smhalf (\etabar \eta^{\rm T}  - \eta \etabar^{\rm T} )]^s
\ef,
\end{equation}
which coincides with \reff{eq.grassmann.sym.IJ}.
The remainder of the proof is as in
\reff{eq.grassmann.sym.IJ}--\reff{eq.grassmann.sym.IJ.last}.
\qed

\section{Proofs of Laplacian-parametrized Cayley identities}  \label{sec.Lap}

In this section we prove Theorems~\ref{thm.cayley.rowlapl.pre}
and \ref{thm.cayley.lapl.pre}.
The proofs use a Grassmann representation of the differential operator,
and are closely patterned after the proofs of the
ordinary and symmetric Cayley identities in
Sections~\ref{sec.grassmann.2.ordinary} and \ref{sec.grassmann.2.symmetric},
respectively.

\subsection{Laplacian-parametrized ordinary Cayley identity}
   \label{sec.Lap.ordinary}

Let us begin by recalling the definitions of the matrices
arising in Theorem~\ref{thm.cayley.rowlapl.pre}:
\begin{subeqnarray}
   (\XRL)_{ij}
   & = &  \cases{ x_{ij}     & if $i \neq j$  \cr
                   \noalign{\vskip 6pt}
                   - \sum\limits_{k \neq i} x_{ik}  & if $i = j$ \cr
                 }
      \\[4mm]
   (\partialRL)_{ij}
   & = &   \cases{ \partial/\partial x_{ij}     & if $i \neq j$  \cr
                   \noalign{\vskip 6pt}
                   0                            & if $i = j$ \cr
                 }
 \label{def.XRL.partialRL}
\end{subeqnarray}
and of course $T = \diag(t_i)$.
In order to maximize the correspondences with the proof
in Section~\ref{sec.grassmann.2.ordinary},
it is convenient to prove instead the dual (and of course equivalent)
result for {\em column}\/-Laplacian matrices:
\begin{subeqnarray}
   (\XCL)_{ij}
   & = &  \cases{ x_{ij}     & if $i \neq j$  \cr
                   \noalign{\vskip 6pt}
                   - \sum\limits_{k \neq i} x_{ki}  & if $i = j$ \cr
                 }
      \\[4mm]
   (\partialCL)_{ij}
   & = &   \cases{ \partial/\partial x_{ij}     & if $i \neq j$  \cr
                   \noalign{\vskip 6pt}
                   0                            & if $i = j$ \cr
                 }
 \label{def.XCL.partialCL}
\end{subeqnarray}
(note that $\partialCL = \partialRL$).
In what follows we shall drop the superscripts ``col-Lap''
in order to lighten the notation,
but it is important to remember the definitions \reff{def.XCL.partialCL}.

\proofof{Theorem~\ref{thm.cayley.rowlapl.pre}}
%
We introduce Grassmann variables $\eta_i, \etabar_i$ ($1 \le i \le n$)
and use the representation
\be
   \det(U + \partial) - \det(\partial)  \;=\;
   \int \! \scrd_n(\eta,\etabar) \; 
   (e^{\etabar^{\rm T} U \eta}-1)
   e^{\etabar^{\rm T} \partial \eta}
   \;,
 \label{eq.detpartialRL}
\ee
It is convenient to introduce the Grassmann quantities
$\Theta = \sum_{\ell} \eta_{\ell}$ and
$\bar{\Theta} = \sum_{\ell} \etabar_{\ell}$,
so that $\etabar^{\rm T} U \eta = \bar{\Theta}^{\rm T} \Theta$.
Since $\Theta^2 = \bar{\Theta}^2 = 0$,
the first exponential in \reff{eq.detpartialRL} is easily expanded,
and we have
\be
   \det(U + \partial) - \det(\partial)  \;=\;
   \int \! \scrd_n(\eta,\etabar) \; 
   (\bar{\Theta}^{\rm T} \Theta) \,
   e^{\etabar^{\rm T} \partial \eta}
   \;.
\label{eq.detpartial.1234lap}
\ee
Let us now apply \reff{eq.detpartial.1234lap} to 
$\det (T + X)^s$ where $s$ is a positive integer,
using the translation formula \reff{eq.translation}:
by (\ref{def.XCL.partialCL}b) we get
$x_{ij} \to x_{ij} + \etabar_i \eta_j$ for $i \neq j$,
and by (\ref{def.XCL.partialCL}a) this induces
$X \to X + \etabar \eta^{\rm T} - \diag(\bar{\Theta} \eta_i)$.
We therefore obtain
\begin{eqnarray}
   [\det(U + \partial) - \det(\partial)]  \: \det(T + X)^s
    && 
\nonumber \\
   && \hspace{-5cm} =\;
   \int \! \scrd_n(\eta,\etabar) \; 
   (\bar{\Theta}^{\rm T} \Theta) \,
\det(T + X + \etabar \eta^{\rm T} - \diag(\bar{\Theta} \eta_i))^s
   \;.
\end{eqnarray}
But now comes an amazing simplification:
because of the prefactor $\bar{\Theta}^{\rm T} \Theta$
and the nilpotency $\Theta^2 = 0$,
all terms in the expansion of the determinant
arising from the term $\diag(\bar{\Theta} \eta_i)$ simply vanish,
so we can drop $\diag(\bar{\Theta} \eta_i)$:
\be
   [\det(U + \partial) - \det(\partial)]  \: \det (T + X)^s
   \;=\;
   \int \! \scrd_n(\eta,\etabar) \;
   (\bar{\Theta}^{\rm T} \Theta) \, \det(T + X + \etabar \eta^{\rm T})
   \;.
\ee
Assuming that $T+X$ is an invertible real or complex matrix,
we can write this as
\begin{eqnarray}
 [\det(U + \partial) - \det(\partial)]  \: \det (T + X)^s && 
\\
\nonumber
   && \hspace{-5cm} =\;
\det(T + X)^s
   \int \! \scrd_n(\eta,\etabar) \; 
   (\etabar^{\rm T} U \eta) \,
   \det [I + (T + X)^{-1} \etabar \eta^{\rm T} ]^s
   \;.
\end{eqnarray}
Let us now change variables from $(\eta,\etabar)$ to
$(\eta',\etabar') \equiv (\eta, (T + X)^{-1} \etabar)$:
we pick up a Jacobian $\det (T + X))^{-1}$,
and dropping primes we have
\begin{eqnarray}
   [\det(U + \partial) - \det(\partial)]  \: \det (T + X)^s && 
\\
\nonumber
   && \hspace{-5cm} =\;
\det(T + X)^{s-1}
   \int \! \scrd_n(\eta,\etabar) \; 
(\etabar^{\rm T} (T + X^{\rm T}) U \eta) \,
\det (I + \etabar \eta^{\rm T})^s
   \;.
\end{eqnarray}
But since $X$ is column-Laplacian, we have $UX=0$ and hence $X^{\rm T} U = 0$,
so that $\etabar^{\rm T} (T + X^{\rm T}) U \eta$
reduces to $\etabar^{\rm T} T U \eta$.\footnote{
   For a row-Laplacian matrix we would have instead chosen
   to multiply $\etabar \eta^{\rm T}$ by $(T + X)^{-1}$
   on the {\em right}\/, and then made the change of variables
   from $(\eta,\etabar)$ to
   $(\eta',\etabar') \equiv ((T + X^{\rm T})^{-1}\eta, \etabar)$,
   picking up a Jacobian $(\det (T + X^{\rm T}))^{-1} = (\det (T + X))^{-1}$
   and obtaining a prefactor
   $\etabar^{\rm T} U (T + X^{\rm T}) \eta$;
   then $U X^{\rm T} = 0$ because $X$ is row-Laplacian.
}
This expresses the left-hand side of the identity
as the desired quantity $\det(T + X)^{s-1}$ times a factor
\be
   P(s,n,T)  \;\equiv\; 
  \int \! \scrd_n(\eta,\etabar) \;
   (\etabar^{\rm T} T U \eta) \,
          \det(I + \etabar \eta^{\rm T})^s
  \;,
\ee
which we now proceed to calculate.
The matrix $I + \etabar \eta^{\rm T}$ is a rank-1 perturbation 
of the identity matrix;  by Lemma~\ref{lemma.lowrank} we have
\begin{subeqnarray}
   \det(I + \etabar \eta^{\rm T})^s
   & = &
   (1 - \etabar^{\rm T} \eta)^{-s}
      \\[1mm]
   & = & \sum_{\ell=0}^\infty
            (-1)^\ell \, {-s \choose \ell} \, (\etabar^{\rm T} \eta)^\ell
   \;.
 \label{eq.grass1.star2_extra1}
\end{subeqnarray}
Now 
\be
   \etabar^{\rm T} T U \eta
   \;=\;
   \sum_i t_i \etabar_i \eta_i  \,+\,
      \sum_i \sum_{j \neq i} t_i \etabar_i \eta_j
   \;,
\ee
but the terms $\etabar_i \eta_j$ with $i \neq j$
cannot contribute to an integral 
in which the rest of the integrand depends on $\etabar$ and $\eta$
only through products $\etabar_i \eta_i$;
so in the integrand we can replace
$\etabar^{\rm T} T U \eta$ with $\etabar^{\rm T} T \eta$.
Since
\be
   \int \! \scrd_n(\eta,\etabar) \; 
\etabar_i \eta_i \,
(\etabar^{\rm T} \eta)^\ell
   \;=\;
   (n-1)! \, \delta_{\ell,n-1}
\ee
for each index $i$, it follows that
\begin{subeqnarray}
   P(s,n,T)
   & = &
\Big( \sum_i t_i \Big)
   (-1)^{n-1} \, {-s \choose n-1} \, (n-1)!
      \\[1mm]
   & = &
\Big( \sum_i t_i \Big) s(s+1) \cdots (s+n-2)
   \;.
 \label{eq.grass1.star3rlapl}
\end{subeqnarray}
This proves \reff{eq.cayley.1.rowlapl.pre.INTRO} when
$X$ is a real or complex matrix,
${\bf t} = (t_i)$ are real or complex values such that $T+X$ is invertible, 
and $s$ is a positive integer; the general validity
of the identity then follows from Proposition~\ref{prop.equiv.symbolic_s}.
\qed

\subsection{Laplacian-parametrized symmetric Cayley identity}
   \label{sec.Lap.symmetric}

Now we prove the corresponding result for symmetric Laplacian matrices.
Let us recall the definitions:
\begin{subeqnarray}
   (\XSL)_{ij}
   & = &   \cases{ x_{ij}     & if $i < j$  \cr
                   \noalign{\vskip 6pt}
                   x_{ji}     & if $i > j$  \cr
                   \noalign{\vskip 6pt}
                   - \sum\limits_{k \neq i} x_{ik}  & if $i = j$ \cr
                 }
   \\[4mm]
   (\partialSL)_{ij}
   & = & \;=\;  \cases{ \partial/\partial x_{ij}     & if $i < j$  \cr
                   \noalign{\vskip 6pt}
                   \partial/\partial x_{ji}     & if $i > j$  \cr
                   \noalign{\vskip 6pt}
                   0                            & if $i = j$ \cr
                 }
 \label{def.XSL.partialSL}
\end{subeqnarray}
Once again we drop the superscripts ``sym-Lap'' to lighten the notation.

\proofof{Theorem~\ref{thm.cayley.lapl.pre}}
Again we introduce Grassmann variables $\eta_i, \etabar_i$ ($1 \le i \le n$)
and apply the representation \reff{eq.detpartial.1234lap}
to $\det (T + X)^s$
where $s$ is a positive integer.
Using the translation formula \reff{eq.translation}
and defining $\Theta$ and $\bar{\Theta}$ as before,
we obtain
\begin{eqnarray}
 [\det(U + \partial) - \det(\partial) ]  \: \det (T + X)^s && 
   \nonumber \\
   && \hspace{-5.5cm} =\;
   \int \! \scrd_n(\eta,\etabar) \; 
   (\bar{\Theta} \Theta) \,
\det(T + X 
+ \etabar \eta^{\rm T}
- \eta \etabar^{\rm T}
- \diag(\etabar_i \Theta)
- \diag(\bar{\Theta} \eta_i)
)^s
   \;.
   \nonumber \\
\end{eqnarray}
We again argue that because of the prefactor $\bar{\Theta}^{\rm T} \Theta$
and the nilpotencies $\Theta^2 = \bar{\Theta}^2 = 0$,
all terms in the expansion of the determinant
arising from the terms $\diag(\etabar_i \Theta)$
and $\diag(\bar{\Theta} \eta_i)$ simply vanish,
so we can drop these two terms:
\begin{eqnarray}
 [\det(U + \partial) - \det(\partial) ]  \: \det (T + X)^s &&
   \nonumber \\    && \hspace{-5.5cm} =\;
   \int \! \scrd_n(\eta,\etabar) \;
   (\bar{\Theta} \Theta) \,
\det(T + X
+ \etabar \eta^{\rm T}
- \eta \etabar^{\rm T}
)^s
   \;.
\end{eqnarray}
Assuming that $T+X$ is an invertible real or complex matrix,
we can write this as
\begin{eqnarray}
 & &
 [\det(U + \partial) - \det(\partial) ]  \: \det (T + X)^s
   \nonumber \\
 & &
 \;=\;
   \det (T + X)^s
   \int \! \scrd_n(\eta,\etabar) \;
   (\etabar^{\rm T} U \eta) \,
\det[I + (T + X)^{-1} (\etabar \eta^{\rm T} - \eta \etabar^{\rm T})]^s
   \;.
\end{eqnarray}
Let us now change variables from $(\eta,\etabar)$ to
$(\eta',\etabar') \equiv 
(
\eta, (T + X)^{-1} \etabar)$;
we pick up a Jacobian $(\det (T + X))^{-1}$,
and dropping primes we have
\begin{eqnarray}
   & &
  [\det(U + \partial) - \det(\partial)] \: \det (T + X)^s
     \nonumber \\
   & &
   \;=\;
\det(T + X)^{s-1}
   \int \! \scrd_n(\eta,\etabar) \; 
(\etabar^{\rm T}
(T + X^{\rm T})
U 
\eta) \,
\det [I + \etabar \eta^{\rm T}
- 
(T + X)^{-1}
\eta \etabar^{\rm T}
(T + X)^{\rm T}
]^s
   \;.
\nonumber \\
\end{eqnarray}
But since $X$ is symmetric Laplacian (hence column-Laplacian),
we have $X^{\rm T} U=0$,
and the prefactor reduces to
$\etabar^{\rm T} T U \eta$.
Now we apply the Corollary \ref{corol.lowranksym} to the
determinant expression, and obtain
\begin{eqnarray}
  & &
  [\det(U + \partial) - \det(\partial)] \: \det (T + X)^s
      \nonumber \\
  & &
  \qquad =\;
\det(T + X)^{s-1}
   \int \! \scrd_n(\eta,\etabar) \; 
(\etabar^{\rm T} T U \eta)
(1 - \etabar^{\rm T} \eta)^{-2s}
   \;.
\end{eqnarray}
This expresses the left-hand side of the identity
as the desired quantity $\det(T + X)^{s-1}$ times a factor
\be
   P(s,n,T)  \;\equiv\; 
\int \! \scrd_n(\eta,\etabar) \;
(\etabar^{\rm T} T U \eta) \,
(1 - \etabar^{\rm T} \eta)^{-2s}
  \;,
\ee
which we now proceed to calculate.
As in the row-Laplacian case, the terms in $\etabar^{\rm T} (TU) \eta$
of the form $\etabar_i \eta_j$ with $i \neq j$ cannot contribute to the
integral, as the rest of the integrand depends on $\eta$ and $\etabar$
only through products $\etabar_i \eta_i$,
so we can replace $\etabar^{\rm T} (TU) \eta$ with $\etabar^{\rm T} T \eta$.
Since
\be
   \int \! \scrd_n(\eta,\etabar) \; 
\etabar_i \eta_i \,
(\etabar^{\rm T} \eta)^\ell
   \;=\;
   (n-1)! \, \delta_{\ell,n-1} 
\ee
for each $i$, it follows that
\begin{subeqnarray}
   P(s,n,T)
   & = &
\Big( \sum_i t_i \Big)
   (-1)^{n-1} \, {-2s \choose n-1} \, (n-1)!
      \\[1mm]
   & = &
\Big( \sum_i t_i \Big)\,
   2s(2s+1) \cdots (2s+n-2)
   \;.
 \label{eq.grass1.star3lapl}
\end{subeqnarray}
This proves \reff{eq.thm.Lap.sym.cayley}
when $X$ is a symmetric real or complex matrix,
${\bf t} = (t_i)$ are real or complex values such that $T+X$ is invertible,
and $s$ is a positive integer;
the general validity
of the identity then follows from Proposition~\ref{prop.equiv.symbolic_s}.
\qed

\section[Proofs of product-parametrized
   and border-parametrized rectangular Cayley identities]
   {Proofs of product-parametrized and border-\hfill\break
   parametrized rectangular Cayley identities}
   \label{sec.border}

In this section we prove the product-parametrized and border-parametrized
rectangular Cayley identities
(Theorems~\ref{thm.productcayley} and \ref{thm.borderedcayley})
and then discuss the close relationship between them.

\subsection{Product-parametrized rectangular Cayley identity}
   \label{subsec.product}

Before beginning the proof of Theorem~\ref{thm.productcayley},
let us observe that the quantity $\det(MA)$ appearing
in the statement of the theorem
[cf.\ \reff{eq.defMinXa2m} for the definition of the matrix $M$
 in terms of $X$ and $B$]
has an alternate expression as follows:

\begin{lemma}
   \label{lemma.productcayley}
Let $A$, $B$, $X$ and $M$ be as in Theorem~\ref{thm.productcayley}.
Then
\be
\det(MA)
\;=\;
\sum_{
\begin{scarray}
L \subseteq [m] \\ |L|=k
\end{scarray} }
\epsilon(I,L) \, (\det(B^{\rm T} A)_{JL}) \, (\det(X A)_{I^c L^c})
\label{eq.defMAthird}
\ee
\end{lemma}

\proof
The definition \reff{eq.defMinXa2m} of $M$ can be rewritten as
\begin{subeqnarray}
   M_{I^c \star}  & = &  X_{I^c \star}  \\[1mm]
   M_{I \star}    & = &  (B^{\rm T})_{J \star}
\end{subeqnarray}
Therefore
\begin{subeqnarray}
   (MA)_{I^c \star}  & = &  (XA)_{I^c \star}  \\[1mm]
   (MA)_{I \star}    & = &  (B^{\rm T} A)_{J \star}
\end{subeqnarray}
for any matrix $A$.  We now apply multi-row Laplace expansion
\reff{eq.multirowlaplace} with row set $I$
and summation variable $L$;  this yields \reff{eq.defMAthird}.
\qed



We are now ready to prove Theorem~\ref{thm.productcayley}.
In order to bring out the ideas behind the proof as clearly as possible,
we will first fully develop the reasoning
proving the ``basic'' identity \reff{eq.cayley.XA2m}
--- which is actually quite simple ---
and then describe the modifications needed to handle the all-minors case.

\proofof{Theorem~\ref{thm.cayley.XA2m}}
Let us apply Corollary~\ref{cor.blockmatrixdet} to the determinant 
$\det(\partial B)$ and then introduce a Grassmann representation for the
resulting block determinant:  we obtain
\be
\det(\partial B) 
\;=\;
\det
\left(
\begin{array}{c|c}
0_m & \de \\
\hline
-B  & I_n
\end{array}
\right)
\;=\;
   \int \! 
\scrd_m(\psi, \psibar) \,
\scrd_n(\eta, \etabar) \,
     \exp\!\left[ 
\psibar^{\rm T} \partial \eta
- \etabar^{\rm T} B \psi
+ \etabar^{\rm T} \eta
           \right]
   \label{eq.grassmann.XA.2m.PRE}
\ee
where $\psi_i, \psibar_i$ ($1 \le i \le m$)
and $\eta_i, \etabar_i$ ($1 \le i \le n$)
are Grassmann variables.
By the translation formula \reff{eq.translation}, we have
\begin{equation}
\det(\partial B) \, f(X)
\;=\;
   \int \! 
\scrd_m(\psi, \psibar) \,
\scrd_n(\eta, \etabar) \,
     \exp\!\left[ 
\etabar^{\rm T} (-B \psi + \eta)
           \right]
f(X \,+\, \psibar \eta^{\rm T})
\end{equation}
for an arbitrary polynomial $f$.
We shall use this formula in the case
$f(X)=\det (X A)^s$
where $s$ is a positive integer, so that
\be
\det(\partial B) \,
\det (X A)^s
\;=\;
   \int \! 
\scrd_m(\psi, \psibar) \,
\scrd_n(\eta, \etabar) \,
     \exp\!\left[ 
\etabar^{\rm T} (-B \psi + \eta)
           \right]
\det\Big[
(X+\psibar \eta^{\rm T}) A
    \Big]^s
 \;.
   \label{eq.grassmann.XA.2m}
\ee
It is convenient to introduce the shorthand
\be
X^{\mathrm{trans}}
\;\equiv\;
X+\psibar \eta^{\rm T}
\ef.
\ee

Suppose now that $XA$ is an invertible real or complex matrix.
Then we have
\begin{equation}
X^{\rm trans} A
\;\equiv\;
(X + \psibar \eta^{\rm T})A
  \;=\;
(XA) [I_m + ((XA)^{-1} \psibar) (\eta^{\rm T} A)]
  \;.
\end{equation}
Let us now change variables from $(\psi, \psibar, \eta, \etabar)$
to $(\psi', \psibar', \eta', \etabar') \equiv
    (\psi, (XA)^{-1} \psibar, \eta, \etabar)$,
with Jacobian $\det (XA)^{-1} = (\det X A)^{-1}$.
In the new variables we have
(dropping now the primes from the notation)
\begin{equation}
X^{\rm trans} A
  \;=\;
(XA) (I_m + \psibar \eta^{\rm T} A)
  \;,
\end{equation}
and the translated determinant is given by
\begin{equation}
\det(X^{\rm trans} A )  \;=\;
\det(XA)
\det(I_m + \psibar \eta^{\rm T} A)
\;,
\end{equation}
so that
\be
   \det(\partial B)  \,  \det(XA)^s
       \nonumber \\
   \;=\;
\det(XA)^{s-1}
   \int \! \scrd_m(\psi,\psibar) \, \scrd_n(\eta,\etabar) \:
   e^{-\etabar^{\rm T} B \psi + \etabar^{\rm T} \! \eta} \,
   \det(I + \psibar \eta^{\rm T} A)^s
   \;.
\ee
Applying Lemma~\ref{lemma.lowrank} to the rightmost determinant
yields
\be
   \det(I + \psibar \eta^{\rm T} A)
   \;=\;
(1- \psibar^{\rm T} A^{\rm T} \eta)^{-1}  \;,
\ee
so that we are left with the Grassmann-integral expression
\be
   \det(\partial B)  \,  \det(XA)^s
   \;=\;
\det(XA)^{s-1}  \!
   \int \! \scrd_m(\psi,\psibar) \, \scrd_n(\eta,\etabar) \:
   e^{-\etabar^{\rm T} B \psi + \etabar^{\rm T} \! \eta} \,
(1- \psibar^{\rm T} A^{\rm T} \eta)^{-s}
   \;.
\label{eq.5486976.rifer}
\ee
We have therefore proven that $\det(\partial B)  \,  \det(XA)^s$
equals the desired quantity $\det(XA)^{s-1}$
multiplied by a factor
\be
   b(s,A,B)  \;=\;
   \int \! \scrd_m(\psi,\psibar) \, \scrd_n(\eta,\etabar) \:
   e^{-\etabar^{\rm T} B \psi + \etabar^{\rm T} \! \eta} \,
(1- \psibar^{\rm T} A^{\rm T} \eta)^{-s}
 \label{def.bsAB}
\ee
that does not involve the variables $X$
(but still involves the parameters $A$ and $B$).
Now, in the expansion of
\be
   (1- \psibar^{\rm T} A^{\rm T} \eta)^{-s}
   \;=\;
   \sum_{k=0}^\infty (-1)^k \binom{-s}{k} \, (\psibar^{\rm T} A^{\rm T} \eta)^k
   \;,
\ee
only the term $k=m$ survives the integration over the variables $\psibar$,
so we can replace $(1- \psibar^{\rm T} A^{\rm T} \eta)^{-s}$
in the integrand of \reff{def.bsAB} by
$(-1)^m \binom{-s}{m} \, (\psibar^{\rm T} A^{\rm T} \eta)^m$.
Moreover, the same reasoning shows that we can replace
$(\psibar^{\rm T} A^{\rm T} \eta)^m$
by $m! \, \exp(\psibar^{\rm T} A^{\rm T} \eta)$.
We are therefore left with a combinatorial prefactor
\be
(-1)^m
\binom{-s}{m} m!
\;=\;
s(s+1)\cdots(s+m-1)
\ee
multiplying the Grassmann integral
\be
  \int \! \scrd_m(\psi,\psibar) \, \scrd_n(\eta,\etabar) \:
    \exp\!\left[ \psibaretabarvert^{\!\! \rm T}
                 \left(
                    \begin{array}{c|c}
                       0_m & A^{\rm T} \\ 
                       \hline
                       -B  & I_n
                    \end{array}
                 \right)
                 \psietavert
           \right]
  \;=\;
   \det\!      \left(
                    \begin{array}{c|c}
                       0_m & A^{\rm T} \\ 
                       \hline
                       -B  & I_n
                    \end{array}
                 \right)
   \,,
\ee
which equals $\det(A^{\rm T} B)$ by Corollary~\ref{cor.blockmatrixdet}.

This proves the ``basic'' identity \reff{eq.cayley.XA2m}
whenever $XA$ is an invertible real or complex matrix
and $s$ is a positive integer.
Now, if $A$ has rank $< m$, then both sides of
\reff{eq.cayley.XA2m} are identically zero;
while if $A$ has rank $m$, then
$XA$ is invertible for a nonempty open set of matrices $X$.
The general validity of the identity \reff{eq.cayley.XA2m} therefore
follows from Proposition~\ref{prop.equiv.symbolic_s}.

Now let us consider the modifications needed to prove
the all-minors identity \reff{eq.cayley.XA2m.2}.
For a while these modifications will run along the same lines as
those used in the proof of the two-matrix rectangular Cayley identity
(Section~\ref{sec.grassmann.2.tmrect}).
Thus, a factor $\epsilon(I,J) \Big( \prod \psibar \psi \Big)_{I^c,J^c}$
gets inserted into the Grassmann integral
\reff{eq.grassmann.XA.2m.PRE}--\reff{eq.grassmann.XA.2m};
after the change of variables (and dropping of primes) it becomes
$\epsilon(I,J) \Big( \prod (XA\psibar) \psi \Big)_{I^c,J^c}$.
So we have, in place of equation (\ref{eq.5486976.rifer}),
the modified expression
\begin{eqnarray}
   & &
   \det[( \partial B )_{IJ}]  \,  \det(XA)^s
   \;=\;
   \epsilon(I,J) \, \det(XA)^{s-1}
   \int \! \scrd_m(\psi,\psibar) \, \scrd_n(\eta,\etabar)
   \,
   \Big( \prod (XA\psibar) \psi \Big)_{I^c,J^c}
        \nonumber \\
   & & \qquad\qquad\qquad \times\;
   \exp\Bigl[ \etabar^{\rm T} \! \eta \,-\, \etabar^{\rm T} B \psi \Bigr]  \,
(1- \psibar^{\rm T} A^{\rm T} \eta)^{-s}
   \;.
\end{eqnarray}
Once again we argue that the integration over variables $\psibar$
allows to replace
\be
(1- \psibar^{\rm T} A^{\rm T} \eta)^{-s}
\quad
\longrightarrow
\quad
(-1)^k
\binom{-s}{k} k!
\;
\exp[\psibar^{\rm T} A^{\rm T} \eta]
\ee
since, in both cases, only the $k$-th term of the expansion survives.
We therefore have
\begin{eqnarray}
   & &
   \det[( \partial B )_{IJ}]  \,  \det(XA)^s
   \;=\;
   \epsilon(I,J) \, \det(XA)^{s-1}
s(s+1)\cdots(s+k-1)
        \nonumber \\
   & &
   \qquad \times\;
   \int \! \scrd_m(\psi,\psibar) \, \scrd_n(\eta,\etabar)
   \,
   \Big( \prod (\psibar A^{\rm T} X^{\rm T}) \psi \Big)_{I^c,J^c}
   \,
   \exp[ \eta^{\rm T} \eta \,-\, \eta^{\rm T} B \psi
                 \,+\, \psibar^{\rm T} A^{\rm T} \eta
       ]
   \;.
   \nonumber \\
\end{eqnarray}
Now we perform the integration over $\eta$ and $\etabar$
using Wick's theorem for ``complex'' fermions in the ``source'' form
[cf.\ \reff{eq.wick.complexfermions.1}], yielding
\begin{eqnarray}
   & &
   \det[( \partial B )_{IJ}]  \,  \det(XA)^s
   \;=\;
   \epsilon(I,J) \, \det(XA)^{s-1}
s(s+1)\cdots(s+k-1)
        \nonumber \\
   & &
   \qquad\qquad\qquad \times\;
   \int \! \scrd_m(\psi,\psibar)
   \,
   \Big( \prod (\psibar A^{\rm T} X^{\rm T}) \psi \Big)_{I^c,J^c}
   \,
   \exp[ \psibar^{\rm T} A^{\rm T} B \psi ]
   \;. \qquad
\end{eqnarray}
Next we perform the integration over $\psi$ and $\psibar$
using Wick's theorem for ``complex'' fermions in the
``correlation function'' form \reff{eq.wick.complexfermions.3a},
yielding\footnote{
   Here we have made in \reff{eq.wick.complexfermions.3a}
   the substitutions $A \to A^{\rm T}B$, $B \to (I)_{J^c \star}$,
   $C \to (A^{\rm T} X^{\rm T})_{\star I^c}$,
   $I \to K$, $J \to L$.
}
\begin{eqnarray}
   & &
   \int \! \scrd_m(\psi,\psibar)
      \, \Big( \prod (\psibar A^{\rm T} X^{\rm T}) \psi \Big)_{I^c,J^c}
      \, \exp[ \psibar^{\rm T} A^{\rm T} B \psi ]
   \nonumber \\[2mm]
   & &
   \qquad = \;
   \sum_{|K| = |L| = n-k} \epsilon(K,L) \,
         (\det I_{J^c L}) \, (\det (A^{\rm T} B)_{K^c L^c}) \,
                              (\det (A^{\rm T} X^{\rm T})_{K I^c})
   \nonumber \\[2mm]
   & &
   \qquad = \;
   \sum_{|K| = n-k} \epsilon(K,J^c) \,
          (\det (A^{\rm T} B)_{K^c J}) \, (\det (A^{\rm T} X^{\rm T})_{K I^c})
   \nonumber \\[2mm]
   & &
   \qquad = \;
   \sum_{|K| = k} \epsilon(K^c,J^c) \,
          (\det (A^{\rm T} B)_{K J}) \, (\det (A^{\rm T} X^{\rm T})_{K^c I^c})
   \;.
 \label{eq.proofproductcayley.star1}
\end{eqnarray}
We now use $\epsilon(K^c,J^c) = \epsilon(J,K)$
and $\epsilon(I,J) \epsilon(J,K) = \epsilon(I,K)$;
it follows that $\epsilon(I,J)$ times \reff{eq.proofproductcayley.star1}
equals \reff{eq.defMAthird}.
\qed

\bigskip

{\bf Remark.}  The $s=1$ special case of the all-minors identity
\reff{eq.cayley.XA2m.2} has an easy elementary proof, which actually
proves a stronger result.
Note first that by the multilinearity of the determinant
$\det(XA)$ in the variables $\{x_{ij}\}$, we have
\be
\det(M A)
\;=\;
\left(
\prod_p (\partial B)_{i_p j_p}
\right)
\det(XA)
\;.
\ee
Moreover, because permuting the indices $j_1,\ldots,j_k$
amounts to permuting the rows of $B^{\rm T}$ and hence permuting
a subset of the rows of $M$, we have
\be
\det(M A)
\;=\;
\sgn(\sigma)
\left(
\prod_p (\partial B)_{i_p j_{\sigma(p)}}
\right)
\det(XA)
  \label{eq.productcayley.remark.sigma}
\ee
for any permutation $\sigma \in \scrs_k$.
Summing this over $\sigma \in \scrs_k$, we obtain
\be
   k! \, \det(M A)  \;=\; [\det (\partial B)_{IJ}] \, \det(XA)
   \;,
\ee
which is nothing other than the $s=1$ case of \reff{eq.cayley.XA2m.2}.
But it is amusing to note that \reff{eq.productcayley.remark.sigma}
holds for {\em each}\/ $\sigma \in \scrs_k$,
not just when summed over $\sigma \in \scrs_k$.

\subsection{Border-parametrized rectangular Cayley identity}
   \label{subsec.border}


\proofof{Theorem~\ref{thm.borderedcayley}}
We introduce Grassmann variables
$\psi_i, \psibar_i$ ($1 \le i \le m$)
and
$\eta_i, \etabar_i$ ($1 \le i \le n-m$).
We use the representation
\begin{equation}
   \det(\widehat{\partial})
   \;=\;
   \int \! 
\scrd_m(\psi, \psibar) \,
\scrd_{n-m}(\eta, \etabar) \,
     \exp\!\left[ \psibar^{\rm T} \partial \psietavert 
                  \,+\,
                  \etabar^{\rm T} B \psietavert
           \right]
\ef.
\end{equation}
By the translation formula \reff{eq.translation}, we have
\begin{equation}
\det(\widehat{\partial}) \, f(X)
\;=\;
   \int \! 
\scrd_m(\psi, \psibar) \,
\scrd_{n-m}(\eta, \etabar) \,
     \exp\!\left[ \etabar^{\rm T} B \psietavert \right] \,
f\Biggl( X \,+\,  \psibar \psietavert^{\!\! \rm T} \Biggr)
\end{equation}
for an arbitrary polynomial $f$.
We shall use this formula in the case $f(X)=\det(\widehat{X})^s$
where $s$ is a positive integer, so that
\be
\det(\widehat{\partial}) \,
\det(\widehat{X})^s
\;=\;
   \int \! 
\scrd_m(\psi, \psibar) \,
\scrd_{n-m}(\eta, \etabar) \,
     \exp\!\left[ \etabar^{\rm T} B \psietavert \right] \,
\det\Biggl[ \widehat{X} \,+\, \psibarzerovert \psietavert^{\!\! \rm T}
    \Biggr]^s
 \;.
   \label{eq.grassmann.XA.IJ}
\ee
It is convenient to introduce the shorthand
\be
(\widehat{X})^{\mathrm{trans}}
\;\equiv\;
\widehat{X} \,+\, \psibarzerovert \psietavert^{\!\! \rm T}
\ee
for the argument of $\det$.

Let us now assume that $\widehat{X}$ is an invertible real or complex
matrix, and change variables from $\displaystyle{\psietavert}$ to
$\displaystyle{\psietaprimevert = \widehat{X}^{-\rm T} \psietavert}$
with Jacobian $(\det \widehat{X})^{-1}$.
Dropping primes from the new variables, we observe that the expression for
the translated matrix can be written as
\begin{equation}
(\widehat{X})^{\mathrm{trans}}
\;=\;
\widehat{X}
\biggl[I \,+\,  \psibarzerovert \psietavert^{\!\! \rm T} \biggr]
   \;,
\end{equation}
so that
\begin{equation}
\det (\widehat{X})^{\mathrm{trans}}
\;=\;
(\det \widehat{X}) 
\det\! \Biggl[I \,+\,  \psibarzerovert \psietavert^{\!\! \rm T} \Biggr]
\ef.
\end{equation}
Applying Lemma~\ref{lemma.lowrank} to the rightmost determinant
yields
\be
\det\! \Biggl[I \,+\,  \psibarzerovert \psietavert^{\!\! \rm T} \Biggr]
  \;=\;
(1- \psibar^{\rm T} \psi)^{-1}
   \;,
\ee
so that we are left with the Grassmann-integral expression
\begin{subeqnarray}
\det(\widehat{\partial}) \,
\det(\widehat{X})^s
&=&
\det(\widehat{X})^{s-1}
   \int \! 
\scrd_m(\psi, \psibar) \,
\scrd_{n-m}(\eta, \etabar) \,
   \nonumber \\
&& \qquad \times
     \exp\!\left[ \etabar^{\rm T} B \widehat{X}^{\rm T} \psietavert \right]
     \: (1- \psibar^{\rm T} \psi)^{-s}
   \\[2mm]
&=&
\det(\widehat{X})^{s-1}
   \int \! 
\scrd_m(\psi, \psibar) \,
\scrd_{n-m}(\eta, \etabar) \,
   \nonumber \\
&& \qquad \times
     \exp\!\left[ \etabar^{\rm T} B X^{\rm T} \psi
                  \,+\, \etabar^{\rm T} B A^{\rm T} \eta
           \right]
     \: (1- \psibar^{\rm T} \psi)^{-s}
   \;. \qquad
\end{subeqnarray}
As the integrand depends on the Grassmann variables only through
combinations of the form $\psibar_i \psi_j$, $\etabar_i \eta_j$ and
$\etabar_i \psi_j$ (i.e.\ there is no $\psibar_i \eta_j$),
we can drop all the terms $\etabar_i \psi_j$,
as these terms would certainly remain unpaired in the expansion.
This removes all the dependence on $X$ in the integrand,
and proves that $\det(\widehat{\partial})$ is a
Bernstein--Sato operator for $\det(\widehat{X})$.
We are left with the determination of the prefactor $b(s)$,
which is given by
\be
   b(s)  \;=\;
   \int \! 
\scrd_m(\psi, \psibar) \,
\scrd_{n-m}(\eta, \etabar) \,
     \exp\!\left[ \etabar^{\rm T} B A^{\rm T} \eta \right]
     \:
(1- \psibar^{\rm T} \psi)^{-s}
   \;.
\ee
Integration over $\eta$ and $\etabar$
gives a factor $\det(B A^{\rm T}) = \det(A B^{\rm T})$,
while the integration over $\psi$ and $\psibar$
is identical to the one performed in the case of the ordinary Cayley identity
[cf.\ \reff{eq.grass1.star2}--\reff{eq.grass1.star3}]
and gives $s(s+1) \cdots (s+m-1)$.
This proves \reff{eq.borderedcayley.1} whenever $\widehat{X}$
is an invertible real or complex matrix,
$B$ is an arbitrary real or complex matrix,
and $s$ is a positive integer.
Now, if $A$ has rank $< n-m$, then both sides of
\reff{eq.borderedcayley.1} are identically zero;
while if $A$ has rank $n-m$, then
$\displaystyle{ \widehat{X} = \left(\!\! \begin{array}{cc} X \\ A \end{array}
                                    \!\! \right) }$
is invertible for a nonempty open set of matrices $X$.
The general validity of the identity therefore
follows from Proposition~\ref{prop.equiv.symbolic_s}.
\qed

\subsection[Relation between product-parametrized and border-parametrized
   identities]{Relation between product-parametrized and border-\hfill\break
   parametrized identities}
   \label{subsec.relation}

Let us begin by recalling the product-parametrized Cayley identity
\reff{eq.cayley.XA2m} and the border-parametrized Cayley identity
\reff{eq.borderedcayley.1}, writing the matrices $A$ and $B$
occurring in them as $A^{(0)}, B^{(0)}$ in the former identity
and $A^{(1)}, B^{(1)}$ in the latter:
\begin{eqnarray}
   \!\!\! \reff{eq.cayley.XA2m}:
   & \! &
   \det(\partial B^{(0)}) \, [\det (X A^{(0)})]^s  \;=\;
     b(s) \, \det(A^{(0)\rm T} B^{(0)}) \, [\det (X A^{(0)})]^{s-1}
     \qquad
        \label{eq.cayley.XA2m.BIS}   \\[2mm]
   \!\!\! \reff{eq.borderedcayley.1}:
   & \! &
   \det\!\left( \! \begin{array}{c}
                       \partial \\
                       \hline
                       B^{(1)}
                   \end{array} \! \right)
   \det\!\left( \! \begin{array}{c}
                       X \\
                       \hline
                       A^{(1)}
                   \end{array} \! \right)^{\! s}
   \;=\;
     b(s) \, \det(A^{(1)} B^{(1) \rm T}) \,
   \det\!\left( \! \begin{array}{c}
                       X \\
                       \hline
                       A^{(1)}
                   \end{array} \! \right)^{\! s-1}
     \qquad
        \label{eq.borderedcayley.1.BIS}
\end{eqnarray}
where $b(s) = s (s+1) \cdots (s+m-1)$.
Here $X$ is an $m \times n$ matrix,
while $A^{(0)}$ and $B^{(0)}$ are $n \times m$ matrices,
and $A^{(1)}$ and $B^{(1)}$ are $(n-m) \times n$ matrices.
Note that $A^{(0)}$ and $B^{(0)}$ must have full rank $m$,
otherwise \reff{eq.cayley.XA2m.BIS} is identically zero;
likewise, $A^{(1)}$ and $B^{(1)}$ must have full rank $n-m$.

We will construct the matrices $A^{(0)}$ and $A^{(1)}$
out of a larger ($n \times n$) matrix $A$, as follows:
Let $A$ be an invertible $n \times n$ matrix, and define
\begin{subeqnarray}
   A^{(0)} & = & A_{\star, [m]}  \;=\; \hbox{first $m$ columns of $A$} \\[1mm]
   A^{(1)} & = & (A^{-1})_{[m]^c,\star}  \;=\;
                                       \hbox{last $n-m$ rows of $A^{-1}$}
 \label{eq.A0.A1}
\end{subeqnarray}
Likewise, let $B$ be an invertible $n \times n$ matrix, and define
$B^{(0)}$ and $B^{(1)}$ by the same procedure.
We then have the following facts:

\begin{lemma}
  \label{lemma.relation}
Let $X$ and $Y$ be $m \times n$ matrices, with $m \le n$;
let $A$ and $B$ be invertible $n \times n$ matrices;
and define matrices $A^{(0)}, A^{(1)}, B^{(0)}, B^{(1)}$ as above. Then:
\begin{itemize}
   \item[(a)]
     $\displaystyle{ \det (X A^{(0)})  \;=\; (\det A) \,
                           \det\!\left( \! \begin{array}{c}
                                                X \\
                                                \hline
                                                A^{(1)}
                                            \end{array} \! \right) }$
   \item[(b)]
     $\displaystyle{ \det (Y B^{(0)})  \;=\; (\det B) \,
                           \det\!\left( \! \begin{array}{c}
                                                Y \\
                                                \hline
                                                B^{(1)}
                                            \end{array} \! \right) }$
   \item[(c)]
     $\displaystyle{ \det(A^{(0)\rm T} B^{(0)})  \;=\;
                     (\det A) \, (\det B) \, \det(A^{(1)} B^{(1) \rm T}) }$
\end{itemize}
\end{lemma}

\proof
(a)  First we expand the left-hand side using the Cauchy--Binet
identity \reff{eq.cauchy-binet}:
\be
  \det (X A^{(0)})
  \;=\;
  \sum_{|L|=m} (\det X_{\star L}) \, (\det A^{(0)}_{L \star})
  \;=\;
  \sum_{|L|=m} (\det X_{\star L}) \, (\det A_{L,[m]})
  \;.
 \label{eq.detXA0}
\ee
Next we expand the right-hand side using multi-row Laplace expansion
\reff{eq.multirowlaplace} with row set $[m]$,
followed by the Jacobi identity \reff{eq.JacobyDet}:
\begin{eqnarray}
   \det\!\left( \! \begin{array}{c}
                       X \\
                       \hline
                       A^{(1)}
                   \end{array} \! \right)
   & = &
   \sum_{|L|=m} \epsilon(L) \, (\det X_{\star L}) \, (\det A^{(1)}_{\star L^c})
        \nonumber \\[1mm]
   & = &
   \sum_{|L|=m} \epsilon(L) \, (\det X_{\star L}) \,
                               (\det A^{- \rm T}_{L^c, [m]^c})
        \nonumber \\[1mm]
   & = &
   \sum_{|L|=m} \epsilon(L) \, (\det X_{\star L}) \,
                \epsilon(L) \, (\det A)^{-1} \, (\det A_{L,[m]})
     \;.
 \label{eq.detXA1}
\end{eqnarray}
Comparing \reff{eq.detXA0} and \reff{eq.detXA1} proves (a);
and (b) is of course identical.

(c) First we expand the left-hand side using Cauchy--Binet:
\be
   \det(A^{(0)\rm T} B^{(0)})
   \;=\;
   \sum_{|L|=m} (\det (A^{(0) \rm T})_{\star L}) \, (\det B^{(0)}_{L \star})
   \;=\;
   \sum_{|L|=m} (\det A_{L,[m]}) \, (\det B_{L,[m]})
   \;.
 \label{eq.detA0TB0}
\ee
Next we expand the right-hand side using Cauchy--Binet
and then using the Jacobi identity twice:
\begin{eqnarray}
   \det(A^{(1)} B^{(1) \rm T})
   & = &
   \sum_{|K|=n-m}  (\det A^{(1)}_{\star K}) \, (\det (B^{(1) \rm T})_{K \star})
        \nonumber \\[1mm]
   & = &
   \sum_{|K|=n-m}  (\det (A^{- \rm T})_{K, [m]^c}) \,
                   (\det (B^{- \rm T})_{K, [m]^c})
        \nonumber \\[1mm]
   & = &
   \sum_{|K|=n-m}  \epsilon(K) \, (\det A)^{-1} \, (\det A_{K^c, [m]}) \,
                   \epsilon(K) \, (\det B)^{-1} \, (\det B_{K^c, [m]}) \,
        \nonumber \\[1mm]
   & = &
   (\det A)^{-1} \, (\det B)^{-1}
      \sum_{|L|=m} (\det A_{L,[m]}) \, (\det B_{L,[m]})
   \;.
  \label{eq.detA1B1T}
\end{eqnarray}
Comparing \reff{eq.detA0TB0} and \reff{eq.detA1B1T} proves (c).
\qed

Using Lemma~\ref{lemma.relation}, we see immediately the
equivalence of \reff{eq.cayley.XA2m.BIS} and \reff{eq.borderedcayley.1.BIS}
whenever $A^{(0)}$ and $A^{(1)}$ are related by \reff{eq.A0.A1}
and likewise for $B^{(0)}$ and $B^{(1)}$.

On the other hand, given any $n \times m$ matrix $A^{(0)}$ of rank $m$,
it can be obviously be completed to yield a nonsingular matrix $A$
(which is invertible at least when the matrix elements take values
in a field).  Likewise, an $(n-m) \times n$ matrix $A^{(1)}$ of rank $n-m$
can be completed to yield a nonsingular matrix $A^{-1}$.
So to each $A^{(0)}$ there corresponds a nonempty set of matrices $A^{(1)}$,
and vice versa.

{\bf Remark.}  It is {\em not}\/ in general true that every
{\em pair}\/ $(A^{(0)}, A^{(1)})$ arises from a matrix $A$.
Consider, for instance, $m=1$ and $n=2$:
an easy calculation shows that for arbitrary $A$
we must have $A^{(1)} A^{(0)} = 0$.

\section{Conjectures on minimality}   \label{sec.conj.minimality}

Let us recall that any pair $Q(s,x,\partial/\partial x)$
and $b(s) \not\equiv 0$ satisfying
\be
   Q(s,x,\partial/\partial x) \, P(x)^s
   \;=\;
   b(s) \, P(x)^{s-1}
\ee
is called a Bernstein--Sato pair for the polynomial $P(x)$.
The minimal (with respect to factorization) monic polynomial $b(s)$
for which there exists such a $Q$
is called the Bernstein--Sato polynomial (or $b$-function) of $P$.
Our Cayley-type identities thus provide Bernstein--Sato pairs
for certain polynomials arising from determinants.
But are our polynomials $b(s)$ minimal?

For the ordinary Cayley identities
(Theorems~\ref{thm.cayley}--\ref{thm.multirectcayley}),
it follows from the general theory of prehomogeneous vector spaces
\cite{Igusa_00,Kimura_03}
that the polynomials $b(s)$ found here are indeed minimal,
i.e.\ that the correct $b$-functions are
\begin{equation}
   b(s)  \;=\;
   \cases{s(s+1) \cdots (s+n-1) & for an $n \times n$ matrix \cr
          \noalign{\vskip 3mm}
          s(s+\smhalf) \cdots \left(s+ {n-1 \over 2} \right)
              & for an $n \times n$ symmetric matrix \cr
          \noalign{\vskip 3mm}
          s(s+2) \cdots (s+2m-2)
              & for a $2m \times 2m$ antisymmetric matrix (pfaffian) \cr
          \noalign{\vskip 3mm}
          (s-\smhalf)s(s+\smhalf) \cdots (s+m-1)
             \hspace*{-4cm}  & \cr
          \noalign{\vskip 1mm}
          \quad & for a $2m \times 2m$ antisymmetric matrix (determinant) \cr
          \noalign{\vskip 3mm}
          s(s+1) \cdots (s+m-1) (s+n-m) \cdots (s+n-1)
             \hspace*{-8cm}  & \cr
          \noalign{\vskip 1mm}
          \quad & for a pair of $m \times n$ rectangular matrices \cr
          \noalign{\vskip 3mm}
          s(s+\smhalf) \cdots \left(s+ {m-1 \over 2} \right)
             \left(s+ {n-m-1 \over 2} \right) \cdots
              \left(s+ {n-2 \over 2} \right)
             \hspace*{-8cm}  & \cr
          \noalign{\vskip 1mm}
          \quad & for an $m \times n$ rectangular matrix (symmetric) \cr
          \noalign{\vskip 3mm}
          s(s+2) \cdots (s+2m-2) (s+2n-2m+1) \cdots (s+2n-1)
             \hspace*{-8cm}  & \cr
          \noalign{\vskip 1mm}
          \quad & for an $m \times n$ rectangular matrix (antisymmetric) \cr
          \noalign{\vskip 3mm}
          \prod\limits_{\alpha=1}^\ell
                 \prod\limits_{j=0}^{n_1 -1} (s+n_\alpha-n_1 +j)
             & for $\ell$ matrices of sizes $n_\alpha \times n_{\alpha+1}$ \cr
         }
 \label{eq.minimality.list}
\end{equation}
Indeed, the polynomials $P$ occurring in these identities all correspond to
relative invariants of prehomogeneous vector spaces:
the ordinary, symmetric and antisymmetric Cayley identities
[\reff{eq.cayley.1}, \reff{eq.symcayley.1} and \reff{eq.antisymcayley.1}]
correspond to cases (1), (2) and (3), respectively,
in Kimura's \cite[Appendix]{Kimura_03} table of the irreducible reduced
prehomogeneous vector spaces;
the one-matrix rectangular symmetric and antisymmetric Cayley identities
[\reff{eq.rectcayley.1} and \reff{eq.antisymrectcayley.1}]
correspond to cases (15) and (13) in the same table;
while the two-matrix and multi-matrix rectangular Cayley identities
[\reff{eq.tmrectcayley.1} and \reff{eq.multirectcayley.1}]
correspond to prehomogeneous vector spaces
associated to equioriented quivers of type~${\sf A}$ \cite{Sugiyama_11}.
Whenever $P$ is a relative invariant of a prehomogeneous vector space,
the general theory \cite{Igusa_00,Kimura_03}
allows the immediate identification
of a suitable operator $Q(\partial/\partial x)$
--- namely, the dual of $P$ itself ---
and provides a proof that the corresponding $b(s)$
satisfies $\deg b = \deg P$
and is indeed (up to a constant factor)
the Bernstein--Sato polynomial of $P$.\footnote{
   See \cite[Corollary~6.1.1 and Theorem~6.1.1]{Igusa_00}
   \cite[Proposition~2.22]{Kimura_03}
   for the first two points,
   and \cite[Theorem~6.3.2]{Igusa_00} for the third.
   We are grateful to Nero Budur for explaining to us
   the connection between our results and
   the theory of prehomogeneous vector spaces,
   and in particular for pointing out that this connection
   provides a general proof of minimality.
}

For the Laplacian-parametrized identities
(Theorems~\ref{thm.Lap.cayley} and \ref{thm.Lap.sym.cayley}),
we conjecture that the polynomials $b(s)$ found here are also minimal,
i.e.\ that the correct $b$-functions are
\begin{equation}
   b(s)  \;=\;
   \cases{s(s+1) \cdots (s+n-2)
                & for a Laplacian-parametrized $n \times n$ matrix \cr
          \noalign{\vskip 3mm}
          s(s+\smhalf) \cdots \left(s+ {n-2 \over 2} \right)
              & for a Laplacian-parametrized $n \times n$ symmetric matrix \cr
         }
 \label{eq.minimality.list.Lap}
\end{equation}
Perhaps these identities can also be interpreted
within the framework of prehomogeneous vector spaces;
or perhaps an alternate proof of minimality can be found.

It is even conceivable that the following general fact
about Bernstein--Sato polynomials is true:

\begin{conjecture}
  \label{conj.bernstein-sato}
Let $P(x_1,\ldots,x_n) \not\equiv 0$ be a {\em homogeneous}\/ polynomial
in $n$ variables with coefficients in a field $K$ of characteristic 0,
and let $b(s)$ be its Bernstein--Sato polynomial.
Then $\deg b \ge \deg P$.
\end{conjecture}

\noindent
Simple examples show that we need not have $\deg b \ge \deg P$
if $P$ is not homogeneous:  for instance, $P(x) = 1-x^2$ has $b(s) = s$.
Moreover, slightly more complicated examples show that
one can have $\deg b > \deg P$ even when $P$ is homogeneous:
for instance, the Bernstein--Sato polynomial of
$P(x_1,x_2) = x_1 x_2 (x_1 + ax_2)$ with $a \neq 0$
is (in our ``shifted'' notation)
$s^2 (s-\smfrac{1}{3})(s+\smfrac{1}{3})$
\cite[Corollary~4.14 and Remark~4.15]{Walther_05}
\cite[5.4]{Saito_07}.\footnote{
   More generally, this is the Bernstein--Sato polynomial
   for a homogeneous polynomial $P(x_1,x_2)$ of degree 3
   with generic (e.g.\ random) coefficients \cite{Leykin_private}.
}
However, no one that we have consulted seems to have any counterexample
to Conjecture~\ref{conj.bernstein-sato}.

Let us remark that a {\em necessary}\/ condition
for a polynomial $b(s)$ to be a Bernstein--Sato polynomial
is that its roots should be {\em rational}\/ numbers $< 1$:
this is the content of a famous theorem of Kashiwara
\cite{Kashiwara_76} \cite[Chapter~6]{Bjork_79}
 \cite[Proposition~2.11]{Leykin_01}.\footnote{
   This is in our ``shifted'' notation \reff{eq.intro.bernstein}.
   In the customary notation, the roots are rational numbers
   that are strictly negative.
   See footnote~\ref{footnote.shifted} above.
}
Our polynomials \reff{eq.minimality.list} and \reff{eq.minimality.list.Lap}
satisfy this condition.

For the diagonal-parametrized Cayley identities
(Theorems~\ref{thm.para.cayley} and \ref{thm.para.sym.cayley}),
a slightly more complicated situation arises.
The polynomials $b(s)$ arising from the basic case $I=J=[n]$ of those theorems,
namely
\begin{equation}
   b(s)  \;=\;
   \cases{s(s+1) \cdots (s+n-1) & for an $n \times n$ matrix \cr
          \noalign{\vskip 3mm}
          s(s+\smhalf) \cdots \left(s+ {n-1 \over 2} \right)
              & for an $n \times n$ symmetric matrix \cr
         }
\end{equation}
are definitely {\em not}\/ minimal.
Indeed, as remarked already in Section~\ref{sec.statement.para},
a lower-order Bernstein--Sato pair can be obtained by taking
$I=J = [n] \setminus \{i_0\}$ for any fixed $i_0 \in [n]$:
\begin{equation}
   b(s)  \;=\;
   \cases{s(s+1) \cdots (s+n-2) & for an $n \times n$ matrix \cr
          \noalign{\vskip 3mm}
          s(s+\smhalf) \cdots \left(s+ {n-2 \over 2} \right)
              & for an $n \times n$ symmetric matrix \cr
         }
\end{equation}
We conjecture that these latter polynomials are indeed minimal,
but we have no proof for general $n$.

It is curious that these polynomials are the \emph{same}\/
as we get for the Laplacian-parametrized Cayley identities.
Furthermore, also the $Q$ operators corresponding to these
(conjecturally minimal) polynomials $b(s)$ of degree $n-1$
are somewhat similar in one respect:
namely, for the diagonal-parametrized case the $Q$ operator
is given by any principal minor of size $n-1$
(i.e.\ $I=J=[n] \setminus \{i_0\}$)
of the relevant matrix $D_{\balpha,\bbeta,s}$ of differential operators,
while for the Laplacian-parametrized case
the $Q$ operator is a polynomial $\det(U+\partial) - \det(\partial)$
that is a \emph{sum}\/ over all minors (not necessarily principal)
of size $n-1$.
We do not know whether this resemblance is indicative of
any deeper connection between these identities.



\bigskip

{\bf Note Added:}
After this paper appeared in preprint form,
Nero Budur (private communication, July 2011)
informed us that he has a proof of our minimality conjectures
for the Laplacian-parametrized and diagonal-parametrized Cayley identities.

\appendix

\section{Grassmann algebra and Gaussian integration}
\label{app.grassmann}

In this appendix we collect some needed information on Grassmann algebra
(= exterior algebra)
and Gaussian integration (both bosonic and fermionic).
We begin by recalling the main properties of
determinants, permanents, pfaffians and hafnians
(Section~\ref{app.grassmann.1}).
We then recall the well-known properties of ``bosonic'' Gaussian integration,
i.e.\ Gaussian integration over $\R^n$ or $\C^n$
(Section~\ref{app.grassmann.2}).
Next we define Grassmann algebra (Section~\ref{app.grassmann.3})
and Grassmann--Berezin (``fermionic'') integration
(Section~\ref{app.grassmann.4}).
Finally, we explain the formulae for fermionic Gaussian integration
(Section~\ref{app.grassmann.5}),
which will play a central role in this paper.
Our presentation in these latter three subsections
is strongly indebted to Abdesselam \cite[Section~2]{Abdesselam_03};
see also Zinn-Justin \cite[Chapter~1]{Zinn-Justin}
for a treatment aimed at physicists.

\subsection{Determinants, permanents, pfaffians and hafnians}
   \label{app.grassmann.1}

{\bf Notation:}
If $A$ is an $m \times n$ matrix,
then for subsets of indices $I \subseteq [m]$ and $J \subseteq [n]$
we denote by $A_{IJ}$ the matrix $A$ restricted to rows in $I$
and columns in $J$, all kept in their original order.
We also use the shorthand notation
$A_{\star J} = A_{[m]\,J}$ when all the rows are kept,
and $A_{I \star} = A_{I\,[n]}$ when all the columns are kept.
Finally, if $A$ is invertible,
we denote by $A^{-{\rm T}}$ the matrix $(A^{-1})^{\rm T} = (A^{\rm T})^{-1}$.
\qed

\subsubsection{Permanent and determinant}

Let $R$ be a commutative ring with identity;
we shall consider matrices with entries in $R$.
In particular, if $A = (a_{ij})_{i,j=1}^n$ is an $n \times n$ matrix
with entries in $R$,
we define its {\em permanent}\/
\be
   \per A  \;=\;  \sum_{\sigma \in \scrs_n}
                  a_{1,\sigma(1)} \cdots a_{n,\sigma(n)}
\ee
and its {\em determinant}\/
\be
   \det A  \;=\;  \sum_{\sigma \in \scrs_n}  \sgn(\sigma) \,
                  a_{1,\sigma(1)} \cdots a_{n,\sigma(n)}
   \;.
\ee
Here the sums range over all permutations $\sigma$ of
$[n] \equiv \{1,\ldots,n\}$,
and $\sgn(\sigma) = (-1)^{\#(\hboxrm{\rm even cycles of }\sigma)}$
is the sign of the permutation $\sigma$.
See \cite{Minc_78} and \cite{Prasolov_94}
for basic information on permanents and determinants, respectively.

%

In this paper we shall need only a few of the most elementary properties of
determinants:

\begin{lemma}[Properties of the determinant]
   \label{lemma.properties.det}
\quad\break
\vspace*{-5mm}
\begin{itemize}
   \item[(a)] $\det I = 1$.
   \item[(b)] $\det(AB) = (\det A)(\det B)$.
%
   \item[(c)] {\bf (Cauchy--Binet formula)\ }
      More generally, let $A$ be an $m \times n$ matrix,
      and let $B$ be an $n \times m$ matrix.  Then
\be
   \det(AB)  \;=\; \sum_{\begin{scarray}
                           I \subseteq [n] \\
                           |I| = m
                         \end{scarray}}
        (\det A_{\star I}) (\det B_{I \star})
   \;.
 \label{eq.cauchy-binet}
\ee 
   \item[(d)]  Let $A$ be an $n \times n$ matrix,
and define the {\bf adjugate matrix} $\adj A$ by
\be
   (\adj A)_{ij}  \;=\;  (-1)^{i+j} \, \det A_{\{j\}^c \{i\}^c}
\ee
(note the transpose between the left-hand and right-hand sides).
Then
\be
   (\adj A) \, A  \;=\;  A \, (\adj A)  \;=\;  (\det A) \, I
   \;.
\ee
In particular, $A$ is invertible in the ring $R^{n \times n}$
if and only if $\det A$ is invertible in the ring $R$
(when $R$ is a field, this means simply that $\det A \neq 0$);
and in this case
\be
   A^{-1}  \;=\;  (\det A)^{-1} \, (\adj A)
\ee
{\bf (Cramer's rule)}.
   \item[(e)] {\bf (Jacobi's identity)\ }
More generally, if $I,J \subseteq [n]$ with $|I|=|J|=k$, then
\be
  \det( (A^{-\rm T})_{IJ} )
  \;=\;
  (\det A)^{-1} \, \epsilon(I,J) \, (\det A_{I^c J^c})
 \label{eq.JacobyDet}
\ee
where $\epsilon(I,J) = (-1)^{\sum_{i \in I} i + \sum_{j \in J} j}$.
   \item[(f)] {\bf (Multi-row Laplace expansion)\ }
For any fixed set of rows $I \subseteq [n]$ with $|I|=k$, we have
\be
   \det A  \;=\;
   \sum_{\begin{scarray}
             J \subseteq [n] \\
             |J| = k
         \end{scarray}}
   \epsilon(I,J) \, (\det A_{IJ}) \, (\det A_{I^c J^c})
   \;.
 \label{eq.multirowlaplace}
\ee
\end{itemize}
\end{lemma}


\subsubsection{Hafnian}

Let $A = (a_{ij})_{i,j=1}^{2m}$ be a $2m \times 2m$
{\em symmetric}\/ matrix
with entries in $R$.
We then define the {\em hafnian}\/ \cite{Caianiello_73}
\be
   \hf A  \;=\;  \sum_{M \in \scrm_{2m}} \, \prod_{ij \in M} a_{ij}
   \;,
\ee
where the sum runs over all {\em perfect matchings}\/
of the set $[2m]$,
i.e.\ all partitions of the set $[2m]$ into $m$ disjoint pairs.
There are $(2m-1)!! = (2m)! / (2^m m!)$ terms in this sum.
We have, for example,
\begin{subeqnarray}
   \hf \left( \begin{array}{cc}
                  a_{11}  &  a_{12}  \\
                  a_{12}  & a_{22}
              \end{array}
       \right)
   & = &
   a_{12}
       \\[2mm]
   \hf \left( \begin{array}{cccc}
                  a_{11}  &  a_{12} & a_{13} & a_{14}  \\
                  a_{12}  & a_{22} & a_{23} & a_{24}  \\
                  a_{13}  & a_{23} & a_{33} & a_{34}  \\
                  a_{14}  & a_{24} & a_{34} & a_{44}
              \end{array}
       \right)
   & = &
   a_{12} a_{34} \,+\, a_{13} a_{24} \,+\, a_{14} a_{23}
\end{subeqnarray}
Note that the diagonal elements of $A$ play no role in the hafnian.

Equivalently, we can identify matchings with a subclass of permutations
by writing each pair $ij \in M$ in the order $i<j$
and then writing these pairs in increasing order of their first elements:
we therefore have
\be
   \hf A  \;=\;  \sum_{\sigma \in \scrs^\star_{2m}}
                  a_{\sigma(1) \sigma(2)} \cdots a_{\sigma(2m-1) \sigma(2m)}
   \;,
\ee
where the sum runs over all permutations $\sigma$ of $[2m]$
satisfying $\sigma(1) < \sigma(3) < \ldots < \sigma(2m-1)$
and $\sigma(2k-1) < \sigma(2k)$ for $k=1,\ldots,m$.

Note that if we were to sum over {\em all}\/ permutations,
we would obtain each term in $\hf A$ exactly $2^m m!$ times.
Therefore, {\em if the ring $R$ contains the rationals}\/,
we can alternatively write
\be
   \hf A  \;=\;
   {1 \over 2^m m!}
   \sum_{\sigma \in \scrs_{2m}}  
                  a_{\sigma(1) \sigma(2)} \cdots a_{\sigma(2m-1) \sigma(2m)}
   \;.
\ee


\subsubsection{Pfaffian}

Finally, let $A = (a_{ij})_{i,j=1}^{2m}$ be a $2m \times 2m$
{\em antisymmetric}\/ matrix (i.e.\ $a_{ij} = -a_{ji}$ and $a_{ii} = 0$)
with entries in $R$.\footnote{
   Such matrices are sometimes called {\em alternating}\/ matrices,
   in order to emphasize that the condition $a_{ii} = 0$ is imposed.
   This latter condition is a consequence of $a_{ij} = -a_{ji}$
   whenever $R$ is an integral domain of characteristic $\neq 2$
   (so that $2x=0$ implies $x=0$),
   but not in general otherwise.
   See e.g.\ \cite[section~XV.9]{Lang_02}.
   In this paper we use the term ``antisymmetric'' to denote
   $a_{ij} = -a_{ji}$ {\em and}\/ $a_{ii} = 0$.
 \label{note_alternating}
}
We then define the {\em pfaffian}\/ by
\be
   \pf A  \;=\;  \sum_{M \in \scrm_{2m}}
       \!\!
       \epsilon(\vec{M},\vec{M}_0)
       \!\!
       \prod_{\begin{scarray}
                 (i,j) \!\in\! \vec{M} \\
                 i < j
              \end{scarray}}
      \!\!\!
        a_{ij}
   \;.
 \label{def.pfaffian.app}
\ee
Here the sum runs once again over all perfect matchings $M$ of the set $[2m]$,
and $\vec{M}$ is an (arbitrarily chosen) oriented version of $M$,
i.e.\ for each unordered pair $ij \in M$ one chooses an ordering $(i,j)$
of the two elements.
The value of the summand in \reff{def.pfaffian.app}
will be independent of the choice of $\vec{M}$
because $\epsilon(\vec{M},\vec{M}_0)$ will be odd
under reorderings of pairs (see below), while $A$ is antisymmetric.
Here $\vec{M}_0$ is some fixed oriented perfect matching of $[2m]$
(we call it the ``reference matching'').
The sign $\epsilon(\vec{M}_1,\vec{M}_2)$ is defined as follows:
If $\vec{M}_1 = \{(i_1,i_2), (i_3,i_4), \ldots, (i_{2m-1},i_{2m})\}$
and $\vec{M}_2 = \{(j_1,j_2), (j_3,j_4), \ldots, (j_{2m-1},j_{2m})\}$,
then $\epsilon(\vec{M}_1,\vec{M}_2)$ is the sign of the permutation
that takes $i_1 \cdots i_{2m}$ into $j_1 \cdots j_{2m}$.
(This is well-defined, i.e.\ independent of the order in which the
 ordered pairs of $\vec{M}_1$ and $\vec{M}_2$ are written,
 because interchanging two pairs is an even permutation.)
This quantity is clearly odd under reorderings of pairs
in $\vec{M}_1$ or $\vec{M}_2$, and has the following properties:
\begin{itemize}
   \item[(a)] $\epsilon(\vec{M}_1,\vec{M}_2) = \epsilon(\vec{M}_2,\vec{M}_1)$;
   \item[(b)] $\epsilon(\vec{M},\vec{M}) = +1$;
   \item[(c)] $\epsilon(\vec{M}_1,\vec{M}_2) \, \epsilon(\vec{M}_2,\vec{M}_3) =
        \epsilon(\vec{M}_1,\vec{M}_3)$;
   \item[(d)] $\epsilon(\vec{M}_1,\vec{M}_2) = -1$
        whenever $\vec{M}_1$ and $\vec{M}_2$
        differ by reversal of the orientation of a single edge;
   \item[(e)] $\epsilon(\vec{M}_1,\vec{M}_2) = -1$
        whenever $\vec{M}_1$ and $\vec{M}_2$
        differ by changing directed edges $(a,b),(c,d)$ in $\vec{M}_1$
        to $(b,c),(d,a)$ in $\vec{M}_2$.
\end{itemize}
Indeed, it is not hard to show that $\epsilon(\vec{M}_1,\vec{M}_2)$
is the {\em unique}\/ map from pairs of oriented perfect matchings
into $\{\pm 1\}$ that has these five properties.\footnote{
   This also implies that $\epsilon(\vec{M}_1,\vec{M}_2)$
   can be given an equivalent (more graph-theoretic) definition as follows:
   Form the union $\vec{M}_1 \cup \vec{M}_2$.
   Ignoring orientations, it is a disjoint union of even-length cycles.
   Looking now at the orientations, let us call a cycle
   {\em even}\/ (resp.\ {\em odd}\/) if it has an even (resp.\ odd)
   number of edges pointing in each of the two directions around the cycle.
   We then have
   $\epsilon(\vec{M}_1,\vec{M}_2) = (-1)^{\#(\hboxrm{\rm even cycles})}$.
}
It follows from (a)--(c) that the oriented perfect matchings
fall into two classes (call them A and B)
such that $\epsilon(\vec{M}_1,\vec{M}_2)$
equals $+1$ if $\vec{M}_1$ and $\vec{M}_2$ belong to the same class
and $-1$ if they belong to different classes.
The choice of reference matching $\vec{M}_0$
really amounts, therefore, to choosing one of the
two equivalence classes of matchings as the reference class,
and thereby fixing the sign of the pfaffian.

The choice of $\vec{M}_0$ can be encoded in an antisymmetric matrix $J$ defined by
\be
   J_{ij}
   \;=\;
   \cases{ 1  & if $(i,j) \in \vec{M}_0$  \cr
           \noalign{\vskip 3pt}
           -1 & if $(j,i) \in \vec{M}_0$  \cr
           \noalign{\vskip 3pt}
           0  & otherwise
         }
\ee
and satisfying $\pf(J) = 1$.
The two most common conventions for the reference matching are
\be
   \vec{M}_0 \:=\: \{(1,2), (3,4), \,\ldots,\, (2m\!-\!1,2m) \} \,,
   \qquad
   J  \;=\;
   \left(
   \begin{array}{cc|cc|c}
    0 & 1 & \multicolumn{3}{c}{} \\
   -1 & 0 & \multicolumn{3}{c}{} \\
   \cline{1-4}
    &  & 0 & 1 & \\
    &  & -1 & 0 & \\
   \cline{3-4}
   \multicolumn{4}{c}{} & \ddots
   \end{array}
   \right)
 \label{def.J.appendix}
\ee
and
\be
   \vec{M}_0 \:=\: \{(1,m+1), \, (2,m+2), \,\ldots,\, (m,2m) \}  \,,
   \qquad
   J  \;=\;
   \left( \!
   \begin{array}{cc}
        0   & I_m  \\
       -I_m & 0
   \end{array}
   \!\right)
   \:.
\ee
In this paper we shall adopt the convention \reff{def.J.appendix}.
We thus have
\begin{subeqnarray}
   \pf \left( \begin{array}{cc}
                  0  &  a_{12}  \\
                  -a_{12}  & 0
              \end{array}
       \right)
   & = &
   a_{12}
       \\[2mm]
   \pf \left( \begin{array}{cccc}
                  0  &  a_{12} & a_{13} & a_{14}  \\
                  -a_{12}  & 0 & a_{23} & a_{24}  \\
                  -a_{13}  & -a_{23} & 0 & a_{34}  \\
                  -a_{14}  & -a_{24} & -a_{34} & 0
              \end{array}
       \right)
   & = &
   a_{12} a_{34} \,-\, a_{13} a_{24} \,+\, a_{14} a_{23}
\end{subeqnarray} 

By identifying matchings $M \in \scrm_{2m}$
with permutations $\sigma \in \scrs^\star_{2m}$
as was done for the hafnian, we can equivalently write
\be
   \pf A  \;=\;  \sgn(\sigma_0)
                \sum_{\sigma \in \scrs^\star_{2m}}  \sgn(\sigma) \,
                  a_{\sigma(1) \sigma(2)} \cdots a_{\sigma(2m-1) \sigma(2m)}
 \label{def.pfA}
\ee
where $\sigma_0 \in \scrs^\star_{2m}$ is the permutation
corresponding to the reference matching $M_0$.
[For our choice \reff{def.J.appendix}, $\sigma_0$ is the identity permutation.]
{\em If the ring $R$ contains the rationals}\/,
we can alternatively write
\be
   \pf A  \;=\;
   {1 \over 2^m m!} \: \sgn(\sigma_0)
   \sum_{\sigma \in \scrs_{2m}}  \sgn(\sigma) \,
                  a_{\sigma(1) \sigma(2)} \cdots a_{\sigma(2m-1) \sigma(2m)}
   \;.
\ee

Let us now recall the following basic properties of pfaffians:

\begin{lemma}[Properties of the pfaffian]
   \label{lemma.properties.pfaffians}
Let $A$ be an antisymmetric $2m \times 2m$ matrix
with elements in a commutative ring $R$.  Then:
\begin{itemize}
   \item[(a)] $\pf J = 1$.
   \item[(b)]  $(\pf A)^2 = \det A$.
   \item[(c)]  $\pf(X A X^{\rm T}) = (\det X) (\pf A)$
      for any $2m \times 2m$ matrix $X$.
   \item[(d)]  {\bf (minor summation formula for pfaffians
        \cite{Ishikawa_95,Ishikawa_06})\ }
      More generally, we have
\be
    \pf(X A X^{\rm T})   \;=\;
    \sum\limits_{\begin{scarray}
                    I \subseteq [2m] \\
                    |I| = 2\ell
                 \end{scarray}}
    (\det X_{\star I}) \, (\pf A_{II})
 \label{eq.app.pfaffian.minorsummation}
\ee
 for any $2\ell \times 2m$ matrix $X$ ($\ell \le m$).
 Here $X_{\star I}$ denotes the submatrix of $X$ with columns $I$
 (and all its rows). 
   \item[(e)]  {\bf (Jacobi's identity for pfaffians)\ }
      If $A$ is invertible, then $\pf(A^{-\rm T}) = (\pf A)^{-1}$
      and more generally
\be
  \pf( (A^{-\rm T})_{II} )
  \;=\;
   \epsilon(I) \, (\pf A)^{-1} \, (\pf A_{I^c I^c})
 \label{eq.JacobyPf}
 \label{eq.app.pfaffian.jacobi}
\ee
for any $I \subseteq [2m]$,
where
    $\epsilon(I) = (-1)^{|I| (|I|-1)/2} (-1)^{\sum_{i \in I} i}$.
%
\end{itemize}
\end{lemma}

See \cite{Stembridge_90,Knuth_96,Fulton_98,Hamel_01,Lang_02,Ishikawa_06,%
Fulmek_10}
for further information on pfaffians.

\medskip

{\bf Remark.}
In this paper we will not in fact use the minor summation formula
for pfaffians;  but we will rederive it using Grassmann--Berezin integration.
See Theorem~\ref{thm.wick.realfermions} and the comments following it.

\subsection{Bosonic Gaussian integration}
   \label{app.grassmann.2}

We shall use the following notation:
If $A = (a_{ij})$ is an $m \times n$ matrix,
and $I = (i_1,\ldots,i_k)$ and $J = (j_1,\ldots,j_\ell)$
are sequences of indices (not necessarily distinct or ordered)
in $[m]$ and $[n]$, respectively,
then we denote by $A_{IJ}$ the $k \times \ell$ matrix defined by
\begin{equation}
   (A_{IJ})_{\alpha\beta}  \;=\;  a_{i_\alpha j_\beta}
   \;.
 \label{eq.AIJ.sequence}
\end{equation}
This generalizes our notation $A_{IJ}$ for {\em subsets}\/
$I \subseteq [m]$ and $J \subseteq [n]$,
where a subset is identified with the sequence of its elements
written in increasing order.
We shall also use the corresponding notation for vectors:
namely, if $\lambda = (\lambda_i)$ is an $n$-vector
and $I = (i_1,\ldots,i_k)$ is a sequence of indices
(not necessarily distinct or ordered) in $[n]$,
then we denote by $\lambda_I$ the $k$-vector defined by
$(\lambda_I)_\alpha = \lambda_{i_\alpha}$.

Let $\varphi = (\varphi_i)_{i=1}^n$ be real variables;
we shall write
\begin{equation}
   \scrd\varphi  \;=\;
   \prod_{i=1}^n {d \varphi_i \over \sqrt{2\pi}}
\end{equation}
for Lebesgue measure on $\R^n$ with a slightly unconventional normalization.
Let $A = (a_{ij})_{i,j=1}^n$ be a real symmetric positive-definite
$n \times n$ matrix.  We then have the following fundamental facts
about Gaussian integration on $\R^n$:

\begin{theorem}[Wick's theorem for real bosons]
 \label{thm.wick.realbosons}
Let $A = (a_{ij})_{i,j=1}^n$ be a real symmetric positive-definite
$n \times n$ matrix.  Then:
\begin{itemize}
   \item[(a)]  For any vector $c = (c_i)_{i=1}^n$ in $\R^n$ (or $\C^n$), we have
\begin{equation}
   \int\! \scrd\varphi \,
   \exp\!\left( -\smhalf \varphi^{\rm T} A \varphi \,+\, c^{\rm T} \varphi
         \right)
   \;=\;
   (\det A)^{-1/2} \,
   \exp\!\left( \smhalf  c^{\rm T} A^{-1} c \right)
   \;.
 \label{eq.wick.realbosons.1}
\end{equation}
   \item[(b)]  For any sequence of indices $I = (i_1,\ldots,i_r)$ in $[n]$,
we have
\begin{equation}
   \int\! \scrd\varphi \: \varphi_{i_1} \cdots \varphi_{i_r} \,
   \exp\!\left( -\smhalf \varphi^{\rm T} A \varphi \right)
   \;=\;
   \cases{ 0                               & if $r$ is odd \cr
           \noalign{\vskip 2mm}
           (\det A)^{-1/2} \, \hf((A^{-1})_{II})  & if $r$ is even \cr
         }
 \label{eq.wick.realbosons.2}
\end{equation}
   \item[(c)]
More generally, for any real or complex $r \times n$ matrix $C$, we have
\begin{equation}
   \int\! \scrd\varphi \,
      \Biggl( \prod_{\alpha=1}^r (C\varphi)_\alpha \Biggr) \,
   \exp\!\left( -\smhalf \varphi^{\rm T} A \varphi \right)
   \;=\;
   \cases{ 0                               & if $r$ is odd \cr
           \noalign{\vskip 2mm}
           (\det A)^{-1/2} \, \hf(C A^{-1} C^{\rm T})  & if $r$ is even \cr
         }
 \label{eq.wick.realbosons.3}
\end{equation}
\end{itemize}
\end{theorem} 

{\bf Historical remarks.}
Physicists call these formulae ``Wick's theorem'' because
Gian-Carlo Wick \cite{Wick_50} proved the analogue of
\reff{eq.wick.realbosons.2} for the correlation functions
of a free quantum field (see e.g.\ \cite{Streater_64}).
These formulae are called ``bosonic'' because the functional-integral
formulation for bosonic quantum fields (see e.g.\ \cite{Zinn-Justin})
leads to ordinary integrals over $\R^n$ or $\C^n$
(or infinite-dimensional generalizations thereof).
By contrast, functional integrals for fermionic quantum fields
lead to Grassmann--Berezin integrals,
to be discussed in Sections~\ref{app.grassmann.3}--\ref{app.grassmann.5}.

The formula \reff{eq.wick.realbosons.2} for the moments of a
mean-zero Gaussian measure
goes back at least to Isserlis \cite{Isserlis_18} in 1918.
We thank Malek Abdesselam for drawing our attention to this reference.
\qed


\medskip
   
Now let $\varphi = (\varphi_i)_{i=1}^n$ be complex variables;
we denote complex conjugation by $\overline{\phantom{\varphi}}$
and shall write
\begin{equation}
   \scrd(\varphi,\varphibar)  \;=\;
   \prod_{i=1}^n {(d \real\varphi_i)(d \imag\varphi_i) \over \pi}
\end{equation}
for Lebesgue measure on $\C^n$ with a slightly unconventional normalization.
Let $A = (a_{ij})_{i,j=1}^n$ be an $n \times n$ complex matrix
(not necessarily symmetric or hermitian)
whose hermitian part $\smhalf (A + A^*)$ is positive-definite.
We then have the following fundamental facts
about Gaussian integration on $\C^n$:

\begin{theorem}[Wick's theorem for complex bosons]
 \label{thm.wick.complexbosons}
Let $A = (a_{ij})_{i,j=1}^n$ be an $n \times n$ complex matrix
whose hermitian part $\smhalf (A + A^*)$ is positive-definite.  Then:
\begin{itemize}
   \item[(a)]  For any vectors $b = (b_i)_{i=1}^n$ and $c = (c_i)_{i=1}^n$
      in $\C^n$, we have
\begin{equation}
   \int\! \scrd(\varphi,\varphibar) \,
   \exp\!\left( - \varphibar^{\rm T} A \varphi \,+\, \bar{b}^{\rm T} \varphi
                                             \,+\, \varphibar^{\rm T} c \right)
   \;=\;
   (\det A)^{-1} \,
   \exp\!\left( \bar{b}^{\rm T} A^{-1} c \right)
   \;.
 \label{eq.wick.complexbosons.1}
\end{equation}
   \item[(b)]  For any sequences of indices $I = (i_1,\ldots,i_r)$
      and $J = (j_1,\ldots,j_s)$ in $[n]$, we have
\begin{equation}
   \int\! \scrd(\varphi,\varphibar) \;
  \varphi_{i_1} \cdots \varphi_{i_r} \,
  \varphibar_{j_1} \cdots \varphibar_{j_s}
   \exp\!\left( - \varphibar^{\rm T} A \varphi \right)
   \;=\;
   \cases{ 0                               & if $r \neq s$ \cr
           \noalign{\vskip 2mm}
           (\det A)^{-1} \, \per((A^{-1})_{IJ})  & if $r=s$ \cr
         }
 \label{eq.wick.complexbosons.2}
\end{equation}
%
   \item[(c)]
More generally, for any complex $r \times n$ matrix $B$
and any complex $n \times s$ matrix $C$, we have
\begin{equation}
   \!\!\!
   \int\! \scrd(\varphi,\varphibar)
                  \, \Biggl( \prod_{\alpha=1}^r (B\varphi)_\alpha \Biggr)
                     \Biggl( \prod_{\beta=1}^s (\varphibar^{\rm T} C)_\beta
                           \Biggr)
   \exp\!\left( - \varphibar^{\rm T} A \varphi \right)
   \;=\;
   \cases{ 0                                        & if $r \neq s$ \cr
           \noalign{\vskip 2mm}
           (\det A)^{-1} \, \per(B A^{-1} C) \!\!\! &  if $r=s$ \cr
         }
 \label{eq.wick.complexbosons.3}
\end{equation}
\end{itemize}
\end{theorem} 

%

\bigskip

{\bf Some final remarks.}
1.  In this article we shall use mainly the ``source'' versions of
Wick's theorem, i.e.\ part (a) of
Theorems~\ref{thm.wick.realbosons} and \ref{thm.wick.complexbosons}
and the corresponding theorems for fermions.
The ``correlation function'' versions, i.e.\ parts (b) and (c),
will be used only in Sections~\ref{sec.grassmann.1} and \ref{subsec.product}
and in the second proofs of
Corollaries~\ref{cor.pbABp.Marco} and \ref{cor.Atheta.Marco}
(Section~\ref{sec.grassmann.2.lemmas}).

2.  We have here presented bosonic Gaussian integration
in an analytic context, i.e.\ integration on $\R^n$ or $\C^n$.
A combinatorial abstraction of bosonic Gaussian integration
can be found in \cite{Abdesselam_feynman}.

\subsection{Grassmann algebra}
   \label{app.grassmann.3}

%
%
%
%
%
%
%

Let $R$ be a commutative ring.
Every textbook on elementary abstract algebra defines
the ring $R[x_1,\ldots,x_n]$ of polynomials
in commuting indeterminates $x_1,\ldots,x_n$
with coefficients in $R$, and studies its properties.
Here we would like briefly to do the same for
the ring $R[\chi_1,\ldots,\chi_n]_{\rm Grass}$ of polynomials
in {\em anticommuting}\/ indeterminates $\chi_1,\ldots,\chi_n$:
we call this ring the {\em Grassmann algebra}\/ over $R$
in generators $\chi_1,\ldots,\chi_n$.
(Of course, readers familiar with exterior algebra
will recognize this as nothing other than the exterior algebra $\Lambda(R^n)$
built from the free $R$-module of dimension $n$.\footnote{
   See e.g.\ \cite[section~6.4]{Cohn_03} or \cite[section~XIX.1]{Lang_02}.
})
To lighten the notation, we shall henceforth omit the subscripts
${}_{\rm Grass}$, since it will always be clear by context
whether we are referring to the Grassmann algebra
or to the ordinary polynomial ring.

Here is the precise definition:

\begin{definition}[Grassmann algebra]
Let $R$ be a commutative ring with identity element\/\footnote{
   Much of the elementary theory works also for coefficient rings
   without identity element.  The main change is that polynomials
   and formal power series $\Phi$ must have constant term $c_0$ in $R$
   (and not in $\Z$) in order to be applied to elements $f \in R[\chi]_+$:
   see the paragraphs immediately after Proposition~\ref{prop.puresoul}.
   But this means that we can consider the most important case,
   namely $\Phi = \exp$, only when $R$ has an identity element.
   For this reason it is convenient simply to make this assumption
   from the beginning.
},
and let $\ch_1,\ldots,\ch_n$ be a collection of letters.
The Grassmann algebra $R[\ch_1,\ldots,\ch_n]$
(or $R[\ch]$ for short)
is the quotient of the ring $R\< \ch_1,\ldots,\ch_n\>$
of noncommutative polynomials in the letters $\ch_1,\ldots,\ch_n$
by the two-sided ideal generated by the expressions
$\ch_i \ch_j+ \ch_j \ch_i$ ($1\le i < j\le n$)
and $\ch_i ^2$ ($1\le i\le n$).
We can consider $R[\ch]$ as a ring and also as an $R$-algebra.
\end{definition}

In other words, the generators $\ch_i$ of $R[\ch]$
satisfy the {\em anticommutation relations}\/
\be
    \ch_i \ch_j+ \ch_j \ch_i \;=\; 0
   \qquad\hbox{for all $i,j \in [n]$}
\label{anti}
\ee
as well as the relations
\be
   \chi_i^2  \;=\; 0
   \qquad\hbox{for all $i \in [n]$} \;.
\label{anti2}
\ee
Please note that the anticommutation relation \reff{anti} for $i=j$
states that $2 \ch_i^2=0$;  but this need not imply $\ch_i^2=0$
if the coefficient ring $R$ does not contain an element $\smfrac{1}{2}$.
For this reason we have explicitly adjoined the relations $\chi_i^2 = 0$.
Of course, if the coefficient ring $R$ contains an element $\smfrac{1}{2}$
(e.g.\ if it contains the rationals),
then this extra relation could be replaced by the cases $i=j$ of
$\ch_i \ch_j+ \ch_j \ch_i = 0$.

The first important property of $R[\ch]$ is the following:

\begin{proposition}
$R[\ch]$ is a free $R$-module with basis given by the
$2^n$ monomials $\chi^I = \ch_{i_1} \cdots \ch_{i_p}$
where $I = \{i_1,\ldots,i_p\} \subseteq [n]$ with $i_1<\ldots<i_p$.
\end{proposition}

It follows that each element $f\in R[\ch]$ can be
written uniquely in the form
\be
f \;=\;
\sum_{I\subseteq [n]} f_I \, \chi^I
\label{fexp}
\ee
with $f_I \in R$.
The term $f_\emptyset$ that contains no factors $\chi_i$
is sometimes termed the {\em body}\/ of $f$,
and the rest $\sum_{I \neq \emptyset} f_I \, \chi^I$
is sometimes termed the {\em soul}\/ of $f$.

Multiplication in the Grassmann algebra is of course $R$-bilinear, i.e.
\be
   \Biggl( \sum_{I\subseteq [n]} f_I \, \chi^I \Biggr) \,
   \Biggl( \sum_{J\subseteq [n]} g_J \, \chi^J \Biggr)
   \;=\;
   \sum_{I,J \subseteq [n]} f_I g_J \, \chi^I \chi^J
\ee
where
\be
   \chi^I \chi^J  \;=\;
   \cases{ \sigma(I,J) \, \chi^{I \cup J}  & if $I \cap J = \emptyset$ \cr
           \noalign{\vskip 2mm}
           0   & if $I \cap J \neq \emptyset$ \cr
         }
\ee
and $\sigma(I,J)$ is the sign of the permutation that rearranges
the sequence $IJ$ into increasing order when $I$ and $J$ are each
written in increasing order.

We define the {\em degree}\/ of a monomial $\chi^I$ in the obvious way,
namely, $\deg(\chi^I) = |I|$.
The Grassmann algebra $R[\ch]$ then possesses a natural $\NN$-grading
\be
   R[\ch] \;=\;  \bigoplus_{p=0}^n R[\ch]_p
\ee
where $R[\ch]_p$ is generated, as an $R$-module,
by the monomials of degree $p$.
A coarser grading is the $Z_2$-grading 
\be
   R[\ch] \;=\; R[\ch]_{\rm even} \oplus R[\ch]_{\rm odd}
\ee
where
\begin{eqnarray}
   R[\ch]_{\rm even} & \eqdef &  \bigoplus_{p \: {\rm even}}  R[\ch]_p  \\[2mm]
   R[\ch]_{\rm odd}  & \eqdef &  \bigoplus_{p \: {\rm odd}}   R[\ch]_p
\end{eqnarray}
Note that $R[\ch]_{\rm even}$ is a subalgebra of $R[\ch]$
(but $R[\ch]_{\rm odd}$ is not).
The {\em parity operator}\/ ${\rm P}$, defined by
\be
   {\rm P} \Biggl( \sum_{I\subseteq [n]} f_I \, \chi^I \Biggr) \;=\;
                   \sum_{I\subseteq [n]} (-1)^{|I|} \, f_I \, \chi^I
   \;,
 \label{def.parity}
\ee
is an involutive automorphism of $R[\ch]$:
it acts as the identity on $R[\ch]_{\rm even}$
and as minus the identity on $R[\ch]_{\rm odd}$.
A nonzero element $f$ that belongs to either $R[\ch]_{\rm even}$
or $R[\ch]_{\rm odd}$ is said to be {\em $Z_2$-homogeneous}\/,
and its {\em parity}\/ is $p(f)\eqdef 0$ in the first case
and $p(f)\eqdef 1$ in the second.
It is easy to see that $Z_2$-homogeneous elements $f,g$ satisfy the
commutation/anticommutation relations
\be
  f g \;=\; (-1)^{p(f)p(g)} g f  \;.
 \label{eq.commanticomm}
\ee
In other words, odd elements anticommute with other odd elements,
while even elements commute with homogeneous elements of both types.
One easily deduces the following two consequences:

\begin{proposition}
An even element $f \in R[\ch]_{\rm even}$
commutes with the entire Grassmann algebra.
(In particular, $R[\ch]_{\rm even}$ is a {\em commutative}\/ ring.)
\end{proposition}

\begin{proposition}
  \label{prop.odd}
An odd element $f \in R[\ch]_{\rm odd}$ is nilpotent of order 2,
i.e.\ $f^2 = 0$.\footnote{
   \proofof{Proposition~\ref{prop.odd}}
   Write $f = \sum\limits_{I\,{\rm odd}} f_I \, \chi^I$; then
   $$ f^2 \;=\; \sum_{I\,{\rm odd}} f_I^2 \, (\chi^I)^2  \,+\,
                \sum_{\begin{scarray}
                         I,J \,{\rm odd} \\
                         I \prec J
                      \end{scarray}}
                 f_I f_J \, (\chi^I \chi^J + \chi^J \chi^I)
   $$
   where $I \prec J$ denotes that the distinct sets $I$ and $J$
   are written in increasing lexicographic order.
   But $(\chi^I)^2 = 0$ and $\chi^I \chi^J + \chi^J \chi^I = 0$
   by \reff{anti}/\reff{anti2} and \reff{eq.commanticomm}.
   \qed
}
\end{proposition}

A further consequence of the relations \reff{anti}/\reff{anti2}
and the finiteness of the number $n$ of generators
is that every element of
\be
   R[\ch]_+  \;\eqdef\;  \bigoplus_{p=1}^n R[\ch]_p
\ee
(the set of elements with no term of degree $0$) is nilpotent:

\begin{proposition}
   \label{prop.puresoul}
A ``pure soul'' element $f \in R[\ch]_+$
is nilpotent of order at most $M = \lfloor n/2 \rfloor + 2$,
i.e.\ $f^{M} = 0$.
\end{proposition}

\proof
Write $f = f_0 + f_1$ with $f_0 \in R[\ch]_{\rm even} \cap R[\ch]_+$
and $f_1 \in R[\ch]_{\rm odd}$.
Then $f_0$ and $f_1$ commute [by \reff{eq.commanticomm}]
and $f_1^n = 0$ for $n \ge 2$ [by Proposition~\ref{prop.odd}],
so we have $f^k = f_0^k + k f_0^{k-1} f_1$ for all $k$.
And it is easy to see, using \reff{anti}/\reff{anti2},
that $f_0^k = 0$ for $k > n/2$.
\qed

If $I = \{i_1,\ldots,i_p\}$ is a subset of $[n]$,
the Grassmann algebra $R[\ch_I] \eqdef R[\ch_{i_1},\ldots,\ch_{i_p}]$
is naturally isomorphic to the subalgebra of
$R[\ch]=R[\ch_1,\ldots,\ch_n]$ generated by $\{\chi_i\}_{i \in I}$;
we shall identify these two algebras
and use the same notation $R[\ch_I]$ for both.
In particular, the degree-zero subalgebra $R[\ch]_0$
is identified with the coefficient ring $R$.

If $\Phi(x) = \sum_{k=0}^N c_k x^k$ is any polynomial in a single
indeterminate $x$ with coefficients in either $R$ or $\Z$,
we can of course apply it to any $f \in R[\chi]$ to obtain
$\Phi(f) = \sum_{k=0}^N c_k f^k \in R[\chi]$.
Moreover, if $\Phi(x) = \sum_{k=0}^\infty c_k x^k$ is a
formal power series with coefficients in $R$ or $\Z$,
we can apply it to any $f \in R[\chi]_+$ because $f$ is nilpotent
and the sum is therefore finite.

But more is true:  suppose we have a formal power series
$\Phi(x) = \sum_{k=0}^\infty (c_k/k!) x^k$ with coefficients in $R$ or $\Z$
(note the factorial denominators!).
If the coefficient ring $R$ contains the rationals as a subring
(as it is usually convenient to assume),
then of course $\Phi(f)$ is well-defined for any $f \in R[\chi]_+$.
But we claim that this expression has an unambiguous meaning
for $f \in R[\chi]_+$
{\em even if $R$ does not contain the rationals}\/.
Indeed, let $f = \sum_{I \neq \emptyset} f_I \, \chi^I$
and consider the expansion of
\be
   f^k  \;=\;
   \left( \sum_{I \neq \emptyset} f_I \, \chi^I \right) ^{\! k}  \;=\;
   \sum_{I_1,\ldots,I_k \neq \emptyset} f_{I_1} \cdots f_{I_k} \,
             \chi^{I_1} \cdots \chi^{I_k}
   \;.
\ee
Whenever two or more of the sets $I_1,\ldots,I_k$ are equal
(or indeed have any elements in common),
we have $\chi^{I_1} \cdots \chi^{I_k} = 0$
by \reff{anti}/\reff{anti2}.
When, by contrast, the sets $I_1,\ldots,I_k$ are all distinct,
then there are $k!$ terms with the same coefficient $f_{I_1} \cdots f_{I_k}$
corresponding to the $k!$ different permutations of the factors
$\chi^{I_1}, \ldots, \chi^{I_k}$,
and these terms are either all equal
(if at most one of the sets $I_\alpha$ is odd)
or else add to zero (if two or more of them are odd).
It follows that
\be
   f^k  \;=\;
   k! \! \sum_{\begin{scarray}
                  I_1 \prec\ldots\prec I_k \\
                  \hbox{\scriptsize at most one $|I_\alpha|$ odd}
               \end{scarray}}
   \!\! f_{I_1} \cdots f_{I_k} \, \chi^{I_1} \cdots \chi^{I_k}
\ee
where $I_1 \prec\ldots\prec I_k$ denotes that the distinct sets
$I_1,\ldots,I_k$ are written in increasing lexicographic order.
So we can define
\be
   \Phi(f)  \;=\;
   \sum_{k=0}^n
   c_k \!\!\!
   \sum_{\begin{scarray}
               I_1 \prec\ldots\prec I_k \\
               \hbox{\scriptsize at most one $|I_\alpha|$ odd}
         \end{scarray}}
   \!\! f_{I_1} \cdots f_{I_k} \, \chi^{I_1} \cdots \chi^{I_k}
 \label{def.Phi.f}
\ee
even if the coefficient ring $R$ does not contain the rationals.

The most important case is $\Phi = \exp$.
Note that $\exp(f+g) = \exp(f) \, \exp(g)$
whenever $f,g \in R[\chi]_+ \cap R[\chi]_{\rm even}$,
but not in general otherwise.
In this paper we will apply the exponential {\em only}\/ to even elements.

We also need to define one other type of composition.
Suppose first that $f \in R\< \ch_1,\ldots,\ch_n\>$
is a noncommutative polynomial in the letters $\ch_1,\ldots,\ch_n$
with coefficients in $R$,
and let $\xi_1,\ldots,\xi_n$ be elements of some ring $R'$
(not necessarily commutative) that contains $R$ within its center.
Then the composition $f(\xi_1,\ldots,\xi_n)$,
obtained by substituting each $\chi_i$ by the corresponding $\xi_i$,
is a well-defined element of $R'$;
furthermore, this composition satisfies the obvious laws
\begin{subeqnarray}
   (f+g)(\xi_1,\ldots,\xi_n)  & = &
      f(\xi_1,\ldots,\xi_n) \,+\, g(\xi_1,\ldots,\xi_n)  \\[1mm]
   (fg)(\xi_1,\ldots,\xi_n)  & = &
      f(\xi_1,\ldots,\xi_n) \, g(\xi_1,\ldots,\xi_n)
 \label{eq.composition_laws}
\end{subeqnarray}
Now suppose, instead, that $f$ belongs to the Grassmann algebra
$R[\ch_1,\ldots,\ch_n] =
 R\< \ch_1,\ldots,\ch_n\> / \{\ch_i \ch_j+ \ch_j \ch_i, \ch_i^2\}$
and that $\xi_1,\ldots,\xi_n$ are elements of $R'$ that satisfy
\begin{subeqnarray}
   \xi_i \xi_j+ \xi_j \xi_i & = &  0   \\[1mm]
   \xi_i^2  & = & 0
 \label{eq.relations.xi}
\end{subeqnarray}
for all $i,j$.
Then the composition $f(\xi_1,\ldots,\xi_n)$ is again
a well-defined element of $R'$,
because the relations \reff{eq.relations.xi}
guarantee that any representative of $f$ in $R\< \ch_1,\ldots,\ch_n\>$
will give, after substituting each $\chi_i$ by $\xi_i$,
the same element of $R'$;
moreover, the laws \reff{eq.composition_laws} continue to hold.
In particular, we can take $R'$ to be another Grassmann algebra over $R$
(which may or may not contain some of the $\chi_i$ as generators,
 and which may or may not contain additional generators)
and $\xi_1,\ldots,\xi_n$ to be arbitrary {\em odd}\/ elements
of this Grassmann algebra;
the relations \reff{eq.relations.xi}
hold by virtue of \reff{eq.commanticomm} and Proposition~\ref{prop.odd},
respectively.
We shall exploit this type of composition
in Propositions~\ref{prop.change} and \ref{prop.translation} below.

\subsection{Grassmann--Berezin (fermionic) integration}
   \label{app.grassmann.4}

Thus far we have simply been recalling standard facts about exterior algebra.
But now we go on to introduce a process called
{\em Grassmann--Berezin integration}\/,
which has become a standard tool of theoretical physicists
over the last 40 years but is still surprisingly little known
among mathematicians.
As we shall see, the term ``integration'' is a misnomer,
because the construction is purely algebraic/combinatorial.
But the term is nevertheless felicitous, because Grassmann--Berezin integration
behaves in many ways {\em analogously}\/ to ordinary integration
over $\R^n$ or $\C^n$, and this analogy is heuristically very fruitful.


We start by defining, for each $i \in [n]$,
the derivation $\partial_i  = \frac{\partial}{\partial\ch_i}$
(acting to the right) as the
$R$-linear map $\partial_i \colon\, R[\ch]\rightarrow R[\ch]$
defined by the following action on monomials $\ch_{i_1}\cdots\ch_{i_p}$
with $i_1<\ldots<i_p$:
\be
   \partial_i \, \ch_{i_1}\cdots\ch_{i_p}
   \;\eqdef\;
   \cases{ (-1)^{\al-1}
              \ch_{i_1}\cdots\ch_{i_{\al-1}}\ch_{i_{\al+1}}\cdots\ch_{i_p}
              &  if $i = i_\al$  \cr
           \noalign{\vskip 2mm}
	   0  &  if $i\notin \{i_1,\ldots,i_p\}$  \cr
         }
 \label{def.partiali}
\ee
(Of course, in the former instance there is a {\em unique}\/ index
 $\alpha \in [p]$ for which $i=i_\al$, so this definition is unambiguous.)
It is then easy to see that \reff{def.partiali} holds also
when the indices $i_1,\ldots,i_p$ are not necessarily ordered,
provided that they are all distinct.
Clearly $\partial_i$ is a map of degree $-1$,
i.e.\ $\partial_i \colon\, R[\ch]_p \rightarrow R[\ch]_{p-1}$;
moreover, it takes values in $R[\ch_{\{i\}^c}]$,
i.e.\ the subalgebra generated by $\{\chi_j\}_{j \neq i}$.
Furthermore, the maps $\partial_i$ satisfy
\begin{eqnarray}
  \partial_i^2  & = & 0   \label{eq.partial2}  \\[1mm]
\partial_i \partial_j + \partial_j \partial_i  & = & 0
   \label{eq.partial.anticommute}
\end{eqnarray}
as well as the modified Leibniz rule
\be
   \partial_i (fg)  \;=\;
   (\partial_i f) g \,+\, ({\rm P}f) (\partial_i g)
 \label{eq.leibniz}
\ee
where ${\rm P}$ is the parity operator \reff{def.parity}.

We now make a surprising definition:
{\em integration is the same as differentiation}\/.
That is, we define
\be
   \int \! d\chi_i \: f  \;=\;  \partial_i f   \;.
\ee
We always write the operator of integration {\em to the left}\/
of the integrand, just as we do for the (completely equivalent)
operator of differentiation.  To lighten the notation
we refrain from repeating the $\int$ sign in iterated integrals,
so that $\int \! d\chi_{i_1} \,\cdots\, d\chi_{i_k} \, f$ is a shorthand for
$\int \! d\chi_{i_1} \,\cdots\, \int \! d\chi_{i_k} \, f$ and hence we have
\be
   \int \! d\chi_{i_1} \,\cdots\, d\chi_{i_k} \: f
   \;=\;  \partial_{i_1} \cdots \partial_{i_k} f   \;.
\ee
Note by \reff{eq.partial.anticommute}
that changing the order of integration changes the sign:
\be
   \int \! d\chi_{i_{\sigma(1)}} \,\cdots\, d\chi_{i_{\sigma(k)}} \: f
   \;=\;
   \sgn(\sigma)
   \int \! d\chi_{i_1} \,\cdots\, d\chi_{i_k} \: f
\ee
for any permutation $\sigma \in \scrs_k$.
For instance, reversing the order of integrations gives
\be
   \int \! d\chi_{i_k} \,\cdots\, d\chi_{i_1} \: f
   \;=\;
   (-1)^{k(k-1)/2}
   \int \! d\chi_{i_1} \,\cdots\, d\chi_{i_k} \: f
   \;.
 \label{eq.reverse}
\ee

In particular, when we integrate $f = \sum_{I \subseteq [n]} f_I \, \chi^I$
with respect to {\em all}\/ the generators, we have
\be
   \int \! d\chi_n \,\cdots\, d\chi_1 \: f
   \;=\;
   f_{[n]}
   \;.
\ee
That is, integration with respect to all the generators
simply picks out the coefficient of the ``top'' monomial $I=[n]$,
{\em provided that we write these integrations in reverse order}\/.\footnote{
   Of course, since the integrations are performed from right to left,
   this actually corresponds to performing $\int \! d\chi_1$ first
   and $\int \! d\chi_n$ last.
}
We adopt the shorthand
\be
   \scrd \chi  \;\eqdef\;  d\chi_n \,\cdots\, d\chi_1
   \;,
 \label{def.scrd.chi}
\ee
and we sometimes write $\scrd_n(\chi)$ if we wish to stress
the number of generators.

An important special case arises when $n$ is even, say $n=2m$,
and the generators $\chi_1,\ldots,\chi_n$
are divided into two sets $\psi_1,\ldots,\psi_m$ and
$\psibar_1,\ldots,\psibar_m$,
where we think of each $\psi_i$ as paired with its corresponding $\psibar_i$.
In this case we adopt the shorthand notation
\be
   \scrd(\psi,\psibar)  \;\eqdef\;
   d\psi_1 \, d\psibar_1 \:\cdots\: d\psi_m \, d\psibar_m
 \label{def.scrd.psipsibar}
\ee
and we sometimes write $\scrd_m(\psi,\psibar)$ if we wish to stress
the number of pairs of generators.
Please note that since each pair $d\psi_i \, d\psibar_i$ is Grassmann-even,
we can also write
\be
   \scrd(\psi,\psibar)  \;\eqdef\;
   \prod_{i=1}^m d\psi_i \, d\psibar_i
 \label{def.scrd.psipsibar.bis}
\ee
where the terms in the product can be taken in any order.
Note also that
\be
   \scrd(\psi,\psibar)  \;=\; (-1)^{m(m-1)/2} \, \scrd\psi \: \scrd\psibar
   \;.
\ee
The notation $\bar{\phantom{\psi}}$
is intended to be suggestive of complex conjugation,
but we stress that it has nothing to do with complex numbers:
it merely denotes the extra combinatorial structure on the index set $[2m]$
that arises from the splitting of $[2m]$ into two sets of cardinality $m$
and the fixing of a bijection between the two sets.\footnote{
   We can view this as dividing the $2m$ individuals
   into $m$ males and $m$ females
   and then pairing those individuals into $m$ heterosexual couples.
}
In particular, the coefficient ring $R$ is still completely arbitrary.
The general case ($\chi$) and special case ($\psi,\psibar$)
are nevertheless known in the physics literature
as ``real fermions'' and ``complex fermions'', respectively.
We shall retain this terminology
but shall always put the adjectives ``real'' and ``complex''
in quotation marks in order to warn the reader that they are potentially
misleading.

Now let $A = (A_{ij})_{i,j=1}^n$ be an $n \times n$ matrix,
and define new Grassmann variables by $\xi_i = \sum_{j=1}^n A_{ij} \ch_j$.
We then have
\begin{eqnarray}
   \int \! d\ch_{n} \cdots d\ch_{1} \; \xi_1 \cdots \xi_n
   & = &
   \int \! d\ch_{n} \cdots d\ch_{1} \;
           \sum_{j_1=1}^n A_{1j_1} \ch_{j_1} \cdots
           \sum_{j_n=1}^n A_{nj_n} \ch_{j_n}
      \nonumber \\
   & = & 
   \sum_\si \ep(\si) \, A_{1\si(1)} \cdots A_{n\si(n)}
      \nonumber \\[1mm]
   & = & 
   \det A  \;.
\end{eqnarray}
This result can be reformulated as follows:

\begin{proposition}[Linear change of variables in Grassmann--Berezin integration]
   \label{prop.change}
Let $f \in R[\xi_1,\ldots, \xi_n]$,
and define $F \in R[\ch_1,\ldots, \ch_n]$ by
\be
   F(\ch_1,\ldots, \ch_n) \;=\;
   f(\xi_1,\ldots, \xi_n) \Bigg| _{\xi_i = \sum_{j=1}^n A_{ij} \ch_j}
\ee
[note that this substitution is well-defined by virtue of the
 discussion in the last paragraph of Section~\ref{app.grassmann.3}].
Then
\be
   \int \! d\ch_{n} \cdots d\ch_{1} \; F 
   \;=\;
   (\det A) \, \int \! d\xi_{n} \cdots d\xi_{1} \; f
   \;.
\ee
\end{proposition}

\noindent
Please note that this is the {\em reverse}\/ of the
change-of-variables formula for ordinary multivariate integrals,
i.e.\ in an ordinary integral the factor $\det(A)$
would appear on the other side.

\bigskip

We will also need:

\begin{proposition}[Fubini theorem for Grassmann--Berezin integration]
   \label{prop.fubini}
\quad\break
Let $I=\{i_1,\ldots,i_p\}$ with $i_1<\cdots<i_p$ be a subset of $[n]$,
and let $I^{\rm c}=\{j_1,\ldots,j_{n-p}\}$ with  $j_1<\cdots<j_{n-p}$.
Then for any elements $f\in R[\ch_I]$ and $g\in R[\ch_{I^{\rm c}}]$
we have
\be
  \int \scrd\ch_I \, \scrd\ch_{I^{\rm c}}\ fg
   \;=\;
  {(-1)}^{p(n-p)}
\lp
\int \scrd\ch_I\ f
\rp
\lp
\int \scrd\ch_{I^{\rm c}}\ g
\rp
  \label{eq.prop.fubini}
\ee
where $\scrd\ch_I$ (resp. $\scrd\ch_{I^{\rm c}}$)
is shorthand for $d\ch_{i_p}\cdots d\ch_{i_1}$
(resp. $d\ch_{j_{n-p}} \cdots d\ch_{j_1}$). 
\end{proposition}

\proof
Expanding $f$ and $g$ in monomials, we see that the only terms
contributing to the integrals on either side are the ``top'' monomials
$\chi^I$ and $\chi^{I^c}$ in $f$ and $g$, respectively;
so we can assume without loss of generality that
$f = \chi^I$ and $g = \chi^{I^c}$.
Now use the fact that integration is the same as differentiation,
and successively apply the operators $\partial_i$ for $i \in I^c$
to the product $fg$ using the Leibniz rule \reff{eq.leibniz}.
The differentiations hit only $g$, and we have ${\rm P}f = (-1)^p f$.
It follows that
\be
  \int \scrd\ch_{I^{\rm c}}\ fg
   \;=\;
  {(-1)}^{p(n-p)}
  f
\lp
\int \scrd\ch_{I^{\rm c}}\ g
\rp
   \;.
\ee
The result then follows by integrating both sides with $\scrd\ch_I$.
\qed

Let us remark that the same formula \reff{eq.prop.fubini}
would hold for {\em any}\/ choice of the orderings
in defining $\scrd\ch_I$ and $\scrd\ch_{I^{\rm c}}$,
provided only that we use the {\em same}\/ orderings
on both sides of the equation.
We have chosen to write
$\scrd\ch_I = d\ch_{i_p}\cdots d\ch_{i_1}$
for compatibility with our convention \reff{def.scrd.chi}
that $\scrd \chi =  d\chi_n \,\cdots\, d\chi_1$.

\bigskip

Finally, the following proposition shows that
a Grassmann--Berezin integral over $d\chi_i$
is invariant under translation by an arbitrary {\em odd}\/
element of the Grassmann algebra
{\em that does not involve the variable $\chi_i$}\/.
More generally, one can consider integration over a set
$I \subseteq [n]$ of generators:

\begin{proposition}[Invariance under translation]
   \label{prop.translation}
Let $I=\{i_1,\ldots,i_p\} \subseteq [n]$
and let $\xi_1,\ldots,\xi_n \in R[\chi_{I^c}]_{\rm odd}$
satisfy $\xi_j = 0$ whenever $j \notin I$.
Then
\be
   \int \! d\chi_{i_p} \,\cdots\, d\chi_{i_1} \: f(\chi+\xi)
   \;=\;
   \int \! d\chi_{i_p} \,\cdots\, d\chi_{i_1} \: f(\chi)
 \label{eq.prop.translation}
\ee
where $f(\chi+\xi)$ denotes the substitution defined
at the end of Section~\ref{app.grassmann.3}.
[Recall that the oddness of the $\xi_i$ is required
 for this substitution to make sense.]
\end{proposition}

\proof
The formula~\reff{eq.prop.translation} can be rewritten as
\be
   \partial_{i_p} \,\cdots\, \partial_{i_1} \: f(\chi+\xi)
   \;=\;
   \partial_{i_p} \,\cdots\, \partial_{i_1} \: f(\chi)
   \;.
\ee
To prove this relation, it suffices to consider the cases
in which $f(\chi) = \chi^J$ for some $J \subseteq [n]$.
Now for any $i \in I$ and $j \in [n]$ we have
\be
   \partial_i (\chi_j + \xi_j)  \;=\;  \partial_i  \chi_j  \;=\; \delta_{ij}
\ee
because  $\xi_j \in R[\chi_{I^c}]$.
Using this relation together with the Leibiniz rule,
we see that $\partial_i (\chi + \xi)^J$ equals the
same object in which $\xi_i$ has been replaced by zero.
Doing this successively for $\partial_{i_1},\ldots,\partial_{i_p}$,
we set $\xi_i$ to zero for all $i \in I$.
But by hypothesis these are the only nonzero $\xi_i$.
\qed

\subsection{Fermionic Gaussian integration}  \label{app.grassmann.5}

We are now ready to state the fundamental formulae (``Wick's theorem'')
for fermionic Gaussian integration, which are analogues of
Theorems~\ref{thm.wick.realbosons} and \ref{thm.wick.complexbosons}
for bosonic Gaussian integration.
The main difference is that fermionic integration is a purely
algebraic/combinatorial construction,
so that it works over an arbitrary (commutative) coefficient ring $R$
and does not require any positive-definiteness condition on the matrix $A$.

We begin with the formula for a pure Gaussian integral,
i.e.\ the integral of the exponential of a quadratic form.
So let $R$ be a commutative ring with identity element,
let $\chi_1,\ldots,\chi_n$ be the generators of a Grassmann algebra,
and let $A = (a_{ij})_{i,j=1}^{n}$ be an $n \times n$
{\em antisymmetric}\/ matrix (i.e.\ $a_{ij} = -a_{ji}$ and $a_{ii} = 0$)
with entries in $R$.\footnote{
   See footnote \ref{note_alternating} above
   concerning our use of the term ``antisymmetric''.
}
We use the notation
\be
   \smfrac{1}{2} \,\ch^{\rm T} A\ch
   \;\eqdef\;
   \sum_{1 \le i< j \le n} \ch_i a_{ij}\ch_j
\ee
and observe that the right-hand side makes sense
even if the coefficient ring $R$ does not contain an element $\smfrac{1}{2}$.
We then have the following formula,
which shows that a Gaussian fermionic integral equals a pfaffian:


\begin{proposition}[Gaussian integral for ``real'' fermions]
   \label{prop.gaussian.realfermions}
Let $A$ be an $n\times n$ antisymmetric matrix with coefficients in $R$.
Then
\be
  \int \! d\ch_{n} \,\cdots\, d\ch_{1} \: e^{\frac{1}{2}\ch^{\rm T} A\ch}
  \;=\;
  \int \! d\ch_{1} \,\cdots\, d\ch_{n} \: e^{-\frac{1}{2}\ch^{\rm T} A\ch}
  \;=\;
  \cases{ \pf A  & if $n$ is even \cr
          \noalign{\vskip 2mm}
          0      & if $n$ is odd \cr
        }
 \label{eq.gaussian.pfaffian}
\ee
\end{proposition}

\proof
We expand the exponential using \reff{def.Phi.f} with $\Phi = \exp$
and then integrate $\int \! d\ch_{n} \,\cdots\, d\ch_{1}$.
When $n$ is odd the integral vanishes,
and when $n$ is even (say, $n=2m$)
the only contribution comes from $k=m$ in \reff{def.Phi.f}, yielding
\begin{subeqnarray}
   & &
   \int \! d\ch_{n} \,\cdots\, d\ch_{1} \:
                 e^{\frac{1}{2}\ch^{\rm T} A\ch}
         \nonumber \\[2mm]
   & & \qquad =\;
   \int \! d\ch_{n} \,\cdots\, d\ch_{1}  \!\!\!
       \sum_{\begin{scarray}
                 I_1 \prec \dots \prec I_m  \\
                 I_k = \{i_k,j_k\} \:\hbox{\scriptsize\rm with } i_k < j_k
             \end{scarray}}
       \!\!\!\!\!\!
       a_{i_1j_1} \cdots a_{i_m j_m} \, \chi^{I_1} \cdots \chi^{I_m}
           \qquad \qquad  \\[2mm]
   & & \qquad =\;
  \sum_{\begin{scarray}
               \si \colon\, \si(2i-1)<\si(2i) \:\hbox{\scriptsize\rm and} \\
               \si(1)<\si(3)< \cdots <\si(2m-1)
              \end{scarray} }
       \!\!\!
\ep(\si)\, a_{\sigma(1)\sigma(2)}\cdots a_{\sigma(2m-1)\sigma(2m)}   \\[2mm]
   & & \qquad =\;
   \pf A  \;.
\end{subeqnarray}
(Note that this holds whether or not the coefficient ring $R$
 contains the rational numbers.)

For the variant in which $A$ is replaced by $-A$
and the order of integration is reversed, it suffices to observe that
\be
   \pf(-A)  \;=\; (-1)^{n/2} \pf(A)
 \label{eq.pfminusA}
\ee
[this is an immediate consequence of the definition \reff{def.pfA}]
and use \reff{eq.reverse} to get
\be
  \int \! d\ch_{1} \cdots d\ch_{n} \: e^{-\frac{1}{2}\ch^{\rm T} A\ch}
  \;=\;
  (-1)^{n(n-1)/2} \pf(-A)
  \;=\;
  (-1)^{n^2/2} \pf(A)
  \;=\;
  \pf(A)
  \quad
\ee
because $n$ is even.
\qed

{\bf Remark.} Equation~\reff{eq.gaussian.pfaffian} shows that
there exist two equally natural conventions for Gaussian integrals
with ``real'' fermions:
either we write our quadratic forms as $e^{\frac{1}{2}\ch^{\rm T} A\ch}$
and our integrals as $\int \! d\ch_{n} \,\cdots\, d\ch_{1}$,
or alternatively
we write our quadratic forms as $e^{-\frac{1}{2}\ch^{\rm T} A\ch}$
and our integrals as $\int \! d\ch_{1} \,\cdots\, d\ch_{n}$.
In this paper we have adopted the first convention [cf.\ \reff{def.scrd.chi}].

\bigskip

In the special case in which $n=2m$
and the generators $\chi_1,\ldots,\chi_n$
are divided into two sets $\psi_1,\ldots,\psi_m$ and
$\psibar_1,\ldots,\psibar_m$, we have

\begin{proposition}[Gaussian integral for ``complex'' fermions]
   \label{prop.gaussian.complexfermions}
$\:$
Let $A$ be an
\linebreak
$m\times m$ matrix with coefficients in $R$.
Then
\be
  \int \! d\psi_1 \, d\psibar_1 \:\cdots\: d\psi_m \, d\psibar_m
          \: e^{\psibar^{\rm T} A\psi}
  \;=\;
  \det A
   \;.
 \label{gaussint.fermionic.complex}
\ee
\end{proposition}

\noindent
Please note that here $A$ is an {\em arbitrary}\/ matrix;
no condition of symmetry or antisymmetry need be imposed on it.

\proof
By the change of variables $(\psi',\psibar') = (A\psi, \psibar)$,
we have from Proposition~\ref{prop.change}
\begin{subeqnarray}
  \int \! \scrd(\psi,\psibar) \: e^{\psibar^{\rm T} A\psi}
  & = &
  (\det A)
  \int \! \scrd(\psi,\psibar) \: e^{\psibar^{\rm T} \psi}
       \\[1mm]
  & = &
  (\det A)
  \int \! \scrd(\psi,\psibar) \: \prod_{i=1}^m (1 + \psibar_i \psi_i)
     \;.
\end{subeqnarray}
When we expand the product,
only the term $\psibar_1 \psi_1 \cdots \psibar_m \psi_m$
has a nonzero integral, and the integral of this term is 1.
\qed


Here is an alternate proof, which treats ``complex'' fermions
as a special case of ``real'' fermions and invokes
Proposition~\ref{prop.gaussian.realfermions}:

\alternateproofof{Proposition~\ref{prop.gaussian.complexfermions}}
Let us write $(\chi_1,\ldots,\chi_{2m}) =$  \hfill\break
  $(\psi_1,\ldots,\psi_m,\psibar_1,\ldots,\psibar_m)$.
Then by \reff{def.scrd.psipsibar.bis} we have
\begin{subeqnarray}
\scrd(\psi,\psibar)
   & = &
   (-1)^m d\psibar_m \, d\psi_m \:\cdots\: d\psibar_1 \, d\psi_1  \\[1mm]
   & = &
   (-1)^{m (m+1)/2} \; d\psibar_m\,\cdots\, d\psibar_1 \:
                       d\psi_m\,\cdots\ d\psi_1    \\[1mm]
   & = & (-1)^{m (m+1)/2} \; d\chi_{2m}\,\cdots\, d\chi_1
\end{subeqnarray}
and
\be
\bar\psi^{\rm T} A\psi  \;=\;  \frac{1}{2} \, \chi^{\rm T} K \chi
   \;,
\ee
where $K$ is the $2m \times 2m$ matrix
\be
  K
  \;=\;
  \left(\! \begin{array}{cc}0 & - A^{\rm T} \\A & 0\end{array}\ \!\!\right)
  \;=\;
  \left(\! \begin{array}{cc}0 & I \\A & 0\end{array}\ \!\!\right)
  \left(\! \begin{array}{cc}0 & I \\-I & 0\end{array}\ \!\!\right)
  \left(\! \begin{array}{cc}0 &  A^{\rm T} \\I & 0\end{array}\ \!\!\right)
  \;.
\ee
Then by Lemma~\ref{lemma.properties.pfaffians}(c) we have
\be
  \pf K 
  \;=\;
  \det\!\left(\! \begin{array}{cc}0 & I \\A & 0\end{array}\ \!\right)
  \,\cdot\,
  \pf\!\left(\! \begin{array}{cc}0 & I \\-I & 0\end{array}\ \!\right)
  \;=\;
  (-1)^m (\det A) \:\cdot\: (-1)^{m (m-1)/2}
  \;,
\ee
so that
\be
 \int \! \scrd(\psi,\psibar) \: e^{\psibar^{\rm T} A\psi} 
 \;=\;
 (-1)^{m (m+1)/2} \pf K 
 \;=\;
 \det A   \;.
\ee
\qed

\bigskip


Let us make one further observation concerning Gaussian integrals
(both ``real'' and ``complex''),
which will be useful in proving certain aspects of Wick's theorem.
If $I = (i_1,\ldots,i_k)$ is an arbitrary sequence of indices in $[n]$
--- not necessarily distinct or ordered ---
then the ``real'' Gaussian integral \reff{eq.gaussian.pfaffian}
can be generalized to
\be
  \int \! d\ch_{i_k} \,\cdots\, d\ch_{i_1} \:
             e^{\frac{1}{2}\ch_I^{\rm T} A_{II}\ch_I}
  \;=\;
  \cases{ \pf A_{II}  & if $k$ is even \cr
          \noalign{\vskip 2mm}
          0      & if $k$ is odd \cr
        }
 \label{eq.gaussian.pfaffian.EXTENDED}
\ee
To prove this, it suffices to observe, first of all,
that if the indices $i_1,\ldots,i_k$ are distinct,
then this formula is merely \reff{eq.gaussian.pfaffian}
after a relabeling of indices:
the point is that by our definition \reff{eq.AIJ.sequence},
the variables arise in the {\em same}\/ order in the
integration measure and in the matrix $A_{II}$.
On the other hand, if the sequence $i_1,\ldots,i_k$ contains
any repeated index, then the left-hand side vanishes because
$d\ch_i \, d\ch_i = 0$ by \reff{eq.partial2},
while $\pf A_{II} = 0$ in such a case because the pfaffian
changes sign under simultaneous permutations of rows and columns
[this is a special case of Lemma~\ref{lemma.properties.pfaffians}(c)].
In a similar way, the ``complex'' Gaussian integral
\reff{gaussint.fermionic.complex} can be generalized to
\be
  \int \! d\psi_{i_1} \, d\psibar_{j_1} \:\cdots\: d\psi_{i_k} \, d\psibar_{j_k}
          \: e^{\psibar^{\rm T}_J A_{JI} \psi_I}
  \;=\;
  \det A_{JI}
 \label{gaussint.fermionic.complex.EXTENDED}
\ee
for arbitrary sequences $I = (i_1,\ldots,i_k)$ and $J = (j_1,\ldots,j_k)$.

\bigskip

We are now ready to state the full Wick's theorem for fermions.
Since the ``source'' version of Wick's theorem for fermions
will involve a {\em fermionic}\/ source $\lambda$
[cf.\ \reff{eq.wick.realfermions.1}],
this means that we will be working (at least when discussing this equation)
in an extended Grassmann algebra $R[\chi,\theta]$ with generators
$\chi_1,\ldots,\chi_n$ and $\theta_1,\ldots,\theta_N$ for some $N \ge 1$,
and the sources $\lambda_i$ will belong to the odd part
of the Grassmann {\em sub}\/algebra $R[\theta]$.

We assume without further ado that $n$ is {\em even}\/, i.e.\ $n=2m$.
Recall also from \reff{eq.AIJ.sequence} the notation $A_{IJ}$
for arbitrary {\em sequences}\/ of indices $I$ and $J$.

\begin{theorem}[Wick's theorem for ``real'' fermions]
 \label{thm.wick.realfermions}
Let $R$ be a commutative ring with identity element,
and let $A = (a_{ij})_{i,j=1}^{2m}$ be a $2m \times 2m$ antisymmetric matrix
with elements in $R$.
Then:
\begin{itemize}
   \item[(a)]
      If the matrix $A$ is invertible, we have
\begin{equation}
   \int\! \scrd\chi \,
   \exp\!\left( \smhalf \chi^{\rm T} A \chi \,+\, \lambda^{\rm T} \chi
         \right)
   \;=\;
   (\pf A) \,
   \exp\!\left( \smhalf  \lambda^{\rm T} A^{-1} \lambda \right)
 \label{eq.wick.realfermions.1}
\end{equation}
     whenever $\lambda = (\lambda_i)_{i=1}^{2m}$
     are odd elements of the Grassmann algebra that do not involve $\chi$,
     i.e.\ $\lambda_i \in R[\theta]_{\rm odd}$.
   \item[(b)]  For any subset $I = \{i_1,\ldots,i_r\} \subseteq [2m]$
with $i_1 < \ldots < i_r$, we have
\begin{equation}
   \int\! \scrd\chi \: \chi_{i_1} \cdots \chi_{i_r} \,
   \exp\!\left( \smhalf \chi^{\rm T} A \chi \right)
   \;=\;
   \cases{ 0                               & if $r$ is odd \cr
           \noalign{\vskip 2mm}
           \epsilon(I) \, \pf A_{I^c I^c}       & if $r$ is even \cr
         }
 \label{eq.wick.realfermions.2a}
\end{equation}
   \item[(c)]  For any sequence of indices $I = (i_1,\ldots,i_r)$ in $[2m]$,
if the matrix $A$ is invertible we have
\begin{equation}
   \int\! \scrd\chi \: \chi_{i_1} \cdots \chi_{i_r} \,
   \exp\!\left( \smhalf \chi^{\rm T} A \chi \right)
   \;=\;
   \cases{ 0                               & if $r$ is odd \cr
           \noalign{\vskip 2mm}
           (\pf A) \, \pf((A^{-{\rm T}})_{II})  & if $r$ is even \cr
         }
 \label{eq.wick.realfermions.2b}
\end{equation}
   \item[(d)]
More generally, for any $r \times 2m$ matrix $C$ with entries in $R$, we have
\begin{equation}
   \int\! \scrd\chi \, (C\chi)_1 \,\cdots\, (C\chi)_r \:
   \exp\!\left( \smhalf \chi^{\rm T} A \chi \right)
   \;=\;
   \cases{ 0                               & if $r$ is odd \cr
           \noalign{\vskip 2mm}
           \sum\limits_{|I| = r}  (\det C_{\star I}) \, \epsilon(I) \,
                 \pf(A_{I^c I^c})  & if $r$ is even \cr
         }
 \label{eq.wick.realfermions.3a}
\end{equation}
where $C_{\star I}$ denotes the submatrix of $C$ with columns in $I$;
and if the matrix $A$ is invertible,
\begin{equation}
   \int\! \scrd\chi \, (C\chi)_1 \,\cdots\, (C\chi)_r \:
   \exp\!\left( \smhalf \chi^{\rm T} A \chi \right)
   \;=\;
   \cases{ 0                               & if $r$ is odd \cr
           \noalign{\vskip 2mm}
           (\pf A) \, \pf(C A^{-{\rm T}} C^{\rm T})  & if $r$ is even \cr
         }
 \label{eq.wick.realfermions.3b}
\end{equation}
\end{itemize}
\end{theorem} 

Note that \reff{eq.wick.realfermions.2a} and \reff{eq.wick.realfermions.2b}
are equivalent by virtue of Jacobi's identity for pfaffians
[cf.\ \reff{eq.app.pfaffian.jacobi}],
and that \reff{eq.wick.realfermions.3a} and \reff{eq.wick.realfermions.3b}
are equivalent by virtue of Jacobi's identity
and the minor summation formula for pfaffians
[cf.\ \reff{eq.app.pfaffian.minorsummation}].
However, we shall give {\em independent}\/ proofs of all these formulae;
as a consequence, our argument provides {\em ab initio}\/ proofs
of these pfaffian identities by means of Grassmann--Berezin integration.


\proofof{Theorem~\ref{thm.wick.realfermions}}
To prove (a), we perform the translation $\chi'= \chi + A^{-1} \lambda$
and use Proposition~\ref{prop.translation}.
We have
\begin{eqnarray}
   \smhalf \chi^{\prime \rm T} A \chi' \,+\, \lambda^{\rm T} \chi'
   & = &
   \smhalf (\chi^{\rm T} + \lambda^{\rm T} A^{-\rm T}) A (\chi + A^{-1} \lambda)
      \,+\, \lambda^{\rm T} (\chi + A^{-1} \lambda)
    \nonumber \\[1mm]
   & = &
   \smhalf \chi^{\rm T} A \chi \,+\,  \smhalf  \lambda^{\rm T} A^{-1} \lambda
\end{eqnarray}
since $A^{-\rm T} = -A^{-1}$
and $\chi^{\rm T} \lambda = -\lambda^{\rm T} \chi$.
Therefore \reff{eq.wick.realfermions.1}
is an immediate consequence of Proposition~\ref{prop.gaussian.realfermions}.

For (b)--(d), let us first remark that if $r$ is odd,
then the integral must vanish, because the expansion of the exponential
can only give an even number of factors of $\chi$,
so the ``top'' monomial in $\chi$ cannot be generated
(recall that the total number of generators is even, i.e.\ $n=2m$).

To prove (b), we compute
\begin{eqnarray}
   \int\! \scrd\chi \: \chi^I \:
   \exp\!\left( \smhalf \chi^{\rm T} A \chi \right)
   & = &
   \int\! \scrd\chi \: \chi^I \:
   \exp\!\left( \smhalf \chi^{\rm T}_{I^c} A_{I^c I^c} \chi_{I^c} \right)
        \nonumber \\[1mm]
   & = &
   \epsilon(I)
   \int\! \scrd\chi_I \, \scrd\chi_{I^c} \: \chi^I \:
   \exp\!\left( \smhalf \chi^{\rm T}_{I^c} A_{I^c I^c} \chi_{I^c} \right)
        \nonumber \\[1mm]
   & = &
   \epsilon(I) \, (-1)^{r(2m-r)}
   \left( \int \! \scrd\chi_I \: \chi^I \right)
   \left( \int \! \scrd\chi_{I^c} \:
   \exp\!\left( \smhalf \chi^{\rm T}_{I^c} A_{I^c I^c} \chi_{I^c} \right)
       \right)
        \nonumber \\[1mm]
   & = &
   \epsilon(I) \, \pf(A_{I^c I^c})
      \;.
\end{eqnarray}
The first equality holds because any factor $\xi_i$ with $i \in I$
arising from the expansion of the exponential would be annihilated
by the prefactor $\chi^I$;
the second equality is simply a reordering of the integration variables;
the third equality is Fubini's theorem (Proposition~\ref{prop.fubini});
and the last equality is simply the evaluation of the two integrals
(using Proposition~\ref{prop.gaussian.realfermions} for the second)
together with the fact that $r$ is even.

(c) We specialize result (a) to the case where the ``sources'' $\lambda_i$
are generators of an extended Grassmann algebra,
then differentiate (or equivalently integrate!)\ 
with respect to the sources $\lambda_{i_1}, \ldots, \lambda_{i_r}$,
and finally replace all the $\lambda_i$ by zero\footnote{
   Setting $\lambda \to 0$ can be interpreted
   as extracting the monomials that do not involve $\lambda$,
   or alternatively as a special case of the substitution discussed
   at the end of Section~\ref{app.grassmann.3} (since $0$ is odd).
}:
\begin{eqnarray}
   \int\! \scrd\chi \: \chi_{i_1} \cdots \chi_{i_r} \,
   \exp\!\left( \smhalf \chi^{\rm T} A \chi \right)
   & = &
   (\pf A) \:
   \left. {\partial \over \partial \lambda_{i_1}} \,\cdots\,
             {\partial \over \partial \lambda_{i_r}} \:
   \exp\!\left( \smhalf  \lambda^{\rm T} A^{-1} \lambda \right)
   \right|_{\lambda=0}
        \nonumber \\[1mm]
   & = &
   (\pf A) \:
   \int\! d \lambda_{i_1} \,\cdots\,d \lambda_{i_r} \:
   \exp\!\left( \smhalf  \lambda^{\rm T}_I (A^{-1})_{II} \lambda_I \right)
        \nonumber \\[1mm]
   & = &
   (\pf A) \; (-1)^{r(r-1)/2}
   \int\! \scrd\lambda_I \:
   \exp\!\left( \smhalf  \lambda^{\rm T}_I (A^{-1})_{II} \lambda_I \right)
        \nonumber \\[1mm]
   & = &
   (\pf A) \; (-1)^{r(r-1)/2} \pf((A^{-1})_{II})
        \nonumber \\[1mm]
   & = &
   (\pf A) \; (-1)^{r(r-1)/2} \: (-1)^{r/2} \: \pf((-A^{-1})_{II})
        \nonumber \\[1mm]
   & = &
   (\pf A) \, \pf((A^{-\rm T})_{II})
 \label{eq.wick.realfermions.proofc}
\end{eqnarray}
Here the third equality involved reordering the integration variables
from increasing to decreasing order;
the fourth equality performed the Gaussian integral
using \reff{eq.gaussian.pfaffian.EXTENDED};
the fifth equality used \reff{eq.pfminusA};
and the final equality used $-A^{-1} = A^{-\rm T}$
and $(-1)^{r^2/2} = 1$ (which holds since $r$ is even).

(d)  We have
\begin{eqnarray}
   &  &
   \int\! \scrd\chi \: (C\chi)_1 \cdots (C\chi)_r \:
   \exp\!\left( \smhalf \chi^{\rm T} A \chi \right)
            \nonumber \\
   &  &
   \qquad \;=\;
   \sum_{i_1,\ldots,i_r}  C_{1 i_1} \cdots C_{r i_r} \,
   \int\! \scrd\chi \: \chi_{i_1} \cdots \chi_{i_r} \:
   \exp\!\left( \smhalf \chi^{\rm T} A \chi \right)
      \;, \qquad
\end{eqnarray}
but the only nonvanishing contributions in the sum
come when $i_1,\ldots,i_r$ are all distinct;
so we can first require $i_1 < \ldots < i_r$
and then sum over permutations, yielding
\be
   \sum\limits_{|I|=r} \, \sum_{\sigma \in \scrs_r}
       C_{1 i_{\sigma(1)}} \cdots C_{r i_{\sigma(r)}}
   \,
  \sgn(\sigma)
   \int\! \scrd\chi \; \chi^I \:
   \exp\!\left( \smhalf \chi^{\rm T} A \chi \right)
   \;=\;
   \sum\limits_{|I| = r}  (\det C_{\star I}) \, \epsilon(I) \, \pf(A_{I^c I^c})
\ee
by part (b).
If $A$ is invertible, we can write
\begin{eqnarray}
   \int\! \scrd\chi \, (C\chi)_1 \,\cdots\, (C\chi)_r \:
   \exp\!\left( \smhalf \chi^{\rm T} A \chi \right)
   & = &
   {\partial \over \partial \lambda_{i_1}} \cdots
        {\partial \over \partial \lambda_{i_r}} \:
   \left.
   \int\! \scrd\chi \:
   \exp\!\left( \smhalf \chi^{\rm T} A \chi  \,+\, \lambda^{\rm T} C \chi\right)
   \right|_{\lambda=0}
          \nonumber \\[1mm]
   & = &
   {\partial \over \partial \lambda_{i_1}} \,\cdots\,
        {\partial \over \partial \lambda_{i_r}} \,
   \left.
   (\pf A) \, \exp\!\left( \smhalf \lambda^{\rm T} C A^{-1} C^{\rm T} \lambda
                    \right)
   \right|_{\lambda=0}
            \nonumber \\
\end{eqnarray}
by using (a).
Then by the same reasoning as in \reff{eq.wick.realfermions.proofc}
we can see that this equals $(\pf A) \, \pf(C A^{-{\rm T}} C^{\rm T})$.
\qed

Next we state Wick's theorem for ``complex'' fermions.
Once again, when discussing the ``source'' version of this theorem
[cf.\ \reff{eq.wick.complexfermions.1}],
we will work in an extended Grassmann algebra $R[\psi,\psibar,\theta]$
in which the sources $\lambda_i,\lambdabar_i$ belong to the odd part
of the Grassmann subalgebra $R[\theta]$.

\begin{theorem}[Wick's theorem for ``complex'' fermions]
 \label{thm.wick.complexfermions}
Let $R$ be a commutative ring with identity element,
and let $A = (a_{ij})_{i,j=1}^{n}$ be a $n \times n$ matrix
with elements in $R$.
Then:
\begin{itemize}
   \item[(a)]
      If the matrix $A$ is invertible, we have
\begin{equation}
   \int\! \scrd(\psi,\psibar) \,
   \exp\!\left( \psibar^{\rm T} A \psi \,+\, \lambdabar^{\rm T} \psi
                                       \,+\, \psibar^{\rm T} \lambda
         \right)
   \;=\;
   (\det A) \,
   \exp\!\left( - \lambdabar^{\rm T} A^{-1} \lambda \right)
 \label{eq.wick.complexfermions.1}
\end{equation}
     whenever $\lambda = (\lambda_i)_{i=1}^{n}$
     and $\lambdabar = (\lambdabar_i)_{i=1}^{n}$
     are odd elements of the Grassmann algebra that do not involve
     $\psi$ and $\psibar$,
     i.e.\ $\lambda_i,\lambdabar_i \in R[\theta]_{\rm odd}$.
   \item[(b)]  For any subsets $I = \{i_1,\ldots,i_r\}$
     and $J = \{j_1,\ldots,j_r\}$ of $[n]$
     having the same cardinality $r$,
     with $i_1 < \ldots < i_r$ and $j_1 < \ldots < j_r$, we have
\begin{equation}
   \int\! \scrd(\psi,\psibar)  \:
      \Biggl( \prod_{\alpha=1}^r \psibar_{i_\alpha} \psi_{j_\alpha} \Biggr)  \:
   \exp\!\left( \psibar^{\rm T} A \psi \right)
   \;=\;
   \epsilon(I,J) \, (\det A_{I^c J^c})
       \;.
 \label{eq.wick.complexfermions.2a}
\end{equation}
   [If there is an unequal number of factors $\psi$ and $\psibar$,
    then the integral is zero.]
   \item[(c)]  For any sequences of indices $I = (i_1,\ldots,i_r)$
     and $J = (j_1,\ldots,j_r)$ in $[n]$ of the same length $r$,
     if the matrix $A$ is invertible we have
\begin{equation}
   \int\! \scrd(\psi,\psibar) \:
      \Biggl( \prod_{\alpha=1}^r \psibar_{i_\alpha} \psi_{j_\alpha} \Biggr)  \:
   \exp\!\left( \psibar^{\rm T} A \psi \right)
   \;=\;
   (\det A) \, \det((A^{-{\rm T}})_{IJ})
      \;.
 \label{eq.wick.complexfermions.2b}
\end{equation}
   [Again, if there is an unequal number of factors $\psi$ and $\psibar$,
    then the integral is zero.]
   \item[(d)]
More generally, for any $r \times n$ matrix $B$
     and $n \times r$ matrix $C$ with entries in $R$, we have
\begin{eqnarray}
   & &
   \int\! \scrd(\psi,\psibar) \,
      \Biggl( \prod_{\alpha=1}^r \left(\bar\psi C\right)_\alpha
                                 \left(B \psi \right)_\alpha
      \Biggr) \:
   \exp\!\left( \psibar^{\rm T} A \psi \right)
          \nonumber \\[1mm]
   & &
   \qquad\qquad
   \;=\;
   \sum\limits_{|I|=|J|=r} \epsilon(I,J) \, (\det B_{\star J})
         (\det A_{I^c J^c}) (\det C_{I \star})
   \;,
   \qquad\qquad
 \label{eq.wick.complexfermions.3a}
\end{eqnarray}
and if the matrix $A$ is invertible,
\begin{equation}
   \int\! \scrd(\psi,\psibar) \,
      \Biggl( \prod_{\alpha=1}^r \left(\bar\psi C\right)_\alpha
                                 \left(B \psi \right)_\alpha
      \Biggr) \:
   \exp\!\left( \psibar^{\rm T} A \psi \right)
   \;=\;
   (\det A) \, \det(B A^{-1} C)
      \;.
 \label{eq.wick.complexfermions.3b}
\end{equation}
\end{itemize}
\end{theorem} 

Note that \reff{eq.wick.complexfermions.2a}
and \reff{eq.wick.complexfermions.2b}
are equivalent by virtue of Jacobi's identity \reff{eq.JacobyDet},
and that \reff{eq.wick.complexfermions.3a} and \reff{eq.wick.complexfermions.3b}
are equivalent by virtue of Jacobi's identity
together with the Cauchy--Binet identity \reff{eq.cauchy-binet}.
However, we shall give {\em independent}\/ proofs of all these formulae;
as a consequence, our argument provides an {\em ab initio}\/ proof
of Jacobi's identity by means of Grassmann--Berezin integration.


The proof of Theorem~\ref{thm.wick.complexfermions}
follows closely the pattern used in the proof
of Theorem~\ref{thm.wick.realfermions},
but with slightly different combinatorics.

\proofof{Theorem~\ref{thm.wick.complexfermions}}
To prove (a), we perform the translations 
$\psi'= \psi - A^{-1} \lambda$ and
$\bar\psi'= \bar\psi - A^{-\rm T} \bar\lambda$
and use Proposition~\ref{prop.translation}.
We have
\begin{eqnarray}
   & &
   \bar\psi^{\prime \rm T} A \psi' \,+\,
       \bar\lambda^{\rm T} \psi' \,+\, \bar\psi^{\prime \rm T} \lambda
      \nonumber \\[1mm]
   & & \qquad =\;
   (\bar\psi^{\rm T} - \bar\lambda^{\rm T} A^{-1}) A (\psi - A^{-1} \lambda)
      \,+\, \bar\lambda^{\rm T} (\psi - A^{-1} \lambda)
      \,+\, (\bar\psi^{\rm T} - \bar\lambda^{\rm T} A^{-1}) \lambda \quad
      \nonumber \\[1mm]
   & & \qquad =\;
   \bar\psi^{\rm T} A \psi \,-\,   \bar\lambda^{\rm T} A^{-1} \lambda  \;.
\end{eqnarray}
Therefore \reff{eq.wick.complexfermions.1}
is an immediate consequence of Proposition~\ref{prop.gaussian.complexfermions}.

To prove (b), we begin by observing as before that
\be
   \int\! \scrd(\psi,\psibar)  \:
   \Biggl( \prod_{\alpha=1}^r \psibar_{i_\alpha} \psi_{j_\alpha} \Biggr) \:
   \exp\!\left( \bar\psi^{\rm T} A \psi \right)
   \;=\;
   \int\! \scrd(\psi,\psibar)  \:
   \Biggl( \prod_{\alpha=1}^r \psibar_{i_\alpha} \psi_{j_\alpha} \Biggr) \:
   \exp\!\left( \bar\psi^{\rm T}_{I^c} A_{I^c J^c} \psi_{J^c} \right)
   \;.
 \label{eq.proof.wick.complexfermions.b1}
\ee
We now proceed to reorder the integration measure.
We have $I = \{i_1,\ldots,i_r\}$ with $i_1 < \ldots < i_r$,
and let us write $I^c = \{i'_1,\ldots,i'_{n-r}\}$
with $i'_1 < \ldots < i'_{n-r}$; and likewise for~$J$.
Then $\scrd(\psi,\psibar)$ is the product of factors
$d\psi_i \, d\psibar_i$ taken in arbitrary order;
we choose the order to be $II^c$ (i.e., $i_1 \cdots i_r i'_1 \cdots i'_{n-r}$).
We now leave the factors $\psibar$ in place,
but reorder the factors $\psi$ to be in order $JJ^c$ rather than $II^c$:
this produces a sign $\epsilon(I,J)$.
Using the notation
\be
   \scrd(\psi_J,\psibar_I)  \;=\;
   d\psi_{j_1} d\psibar_{i_1} \,\cdots\, d\psi_{j_r} d\psibar_{i_r}
\ee
and likewise for the complementary sets, we have proven that
\be
  \scrd(\psi,\psibar)  \;=\;
  \epsilon(I,J) \: \scrd(\psi_J,\psibar_I) \: \scrd(\psi_{J^c},\psibar_{I^c})
  \;.
\ee
This in turn can be trivially rewritten as
$\epsilon(I,J) \: \scrd(\psi_{J^c},\psibar_{I^c})  \: \scrd(\psi_J,\psibar_I)$
since the two factors $\scrd$ are Grassmann-even.
We can now apply this measure to the integrand in
\reff{eq.proof.wick.complexfermions.b1};
integrating the monomial against $\scrd(\psi_J,\psibar_I)$ gives 1,
and integrating the exponential against $\scrd(\psi_{J^c},\psibar_{I^c})$
gives $\det A_{I^c J^c}$ by Proposition~\ref{prop.gaussian.complexfermions}.

(c)  As in the case of ``real'' fermions,
we specialize result (a) to the case where the ``sources''
$\lambda_i,\lambdabar_i$ are generators of an extended Grassmann algebra,
and then differentiate with respect to them:
\begin{eqnarray}
   & &
   \int\! \scrd(\psi,\psibar)  \:
      \Biggl( \prod_{\alpha=1}^r \psibar_{i_\alpha} \psi_{j_\alpha} \Biggr)  \:
   \exp\!\left( \bar\psi^{\rm T} A \psi \right) 
    \nonumber\\[1mm]
   & & \qquad =\;
   (\det A) \:
   \left.  \Biggl( \prod_{\alpha=1}^r
                   - \frac{\partial}{\partial \lambda_{i_\alpha}}
                     \frac{\partial}{\partial \lambdabar_{j_\alpha}} \Biggr) \:
   \exp\!\left( -  \bar\lambda^{\rm T} A^{-1} \lambda \right)
   \right|_{\bar\lambda=\lambda=0}
        \nonumber \\[2mm]
   & & \qquad =\;
   (\det A) \:
   (-1)^r \int \! \scrd(\lambda_I,\lambdabar_J)
   \exp\!\left( -  \lambdabar^{\rm T}_J (A^{-1})_{JI} \lambda_I \right)
        \nonumber \\[2mm]
   & & \qquad =\;
   (\det A) \: (-1)^r \: \det((-A^{-1})_{JI})
        \nonumber \\[2mm]
   & & \qquad =\;
   (\det A) \, \det((A^{-\rm T})_{IJ})
  \;.
 \label{eq.wick.complexfermions.proofc}
\end{eqnarray}
Here the minus sign in the first equality comes from the fact
that differentiation of the source term $\psibar^{\rm T} \lambda$
with respect to $\lambda_i$ yields $-\psibar_i$
according to \reff{eq.leibniz};
the second equality says that differentiation is the same as integration;
and the third equality uses \reff{gaussint.fermionic.complex.EXTENDED}.

(d)  We have
\begin{eqnarray}
 & &
 \!\!\!\!\!
 \int\! \scrd(\psi,\psibar) \,
      \Biggl( \prod_{\alpha=1}^r \left(\bar\psi C\right)_\alpha
                                 \left(B \psi \right)_\alpha
      \Biggr) \:
   \exp\!\left( \bar\psi^{\rm T} A \psi \right)
       \nonumber \\
 & &
 \!\!\!\!\!
 \quad =\;
 \sum_{i_1,\ldots,i_r}  \sum_{j_1,\ldots,j_r}
 \left( \prod_{\alpha=1}^r C_{i_\alpha, \alpha} B_{\alpha, j_\alpha} \right)
    \int\! \scrd(\psi,\psibar) \,
      \Biggl( \prod_{\alpha=1}^r \psibar_{i_\alpha} \psi_{j_\alpha} \Biggr)  \:
   \exp\!\left( \bar\psi^{\rm T} A \psi \right)   ,
    \qquad\qquad
\end{eqnarray}
but the only nonvanishing contributions in the sum
come when $i_1,\ldots,i_r$, and also $j_1,\ldots,j_r$, are all distinct;
so we can first require $i_1 < \ldots < i_r$, and $j_1 < \ldots < j_r$
and then sum over permutations, yielding
\begin{eqnarray}
 & &
   \sum_{|I|=r} \, \sum_{|J|=r} \, \sum_{\sigma,\tau \in \scrs_r} 
       C_{i_{\sigma(1)},1} \cdots C_{i_{\sigma(r)},r} \,
       B_{1,j_{\tau(1)}} \cdots B_{r,j_{\tau(r)}}
    \; \times
   \nonumber  \\ 
 & &
 \qquad\qquad
       \sgn(\sigma) \, \sgn(\tau) \,
       \int\! \scrd(\psi,\psibar)  \:
      \Biggl( \prod_{\alpha=1}^r \psibar_{i_\alpha} \psi_{j_\alpha} \Biggr)  \:
   \exp\!\left( \bar\psi^{\rm T} A \psi \right)
   \nonumber  \\ 
 & &
 \qquad =\;
   \sum\limits_{|I| = r} \, \sum\limits_{|J|=r} (\det B_{\star J}) \, \epsilon(I,J) \, (\det A_{I^c J^c}) (\det C_{I \star}) 
   \;.
\end{eqnarray}
by part (b).
If $A$ is invertible, we can write
\begin{eqnarray}
   & &
   \int\! \scrd(\psi,\psibar) \,
      \Biggl( \prod_{\alpha=1}^r \left(\bar\psi C\right)_\alpha
                                 \left(B \psi \right)_\alpha
      \Biggr) \:
          \exp\!\left( \bar\psi^{\rm T} A \psi \right) 
        \nonumber \\
   & & \qquad =\;
  \Biggl(\prod_{i=1}^r - \frac{\partial}{ \partial \lambda_{i}} \frac{\partial }{ \partial \bar\lambda_{i}} \Biggr)\:
   \left.
     \int\! \scrd(\psi,\psibar) \,
   \exp\!\left( \bar\psi^{\rm T} A \psi \,+\, \bar\lambda^{\rm T} B \psi
                                       \,+\, \bar\psi^{\rm T} C \lambda
         \right)   \right|_{\bar\lambda=\lambda=0}
      \nonumber \\[1mm]
   & & \qquad =\;
   \Biggl(\prod_{i=1}^r - \frac{\partial}{ \partial \lambda_{i}} \frac{\partial }{ \partial \bar\lambda_{i}} \Biggr)\:
   \left.
  (\det A) \, \exp\!\left( - \bar\lambda^{\rm T} B A^{-1} C \lambda
                    \right)
   \right|_{\bar\lambda=\lambda=0}
\end{eqnarray}
by using (a).
Then by the same reasoning as in \reff{eq.wick.complexfermions.proofc}
we see that this equals $(\det A) \, \det(B A^{-1} C)$.
\qed

\subsection{Summary of bosonic and fermionic Gaussian integration}

The four types of Gaussian integration and their basic characteristics
are summarized in Table~\ref{table.Gaussian.summary}.
Let us stress the following facts:
\begin{itemize}
   \item Fermionic integration is a purely algebraic/combinatorial operation:
      no analytic conditions on the matrix $A$ (such as positive-definiteness)
      are needed, and in fact fermionic integration makes sense over an
      arbitrary (commutative) coefficient ring~$R$.
      Bosonic integration, by contrast, is an analytic operation
      (at least as we have defined it here) and requires positive-definiteness
      of $A$ or its hermitian part.
   \item For complex bosons and ``complex'' fermions, no symmetry or
      antisymmetry conditions are imposed on the matrix $A$.
      By contrast, for real bosons and ``real'' fermions,
      the matrix $A$ must be symmetric or antisymmetric, respectively.
\end{itemize}

%
%
\begin{table}[t]
\hspace*{-0.8cm}
\small
\begin{tabular}{|c||c|c||c|c|}
\hline
   Type of   & Combinatorial & Analytic   & Value of  & Contractions in  \\
   Variables & conditions    & conditions & Gaussian int. &
                                              $2k$-point correlation \\
\hline\hline
   Real bosons     &  $A$ symmetric  &
     $A$ pos.-def.  &  $(\det A)^{-1/2}$  & $2k \times 2k$ hafnian  \\
\hline
   Complex bosons  &  none  &
       $A + A^*$ pos.-def. &
       $(\det A)^{-1}$  &  $k \times k$ permanent \\
\hline
   ``Real'' fermions   &  $A$ antisymm.  &  none &
       $\pf A$  &  $2k \times 2k$ pfaffian \\
\hline
   ``Complex'' fermions   &  none  &  none  &
       $\det A$ & $k \times k$ determinant \\
\hline
\end{tabular}
\caption{
   Summary of the four types of Gaussian integration.
   Here $A$ is the matrix appearing in the quadratic form
   in the exponential.
   See text for details.
}
  \label{table.Gaussian.summary}
\end{table}

\section{Some useful identities}  \label{app.identities}

In this appendix we collect some auxiliary results that
will be used at various places in this paper:
identities for sums of products of binomial coefficients
(Section~\ref{app.binomial})
and for determinants and pfaffians (Section~\ref{app.det}),
lemmas on matrix factorization (Section~\ref{app.matrix.decomp}),
and an identity that we call the ``dilation-translation formula''
(Section~\ref{app.gen.transl}).

\subsection{Binomial identities}  \label{app.binomial}

We collect here a few identities for sums of products of
binomial coefficients that will be needed in the proof of
the rectangular Cayley identities
(Sections~\ref{sec.grassmann.2.tmrect} and ~\ref{sec.grassmann.2.rect}).
We use the standard convention \cite{Graham_94}
for the definition of binomial coefficients:
\begin{equation}
\label{eq.binomdef}
\binom{r}{k}  \;=\;
\cases{
   r (r-1) \cdots (r-k+1)/k! & for integer $k > 0$  \cr
   \noalign{\vskip 0.5mm}
   1                         & for $k=0$   \cr
   \noalign{\vskip 0.5mm}
   0                         & for integer $k<0$  \cr
}
\end{equation}
where $r$ is an indeterminate and $k$ is always an integer.
Multinomial coefficients
\begin{equation}
\binom{a_1 + \cdots + a_k}{a_1, \ldots, a_k}
=
\frac{(a_1 + \cdots + a_k)!}{a_1! \cdots a_k!}
\end{equation}
will, by contrast,
be used only when all the $a_i$ are nonnegative integers.


We now state some easy combinatorial lemmas involving binomial coefficients.
They are all either contained in the textbook of
Graham--Knuth--Patashnik \cite[Chapter 5]{Graham_94}
or derived in the same fashion.

\begin{lemma}[parallel summation]   
\begin{equation}
\sum_{k \le m} \binom{r+k}{k}  \;=\; \binom{r+m+1}{m}
\qquad\hbox{\rm for $m$ integer}
\ef.
 \label{eq.lemma.A1}
\end{equation}
\end{lemma}

\begin{lemma}[Chu--Vandermonde convolution]
  \label{lemma.vandermonde}
\begin{equation}
\sum_k \binom{w}{k\vphantom{p}} \binom{m}{p-k} = \binom{w+m}{p}
\ef.
\end{equation}
\end{lemma}

\proof
For positive integers $w$ and $m$,
this is the number of ways of selecting $p$ people
out of a group of $w$ women and $m$ men.
Since, for any fixed integer $p$,
both sides are polynomials in $w$ and $m$,
the identity holds as a polynomial identity.
\qed

\begin{lemma}
\begin{equation}
\sum_{\begin{scarray}
        k,h,l \geq 0 \\
        k+h+l = m
      \end{scarray}}
(-1)^k
\binom{a\vphantom{b}}{h} \binom{b}{k} \binom{h+l}{h}
\;=\;
\binom{a-b+m}{m}
\ef.
 \label{eq.lemma.A3}
\end{equation}
\end{lemma}

\proof
The left-hand side equals
\begin{subeqnarray}
\sum_{k=0}^{m} \sum_{h=0}^{m-k}
   (-1)^k \binom{a\vphantom{b}}{h} \binom{b}{k} \binom{m-k}{m-k-h}
& = &
\sum_{k=0}^{m} (-1)^k \binom{a+m-k}{m-k} \binom{b}{k}
  \qquad
   \\
& = &
\sum_{k=0}^{m} (-1)^m \binom{-a-1}{m-k} \binom{b}{k}
   \\
& = &
(-1)^m \binom{b-a-1}{m}
   \\
& = &
\binom{a-b+m}{m}
\end{subeqnarray}
where the first and third equalities use the Chu--Vandermonde convolution.
\qed

\begin{lemma}
\label{lemma2}
\begin{equation}
\sum_{\begin{scarray}
         k,h,l \geq 0 \\
         k+h+l \leq m
      \end{scarray}}
(-1)^k \binom{a\vphantom{b}}{h} \binom{b}{k} \binom{h+l}{h}
\;=\;
\binom{a-b+m+1}{m}
\ef.
\end{equation}
\end{lemma}

\proof
Sum \reff{eq.lemma.A3} and use \reff{eq.lemma.A1}.
\qed

\subsection{Determinant and pfaffian identities}
\label{app.det}

We collect here some identities for determinants and pfaffians
that will be needed in Section~\ref{sec.grassmann.2}:
these concern the determinant or pfaffian of a partitioned matrix
and the change of determinant under low-rank perturbation.
These identities are well-known, but for completeness
we will give compact proofs using Grassmann--Berezin integration.
We will also give a (possibly new) ``fermionic'' analogue of the
low-rank-perturbation formula,
which will play a crucial role throughout Section~\ref{sec.grassmann.2}.

Let us begin with the formula for the determinant of a partitioned matrix,
due to Schur \cite[Hilfssatz, pp.~216--217]{Schur_17}:

\begin{proposition}[Schur's formula for the determinant of a partitioned matrix]
  \label{prop.schur}
   Consider a partitioned matrix of the form
\be
   M \;=\; \left( \begin{array}{c|c}
                    A & B \\
                    \hline
                    C & D
               \end{array} \right)
\ee
where $A,B,C,D$ are matrices of sizes
$m \times m$, $m \times n$, $n \times m$ and $n \times n$, respectively,
with elements in a commutative ring with identity.
\begin{itemize}
   \item[(a)]  If $A$ is invertible, then
       $\det M = (\det A) \det(D - C A^{-1} B)$.
   \item[(b)]  If $D$ is invertible, then
       $\det M = (\det D) \det(A - B D^{-1} C)$.
\end{itemize}
\end{proposition}

The matrix $D - C A^{-1} B$ is called the {\em Schur complement}\/
of $A$ in the partitioned matrix $M$;
see \cite{Cottle_74,Ouellette_81,Carlson_86,Zhang_05} for reviews.
One well-known proof of Schur's formula is based on the identity
\be
   \left( \begin{array}{c|c}
                    A & B \\
                    \hline
                    C & D
          \end{array} \right)
   \;=\;
   \left( \begin{array}{c|c}
                    A & 0 \\
                    \hline
                    C & I_n
          \end{array} \right)
   \left( \begin{array}{c|c}
                    I_m & A^{-1} B \\
                    \hline
                    0 & D - C A^{-1} B
          \end{array} \right)
   \;.
\ee
Let us give a quick proof of Schur's formula
using Grassmann--Berezin integration:

\proofof{Proposition~\ref{prop.schur}}
Let $R$ be the commutative ring with identity in which the elements of
$A,B,C,D$ take values.
We introduce Grassmann variables $\psi_i, \psibar_i$ ($1 \le i \le m$) and
$\eta_j, \etabar_j$ ($1 \le j \le n$)
and work in the Grassmann algebra $R[\psi,\psibar,\eta,\etabar]$.
We have
\begin{eqnarray}
   \det \left( \begin{array}{c|c}
                    A & B \\
                    \hline
                    C & D
               \end{array} \right)
& = &
\int \! \scrd_m(\psi,\psibar) \, \scrd_n(\eta,\etabar)
\:
\exp \!\left[
(\psibar^{\rm T}, \etabar^{\rm T})
\left( \begin{array}{c|c}
                    A & B \\
                    \hline
                    C & D
               \end{array} \right)
\left( \! \begin{array}{c}
\psi \\ \eta
               \end{array} \! \right)
\right]
\nonumber
\\
& = &
\int \! \scrd_m(\psi,\psibar) \, \scrd_n(\eta,\etabar)
\:
\exp \!\left[
   \psibar^{\rm T} A \psi  +  \psibar^{\rm T} B \eta +
   \etabar^{\rm T} C \psi  +  \etabar^{\rm T} D \eta
\right]  \;.
   \nonumber \\
\end{eqnarray}
If $A$ is invertible, we can perform the integration over $\psi,\psibar$
using Wick's theorem for ``complex'' fermions
(Theorem~\ref{thm.wick.complexfermions}), yielding
\be
   (\det A) \int \! \scrd_n(\eta,\etabar) \,
\exp \!\left( \etabar^{\rm T} D \eta - \etabar^{\rm T} C A^{-1} B \eta \right)
   \;.
\ee
Then performing the integration over $\eta,\etabar$ yields
$(\det A) \det(D - C A^{-1} B)$.
This proves (a);  and the proof of (b) is identical.
\qed

Here are some important special cases of Proposition~\ref{prop.schur}:

\begin{corollary}
   \label{cor.blockmatrixdet}
Let $U,V$ be $m \times n$ matrices with elements in a commutative ring
with identity.  Then
\be
   \det(U V^{\rm T})  \;=\;
   \det \left( \begin{array}{c|c}
                    0_m & U \\
                    \hline
                    -V^{\rm T}  & I_n
               \end{array} \right)
   \;,
\ee
where $0_m$ is the $m \times m$ zero matrix
and $I_n$ is the $n \times n$ identity matrix.
\end{corollary}

\begin{corollary}[matrix determinant lemma]
   \label{cor.sylvester}
Let $A$ be an invertible $m \times m$ matrix,
let $W$ be an invertible $n \times n$ matrix,
and let $U,V$ be $m \times n$ matrices,
all with elements in a commutative ring with identity.
Then
\be
   \det(A + UWV^{\rm T})  \;=\;
   (\det A) \, (\det W) \, \det(W^{-1} + V^{\rm T} A^{-1} U)
   \;.
 \label{eq.sylvester}
\ee
In particular, if we take $W = I_n$,
then
\be
   \det(A + UV^{\rm T})  \;=\;
   (\det A) \, \det(I_n + V^{\rm T} A^{-1} U)
   \;.
 \label{eq.sylvester.1}
\ee
If in addition we take $A = I_m$, then
\be
   \det(I_m + UV^{\rm T})  \;=\;  \det(I_n + V^{\rm T} U)
   \;.
 \label{eq.sylvester.2}
\ee
\end{corollary}

Corollary~\ref{cor.sylvester} is sometimes known as the
``matrix determinant lemma'',
and the special case \reff{eq.sylvester.2} is sometimes known as
``Sylvester's theorem for determinants''.
When $n \ll m$, we can interpret \reff{eq.sylvester}--\reff{eq.sylvester.2}
as formulae for the change of determinant under a low-rank perturbation:
see Lemma~\ref{lemma.lowrank} below for an explicit statement.

Analogues of Proposition~\ref{prop.schur}
and Corollary~\ref{cor.blockmatrixdet} exist also for pfaffians:

\begin{proposition}[Pfaffian of a partitioned matrix]
  \label{prop.schur.pfaffian}
   Consider a partitioned matrix of the form
\be
   M \;=\; \left( \begin{array}{c|c}
                    A & B \\
                    \hline
                    -B^{\rm T} & D
               \end{array} \right)
\ee
where $A,B,D$ are matrices of sizes
$2m \times 2m$, $2m \times 2n$ and $2n \times 2n$, respectively,
with elements in a commutative ring with identity,
and $A$ and $D$ are antisymmetric.
\begin{itemize}
   \item[(a)]  If $A$ is invertible, then
       $\pf M = (\pf A) \pf(D + B^{\rm T} A^{-1} B)$.
   \item[(b)]  If $D$ is invertible, then
       $\pf M = (\pf D) \pf(A + B D^{-1} B^{\rm T})$.
\end{itemize}
\end{proposition}

\noindent
Note that the matrices $D + B^{\rm T} A^{-1} B$
and $A + B D^{-1} B^{\rm T}$ appearing here
are just the usual Schur complements $D - C A^{-1} B$ and $A - B D^{-1} C$
specialized to $C = -B^{\rm T}$.

\begin{corollary}
   \label{cor.blockmatrixpf}
Let $U$ be a $2m \times 2n$ matrix with elements in a commutative ring
with identity.  Then
\be
   \pf(U J_{2n} U^{\rm T})  \;=\;
   (-1)^m \pf(-U J_{2n} U^{\rm T})  \;=\;
   (-1)^m
   \pf \left( \begin{array}{c|c}
                    0_{2m} & U \\
                    \hline
                    -U^{\rm T}  & J_{2n}
               \end{array} \right)
\ee
where $J_{2n}$ is the standard $2n \times 2n$ symplectic form
\reff{def.J.appendix}.
\end{corollary}

\proofof{Proposition~\ref{prop.schur.pfaffian}}
This time we introduce ``real'' Grassmann variables
$\theta_i$ ($1 \le i \le 2m$) and $\lambda_i$ ($1 \le i \le 2n$).
We have
\be
   \pf M
   \;=\;
   \int \! \scrd_{2m}(\theta) \, \scrd_{2n}(\lambda)
   \:
   \exp \!\left[
      \smhalf \theta^{\rm T} A \theta +  \theta^{\rm T} B \lambda
      + \smhalf \lambda^{\rm T} D \lambda \right]
   \;.
\ee
If $A$ is invertible, we can perform the integration over $\theta$
using Wick's theorem for ``real'' fermions
(Theorem~\ref{thm.wick.realfermions}), yielding
\be
   (\pf A) \int \! \scrd_{2n}(\lambda) \,
\exp \left( \smhalf \lambda^{\rm T} D \lambda
            + \smhalf \lambda^{\rm T} B^{\rm T} A^{-1} B \lambda \right)
   \;.
\ee
Then performing the integration over $\lambda$ yields
$(\pf A) \pf(D + B^{\rm T} A^{-1} B)$.
This proves (a);  and the proof of (b) is identical.
\qed

\bigskip

We next wish to prove an analogue of Corollary~\ref{cor.sylvester}
when the entries in the various matrices belong,
not to a commutative ring, but to a Grassmann algebra.
More precisely, the entries in $A$ and $W$ will be {\em even}\/ elements
of the Grassmann algebra, while the entries in $U$ and $V$
will be {\em odd}\/ elements of the Grassmann algebra.

\begin{proposition}[fermionic matrix determinant lemma]
   \label{prop.sylvester.fermionic}
Let $\scrg$ be a Grassmann algebra over a commutative ring with identity;
let $A$ be an invertible $m \times m$ matrix
and $W$ an invertible $n \times n$ matrix,
whose elements belong to $\scrg_{\rm even}$;
and let $U,V$ be $m \times n$ matrices
whose elements belong to $\scrg_{\rm odd}$.
Then\footnote{
   Note that all the matrix elements in all these determinants
   belong to the {\em commutative}\/ ring $\scrg_{\rm even}$;
   therefore, these determinants are unambiguously defined.
}
\be
   \det(A + UWV^{\rm T})  \;=\;
   (\det A) \, (\det W)^{-1} \, \det(W^{-1} + V^{\rm T} A^{-1} U)^{-1}
   \;.
 \label{eq.sylvester.fermionic}
\ee
In particular, if we take $W = I_n$,
then
\be
   \det(A + UV^{\rm T})  \;=\;
   (\det A) \, \det(I_n + V^{\rm T} A^{-1} U)^{-1}
   \;.
 \label{eq.sylvester.fermionic.1}
\ee
If in addition we take $A = I_m$, then
\be
   \det(I_m + UV^{\rm T})  \;=\;  \det(I_n + V^{\rm T} U)^{-1}
   \;.
 \label{eq.sylvester.fermionic.2}
\ee
\end{proposition}

Note that \reff{eq.sylvester.fermionic}--\reff{eq.sylvester.fermionic.2}
differ from \reff{eq.sylvester}--\reff{eq.sylvester.2}
by replacing some determinants with their inverse.

\proof
Let us first observe that since $A$ and $W$ are invertible,
we can rewrite \reff{eq.sylvester.fermionic} as
\be
   \det(I_m + A^{-1} U WV^{\rm T})  \;=\;
   \det(I_n + W V^{\rm T} A^{-1} U)^{-1}
   \;,
\ee
which is equivalent to
\be
   \det(I_m + \widetilde{U} \widetilde{V}^{\rm T})
   \;=\;  \det(I_n + \widetilde{V}^{\rm T} \widetilde{U})^{-1}
\ee
under the (invertible) change of variables
$\widetilde{U} = A^{-1} U$, $\widetilde{V} = V W^{\rm T}$.
So it suffices to prove \reff{eq.sylvester.fermionic.2}.

Let us first prove \reff{eq.sylvester.fermionic.2} when the ring $R$ is $\R$.
We augment the Grassmann algebra by introducing
generators $\eta_j, \etabar_j$ ($1 \le j \le n$),
and we also introduce bosonic variables
$\varphi_i, \varphibar_i$ ($1 \le i \le m$).
We wish to consider the mixed bosonic-fermionic integral
\be
   \int \! \scrd_m(\varphi,\varphibar) \, \scrd_n(\eta,\etabar) \,
\exp \left[
   -\varphibar^{\rm T} \varphi  +  \varphibar^{\rm T} U \eta +
   \etabar^{\rm T} V^{\rm T} \varphi  +  \etabar^{\rm T} \eta
\right]
   \;.
 \label{eq.sylvester.mixedint}
\ee
Please note that the quantities $U \eta$, $\etabar^{\rm T} V^{\rm T}$
and $\etabar^{\rm T} \eta$ all belong to the even ``pure soul'' part 
of the augmented Grassmann algebra, and in particular are nilpotent;
therefore the integrand can be interpreted as
\be
   \exp(-\varphibar^{\rm T} \varphi) \,
   \sum_{k=0}^\infty {1 \over k!} \,
     (\varphibar^{\rm T} U \eta + \etabar^{\rm T} V^{\rm T} \varphi
                                + \etabar^{\rm T} \eta)^k
\ee
where the sum over $k$ is in fact {\em finite}\/.
The Gaussian integration over $\varphi,\varphibar$
can thus be interpreted as separate Gaussian integrations
for the coefficients (which belong to $\R$)
of each monomial in the augmented Grassmann algebra;
this makes perfect analytic sense.
Let us now evaluate \reff{eq.sylvester.mixedint} in two ways:
Integrating first over $\varphi,\varphibar$
and then over $\eta,\etabar$, we get
\be
   \int \! \scrd_n(\eta,\etabar) \,
   \exp \left[ \etabar^{\rm T} \eta + \etabar^{\rm T} V^{\rm T} U \eta
        \right]
   \;=\;
   \det(I_n + V^{\rm T} U)
   \;.
\ee
On the other hand, integrating first over $\eta,\etabar$
and then over $\varphi,\varphibar$, we get
\be
   \int \! \scrd_m(\varphi,\varphibar) \,
   \exp \left[ -\varphibar^{\rm T} \varphi 
               -\varphibar^{\rm T} U V^{\rm T} \varphi
        \right]
   \;=\;
   \det(I_m + U V^{\rm T})^{-1}
   \;.
\ee
This proves \reff{eq.sylvester.fermionic.2} when the ring $R$ is $\R$.

Let us finally give an abstract argument showing that
\reff{eq.sylvester.fermionic.2} holds for arbitrary commutative rings $R$.
The point is that the matrix elements of $U$ and $V$ belong to
a Grassmann algebra over some finite set of generators $\chi_1,\ldots,\chi_N$,
and hence can be written as
\begin{subeqnarray}
   U_{ij}  & = & \sum_{K\,{\rm odd}} \alpha_{ij;K} \, \chi^K  \\[2mm]
   V_{ij}  & = & \sum_{K\,{\rm odd}} \beta_{ij;K} \, \chi^K
\end{subeqnarray}
for some coefficients $\alpha_{ij;K}, \beta_{ij;K} \in R$.
Now, both sides of \reff{eq.sylvester.fermionic.2}
are of the form $\sum\limits_{L\,{\rm even}} \gamma_L \, \chi^L$
where the coefficients $\gamma_L$ are polynomials in
$\{\alpha_{ij;K}, \beta_{ij;K}\}$ with integer coefficients.\footnote{
   The determinant on the right-hand side is an element of the
   Grassmann algebra whose ``body'' term is 1;
   therefore, it is invertible in the Grassmann algebra
   and the coefficients of its inverse are polynomials
   (with integer coefficients) in its own coefficients.
}
But we have just shown that these two polynomials coincide
whenever $\{\alpha_{ij;K}, \beta_{ij;K}\}$
are replaced by any set of specific values in $\R$.
Therefore, they must coincide
{\em as polynomials in the indeterminates $\{\alpha_{ij;K}, \beta_{ij;K}\}$}\/.
But this implies that they are equal when $\{\alpha_{ij;K}, \beta_{ij;K}\}$
are replaced by specific values in {\em any}\/ commutative ring $R$.
\qed

It is convenient, for our applications, to rephrase
\reff{eq.sylvester.2} and \reff{eq.sylvester.fermionic.2}
as explicit formulae for the change of determinant
under a low-rank perturbation.
So let $u_1,\ldots,u_n$ and $v_1,\ldots,v_n$ be vectors of length $m$
whose entries are elements of a Grassmann algebra
and are either all Grassmann-even or all Grassmann-odd.
We call these cases $\epsilon=+1$ and $\epsilon=-1$, respectively.
We then have the following formula for the determinant of
a rank-$n$ perturbation of the identity matrix:

\begin{lemma}[low-rank perturbation lemma]
   \label{lemma.lowrank}
Let $u_1,\ldots,u_n$ and $v_1,\ldots,v_n$ be as above,
and define the $n \times n$ matrix $M$ by
\begin{equation}
M_{\alpha \beta} \;=\; v_{\alpha} \cdot u_{\beta}
   \;\equiv\; \sum_{i=1}^m (v_{\alpha})_i (u_{\beta})_i
\ef.
\end{equation}
Then we have
\be
\det \Big(I_m+ \sum_{\alpha=1}^n u_{\alpha} v_{\alpha}^{\rm T}  \Big)
 \;=\;
 \det(I_n + M)^{\epsilon}
 \label{eq.lemma.lowrank}
\ee
where $\epsilon= \pm 1$ is as above.
\end{lemma}



Here is one special case of the low-rank perturbation lemma
that will be useful in
Sections~\ref{sec.grassmann.2.symmetric} and \ref{sec.Lap.symmetric}
in treating symmetric Cayley identities:

\begin{corollary}
   \label{corol.lowranksym}
Let $\eta = (\eta_1,\ldots,\eta_m)$ and
$\etabar = (\etabar_1,\ldots,\etabar_m)$
be Grassmann variables,
and let $A$ be an invertible $m \times m$ symmetric matrix
whose elements are Grassmann-even (hence commute with everything).
Then
\be
\det( I_m + \etabar \eta^{\rm T} - A^{-1} \eta \etabar^{\rm T} A)
 \;=\;
(1- \etabar^{\rm T} \eta)^{-2}
   \;.
\ee
\end{corollary}

\proof
Applying Lemma~\ref{lemma.lowrank} with $u_1 = \etabar$, $v_1 = \eta$,
$u_2 = A^{-1} \eta$ and $v_2 = -A^{\rm T} \etabar = -A\etabar$ gives
\be
\det( I_m + \etabar \eta^{\rm T} - A^{-1} \eta \etabar^{\rm T} A)
 \;=\;
{\det}^{-1}  \!
\left(
\begin{array}{cc}
1+ \eta^{\rm T} \etabar         & \eta^{\rm T} A^{-1} \eta \\
-\etabar^{\rm T} A \etabar      & 1 - \etabar^{\rm T} \eta
\end{array}
\right)
\ee
Since $A$ and $A^{-1}$ are symmetric,
the off-diagonal elements vanish, which gives the result.
\qed

Another special case of the low-rank perturbation lemma
will arise in Section~\ref{sec.grassmann.2.multi}:
here the $m \times m$ matrix $I_m + U V^{\rm T}$ occurring on
the left-hand side on \reff{eq.lemma.lowrank}
will be written as a product of {\em rectangular}\/ matrices,
each of which is a rank-one perturbation of the corresponding
rectangular pseudo-identity matrix.
(The $m \times n$ pseudo-identity matrix $\widehat{I}_{mn}$
 has matrix elements $(\widehat{I}_{mn})_{ij} = \delta_{ij}$.)
Direct application of the low-rank perturbation lemma
to such a product matrix yields a rather messy result,
but after some row operations we can obtain a fairly neat alternative formula:

\begin{corollary}
   \label{cor.lowrankrect}
Fix integers $\ell \ge 1$ and $n_1,\ldots,n_{\ell} \ge 1$
with $n_\alpha \ge n_1$ for $2 \le \alpha \le \ell$,
and write $n_{\ell+1} = n_1$.
Let $x_1,\ldots,x_\ell$ and $y_1,\ldots,y_\ell$ be vectors,
where $x_\alpha$ is of length $n_\alpha$
and $y_\alpha$ is of length $n_{\alpha+1}$,
whose entries are elements of a Grassmann algebra
and are either all commuting ($\epsilon=+1$)
or else all anticommuting ($\epsilon=-1$).
Then we have
\be
   \det\Biggl( \prod_{\alpha=1}^\ell
            (\widehat{I}_{n_\alpha n_{\alpha+1}} - x_\alpha y_\alpha^{\rm T})
       \Biggr)
   \;=\;
   (\det N)^\epsilon
\label{eq.cor.lowrankrect}
\ee
where the product is read from left ($\alpha=1$) to right ($\alpha=\ell$),
and the $\ell \times \ell$ matrix $N$ is defined by
\be
   N_{\alpha \beta}
   \;=\;
   \cases{ \sum\limits_{i=n_1+1}^{n_{\alpha+1,\beta}} y^\alpha_i x^\beta_i
                   & \hbox{\rm if } $\alpha < \beta$  \cr
           \noalign{\vskip 6pt}
           \delta_{\alpha\beta} - \sum\limits_{i=1}^{n_1} y^\alpha_i x^\beta_i
                   & \hbox{\rm if } $\alpha \ge \beta$  \cr
         }
\label{eq.cor.lowrankrect.defN}
\ee
where $n_{\alpha,\beta} = \min\limits_{\alpha \le \gamma \le \beta} n_\gamma$.
\end{corollary}

\proof
Note first that 
$\widehat{I}_{n_\alpha n_{\alpha+1}}
 \widehat{I}_{n_{\alpha+1} n_{\alpha+2}} \cdots
 \widehat{I}_{n_{\beta-1} n_{\beta}}$
is an $n_\alpha \times n_{\beta}$ matrix whose $ij$ element is
1 if $i=j\leq n_{\alpha,\beta}$ and 0 otherwise.
In particular,
$\widehat{I}_{n_1 n_2}
 \widehat{I}_{n_2 n_3} \cdots
 \widehat{I}_{n_{\alpha-1} n_{\alpha}} 
= \widehat{I}_{n_{1} n_{\alpha}}$.
For $v$ a vector of length $m \geq n_1$, define $\bar{v}$ as the vector
restricted to the first $n_1$ components.
So we can expand the matrix on the left-hand side of
\reff{eq.cor.lowrankrect} as
\begin{eqnarray}
&&
\!\!\!\!\!
\prod_{\alpha=1}^\ell
(\widehat{I}_{n_\alpha n_{\alpha+1}} - x_\alpha y_\alpha^{\rm T})
\nonumber
\\
&&
\;=\;
I_{n_1} - \sum_{\alpha=1}^{\ell}
\widehat{I}_{n_1 n_2}
 \widehat{I}_{n_2 n_3} \cdots
 \widehat{I}_{n_{\alpha-1} n_{\alpha}} 
x_{\alpha} 
y_{\alpha}^{\rm T}
(\widehat{I}_{n_{\alpha+1} n_{\alpha+2}} - x_{\alpha+1} y_{\alpha+1}^{\rm T})
\cdots
(\widehat{I}_{n_{\ell} n_{\ell+1}} - x_{\ell} y_{\ell}^{\rm T})
\nonumber
\\
&&
\;=\;
I_{n_1} - \sum_{\alpha=1}^{\ell}
\bar{x}_{\alpha} 
y_{\alpha}^{\rm T}
(\widehat{I}_{n_{\alpha+1} n_{\alpha+2}} - x_{\alpha+1} y_{\alpha+1}^{\rm T})
\cdots
(\widehat{I}_{n_{\ell} n_{\ell+1}} - x_{\ell} y_{\ell}^{\rm T})
   \;.
\end{eqnarray}
In this form, we are ready to apply Lemma~\ref{lemma.lowrank} with vectors 
$u_{\alpha} = -\bar{x}_{\alpha}$ and
$v_{\alpha}^{\rm T} = y_{\alpha}^{\rm T}
(\widehat{I}_{n_{\alpha+1} n_{\alpha+2}} - x_{\alpha+1} y_{\alpha+1}^{\rm T})
\cdots
(\widehat{I}_{n_{\ell} n_{\ell+1}} - x_{\ell} y_{\ell}^{\rm T})$
for $\alpha =1, \ldots, \ell$. This gives
\be
\label{eq.cor.lowrankrect.p1}
   \det\Biggl( \prod_{\alpha=1}^\ell
               (\widehat{I}_{n_\alpha n_{\alpha+1}} - x_\alpha y_\alpha^{\rm T})
       \Biggr)
   \;=\;
   (\det N^{(0)})^\epsilon
\ee
with
$N^{(0)}_{\alpha \beta} = \delta_{\alpha \beta} + v_{\alpha}^{\rm T} u_{\beta}$.
Now observe that
\be
  \label{eq.54876978}
v_{\alpha}^{\rm T}
\;=\;
\bar{y}_{\alpha}^{\rm T}
- \sum_{\beta > \alpha}
\big(
y_{\alpha}^{\rm T}
\widehat{I}_{n_{\alpha+1} n_{\alpha+2}}
\cdots
\widehat{I}_{n_{\beta-1} n_{\beta}}
x_{\beta}
\big)
\,
v_{\beta}^{\rm T}
\ee
and call
$c_{\alpha \beta} = y_{\alpha}^{\rm T}
\widehat{I}_{n_{\alpha+1} n_{\alpha+2}}
\cdots
\widehat{I}_{n_{\beta-1} n_{\beta}}
x_{\beta}
= \sum_{i=1}^{n_{\alpha+1, \beta}}
y^{\alpha}_i x^{\beta}_i
$.
Thus, defining the $\ell \times \ell$ upper-triangular matrix $\widehat{C}$ as
\be
\widehat{C}_{\alpha \beta} \;=\;
\cases{ 
c_{\alpha \beta}      & if $\alpha < \beta$ \cr
\noalign{\vskip 4pt}
\delta_{\alpha \beta} & if $\alpha = \beta$ \cr
\noalign{\vskip 4pt}
0                     & if $\alpha > \beta$ \cr
}
\ee
we have
\be
   \sum_\beta  \widehat{C}_{\alpha \beta} v_{\beta}^{\rm T}
   \;=\;
   \bar{y}_{\alpha}^{\rm T}
   \;.
  \label{eq.54876978.BIS}
\ee
Clearly $\det \widehat{C} = 1$,
so if we define $N = C N^{(0)}$ we have $\det N = \det N^{(0)}$.
It is easy to see, using \reff{eq.54876978.BIS},
that $N$ is exactly the matrix given in \reff{eq.cor.lowrankrect.defN}.
\qed

Please note that all the entries of the matrix $N$
are polynomials of degree {\em at most two}\/ in the variables $x$ and $y$ ---
unlike the matrix $N^{(0)}$ coming from the bare application
of the low-rank perturbation lemma,
which contains terms of degree as high as $2\ell$.

\bigskip

Corollary~\ref{cor.blockmatrixdet} is in fact the case $\ell=2$
of a more general lemma that holds for $\ell \geq 2$,
and will be needed in Section~\ref{sec.grassmann.2.multi}:

\begin{lemma}
   \label{lem.blockmatrixdet.multi}
Fix integers $\ell \geq 2$ and $n_1,\ldots,n_\ell \ge 1$,
and write $n_{\ell+1} = n_1$.
Let $U_1, \ldots, U_\ell$ be 
matrices with elements in a commutative ring
with identity, $U_\alpha$ being of dimension 
$n_{\alpha} \times n_{\alpha+1}$.
Define
\be
M_\ell(U_1, \ldots, U_\ell)
:=
\left(
\begin{array}{ccccc}
0_{n_1}    & -U_1 & 0        & \cdots & 0             \\ 0          & I_{n_2}  & -U_2 &        & 0             \\
0          & 0        & I_{n_3}  & \ddots & \vdots        \\
\vdots     &          &          & \ddots & -U_{\ell-1} \\
U_\ell & 0        & \cdots   & 0      & I_{n_{\ell}}  
\end{array}
\right)
\;,
\ee
where $0_n$ is the $n \times n$ zero matrix
and $I_n$ is the $n \times n$ identity matrix.
Then
\be
   \det M_\ell(U_1, \ldots, U_\ell)
\;=\;
   \det(U_1 \cdots U_\ell)
   \;.
\ee
\end{lemma}
\proof
We prove this by induction on $\ell$.
The case $\ell=2$ is already proven by Corollary \ref{cor.blockmatrixdet}.
For $\ell \geq 2$ we use the Grassmann representation of $\det(M_\ell)$:
\begin{eqnarray}
& & \!\!\!
\det M_\ell(U_1, \ldots, U_\ell)
  \;=\;
\int \! \scrd_{n_1}(\psi^1,\psibar^1) \cdots
\scrd_{n_{\ell}}(\psi^{\ell},\psibar^{\ell})
   \nonumber \\
& & \qquad\qquad\qquad \times\;
\exp \!\left[
\sum_{\alpha=2}^{\ell}
( \psibar^{\alpha} \psi^{\alpha} 
- \psibar^{\alpha-1} U_{\alpha-1} \psi^{\alpha} )
+ \psibar^{\ell} U_\ell \psi^{1}
\right]
   \;.
   \qquad
  \label{eq.proof.lem.blockmatrixdet.multi}
\end{eqnarray}
Now perform the integration over $(\psi^{\ell},\psibar^{\ell})$:
highlighting the factors in the integrand that involve these fields,
we see that
\be
\int \! \scrd_{n_{\ell}}(\psi^{\ell},\psibar^{\ell})
\exp \!\left[
\psibar^{\ell} \psi^{\ell} 
- \psibar^{\ell-1} U_{\ell-1} \psi^{\ell}
+ \psibar^{\ell} U_\ell \psi^{1}
\right]
  \;=\;
\exp \!\left(
\psibar^{\ell-1} U_{\ell-1} U_\ell \psi^{1}
\right)
\;.
\ee
Comparing this with \reff{eq.proof.lem.blockmatrixdet.multi}, we see that
\be
   \det M_\ell(U_1, \ldots, U_\ell)
\;=\;
   \det M_{\ell-1}(U_1, \ldots, U_{\ell-2}, U_{\ell-1} U_\ell)
\;,
\ee
which provides the required inductive step.
\qed

{\bf Remark.}  Just as Corollary \ref{cor.blockmatrixdet} is a
specialization of the more general Proposition \ref{prop.schur},
so Lemma \ref{lem.blockmatrixdet.multi} has a similar
generalization, proven through an identical procedure
(of which the details are left to the reader), namely:

\begin{lemma}
   \label{lem.blockmatrixdet.multi.gen}
Fix integers $\ell \geq 2$ and $n_1,\ldots,n_\ell \ge 1$,
and write $n_{\ell+1} = n_1$.
Let $B_1,\ldots,B_\ell$ be 
matrices with elements in a commutative ring
with identity, $B_\alpha$ being of dimension 
$n_{\alpha} \times n_{\alpha+1}$.
Let $A_1,\ldots,A_\ell$ be square
matrices with elements in the same commutative ring,
$A_\alpha$ being of dimension $n_{\alpha} \times n_{\alpha}$.
Assume that $A_2,\ldots,A_\ell$ are invertible.
Define
\be
M(A_1, \ldots, A_\ell; B_1, \ldots, B_\ell)
:=
\left(
\begin{array}{ccccc}
A_1    & -B_1 & 0        & \cdots & 0             \\
0          & A_2  & -B_2 &        & 0             \\
0          & 0        & A_3  & \ddots & \vdots        \\
\vdots     &          &          & \ddots & -B_{\ell-1} \\
B_\ell & 0        & \cdots   & 0      & A_\ell  
\end{array}
\right)
\;.
\ee
Then
\be
   \det M(A_1, \ldots, A_\ell; B_1, \ldots, B_\ell)
   \;=\;
   \det( A_1 + B_1 A_2^{-1} B_2 \cdots A_\ell^{-1} B_\ell)
   \;
   \prod\limits_{j=2}^\ell \det A_j
   \;.
\ee
If also $A_1$ is invertible, then we can obtain
an expression with a form of cyclic symmetry:
\be
   \det M(A_1, \ldots, A_\ell; B_1, \ldots, B_\ell)
   \;=\;
   \det( I_{n_1} + A_1^{-1} B_1 A_2^{-1} B_2 \cdots A_\ell^{-1} B_\ell)
   \;
   \prod\limits_{j=1}^\ell \det A_j
   \;.
\ee
\end{lemma}

\subsection{Matrix factorization lemmas}   \label{app.matrix.decomp}

In Sections~\ref{sec.grassmann.2.antisymm},
 \ref{sec.grassmann.2.tmrect}, \ref{sec.grassmann.2.rect},
 \ref{sec.grassmann.2.antisymrect} and \ref{sec.grassmann.2.multi}
we shall need some matrix factorization lemmas having the general form
\begin{quote}
   For any matrix $X$ of the form ..... there exists a matrix $A$
   of the form ..... such that $\Phi(X,A) = 0$
   [where $\Phi$ denotes a specified collection of polynomial
    or rational functions].
\end{quote}
or the multi-matrix generalization thereof:
\begin{quote}
   For any matrices $X,Y,\ldots$ of the form ..... there exist
   matrices $A,B,\ldots$ of the form .....
   such that $\Phi(X,Y,\ldots,A,B,\ldots) = 0$.
\end{quote}
The prototype for such matrix decomposition lemmas is
the well-known Cholesky factorization \cite[Theorem~4.2.5]{Golub_96}:

\begin{lemma}[Cholesky factorization]
   \label{lemma.cholesky}
Let $X$ be a real symmetric positive-definite $n \times n$ matrix.
Then there exists a unique lower-triangular real matrix $A$
with strictly positive diagonal entries, such that $X = A A^{\rm T}$.
\end{lemma}

\noindent
A similar but less well-known result is the following factorization
for antisymmetric matrices
\cite{Weyl_50,Murnaghan_31,Bunch_82,Benner_00,Xu_03}:

\begin{lemma}[$\bm{A J A^{\rm T}}$ factorization of an antisymmetric matrix]
    \label{lemma.decomp.1.FULL}
Let $X$ be a (real or complex) antisymmetric $2m \times 2m$ matrix.
Then there exists a (real or complex, respectively)
$2m \times 2m$ matrix $A$ such that $X = A J A^{\rm T}$,
where $J$ is defined in \reff{def.J.appendix}.
In particular, if $X$ is nonsingular, then $A \in GL(2m)$.
[The form of $A$ can be further restricted in various ways,
 but we shall not need this.]
\end{lemma}

In our applications we shall not need the uniqueness of $A$,
but merely its existence.
Nor shall we need any particular structure of $A$ (e.g.\ triangularity)
beyond lying in $GL(n)$ or $O(n)$ or $Sp(2n)$ as the case may be.
Finally, and most importantly, we shall not need the
existence of $A$ for {\em all}\/ matrices $X$ of a given type,
but only for those in some nonempty open set
(for instance, a small neighborhood of the identity matrix).
We shall therefore give easy existence proofs
using the implicit function theorem.  It is an interesting open question
whether our decompositions actually extend to arbitrary matrices $X$
in the given classes.

More precisely, we shall need the following decomposition lemmas
in addition to Lemma~\ref{lemma.decomp.1.FULL}.
The matrices $\widehat{I}_{mn}$ are defined in \reff{def.Ihat}.

\begin{lemma}
    \label{lemma.decomp.2}
Let $X$ and $Y$ be a (real or complex) $m \times n$ matrices ($m \le n$)
of rank $m$ that are sufficiently close to the matrix $\widehat{I}_{mn}$.
Then there exist matrices $P, R \in GL(m)$ and $Q \in GL(n)$
such that $X=P \widehat{I}_{mn} Q$ and $Y=R \widehat{I}_{mn} Q^{-\rm T}$.
\end{lemma}

\begin{lemma}
    \label{lemma.decomp.3}
Let $X$ be a (real or complex) $m \times n$ matrix ($m \le n$)
of rank $m$ that is sufficiently close to the matrix $\widehat{I}_{mn}$.
Then there exist matrices $P \in GL(m)$ and $Q \in O(n)$
such that $X=P \widehat{I}_{mn} Q$.
\end{lemma}

\begin{lemma}
    \label{lemma.decomp.4}
Let $X$ be a (real or complex) $2m \times 2n$ matrix ($m \le n$)
of rank $2m$ that is sufficiently close to the matrix $\widehat{I}_{2m,2n}$.
Then there exist matrices $P \in GL(2m)$ and $Q \in Sp(2n)$
such that $X=P \widehat{I}_{2m,2n} Q$.
\end{lemma}

\begin{lemma}
    \label{lemma.decomp.5}
Let $\ell \ge 1$ and $n_1,\ldots,n_{\ell+1} \ge 1$;
and let $\{X_\alpha\}_{1 \le \alpha \le \ell}$
be (real or complex) $n_\alpha \times n_{\alpha+1}$ matrices
of rank $\min(n_\alpha, n_{\alpha+1})$
that are sufficiently close to the matrix
$\widehat{I}_{n_\alpha n_{\alpha+1}}$.
Then there exist matrices $\{P_\alpha\}_{1 \le \alpha \le \ell+1}$
with $P_\alpha \in GL(n_\alpha)$ such that
$X_\alpha =
 P_{\alpha} \, \widehat{I}_{n_{\alpha} n_{\alpha+1}} \, P_{\alpha+1}^{-1}$.
\end{lemma}

\noindent
Lemmas~\ref{lemma.decomp.1.FULL}--\ref{lemma.decomp.5} will be needed in
Sections~\ref{sec.grassmann.2.antisymm},
 \ref{sec.grassmann.2.tmrect}, \ref{sec.grassmann.2.rect},
 \ref{sec.grassmann.2.antisymrect} and \ref{sec.grassmann.2.multi},
respectively.
In addition, Cholesky factorization could be used in
Section~\ref{sec.grassmann.2.symmetric} but we were able to avoid it;
see the Remark at the end of that section.

The proofs of these lemmas will all follow the same pattern.
First we find an explicit pair $X_0,A_0$ (or the multi-matrix generalization)
with the needed properties.  Then we linearize the functions $\Phi$
in a neighborhood of $(X_0,A_0)$, and we show that the tangent space
for $A$ at $A_0$ is mapped onto the full tangent space for $X$ at $X_0$.
The required existence of $A$ for $X$ in a neighborhood of $X_0$
then follows from the implicit function theorem.\footnote{
   We will use the implicit function theorem in the following form:
   Let $U \subseteq \R^N$ and $V \subseteq \R^p$ be open sets,
   and let $f \colon\, U \times V \to \R^N$ be a $C^k$ function ($k \ge 1$).
   Let $u_0 \in U$ and $v_0 \in V$ satisfy $f(u_0,v_0) = 0$,
   with $(\partial f/\partial u)(u_0,v_0)$ nonsingular.
   Then there exist neighborhoods $U' \ni u_0$ and $V' \ni v_0$
   such that for all $v \in V'$ there exists a unique $u \in U'$
   satisfying $f(u,v) = 0$;  moreover, the map $v \mapsto u$ is $C^k$.
}
For completeness we will also show how Lemmas~\ref{lemma.cholesky}
and \ref{lemma.decomp.1.FULL}, which are known to hold globally,
have simple proofs in this ``infinitesimal'' setting.
Nearly all these proofs will be easy;
only the last (Lemma~\ref{lemma.decomp.5}) turns out to be slightly tricky.

\proofof{Lemma~\ref{lemma.cholesky} for $X$ near $I$}
Linearizing $X = A A^{\rm T}$ in a neighborhood of $(X_0,A_0) = (I,I)$
by writing $X = I + X'$ and $A = I + A'$, we have $X' = A' + (A')^{\rm T}$
(plus higher-order corrections that we always drop).
Then an explicit solution is given by the lower-triangular matrix
\be
   A'_{ij}  \;=\;
   \cases{ X'_{ij}/2   &  if $i=j$  \cr
           \noalign{\vskip 4pt}
           X'_{ij}     &  if $i>j$  \cr
           \noalign{\vskip 4pt}
           0           &  if $i<j$  \cr
         }
\ee
The rest follows from the implicit function theorem.\footnote{
   For completeness let us make explicit how the
   implicit function theorem is used in this case;
   the analogous reasoning for the remaining lemmas can be supplied
   by the reader.

   If $M = (m_{ij})_{i,j=1}^n$ is an $n \times n$ matrix,
   let us write $[M]_{\rm LT} = (m_{ij})_{1 \le i \le j \le n}$
   to denote its lower-triangular part.
   We then use the implicit function theorem (see the preceding footnote)
   as follows:
   Let $u = A =$ a generic lower-triangular matrix
     $(a_{ij})_{1 \le i \le j \le n}$,
   $v = [X]_{\rm LT}$,
   $u_0 = v_0 = [I]_{\rm LT}$,
   $f(u,v) = [A A^{\rm T} - X]_{\rm LT}$.
   (Since $A A^{\rm T} - X$ is manifestly symmetric,
    it vanishes if and only if its lower-triangular part does.)
   Then $(\partial f/\partial u)(u_0,v_0)$ is the linear map
   $A' \mapsto A' + (A')^{\rm T}$; in other words, we have
   $$
      {\partial f_{ij} \over \partial u_{kl}} (u_0,v_0)
      \;=\;
      \delta_{(ij),(kl)} \, (1 + \delta_{kl})
     \;,
   $$
   which is a diagonal matrix with nonzero entries (namely, 1 and 2),
   hence nonsingular.
}
\qed

\proofof{Lemma~\ref{lemma.decomp.1.FULL} for $X$ near $J$}
Linearizing $X = A J A^{\rm T}$ in a neighborhood of $(X_0,A_0) = (J,I)$
by writing $X = J + X'$ and $A = I - A' J$, we have $X' = A' - (A')^{\rm T}$.
Then one explicit solution is given by the strictly lower-triangular matrix
\be
   A'_{ij}  \;=\;
   \cases{ X'_{ij}     &  if $i>j$  \cr
           \noalign{\vskip 4pt}
           0           &  if $i \le j$  \cr
         }
\ee
Another explicit solution is given by the antisymmetric matrix $A' = X'/2$.
\qed

\proofof{Lemma~\ref{lemma.decomp.2}}
Linearizing $X=P \widehat{I}_{mn} Q$ and $Y=R \widehat{I}_{mn} Q^{-\rm T}$
in a neighborhood of
$(X_0,Y_0,P_0,R_0,Q_0) = (\widehat{I}_{mn}, \widehat{I}_{mn}, I_m, I_m, I_n)$
by writing $X = \widehat{I}_{mn} + X'$ and so forth, we obtain
$X' = P' \widehat{I}_{mn} + \widehat{I}_{mn} Q'$
and $Y' = R' \widehat{I}_{mn} - \widehat{I}_{mn} (Q')^{\rm T}$.
In terms of the block decompositions
$\widehat{I}_{mn} \equiv (I_m , 0_{m \times (n-m)})$,
$X' = (X'_1 , X'_2)$, $Y' = (Y'_1 , Y'_2)$
and
$Q' = \left( \!\! \begin{array}{cc}
                    Q'_{11}  &  Q'_{12}  \\[1mm]
                    Q'_{21}  &  Q'_{22}
                \end{array}
      \!\! \right)$,
we have
\begin{subeqnarray}
     X'_1  & = &  P' + Q'_{11}   \\
     X'_2  & = &  Q'_{12}   \\
     Y'_1  & = &  R' - (Q'_{11})^{\rm T}  \\
     Y'_2  & = &  - (Q'_{21})^{\rm T}
\end{subeqnarray}
We can choose $Q'_{11}$ and $Q'_{22}$ arbitrarily;
then the remaining unknowns $P', R', Q'_{12}, Q'_{21}$
are uniquely determined.
%
\qed

\medskip

We shall actually prove the following generalization of
Lemma~\ref{lemma.decomp.3}:

\addtocounter{defin}{-3}
\begin{lemma}
\hspace*{-3.2mm} ${}^{\bf\prime}$ \hspace{1mm}
Let $X$ be a (real or complex) $m \times n$ matrix ($m \le n$)
of rank $m$ that is sufficiently close to the matrix $\widehat{I}_{mn}$,
and let $Y$ be a (real or complex) $n \times n$ symmetric matrix
that is sufficiently close to the identity matrix $I_n$.
Then there exist matrices $P \in GL(m)$ and $Q \in GL(n)$
such that $X=P \widehat{I}_{mn} Q$ and $Y = Q^{\rm T} Q$.
\end{lemma}

\noindent
When $Y = I_n$ this reduces to Lemma~\ref{lemma.decomp.3}.

\proofof{Lemma~\ref{lemma.decomp.3}${}'$}
Linearizing $X=P \widehat{I}_{mn} Q$ and $Y = Q^{\rm T} Q$
in a neighborhood of $(X_0,Y_0,P_0,Q_0) = (\widehat{I}_{mn}, I_n, I_m, I_n)$
by writing $X = \widehat{I}_{mn} + X'$ and so forth,
we obtain $X' = P' \widehat{I}_{mn} + \widehat{I}_{mn} Q'$
and $Y' = Q' + (Q')^{\rm T}$.
In terms of the block decompositions
$\widehat{I}_{mn} \equiv (I_m , 0_{m \times (n-m)})$,
$X' = (X'_1 , X'_2)$,
$Y' = \left( \!\! \begin{array}{cc}
                    Y'_{11}  &  Y'_{12}  \\[1mm]
                    (Y'_{12})^{\rm T}  &  Y'_{22}
                \end{array}
      \!\! \right)$
with $Y'_{11}$ and $Y'_{22}$ symmetric,
and
$Q' = \left( \!\! \begin{array}{cc}
                    Q'_{11}  &  Q'_{12}  \\[1mm]
                    Q'_{21}  &  Q'_{22}
                \end{array}
      \!\! \right)$,
we have
\begin{subeqnarray}
     X'_1  & = &  P' + Q'_{11}   \\
     X'_2  & = &  Q'_{12}   \\
     Y'_{11}  & = &  Q'_{11} + (Q'_{11})^{\rm T}   \\
     Y'_{22}  & = &  Q'_{22} + (Q'_{22})^{\rm T}   \\
     Y'_{12}  & = &  Q'_{12} + (Q'_{21})^{\rm T}
\end{subeqnarray}
We can choose arbitrarily the antisymmetric parts of $Q'_{11}$ and $Q'_{22}$
(e.g.\ by taking $Q'_{11}$ and $Q'_{22}$ lower-triangular,
 or alternatively by taking $Q'_{11}$ and $Q'_{22}$ symmetric);
then the remaining unknowns are uniquely determined.
\qed

\medskip

Similarly, let us prove the following generalization of
Lemma~\ref{lemma.decomp.4}:

\begin{lemma}
\hspace*{-3.2mm} ${}^{\bf\prime}$ \hspace{1mm}
Let $X$ be a (real or complex) $2m \times 2n$ matrix ($m \le n$)
of rank $2m$ that is sufficiently close to the matrix $\widehat{I}_{2m,2n}$,
and let $Y$ be a (real or complex) $2n \times 2n$ antisymmetric matrix
that is sufficiently close to $J_{2n}$.
Then there exist matrices $P \in GL(2m)$ and $Q \in GL(2n)$
such that $X=P \widehat{I}_{2m,2n} Q$ and $Y = Q^{\rm T} J Q$.
\end{lemma}

\noindent
When $Y = J_{2n}$ this reduces to Lemma~\ref{lemma.decomp.4}.

\proofof{Lemma~\ref{lemma.decomp.4}${}'$}
Linearizing $X=P \widehat{I}_{2m,2n} Q$ and $Y = Q^{\rm T} J Q$
in a neighborhood of 
$(X_0,Y_0,P_0,Q_0) = (\widehat{I}_{2m,2n}, J_{2n}, I_{2m}, I_{2n})$
by writing $X = \widehat{I}_{2m,2n} + X'$,
$Y = J_{2n} + Y'$, $P = I_{2m} + P'$ and $Q = I_{2n} - J_{2n} Q'$,
we obtain $X' = P' \widehat{I}_{2m,2n} - \widehat{I}_{2m,2n} J_{2n} Q'$
and $Y = Q' - (Q')^{\rm T}$.
In terms of the block decompositions
$\widehat{I}_{mn} \equiv (I_m , 0_{m \times (n-m)})$,
$X' = (X'_1 , X'_2)$,
$Y' = \left( \!\! \begin{array}{cc}
                    Y'_{11}  &  Y'_{12}  \\[1mm]
                    -(Y'_{12})^{\rm T}  &  Y'_{22}
                \end{array}
      \!\! \right)$
with $Y'_{11}$ and $Y'_{22}$ antisymmetric,
and
$Q' = \left( \!\! \begin{array}{cc}
                    Q'_{11}  &  Q'_{12}  \\[1mm]
                    Q'_{21}  &  Q'_{22}
                \end{array}
      \!\! \right)$,
we have
\begin{subeqnarray}
     X'_1  & = &  P' - J_{2m} Q'_{11}   \\
     X'_2  & = &  - J_{2m} Q'_{12}   \\
     Y'_{11}  & = &  Q'_{11} - (Q'_{11})^{\rm T}   \\
     Y'_{22}  & = &  Q'_{22} - (Q'_{22})^{\rm T}   \\
     Y'_{12}  & = &  Q'_{12} - (Q'_{21})^{\rm T}
\end{subeqnarray}
We can choose arbitrarily the symmetric parts of $Q'_{11}$ and $Q'_{22}$;
then the remaining unknowns are uniquely determined.
\qed

In preparation for the proof of Lemma~\ref{lemma.decomp.5},
it is convenient to introduce some simple notation
for decomposing rectangular matrices.
Given an $m \times n$ matrix $Y$,
we define the strictly lower-triangular $m \times m$ matrix $\scrl(Y)$
by taking the strictly lower-triangular part of $Y$
and either deleting the last $n-m$ columns (if $m < n$)
or appending $m-n$ columns of zeros (if $m > n$).
Likewise, we define the upper-triangular $n \times n$ matrix $\scru(Y)$
by taking the upper-triangular part of $Y$
and either deleting the last $m-n$ rows (if $m > n$)
or appending $n-m$ rows of zeros (if $m < n$).
It follows immediately from these definitions that
\be
   Y  \;=\; \scrl(Y) \widehat{I}_{mn} \,+\, \widehat{I}_{mn} \scru(Y)
   \;.
\ee

%

\proofof{Lemma~\ref{lemma.decomp.5}}
Linearizing
$X_\alpha =
 P_{\alpha} \, \widehat{I}_{n_{\alpha} n_{\alpha+1}} \, P_{\alpha+1}^{-1}$
in a neighborhood of $X_\alpha = \widehat{I}_{n_{\alpha} n_{\alpha+1}}$
and $P_\alpha = I_{n_\alpha}$
by writing
$X_\alpha = \widehat{I}_{n_{\alpha} n_{\alpha+1}} + X'_\alpha$
and $P_\alpha = I_{n_\alpha} + P'_\alpha$,
we obtain $X'_\alpha = P'_\alpha \widehat{I}_{n_{\alpha} n_{\alpha+1}}
   - \widehat{I}_{n_{\alpha} n_{\alpha+1}} P'_{\alpha+1}$.

Let us decompose each $P'_\alpha$ as a sum of its
strictly lower-triangular part $L_\alpha$
and its upper-triangular part $U_\alpha$.
We therefore need to solve the equations
\be
   X'_\alpha 
   \;=\;
   (L_\alpha + U_\alpha) \widehat{I}_{n_{\alpha} n_{\alpha+1}}
   \,-\,
   \widehat{I}_{n_{\alpha} n_{\alpha+1}} (L_{\alpha+1} + U_{\alpha+1})
   \;,
 \label{eq.proof.lemma.decomp.5.a}
\ee
where $X_1,\ldots,X_\ell$ are considered as parameters
and $L_1,\ldots,L_{\ell+1},U_1,\ldots,U_{\ell+1}$ are considered as unknowns.
But let us prove a bit more,
namely that the matrices $U_1$ and $L_{\ell+1}$
can be considered as parameters (i.e.\ can be chosen arbitrarily).
Note first that the system \reff{eq.proof.lemma.decomp.5.a} can be solved by
\begin{subeqnarray}
   L_\alpha  & = &
     \scrl\bigl(
         X'_\alpha - U_\alpha \widehat{I}_{n_{\alpha} n_{\alpha+1}}
                   + \widehat{I}_{n_{\alpha} n_{\alpha+1}} L_{\alpha+1}
          \bigr)
       \\[2mm]
   U_{\alpha+1}  & = &
     \scru\bigl(
         X'_\alpha - U_\alpha \widehat{I}_{n_{\alpha} n_{\alpha+1}}
                   + \widehat{I}_{n_{\alpha} n_{\alpha+1}} L_{\alpha+1}
          \bigr)
 \label{eq.proof.lemma.decomp.5.b}
\end{subeqnarray}
in the sense that any solution of \reff{eq.proof.lemma.decomp.5.b}
provides a solution of \reff{eq.proof.lemma.decomp.5.a}.
The equations \reff{eq.proof.lemma.decomp.5.b}
appear at first glance to be entangled,
i.e.\ $L_\alpha$ depends on $U_\alpha$ and vice versa.
But this is only an appearance, because the operators $\scrl$  and $\scru$
``see'', respectively, only the strictly-lower-triangular and
upper-triangular parts of the matrix on which they act.
Therefore, the system \reff{eq.proof.lemma.decomp.5.b} can be rewritten as
\begin{subeqnarray}
   L_\alpha  & = &
     \scrl\bigl(
         X'_\alpha + \widehat{I}_{n_{\alpha} n_{\alpha+1}} L_{\alpha+1}
          \bigr)
       \\[2mm]
   U_{\alpha+1}  & = &
     \scru\bigl(
         X'_\alpha - U_\alpha \widehat{I}_{n_{\alpha} n_{\alpha+1}}
          \bigr)
 \label{eq.proof.lemma.decomp.5.c}
\end{subeqnarray}
But these latter equations can manifestly be solved sequentially
for $L_\ell,\ldots,L_1$ (given $L_{\ell+1}$)
and for $U_2,\ldots,U_{\ell+1}$ (given $U_1$).
\qed

{\bf Remark.}
It is instructive to count parameters and variables in
Lemma~\ref{lemma.decomp.5}.
The parameters are $X_1,\ldots,X_\ell$,
and their number is $N_p = \sum_{\alpha=1}^\ell n_\alpha n_{\alpha+1}$.
(The matrices $U_1$ and $L_{\ell+1}$,
 which can be chosen arbitrarily, do not count as extra parameters
 because they merely redefine the matrices $X_1$ and $X_\ell$, respectively.)
The variables are $L_1,\ldots,L_{\ell}$ and $U_2,\ldots,U_{\ell+1}$,
and their number is
\be
   N_v  \;=\;  {n_1 (n_1-1) \over 2} \,+\, n_2^2 \,+\, \ldots \,+\, n_\ell^2
               \,+\, {n_{\ell+1} (n_{\ell+1}+1) \over 2}
   \;.
\ee
Therefore
\be
   N_v - N_p  \;=\;  \sum_{\alpha=1}^\ell
    {(n_\alpha - n_{\alpha+1}) (n_\alpha - n_{\alpha+1} - 1)  \over 2}
   \;\ge\;  0
   \;.
\ee
The $N_v - N_p$ extra variables were fixed by our choice of
the operators $\scrl$ and $\scru$:
we decided to append $m-n$ columns of zeros to $\scrl$ when $m > n$,
and $n-m$ rows of zeros to $\scru$ when $m < n$,
but we could equally well have inserted an arbitrary
strictly-lower-triangular matrix of size $m-n$ into $\scrl$
and an arbitrary upper-triangular matrix of size $n-m$ into $\scru$.
This gives $(m-n)(m-n-1)/2$ additional variables in both cases,
which precisely accounts for $N_v - N_p$.
\qed

\subsection{Dilation-translation formula}  \label{app.gen.transl}

In Sections~\ref{sec.param.1} and \ref{sec.param.2}
we will need the following well-known generalization
of the translation formula \reff{eq.translation}:

\begin{lemma}[Dilation-translation formula]
   \label{lemma.diltrans}
Let $P(z)$ be a polynomial in a single indeterminate $z$,
with coefficients in a commutative ring $R$ containing the rationals,
and let $a$ and $b$ be indeterminates.  Then
\be
   \exp\left( (a+bz) {\partial \over \partial z} \right) \, P(z)
   \;=\;
   P\Bigl( e^b z \,+\, {e^b - 1 \over b} a \Bigr)
 \label{eq.transl2}
\ee
as an identity in the ring $R[z,a][[b]]$ of formal power series in $b$
whose coefficients are polynomials in $z$ and $a$.
[Here the exponential is defined by its Taylor series,
 as are $e^b$ and $(e^b - 1)/b$.]

In particular, this identity can be evaluated at any {\em nilpotent}\/
element $b \in R$, as both sides then reduce to finite sums.
\end{lemma}

{\bf Remark.}  If $b$ is nilpotent of order 2 (i.e.\ $b^2 = 0$),
then the formula simplifies further to
\be
   \exp\left( (a+bz) {\partial \over \partial z} \right) \, P(z)
   \;=\;
   P\Bigl( (1+b) z \,+\, (1 + \smfrac{b}{2}) a \Bigr)
\ee
In our applications in Sections~\ref{sec.param.1} and \ref{sec.param.2}
we will have $b^2 = 0$ and also $ba = 0$, in which case the
identity holds even when the ring $R$ does not contain the rationals.
\qed

\proofof{Lemma~\ref{lemma.diltrans}}
When $a=0$, the formula \reff{eq.transl2} states the well-known fact that
the operator $z \, \partial/\partial z$ generates dilations;
it is easily checked by applying both sides to $z^n$.
The general case is handled by the change of variables $w = z + a/b$.
\qed

{\bf Remark.}
The formula \reff{eq.transl2} is a special case of a more general
formula for operators of the form $\exp[ t g(z) \, \partial/\partial z]$:
\be
   \exp\left( t g(z) {\partial \over \partial z} \right) \, P(z)
   \;=\;
   P\Bigl( \widetilde{z}(t;z) \Bigr)
\ee
where $\widetilde{z}(t;z)$ is the solution of the differential
equation $d \widetilde{z}(t;z)/dt = g(\widetilde{z}(t;z))$
with initial condition $\widetilde{z}(0;z) = z$.
\qed

\section*{Acknowledgments}

We wish to thank Marco Polin for collaborating in early stages of
this work, and for many helpful comments.
We also wish to thank Anton Leykin for patiently answering
our many naive questions about Bernstein--Sato polynomials;
Malek Abdesselam for supplying historical references;
and Alex Scott for valuable conversations.
Finally, we are grateful to Nero Budur for explaining to us
the connection between our results and
the theory of prehomogeneous vector spaces,
and for drawing our attention to the recent work of Sugiyama \cite{Sugiyama_11}.

We wish to thank the Isaac Newton Institute for Mathematical Sciences,
University of Cambridge, for generous support during the programme on
Combinatorics and Statistical Mechanics (January--June 2008).
One of us (A.D.S.)\ also thanks the
Institut Henri Poincar\'e -- Centre Emile Borel
for hospitality during the programmes on
Interacting Particle Systems, Statistical Mechanics and Probability Theory
(September--December 2008)
and Statistical Physics, Combinatorics and Probability
(September--December 2009).
Finally, we thank the Laboratoire de Physique Th\'eorique 
of the \'Ecole Normale Sup\'erieure (Paris)
for hospitality during May--June 2010 and April--June 2011.

This research was supported in part by
U.S.\ National Science Foundation grants PHY--0099393 and PHY--0424082.


\begin{thebibliography}{199}

\bibitem{Abdesselam_feynman}  A. Abdesselam,
   Feynman diagrams in algebraic combinatorics,
    S\'em. Lothar. Combin. B49 (2003)
   article B49c,
   math.CO/0212121 at arXiv.org.

\bibitem{Abdesselam_03} A. Abdesselam,
    Grassmann--Berezin calculus and theorems of the matrix-tree type,
     Adv. Appl. Math. 33 (2004) 51--70,
    math.CO/0306396 at arXiv.org.

\bibitem{cayley_history}  A. Abdesselam, T. Crilly and A.D. Sokal,
   The tangled history of the ``Cayley'' identity
   $\det(\partial) (\det X)^s = s(s+1) \cdots (s+n-1) (\det X)^{s-1}$,
   in preparation.

   
\bibitem{Aslaksen_96}  H. Aslaksen, Quaternionic determinants,
    Math. Intelligencer 18 
    (1996) 57--65.

\bibitem{Atiyah_70}  M.F. Atiyah, Resolution of singularities and
   division of distributions,  Comm. Pure Appl. Math. 23  (1970)
   145--150.

\bibitem{Benner_00}  P. Benner, R. Byers, H. Fassbender, V. Mehrmann,
   D. Watkins, Cholesky-like factorizations of skew-symmetric matrices,
    Electron. Trans. Numer. Anal. 11 (2000) 85--93;
   addendum  11 (2000) 93A.

\bibitem{Bernstein_72}  I.N. Bern\v{s}te\u{\i}n,
   The analytic continuation of generalized functions with respect to a
   parameter,
    Funkcional. Anal. i Prilo\v{z}en. 6 
    (1972) 26--40
   [=  Funct. Anal. Appl. 6  (1972) 273--285].

\bibitem{Bernstein_69}  I.N. Bern\v{s}te\u{\i}n, S.I. Gel'fand,
   Meromorphy of the function $P^{\lambda}$,
    Funkcional. Anal. i Prilo\v{z}en. 3 
    (1969) 84--85
   [=  Funct. Anal. Appl. 3  (1969) 68--69].

\bibitem{Bjork_79}  J.-E. Bj\"ork,  Rings of Differential Operators,
   North-Holland, Amsterdam--Oxford--New York, 1979.

\bibitem{Blekher_82}  P.M. Blekher,
   Integration of functions in a space with a complex number of dimensions,
    Teoret. Mat. Fiz. 50 (1982) 370--382 
   [English translation:  Theor. Math. Phys. 50 (1982) 243--251].

\bibitem{Blind_97}  B. Blind,
   Distributions homog\`enes sur une alg\`ebre de Jordan,
    Bull. Soc. Math. France 125  (1997) 493--528.

\bibitem{Bunch_82}  J.R. Bunch, A note on the stable decomposition of
   skew-symmetric matrices,  Math. Comp. 38   (1982) 475--479.

\bibitem{Caianiello_73}  E.R. Caianiello,  Combinatorics and
   Renormalization in Quantum Field Theory, 
   Benjamin, Reading, Mass.--London--Amsterdam, 1973.

\bibitem{Canfield_81}  E.R. Canfield and S.G. Williamson,
   Hook length products and Cayley operators of classical invariant theory,
    Lin. Multilin. Alg. 9  (1981) 289--297.

\bibitem{Capelli_1882}  A. Capelli, Fondamenti di una teoria generale
   delle forme algebriche,  Atti Reale Accad. Lincei, Mem. Classe
   Sci. Fis. Mat. Nat. (serie 3)  12  (1882) 529--598.

\bibitem{Capelli_1887}  A. Capelli, Ueber die Zur\"uckf\"uhrung der
   Cayley'schen Operation $\Omega$ auf gew\"ohnliche Polar-Operationen,
    Math. Annalen 29  (1887) 331--338.

\bibitem{Capelli_1888}  A. Capelli, Ricerca delle operazioni invariantive
   fra pi\`u serie di variabili permutabili con ogni altra operazione
   invariantiva fra le stesse serie,
    Atti Reale Accad. Sci. Fis. Mat. Napoli (serie 2)  1 (1888)
   1--17.

\bibitem{Capelli_1890}  A. Capelli, Sur les op\'erations dans
   la th\'eorie des formes alg\'ebriques,
    Math. Annalen 37  (1890) 1--37.

\bibitem{Capelli_02}  A. Capelli,  Lezioni sulla Teoria delle Forme
   Algebriche, Pellerano, Napoli, 1902.

\bibitem{Caracciolo_04}  S. Caracciolo, J.L. Jacobsen, H. Saleur,
   A.D. Sokal, A. Sportiello, Fermionic field theory for trees and forests,
    Phys. Rev. Lett. 93  (2004) 080601,
   cond-mat/0403271 at arXiv.org.

\bibitem{CSS_hyperforests}  S. Caracciolo, A.D. Sokal, A. Sportiello,
   Grassmann integral representation for spanning hyperforests,
    J. Phys. A: Math. Theor. 40  (2007) 13799--13835,
   arXiv:0706.1509 [math-ph] at arXiv.org.

\bibitem{CSS_Capelli_08}  S. Caracciolo, A.D. Sokal, A. Sportiello,
   Noncommutative determinants, Cauchy--Binet formulae,
   and Capelli-type identities.  I.~Generalizations of the Capelli and
   Turnbull identities, 
    Electron. J. Combin. 16(1) (2009) \#R103,
   arXiv:0809.3516 [math.CO] at arXiv.org.

\bibitem{Caracciolo_et_al_in_prep}  S. Caracciolo, A.D. Sokal, A. Sportiello, 
    work in progress.

\bibitem{Carlson_86}  D. Carlson, What are Schur complements, anyway?,
    Lin. Alg. Appl. 74 (1986) 257--275.

\bibitem{Cayley_1846}  A. Cayley, On linear transformations,
    Cambridge and Dublin Math. J. 1 (1846) 104--122.
   [Also in  The Collected Mathematical Papers
    of Arthur Cayley, Cambridge University Press,
    Cambridge, 1889--1897, vol.~1, pp.~95--112.]

\bibitem{Cayley_collected}  A. Cayley,  The Collected Mathematical Papers
    of Arthur Cayley, 13 vols., Cambridge University Press,
    Cambridge, 1889--1897.
    [Also republished by Johnson Reprint Corp., New York, 1963.]

\bibitem{Chaiken_82}  S. Chaiken, A combinatorial proof of the
   all minors matrix-tree theorem,  SIAM J. Alg. Disc. Meth.
    3  (1982) 319--329.

\bibitem{Clebsch_1861}  A. Clebsch, Ueber symbolische Darstellung
   algebraischer Formen,  J. Reine Angew. Math. 59  (1861) 1--62.

\bibitem{Clebsch_1872}  A. Clebsch,  Theorie der bin\"aren
   algebraischen Formen, B.G.~Tuebner, Leipzig, 1872.

\bibitem{Cohn_03}  P.M. Cohn,  Basic Algebra: Groups, Rings and Fields,
   Springer-Verlag, London--Berlin--Heidelberg, 2003.

\bibitem{Cottle_74}  R.W. Cottle, Manifestations of the Schur complement,
    Lin. Alg. Appl. 8  (1974) 189--211.

\bibitem{Coutinho_95}  S.C. Coutinho,
    A Primer of Algebraic $D$-Modules,
   London Mathematical Society Student Texts \#33,
   Cambridge University Press, Cambridge, 1995.

\bibitem{Creutz_78}  M. Creutz, On invariant integration over $SU(N)$,
    J. Math. Phys. 19  (1978) 2043--2046.

\bibitem{Dolgachev_03}  I. Dolgachev,  Lectures on Invariant Theory,
   London Mathematical Society Lecture Note Series \#296,
   Cambridge University Press, Cambridge, 2003.

\bibitem{Etingof_99} P. Etingof, Note on dimensional regularization,
   in P.~Deligne  et al., (Eds.),
    Quantum Fields and Strings: A Course for Mathematicians,
   American Mathematical Society, Providence RI, 1999,
   vol.~1, pp.~597--607.

\bibitem{Faraut_94}  J. Faraut, A. Kor\'anyi,
    Analysis on Symmetric Cones,
   Oxford University Press, Oxford--New York, 1994, Chapter VII.

\bibitem{Foata_94}  D. Foata, D. Zeilberger, Combinatorial proofs of
   Capelli's and Turnbull's identities from classical invariant theory,
    Electron. J. Combin.  1  (1994) \#R1.

\bibitem{Fulmek_10}  M. Fulmek, Graphical condensation, overlapping Pfaffians
   and superpositions of matchings,
    Electron. J. Combin.  17  (2010) \#R83.

\bibitem{Fulton_91}  W. Fulton, J. Harris,
    Representation Theory: A First Course,
   Springer-Verlag, New York--Heidelberg, 1991, Appendix F.

\bibitem{Fulton_98}  W. Fulton, P. Pragacz,
    Schubert Varieties and Degeneracy Loci,
   Lecture Notes in Mathematics \#1689,
   Springer-Verlag, Berlin--Heidelberg, 1998, Appendix D.

\bibitem{Garding_48}  L. G{\aa}rding, Extension of a formula by Cayley
   to symmetric determinants,  Proc. Edinburgh Math. Soc. 
   8 (1948) 73--75.

\bibitem{Gelfand_54}  I. Gelfand, Some aspects of functional analysis and
   algebra, in  Proceedings of the International Congress of
   Mathematicians, Amsterdam, 1954, 
   Noordhoff, Groningen / North-Holland, Amsterdam, 1957, vol.~1, pp.~253--276.

\bibitem{Gelfand_64}  I.M. Gel'fand, G.E. Shilov,
    Generalized Functions, vol. 1,
   Academic Press, New York--London, 1964.

\bibitem{Golub_96} G.H. Golub, C.F. Van Loan,  Matrix
   Computations, 3rd edition, The Johns Hopkins University Press,
   Baltimore, 1996.

\bibitem{Gordan_1872}  P. Gordan, Ueber Combinanten,
    Math. Annalen 5  (1872) 95--122.

\bibitem{Grace_03} J.H. Grace, A. Young,  The Algebra of Invariants,
   Cambridge University Press, Cambridge, 1903, pp.~259--260.

\bibitem{Graham_94}  R.L. Graham, D.E. Knuth, O. Patashnik,
   Concrete Mathematics: A Foundation for Computer Science,
  2nd ed., Addison-Wesley, Reading, Mass., 1994.

\bibitem{Hamel_01}  A.M. Hamel, Pfaffian identities: a combinatorial approach,
   J. Combin. Theory A  94  (2001) 205--217.

\bibitem{Hironaka_64}  H. Hironaka, Resolution of singularities of an
   algebraic variety over a field of characteristic zero,
    Ann. Math. 79  (1964) 109--326.

\bibitem{Howe_89}  R. Howe, Remarks on classical invariant theory,
    Trans. Amer. Math. Soc. 313  (1989) 539--570
   and erratum  318 (1990) 823 .

\bibitem{Howe_91}  R. Howe, T. Umeda, The Capelli identity,
   the double commutant theorem, and multiplicity-free actions,
    Math. Ann. 290  (1991) 565--619.

\bibitem{Igusa_00}  J.-i. Igusa, An Introduction to the Theory of Local Zeta
   Functions
   (AMS/IP Studies in Advanced Mathematics \#14),
   American Mathematical Society, Providence RI, 2000.

\bibitem{Ishikawa_95}  M. Ishikawa, M. Wakayama,
   Minor summation formula of Pfaffians,
    Lin. Multilin. Alg. 39 (1995) 285--305 . 

\bibitem{Ishikawa_06}  M. Ishikawa, M. Wakayama,
   Applications of minor summation formula. III.~Pl\"ucker relations,
   lattice paths and Pfaffian identities,
    J. Combin. Theory A 113  (2006) 113--155. 

\bibitem{Isserlis_18}  L. Isserlis, On a formula for the product-moment
   coefficient of any order of a normal frequency distribution
   in any number of variables,
    Biometrika 12  (1918) 134--139.

\bibitem{Kashiwara_76}  M. Kashiwara, $B$-functions and holonomic systems.
   Rationality of roots of $B$-functions,
    Invent. Math. 38 (1976/77) 33--53.

\bibitem{Khekalo_01}  S.P. Kh\`ekalo, Riesz potentials in the space of
   rectangular matrices, and the iso-Huygens deformation of the
   Cayley-Laplace operator,
    Dokl. Akad. Nauk 376  (2001) 168--170
   [English translation in  Dokl. Math. 63  (2001) 35--37].

\bibitem{Khekalo_05}  S.P. Kh\`ekalo, The Cayley-Laplace differential
   operator on the space of rectangular matrices,
     Izv. Ross. Akad. Nauk Ser. Mat. 69  (2005) 195--224
   [English  translation in  Izv. Math. 69  (2005) 191--219].

\bibitem{Kimura_82}  T. Kimura, The $b$-functions and holonomy diagrams of
   irreducible regular prehomogeneous vector spaces,
   Nagoya Math. J. 85 (1982) 1--80.

\bibitem{Kimura_03}  T. Kimura, Introduction to Prehomogeneous Vector Spaces
   (Translations of Mathematical Monographs \#215),
   American Mathematical Society, Providence, RI, 2003.

\bibitem{Kinoshita_02}  K. Kinoshita, M. Wakayama,
   Explicit Capelli identities for skew symmetric matrices,
    Proc. Edinburgh Math. Soc. 45  (2002) 449--465.

\bibitem{Knuth_96}  D.E. Knuth, Overlapping pfaffians,
    Electron. J. Combin.  3, no.~2,  (1996) \#R5.

\bibitem{Kostant_91}  B. Kostant, S. Sahi, The Capelli identity,
   tube domains, and the generalized Laplace transform,
    Adv. Math. 87  (1991) 71--92.

\bibitem{Krause_00}  G.R. Krause, T.H. Lenagan,
    Growth of Algebras and Gelfand--Kirillov Dimension,
   rev.\ ed.,
   American Mathematical Society, Providence RI, 2000.

\bibitem{Lang_02}  S. Lang,  Algebra, revised 3rd ed.,
   Springer-Verlag, New York, 2002.

\bibitem{Leykin_01}  A. Leykin, Constructibility of the set of polynomials
   with a fixed Bernstein-Sato polynomial: an algorithmic approach,
    J. Symbolic Comput. 32  (2001) 663--675. 

\bibitem{Leykin_private}  A. Leykin, private communication (2006).

\bibitem{Minc_78}  H. Minc,  Permanents,
   Encyclopedia of Mathematics and its Applications \#6,
   Addison-Wesley, Reading MA, 1978.

\bibitem{Moon_94}  J.W. Moon, Some determinant expansions and the
   matrix-tree theorem,  Discrete Math. 124  (1994) 163--171.

\bibitem{Muir}  T. Muir,  The Theory of Determinants in the
   Historical Order of Development, 4 vols.,
   Macmillan, London, 1906--23.

\bibitem{Murnaghan_31}  F.D. Murnaghan, A. Wintner,
   A canonical form for real matrices under orthogonal transformations,
    Proc. Nat. Acad. Sci. USA 17  (1931) 417--420.

\bibitem{Muro_99}  M. Muro, Singular invariant hyperfunctions on the
    space of real symmetric matrices,
    Tohoku Math. J. 51  (1999) 329--364.
 
\bibitem{Muro_01}  M. Muro, Singular invariant hyperfunctions on the
    space of complex and quaternion Hermitian matrices,
     J. Math. Soc. Japan 53  (2001) 589--602.

\bibitem{Muro_03}  M. Muro, Singular invariant hyperfunctions on the
    square matrix space and the alternating matrix space,
     Nagoya Math. J. 169  (2003) 19--75.

\bibitem{Olver_99}  P.J. Olver,  Classical Invariant Theory,
   London Mathematical Society Student Texts \#44,
   Cambridge University Press, Cambridge, 1999.

\bibitem{Ouellette_81}  D.V. Ouellette, Schur complements and statistics,
    Lin. Alg. Appl. 36  (1981) 187--295.

\bibitem{Ozeki_80}  I. Ozeki, On the microlocal structure of a regular
   prehomogeneous vector space associated with $GL(8)$,
   Proc. Japan Acad. Ser. A Math. Sci. 56 (1980) 18--21. 

\bibitem{Ozeki_90}  I. Ozeki, On the microlocal structure of the regular
    prehomogeneous vector space associated with $SL(5) \times GL(4)$,
    Publ. Res. Inst. Math. Sci. Kyoto Univ. 26 (1990) 539--584. 

\bibitem{Prasolov_94}  V.V. Prasolov,  Problems and Theorems in
   Linear Algebra,
   Translations of Mathematical Monographs \#134,
   American Mathematical Society, Providence, RI, 1994.

\bibitem{Rais_72}  M. Ra\"{\i}s, Distributions homog\`enes sur des
   espaces de matrices,
    M\'em. Soc. Math. France 30  (1972) 3--109.

\bibitem{Ricci_86}  F. Ricci, E.M. Stein,
   Homogeneous distributions on spaces of Hermitean matrices,
    J. Reine Angew. Math. 368  (1986) 142--164.

\bibitem{Rubenthaler_87}  H. Rubenthaler, G. Schiffmann,
   Op\'erateurs diff\'erentiels de Shimura et espaces pr\'ehomog\`enes,
    Invent. Math. 90  (1987) 409--442.

\bibitem{Rubin_06}  B. Rubin, Riesz potentials and integral geometry
    in the space of rectangular matrices,
     Adv. Math. 205  (2006) 549--598.

\bibitem{Saito_07}  M. Saito, Multiplier ideals, $b$-function,
    and spectrum of a hypersurface singularity,
     Compositio Math. 143  (2007) 1050--1068.

\bibitem{Sato_06}  F. Sato, K. Sugiyama,
   Multiplicity one property and the decomposition of $b$-functions,
   Internat. J. Math. 17 (2006) 195--229.

\bibitem{Sato_80}  M. Sato, M. Kashiwara, T. Kimura, T. Oshima,
   Micro-local analysis of prehomogeneous vector spaces,
   Inventiones Math. 62 (1980) 117--179.

\bibitem{Sato_74}  M. Sato, T. Shintani,
    On zeta functions associated with prehomogeneous vector spaces,
     Ann. Math. 100  (1974) 131--170.

\bibitem{Schur_17}  J. [Isaai] Schur,
   \"Uber Potenzreihen, die im Innern des Einheitskreises beschr\"ankt sind,
    J. Reine Angew. Math. 147  (1917) 205--232
   [English translation:  On power series which are bounded in the interior
    of the unit circle. I., in I.~Gohberg, ed.,
     I.~Schur Methods in Operator Theory and Signal Processing,
    Birkh\"auser, Basel, 1986, pp.~31--59.]

\bibitem{Schur_68}  I. Schur,  Vorlesungen \"uber Invariantentheorie,
   Springer-Verlag, Berlin--Heidelberg--New York, 1968.

\bibitem{Shimura_84}  G. Shimura, On differential operators attached to
   certain representations of classical groups,
    Invent. Math. 77  (1984) 463--488.

\bibitem{Sokal_bcc2005}  A.D. Sokal,
The multivariate Tutte polynomial
(alias Potts model) for graphs and matroids,
in  Surveys in Combinatorics, 2005,
edited by B.S. Webb,
Cambridge University Press, Cambridge--New York, 2005, pp.~173--226,
math.CO/0503607 at arXiv.org.

\bibitem{Speer_76}  E.R. Speer,  Dimensional and analytic renormalization,
   in G. Velo, A.S. Wightman, eds.,  Renormalization Theory,
   Proceedings of the NATO Advanced Study Institute, Erice, 1975,
   NATO Advanced Study Institute Series C: Math. and Phys. Sci., Vol. 23,
   Reidel, Dordrecht, 1976, pp. 25--93.

\bibitem{Stanley_97}  R.P. Stanley,
    Enumerative Combinatorics, vol.~1,
   Cambridge University Press, Cambridge--New York, 1997.

\bibitem{Stein_54}  P. Stein, An extension of a formula of Cayley,
    Proc. Edinburgh Math. Soc. 9  (1954) 91--99.

\bibitem{Stein_82}  P.R. Stein, On an identity from classical invariant theory,
    Lin. Multilin. Alg. 11  (1982) 39--44.

\bibitem{Stembridge_90}  J.R. Stembridge,
   Nonintersecting paths, pfaffians, and plane partitions,
    Adv. Math. 83 (1990) 96--131,
   Section~2 and Lemma~4.2.

\bibitem{Streater_64}  R.F. Streater, A.S. Wightman,
    PCT, Spin and Statistics, and All That,
   Benjamin/Cummings, Reading, Massachusetts, 1964, $2^{nd}$ ed.~1978.

\bibitem{Sugiyama_02}  K. Sugiyama, $b$-functions of regular simple
   prehomogeneous vector spaces, Tsukuba J. Math. 26 (2002) 407--451.

\bibitem{Sugiyama_05}  K. Sugiyama, $b$-function of a prehomogeneous vector
   space with no regular component,
   Comment. Math. Univ. St. Pauli 54 (2005) 99--119. 

\bibitem{Sugiyama_11}  K. Sugiyama, $b$-functions associated with quivers
   of type~${\sf A}$, Transformation Groups 16 (2011), 1183--1222,
   arXiv:1005.3596 [math.RT] at arXiv.org.

\bibitem{Tamm_1}  U. Tamm, The determinant of a Hessenberg matrix
   and some applications in discrete mathematics,
   \url{http://www.mathematik.uni-bielefeld.de/ahlswede/pub/tamm/hessen.ps}
   (undated).

\bibitem{Tamm_2}  U. Tamm, Majorization in lattice path enumeration
   and creating order,
   University of Bielefeld preprint 00-108,
   \url{http://www.math.uni-bielefeld.de/sfb343/preprints/pr00108.ps.gz}
   (2000).

\bibitem{Turnbull_28}  H.W. Turnbull,  The Theory of Determinants,
   Matrices, and Invariants,
   Blackie \& Son, London--Glasgow, 1928, pp.~114--116.
   2nd edition: Blackie \& Son, London--Glasgow, 1945.
   3rd edition: Dover, New York, 1960.

\bibitem{Turnbull_48}  H.W. Turnbull, Symmetric determinants and the
   Cayley and Capelli operators,  Proc. Edinburgh Math. Soc. 
   8 (1948) 76--86.

\bibitem{Turnbull_49}  H.W. Turnbull, Note upon the generalized Cayleyean
   operator,  Canad. J. Math. 1  (1949) 48--56.

\bibitem{Ukai_03}  K. Ukai, $b$-functions of prehomogeneous vector spaces
   of Dynkin--Kostant type for exceptional groups,
   Compositio Math. 135 (2003) 49--101.
   
\bibitem{Umeda_98}  T. Umeda, The Capelli identities, a century after,
    Amer. Math. Soc. Transl. 
    183  (1998) 51--78.

\bibitem{Vivanti_1890}  G. Vivanti, Alcune formole relative all'operazione
   $\Omega$,  Rendiconti Circ. Mat. Palermo 4 (1890) 261--268.

\bibitem{Wakatsuki_04}  S. Wakatsuki, $b$-functions of regular 2-simple
   prehomogeneous vector spaces associated to the symplectic group and
   the orthogonal group,
   Comment. Math. Univ. St. Pauli 53 (2004) 121--137.

\bibitem{Wallace_53}  A.H. Wallace, A note on the Capelli operators
   associated with a symmetric matrix,  Proc. Edinburgh Math. Soc. 
   9  (1953) 7--12.

\bibitem{Walther_05}  U. Walther, Bernstein--Sato polynomial versus
   cohomology of the Milnor fiber for generic hyperplane arrangements,
    Compositio Math. 141  (2005) 121--145.

\bibitem{Weitzenbock_23}  R. Weitzenb\"ock,  Invariantentheorie,
   P.~Noordhoff, Groningen, 1923, pp.~15--16.

\bibitem{Weyl_46}  H. Weyl,  The Classical Groups, Their Invariants
    and Representations, 2nd ed.\ 
    Princeton University Press, Princeton NJ, 1946.

\bibitem{Weyl_50}  H. Weyl,  The Theory of Groups and Quantum Mechanics,
   Dover, New York, 1950, Appendix~3.

\bibitem{Wick_50}  G.C. Wick, The evaluation of the collision matrix,
    Phys. Rev. 80  (1950) 268--272.

\bibitem{Williamson_87}  S.G. Williamson, Generic common minor expansions,
    Lin. Multilin. Alg. 20  (1987) 253--279, Section~4.

\bibitem{Xu_03}  H. Xu, An SVD-like matrix decomposition and its applications,
    Lin. Alg. Appl. 368  (2003) 1--24.

\bibitem{Zhang_05}  F. Zhang, ed.,  The Schur Complement and its
   Applications, Numerical Methods and Algorithms \#4,
   Springer-Verlag, New York, 2005.

\bibitem{Zinn-Justin} J. Zinn-Justin,
     Quantum Field Theory and Critical Phenomena, 4th ed.,
     Clarendon Press, Oxford, 2002.

\end{thebibliography}
\end{document}